\documentclass[a4paper,english,11pt]{article}

\usepackage[a4paper,left=2.5cm,right=2.5cm,top=2.5cm,bottom=2.5cm]{geometry}
\usepackage{bbm}
\usepackage{latexsym,amsthm,amsmath,amscd,amssymb,color,url,amsfonts}
\usepackage[utf8]{inputenc}
\usepackage[T1]{fontenc} 
\usepackage{mathrsfs}
\usepackage[english]{babel}
\usepackage{libertine}
\usepackage{caption} 
\usepackage{graphicx}
\usepackage{MnSymbol}
\usepackage{hyperref}
\usepackage{setspace}
\usepackage{chngpage}
\usepackage{pgf,tikz}
\usetikzlibrary{arrows}
\usepackage[all]{xy}
\usepackage{appendix}
\usepackage{wrapfig}
\usepackage{array,booktabs,arydshln,xcolor}
\usetikzlibrary{positioning}
\usetikzlibrary{arrows}
\usepackage{multirow}
\usepackage{tikz-cd}
\usepackage{multicol}
\usepackage{geometry} 
\tikzcdset{arrow style=tikz,diagrams={>=stealth}}

\newtheorem{lemma}{Lemma}[section]
\newtheorem{thm}[lemma]{Theorem}
\newtheorem{thmintro}{Theorem}[section]

\newtheorem{propointro}[thmintro]{Proposition}
\newtheorem*{thm*}{Theorem}

\newtheorem{rem}[lemma]{Remark}
\newtheorem*{rem*}{Remark}
\newtheorem{propo}[lemma]{Proposition}
\newtheorem{propr}[lemma]{Property}

\newtheorem{cor}[lemma]{Corollary}
\newtheorem{conj}[lemma]{Conjecture}

\newtheorem{example}[lemma]{Example}

\newtheorem{defn}[lemma]{Definition}

\newtheorem{nota}[lemma]{Notation}

\title{Geometric rigidity of quasi-isometries in horospherical products}

\author{\textsc{Tom Ferragut}}

\date{\today}

\begin{document}

\maketitle

\begin{abstract}
%

We prove that quasi-isometries of horospherical products of hyperbolic spaces are geometrically rigid in the sense that they are uniformly close to product maps, this is a generalization of the result obtained by Eskin, Fisher and Whyte in \cite{EFW1}. Our work covers the case of solvable Lie groups of the form $\mathbb{R}\ltimes (N_1\times N_2)$, where $N_1$ and $N_2$ are nilpotent Lie groups, and where the action on $\mathbb{R}$ contracts the metric on $N_1$ while extending it on $N_2$. We obtain new quasi-isometric invariants and classifications for these spaces.
\end{abstract}

\tableofcontents

\section*{Introduction}

\addtocontents{toc}{\setcounter{tocdepth}{-1}}

Let $(X,d_X)$ and $(Y,d_Y)$ be two Gromov hyperbolic spaces. Their \textbf{horospherical product}, denoted by $X\bowtie Y$ is constructed by combining $X$ and $Y$, and lies in the direct product $X\times Y$. It has no longer negative curvature, however its geometry is still very rigid (see Section \ref{SecDefHoro} for the definition). This way of combining two hyperbolic spaces appeared to unify the construction of metric spaces such as Diestel-Leader graphs, treebolic spaces and Sol geometries, which are the horospherical products constructed out of a regular infinite tree or the hyperbolic plane $\mathbb{H}_2$.


\subsection*{Quasi-isometric classification and existing rigidity results}

In \cite{Gromov}, a mainstay of geometric group theory, Gromov points out the importance of quasi-isometric invariants in groups. The quasi-isometric classification of groups, or metric spaces, has since been a wide and prolific research domain (see \cite{Kapo} for a nice survey on this topic). For the family of solvable groups, there are still a lot of open cases. 

The first result was obtained in \cite{FB1} where Farb and Mosher provided a quasi-isometric classification of solvable Baumslag-Solitar groups $\mathrm{BS}(1,n)$. Then Eskin, Fisher and Whyte obtained the quasi-isometric classification of lamplighter groups and Sol geometries in \cite{EFW1} and \cite{EFW2}. In both the works \cite{FB1} and \cite{EFW1}, the horospherical product construction of their respective groups is crucial in their proofs.

The paper \cite{EFW1} also permitted to answer a question asked by Woess in \cite{Woess3} about the existence of vertex-transitive graphs not quasi-isometric to any Cayley graph. Eskin, Fisher and Whyte showed that when $m$ and $n$ are coprime integers, the Diestel-Leader graphs $T_m\bowtie T_n$ are such graphs.
\\\\Throughout \cite{Peng}, \cite{Peng2} and \cite{Dymarz}, using similar methods as in \cite{EFW1} and \cite{EFW2}, Peng and Dymarz generalized the description of the quasi-isometries for Lie groups of the form $\mathbb{R}\ltimes \mathbb{R}^p$. In \cite{Peng} and \cite{Peng2}, Peng proved that a subgroup of finite index of the quasi-isometry group of Lie groups of the form $\mathbb{R}^m\ltimes\mathbb{R}^n$ is a product of groups of bilipschitz maps.

\subsection*{Statement of results}\label{SecRePers}

The main goal of our work is to generalize the methods and techniques developed by Eskin, Fisher and Whyte to a wider set of horospherical products $X\bowtie Y$. In order to do that, the spaces $X$ and $Y$ are endowed with appropriate measures (see Definition \ref{DefAdmissible2}). Once endowed with suitable measures, $X$ and $Y$ are called horopointed admissible spaces.
\\\\To be more precise let $X$ (respectively $X'$, $Y$, $Y'$) be a horopointed admissible space with exponential growth parameter $m$ (respectively $m'$, $n$, $n'$). When $X$ is a regular tree, the parameter $m$ is related to the degree of $X$. When $X$ is a negatively curved Lie group $\mathbb{R}\ltimes_A N$, the parameter $m$ is $\mathrm{tr}(A)$, the trace of $A$. 

Let $\Phi:X\bowtie Y\to X'\bowtie Y'$ be a quasi-isometry. The map $\Phi$ is called a \textbf{product map} if and only if there exist two maps $\Phi^X : X\to X$ (or $\Phi^X : X\to Y$) and $\Phi^Y:Y\to Y$ (or $\Phi^Y : Y\to X$) such that for all $(x,y)\in X\bowtie Y$ we have either:
\begin{align*}
\Phi(x,y)=\left(\Phi^X(x),\Phi^Y(y)\right)\quad\text{or}\quad \Phi(x,y)=\left(\Phi^Y(y),\Phi^X(x)\right).
\end{align*}
Our main theorem states that, when $m> n$ and $m'>n'$, any quasi-isometry $\Phi:X\bowtie Y\to X'\bowtie Y'$ is close to a product map.

\begin{thmintro}[Geometric rigidity]\label{ThmDintro}\ 
\\Let $X$, $X'$, $Y$ and $Y'$ be horo-pointed admissible measured metric spaces with $m> n$ and $m'>n'$ and let $\Phi:X\bowtie Y\to X'\bowtie Y'$ be a quasi-isometry. Then there exist two quasi-isometries $\Phi^X:X\to X'$ and $\Phi^Y:Y\to Y'$ such that:
\begin{align*}
d_{\bowtie}\left(\Phi,\left(\Phi^X,\Phi^Y\right)\right)<+\infty.
\end{align*}
\end{thmintro}

This is a generalization of Theorems $2.1$ and $2.3$ of \cite{EFW1}. While completing the proof of this result, we obtained a first quasi-isometry invariant in horospherical products.

\begin{thmintro}
When $m>n$, the parameter $\frac{m}{n}$ is a quasi-isometry invariant. 
\end{thmintro}

Let $\mathbb{R}\ltimes_{A_1}N_1$ and $\mathbb{R}\ltimes_{A_2}N_2$ be two simply connected, negatively curved, solvable Lie groups (also called \textit{Heintze groups}). In Chapter \ref{SecExample} we show that this couple of Heintze groups is admissible, and that the condition $m> n$ is equivalent to $tr(A_1)> tr(A_2)$. We obtain a necessary condition for the existence of a quasi-isometry on solvable Lie groups. The horospherical product of these two Heintze groups is isomorphic to
\begin{align*}
G:=\mathbb{R}\ltimes_{\mathrm{Diag}(A_1,-A_2)} (N_1\times N_2),
\end{align*}
defined by the diagonal action of $\mathbb{R}$, $t\mapsto (\exp(tA_1),\exp(-tA_2))$ on $N_1\times N_2$. 

We say that $G$ is \textit{Carnot-Sol type} if $N_1$ and $N_2$ are Carnot groups and if $A_1$ and $A_2$ are multiples of Carnot derivations of $N_1$ and $N_2$ respectively. In the literature (see \cite{Pansu} for example), Carnot type stands for Lie groups with $N_2=\{1\}$. Here we extend the denominations to non-hyperbolic Lie groups.
  
Using the previous quasi-isometry invariants we obtain the following quasi-isometry classification.

\begin{thmintro}
Let $G=\mathbb{R}\ltimes_{\mathrm{Diag}(A_1,-A_2)} (N_1\times N_2)$ and $G'=\mathbb{R}\ltimes_{\mathrm{Diag}(A_1',-A_2')} (N_1'\times N_2')$ be Carnot-Sol type, non-unimodular Lie groups, then
\begin{equation}
G\text{ and }G'\text{ are quasi-isometric} \quad \Leftrightarrow \quad G\text{ and }G'\text{ are isomorphic}.
\end{equation}
\end{thmintro}

The case where $N_2= \{1\}$ is treated in Corollary 12.4 of \cite{Pansu}. 

Recall that a group $G$ is called \textit{metabelian} if $[G,G]$ is abelian (when both $N_1$ and $N_2$ are euclidean spaces). In this case, a similar quasi-isometry classification is deduced from the work of Peng \cite{Peng} and \cite{Peng2}. Both the quasi-isometry classification for the metabelian groups and for Carnot-Sol type groups are special cases of the conjecture $19.113$ of \cite{COR} that we recall.

\begin{conj}
Let $S$ and $S'$ be completely solvable Lie groups. Then $S$ and $S'$ are quasi-isometric if and only if they are isomorphic.
\end{conj}

Classifying completely solvable Lie groups up to quasi-isometry would yield the quasi-isometry classification of all connected Lie groups, see \cite{COR2}. 

For $i\in\{1,2\}$, let $N_i$ and $N_i'$ be two simply connected, nilpotent groups and let $A_i\in \mathrm{Lie}(N_i)$ and $A_i'\in \mathrm{Lie}(N_i')$ be derivations. Let $G:= \mathbb{R}\ltimes_{\mathrm{Diag}(A_1,-A_2)} (N_1\times N_2)$ and $G':=\mathbb{R}\ltimes_{\mathrm{Diag}(A_1',-A_2')} (N_1'\times N_2')$. In this general setting of horospherical products of Heintze groups we have the following necessary conditions for being quasi-isometric. 

\begin{propointro}
Let us assume that $\mathrm{tr}(A_1)>\mathrm{tr}(A_2)$ and $\mathrm{tr}(A_1')>\mathrm{tr}(A_2')$. If $G$ and $G'$ are quasi-isometric, then we have that for $i\in\{1,2\}$
\begin{enumerate}
\item $N_i$ and $N_i'$ are bilipschitz;
\item $A_i$ and $\frac{\mathrm{tr}(A_1)}{\mathrm{tr}(A_1')}A_i'$ share the same characteristic polynomial.
\end{enumerate}
\end{propointro}

With the  same setting, using the geometric rigidity on self quasi-isometries of this family of solvable Lie groups, we provide a characterisation of their quasi-isometry group.

Recall that for $F$ a metric space, $\mathrm{QI}(F)\big/\mathord{\sim}$ is the group of self quasi-isometries of $F$, up to finite distance. (This equivalence relation is required since a quasi-isometry only has a coarse inverse.) Recall also that $\mathrm{Bilip}(F)$ stands for the group of self bi-Lipschitz maps of $F$. Then we have: 
\begin{thmintro}\label{ThmIntroQigroup}
If $~tr(A_1)\neq tr(A_2)$:
\begin{align}
\mathrm{QI}\left(\mathbb{R}\ltimes_{\mathrm{Diag}(A_1,-A_2)} (N_1\times N_2)\right)\big/\mathord{\sim}=\mathrm{Bilip}\left(N_1\right)\times \mathrm{Bilip}.\left(N_2\right)
\end{align}
\end{thmintro}
Here we choose the horospherical product metric on $\mathbb{R}\ltimes_{\mathrm{Diag}(A_1,-A_2)} (N_1\times N_2)$. 

In the course of this proof we also obtain that any self quasi-isometry of $\mathbb{R}\ltimes_{\mathrm{Diag}(A_1,-A_2)} (N_1\times N_2)$ is a rough isometry. Le Donne, Pallier and Xie proved in \cite{LDPX} that when you change the left-invariant Riemannian metric of one of these solvable Lie groups, the identity map is a rough similarity. Hence self quasi-isometries are rough isometries with respect to any left-invariant distances.
%
%

\subsection*{Outline of the proof}

Let $X$ and $Y$ be two Gromov hyperbolic spaces, and let $\beta_X:X\to\mathbb{R} $ and $\beta_Y:Y\to \mathbb{R}$ be two Busemann functions. We call height functions $h_X$ and $h_Y$ the opposite of the Busemann functions. The horospherical product of $X$ and $Y$, denoted by $X\bowtie Y$, is defined as the set of points in $X\times Y$ such that the two Busemann functions (or the height functions) add up to zero.
\begin{equation}
X\bowtie Y:=\lbrace (x,y)\in X\times Y\mid \beta_X(x)+\beta_Y(y)=0 \rbrace.\nonumber
\end{equation}
A Busemann function is associated with a unique point on the boundary. We call any geodesic ray in the equivalence class of this point a \textit{vertical} geodesic ray. 
\\\\In order to generalize the proof of Eskin, Fisher and Whyte developed in \cite{EFW1} and \cite{EFW2}, the horospherical products have to be equipped with appropriate measures presented in Definition \ref{DefAdmissible2}.

Briefly speaking, for the measured space $(X,\mu^X)$, the measure $\mu^X$ must verify three assumptions. Assumption $(E1)$ allows us to disintegrate $\mu^X$ on its horospheres, assumption $(E2)$ provides us with a bounded geometry on horospheres and $(E3)$ ensures an exponential contraction (of exponent $m$) of the horospheres' measures in the upward vertical direction.
\\\\Let $X$ (respectively $X'$, $Y$, $Y'$) be a horopointed admissible space with exponential growth parameter $m$ (respectively $m'$, $n$, $n'$).

Most of this paper focuses on proving Theorem \ref{ThmDintro}. To do so we will use three major tools:

\begin{enumerate}
\item[$\bullet$] The \textit{coarse vertical quadrilaterals}, which are realised by four points (the vertices) whose neighbourhoods are linked by vertical geodesics (the edges). In Lemma \ref{LemmaTetrahedron}, we show that coarse vertical quadrilaterals are rigid: two of the four points almost share the same $X$-coordinate and the two other almost share the same $Y$-coordinate. 
\item[$\bullet$] \textit{Box Tilings} of different scales for $X\bowtie Y$, suitable for the vertical flow. The boxes correspond to euclidean rectangular cuboid in the Sol geometry.
\item[$\bullet$] The \textit{coarse differentiation}: given a quasi-isometry $\Phi:X\bowtie Y\to X'\bowtie Y'$, there exists a suitable scale $R$ for the box tilling of $X\bowtie Y$. Suitable here means that the image by $\Phi$ of most vertical geodesic segments of length $R$ are close to a vertical geodesic segment.
\end{enumerate}

With these tools, the proof can be summarized as follows. Let $\Phi:X\bowtie Y\to X'\bowtie Y'$ be a quasi-isometry.

\begin{description}
\item[Step 1] By the \textit{coarse differentiation}, there exists a scale $R$ such that in the box tilling at scale $R$ of $X\bowtie Y$, the quasi-isometry $\Phi$ mostly preserve the vertical direction on most of the boxes at scale $R$. It means that on most of the boxes, most vertical geodesic segments are sent close to a vertical geodesic segment by $\Phi$.
\item[Step 2] Then in most of the boxes at scale $R$, most of the vertical quadrilateral are sent close to vertical quadrilateral by $\Phi$. Therefore, by the rigidity property of these configurations, on most of the boxes $\mathcal{B}$ the quasi-isometry $\Phi$ is close to a product map $\hat{\Phi}_{|\mathcal{B}}=\left(\hat{\Phi}^X,\hat{\Phi}^Y\right)$ or $\left(\hat{\Phi}^Y,\hat{\Phi}^X\right)$.
\item[Step 3] If $m>n$ and $m'>n'$ then all product maps have the form $\hat{\Phi}_{\mathcal{B}}=\left(\Phi^X,\Phi^Y\right)$. Therefore by \textit{gluing} them together, we show that there exists $L>>R$ such that on \textbf{all} boxes at scale $L$, $\Phi$ is close to a product map $\hat{\Phi}=\left(\Phi^X,\Phi^Y\right)$.
\item[Step 4] We show that $\Phi$ quasi-respect the height, then we use this last result on $\Phi^{-1}$ to show that $\Phi$ send \textbf{all} vertical geodesics close to vertical geodesics. Therefore all vertical quadrilateral configurations are preserved by $\Phi$, hence $\Phi$ itself is close to a product map on all $X\bowtie Y$.
\end{description} 

A major technical issue in this proof is to manage the notion of "almost all" vertical geodesic segments having a certain property. The disintegrable measure $\mu$ of assumption $(E1)$ is not suited for this role since it concentrates the measure of a box on its bottom part. Therefore we introduce another disintegrable measure $\lambda$, constructed from $\mu$, which (almost) equally weights the level-sets of the height function $h$ in boxes.   

Such a measure $\lambda^X$ on $X$, together with a similar measure $\lambda^Y$ on $Y$, allows us to define a suitable measure (later denoted by $\eta$) on the family of vertical geodesics contained in a box $\mathcal{B}\subset X\bowtie Y$.
\\\\The geometric rigidity has useful consequences when we understand the boundaries of $X$ and $Y$. In this case, Theorem \ref{ThmDintro} leads to a description of the quasi-isometry-group of $X\bowtie Y$. In the last section of this paper, we detail such a description for the horospherical product of two Heintze groups.

\subsection*{Organization of the paper}

This work, about the geometric rigidity of quasi-isometries between two horospherical products, is organized as followed.
\begin{enumerate}
\item[$\bullet$] In Section \ref{ChapCoarsDifTetra} we display the coarse differentiation in our context, and we discuss particular quadrilateral configurations of $X\bowtie Y$.
\item[$\bullet$] Section \ref{ChapMeasure} focuses on developing all the measure theoretical tools required to achieve the rigidity results.
\item[$\bullet$] Then, in Section \ref{ChapBigProof}, we follow the structure of the proof proposed by Eskin, Fisher and Whyte in \cite{EFW1}, invoking technical tools of previous chapters when required.
\item[$\bullet$] In the last section we present an application of our theorem by providing new quasi-isometric classifications for some families of solvable Lie groups. We also provide a description of the quasi-isometry group of a wider family of solvable Lie groups.
\end{enumerate}

\newpage

\section*{Acknowledgement}

This work was supported by the University of Montpellier, as well as by the Japan Society for the Promotion of Science (JSPS) and the ERC of Enrico Le Donne.

I am deeply grateful to my advisors, Jeremie Brieussel and Constantin Vernicos, for their guidance, for their many reviews and comments. I also wish to thank Xiangdong Xie for the early and stimulating discussions regarding the paper of Eskin, Fisher and Whyte. My thanks extend to Ryokichi Tanaka and Gabriel Pallier for their detailed explanations of the implications of geometric rigidity in the context of Lie groups. Lastly, I am particularly grateful to the reviewer for their meticulous reading, detailed feedback, and thoughtful suggestions, which greatly improved the quality of this paper.

\addtocontents{toc}{\setcounter{tocdepth}{2}} 

\newpage

\section{Context}

\subsection{Gromov hyperbolic, Busemann spaces}

Let $\delta>0$, and let $(X,d_X)$ and $(Y,d_Y)$ be two $\delta$-hyperbolic spaces (See \cite{GDLH} or chap.III H. p.399 of \cite{BH} for more details on Gromov hyperbolic spaces). We present here the context in which we will construct our horospherical product. We require that $X$ and $Y$ are both \textbf{proper}, \textbf{geodesically complete}, \textbf{Busemann spaces}.
\begin{itemize}
\item A metric space is called proper if all closed metric balls are compact. 
\item A \textbf{geodesic line}, respectively \textbf{ray}, \textbf{segment}, of $X$ is the isometric image of a Euclidean line, respectively half Euclidean line, interval, in $X$. We denote by $[x_1,x_2]$ a geodesic segment linking $x_1\in X$ to $x_2\in X$.
\item A metric space $X$ is called geodesically complete if all geodesics are infinitely extendable.
\item A metric space is called Busemann if the distance between any couple of geodesics parametrized by arclength is a convex function. (See Chap.8 and Chap.12 of \cite{Papa} for more details on Busemann spaces.)
\end{itemize}
An important property of Gromov hyperbolic spaces is that they admit a nice compactification thanks to their \textbf{Gromov boundary}. We call two geodesic rays of $X$ equivalent if their images are at finite Hausdorff distance. Let $w\in X$ be a base point. We define $\partial_w X $, the Gromov boundary of $ X $, as the set of equivalence classes of geodesic rays starting from $ w $. While $ \partial_w X $ as a set depends on the choice of the base point $ w $, it is topologically independent of $ w $ under the cone topology. We denote the Gromov boundary simply by $ \partial X $ when the choice of $ w $ does not matter topologically. The cone topology on $X \cup \partial X $ restricts to the natural topology on $ X $, and with this topology, $ X \cup \partial X $ is compact (see \cite{BH} for further details on the cone topology). In this context, the Gromov boundary coincides with the visual boundary.
\\\\Let us fix a point $a\in\partial X$ on the boundary. We call \textbf{vertical geodesic ray}, respectively \textbf{vertical geodesic line}, any geodesic ray in the equivalence class $a$, respectively any geodesic line with one of its half-lines in $a$. The study of these specific geodesic rays is central in this work. 
\\The Busemann assumption removes some technical difficulties in a significant number of proofs in this work. If $X$ is a Busemann space in addition to being Gromov hyperbolic, for all $x\in X$ there exists a unique vertical geodesic ray, denoted by $V_x$, starting at $x$. In  fact the distance between two vertical geodesics starting at $x$ is a  convex and bounded function, hence decreasing and therefore constant equal to $0$.
\\The construction of the \textit{horospherical product} of two Gromov hyperbolic space $X$ and $Y$ requires the so called \textbf{Busemann functions}. Their definition is simplified by the Busemann assumption. Let us consider $\partial X$, the Gromov boundary of $X$ (which, in this setting, is the same as the visual boundary). Both the boundary $\partial X$ and $X\cup \partial X$, endowed with the natural Hausdorff topology, are compact. Then, given $a\in\partial X$ a point on the boundary, and $w\in X$ a base point, we define a Busemann function $\beta_{(a,w)}$ with respect to $a$ and $w$. Let $V_w$ be the unique vertical geodesic ray starting from $w$.
\begin{equation}
\forall\ x \in X,\ \beta_{(a,w)}(x):= \limsup\limits_{t\rightarrow + \infty}(d(x,V_w(t))-t).\nonumber
\end{equation}
In all our results, $X$ and $Y$ will be proper, geodesically complete, Gromov hyperbolic, Busemann spaces, with some additional assumption from time to time.

\subsection{Horospherical products}\label{SecDefHoro}

\begin{figure}
\begin{center}
\includegraphics[scale=1.1]{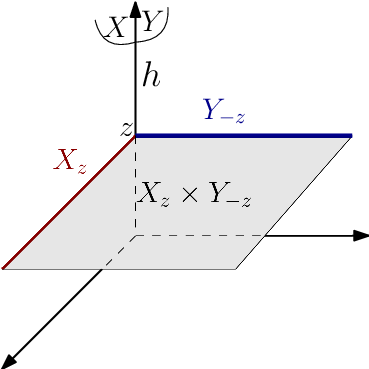} 
\end{center}
\caption{Horospherical product $X\bowtie Y$.}
\label{FigDefHoroProd}
\end{figure}

Let $a^X\in\partial X,a^Y\in\partial Y$ be points on the boundaries and let $w^X\in X,w^Y\in Y$ be base points. Let us denote by $h^X:=-\beta_{\left(a^X,w^X\right)}$ and $h^Y:=-\beta_{\left(a^Y,w^Y\right)}$ the two corresponding height functions. The \textbf{horospherical product} of $X$ and $Y$, relatively to $\left(a^X,w^X\right)$ and $\left(a^Y,w^Y\right)$ , denoted by $X\bowtie Y$ is defined by:
\begin{align*}
X\bowtie Y :=\left\lbrace(x,y)\in X\times Y\mid h^X(x)+h^Y(y)=0\right\rbrace.
\end{align*}
The set $X\bowtie Y$, can be seen as a diagonal in $X\times Y$. It is constructed by gluing $X$ with an upside down copy of $Y$ along their respective horospheres. This construction, illustrated in Figure \ref{FigDefHoroProd}, can also be seen as the union of the direct products between opposite horospheres in $X$ and $Y$
\begin{align*}
X\bowtie Y =\bigsqcup\limits_{z\in\mathbb{R}}X_z\times Y_{-z}.
\end{align*}

From now on, with a slight abuse, we omit the reference to the base points and points on the boundaries in the construction of the horospherical product. 

\begin{figure}[ht]
\begin{center}
\includegraphics[scale=0.90]{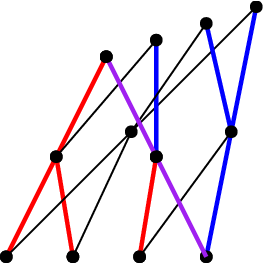}
\captionof{figure}{Small neighbourhood in $T_3\bowtie T_3$.}
\label{FigCayleyGrapheLamplighter}
\end{center}
\end{figure}

To study the geometry of a horospherical product $X\bowtie Y$, we make additional assumptions on $X$ and $Y$. We require them to be Gromov hyperbolic, Busemann, \textit{geodesically complete} and \textit{proper} metric spaces.
\begin{enumerate}
\item $X$ is \textit{geodesically complete} if and only if all geodesic segments of $X$ can be extended into a geodesic bi-infinite line.
\item $X$ is \textit{proper} if and only if all closed metric balls of $X$ are compact.
\end{enumerate}
If $X$ and $Y$ satisfy these two additional conditions, the horospherical product $X\bowtie Y$ is connected (see Property 3.11 of \cite{TF}).
\begin{example}\label{ExampleTrivial}
Let $X$ be a Gromov hyperbolic, Busemann, geodesically complete and proper metric space. Then $X\bowtie \mathbb{R}$ is isometric to $X$. In particular, if $~V^Y$ is a vertical geodesic line of $~Y$, $X\bowtie V^Y$ is an isometric embedding of $X$ in $X\bowtie Y$.
\end{example}
The three (non-trivial) first examples of horospherical products appeared independently in the literature. They correspond to the case where $X$ and $Y$ are either a regular infinite tree $T_m$ of degree $m$ or the hyperbolic plane $\mathbb{H}^2$.
\begin{enumerate}
\item $T_m\bowtie T_n$ is the Diestel-Leader graph $\mathrm{DL}(m,n)$. When $m=n$, this horospherical product is a Cayley graph of the lamplighter group $\mathbb{Z}\wr\mathbb{Z}_m$. See Figure \ref{FigCayleyGrapheLamplighter} for a subset of $T_3\bowtie T_3$.
\item $\mathbb{H}^{2,m}\bowtie \mathbb{H}^{2,n}$ is the Lie group $\mathbb{R}\ltimes_{(m,n)} \mathbb{R}^2=\mathrm{Sol}(m,n)$, one of the eight Thurston geometries when $m=n$. By $\mathbb{H}^{2,m}$ we mean the manifold $\mathbb{R}^2$ endowed with the infinitesimal Riemannian metric $\mathrm{d}s^2=e^{-2mz}\mathrm{d}x^2+\mathrm{d}z^2$. The action associated to the aforementioned semi-direct product is described by $(z,(x,y))\mapsto \left(e^{mz}x,e^{-nz}y\right)$.
\item $T_m \bowtie\mathbb{H}_2$ is a Cayley $2$-complex of the Baumslag-Solitar group $\mathrm{BS}(1,m)$.
\end{enumerate}
The awareness of them being identically constructed from Gromov hyperbolic spaces came later, and a survey on these three examples is provided by Wolfgang Woess in \cite{Woess}. 

An other approach, is to consider the hyperbolic plane $\mathbb{H}^{2,m}$ as the affine Lie group $\mathbb{R}\ltimes_m \mathbb{R}$ with action by multiplication $(z,x)\mapsto e^{mz}x$, and the Sol geometry Sol$(m,n)$ as the Lie group $\mathbb{R}\ltimes_{(m,n)} \mathbb{R}^2$. In this context we have that $(\mathbb{R}\ltimes_m \mathbb{R})\bowtie (\mathbb{R}\ltimes_n \mathbb{R}) = \mathbb{R}\ltimes_{(m,n)} \mathbb{R}^{2}$. The natural next step, is to consider which Lie group can be taken as a component in a horospherical product. 
\\\\A \textbf{Heintze group} is a Lie group of the form $\mathbb{R}\ltimes_A N$ with $N$ a nilpotent Lie group, with $A$ the derivation of the Lie algebra and where all eigenvalues of $A$ have positive real part. Heintze proved in \cite{Heintze} that any simply connected, negatively curved solvable Lie group is isomorphic to a Heintze group.
\\\\Moreover, a Busemann metric space is simply connected, hence any Gromov hyperbolic, Busemann Lie group is isomorphic to a Heintze group. Consequently, Heintze groups are natural candidates for the two components from which a horospherical product is constructed. Let $\mathbb{R}\ltimes_{A_1}N_1$ and $\mathbb{R}\ltimes_{A_2}N_2$ be two Heintze groups, we have
\begin{align*}
\left(\mathbb{R}\ltimes_{A_1}N_1\right)\bowtie\left(\mathbb{R}\ltimes_{A_1}N_1\right)=\mathbb{R}\ltimes_{\mathrm{Diag}(A_1,-A_2)} (N_1\times N_2),
\end{align*}
where $\mathrm{Diag}(A_1,-A_2)$ is the block diagonal matrix containing $A_1$ and $-A_2$ on its diagonal. 
\\\\In his paper \cite{Xie}, Xie classifies the subfamily of all negatively curved Lie groups $\mathbb{R}\ltimes \mathbb{R}^n$ up to quasi-isometry. In Chapter \ref{SecExample}, we provide a description of the quasi-isometry group of the horospherical product of two Heintze groups, namely the solvable Lie groups $\mathbb{R}\ltimes_{\mathrm{Diag}(A_1,-A_2)} (N_1\times N_2)$.

\subsection{Settings}

In this chapter we recall some material about horospherical products.

In order to lighten the notations, we will not fully describe the multiplicative and additive constants involved in inequalities. We will use the following notations instead.

\begin{nota}
Let $A,B\in\mathbb{R}$ and $e$ a parameter (set, real numbers, ...). Let us denote:
\begin{enumerate}
\item $A\preceq_e B$ if and only if there exists a constant $M(e)$ depending only on $e$ such that $A\leq M(e) B$.
\item $A\asymp_e B$ if and only if $B\preceq_e A\preceq_e B$.
\end{enumerate} 
\end{nota}

If the constant $M$ is a specific integer such as $2$, we will simply denote $A\preceq B$, and similarly $A\succeq B$, $A \asymp B$. The notation $\preceq_e$ might also appear for parameters in several results of this paper. In this context it means that there exists a constant depending only on $e$ such that the implied result holds.

A metric space is called geodesically complete if all its geodesic segments can be extended into geodesic lines, therefore when the space is also Gromov hyperbolic and Busemann space, with respect to $a\in\partial X$, any point is included in a vertical geodesic line (not necessarily unique). 

We define the relative distance between two points $x_1$ and $x_2$ of $X$ as:
\begin{equation}
d_r\left(x_1,x_2\right)=d\left(x_1,x_2\right)-\Delta h\left(x_1,x_2\right).\nonumber
\end{equation}
It can be understood as the distance along a level-set of the Busemann function. Let us recall Lemma 4.7 of \cite{TF}.
\begin{lemma}\label{LemmaBackward}\ 
\\Let $X$ be a proper, $\delta$-hyperbolic, Busemann space. Let $V_1$ and $V_2$ be two vertical geodesics of $H$. Let $t_1,t_2\in\mathbb{R}$ and let us denote $D:=\frac{1}{2}d_r\big(V_1(t_1),V_2(t_2)\big)$. Then for all $t\in \left[0,D\right]$
\begin{equation}
\left|d_r\big(V_1(t_1+D-t),V_2(t_1+D-t)\big)-2t\right|\leq 288\delta.
\end{equation}
\end{lemma}

\begin{cor}\label{CoroBackward}
Let $V_1$, $V_2$ be two vertical geodesics of $X$. Then there exists a height $h_{\mathrm{div}}(V_1,V_2)\in\mathbb{R}$ from which $V_1$ and $V_2$ diverge from each other:
\begin{enumerate}
\item $\forall t\geq h_{\mathrm{div}}(V_1,V_2)$,~~  $d\big(V_1(t),V_2(t)\big)\preceq_{\delta} 1$;
\item $\forall t\leq h_{\mathrm{div}}(V_1,V_2)$,~~  $\left|d\big(V_1(t),V_2(t)\big)-2\left(h_{\mathrm{div}}(V_1,V_2)-t\right)\right|\preceq_{\delta} 1$.
\end{enumerate}
\end{cor}

\begin{figure}
\begin{center}
\includegraphics[scale=1]{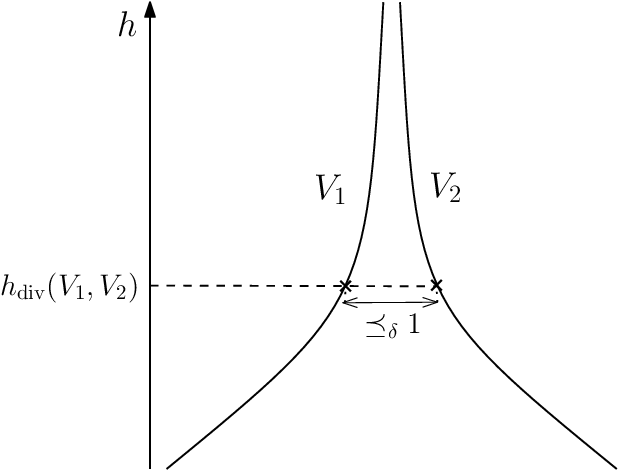} 
\end{center}
\caption{Figure of Lemma \ref{CoroBackward}.}
\label{FigCorBackward}
\end{figure}

This corollary is illustrated in Figure \ref{FigCorBackward}. We also have a more quantitative version.

\begin{lemma}[Lemma 4.3 of \cite{TF}]\label{LEM0}
Let $H$ be a $\delta$-hyperbolic and Busemann metric space, let $x$ and $y$ be two elements of $H$ such that $h(x)\leq h(y)$, and let $\alpha$ be a geodesic linking $x$ to $y$. Let us denote $z=\alpha\left(\Delta h(x,y)+\frac{1}{2}d_r(x,y)\right)$, $x_1:=V_x\left(h(y)+\frac{1}{2}d_r(x,y)\right)$ the point of $V_x$ at height $h(y)+\frac{1}{2}d_r(x,y)$ and $y_1:=V_y\left(h(y)+\frac{1}{2}d_r(x,y)\right)$ the point of $V_y$ at the same height $h(y)+\frac{1}{2}d_r(x,y)$. Then we have:
\begin{enumerate}
\item $\sup\limits_{p\in \alpha}\big(h(p)\big)\geq h(y)+\frac{1}{2}d_r(x,y)-96\delta$;
\item $d\left(z,x_1\right)\leq 144\delta$;
\item $d\left(z,y_1\right)\leq 144\delta$;
\item $d\left(x_1,y_1\right)\leq 288\delta$.
\end{enumerate}
\end{lemma}

We list here some notations we will use in later sections.

\begin{nota}\label{NotaBeginHyper}
Let $X$ be a proper, geodesically complete, $\delta$-hyperbolic, Busemann space.
\begin{enumerate}
\item Let us denote the $r$-neighbourhood of $U$ for all $U\subset X$ and for all $r\geq 0$ by
\begin{equation}
\mathcal{N}_r(U):=\lbrace x\in X\mid d(x,U)\leq r\rbrace.
\end{equation}
\item For all $x\in X$ let us denote by $V_x$ the unique vertical geodesic ray such that $V_x(0)=x$. 
\item For a subset $A\subset X$, let us denote
\begin{align}
h^-(A):=\inf\limits_{x\in A}\big(h(x)\big)\quad;\quad h^+(A):=\sup\limits_{x\in A}\big(h(x)\big).
\end{align}
\item For a subset $A\subset X$ and a height $z\in \mathbb{R}$, we denote the slice of $A$ at the height $z$ by $A_z:=A\cap h^{-1}(z)$. Therefore the horospheres of $X$ are denoted by $X_z$ for $z\in \mathbb{R}$.
\item Given a point $p\in X$ and a radius $r\in\mathbb{R}^+$, let us denote the ball of radius $r$ included in the horosphere $X_{h(p)}$ by $D_r(p):=\left\lbrace x\in X \mid h(x)=h(p)\text{ and } d(x,p)\leq r \right\rbrace=B(p,r)\cap X_{h(p)}$.
\item $\forall z\in\mathbb{R}$, $\forall U\subset X_z$, $\forall r>0$, the $r$-interior of U in $X_z$ is defined by
\begin{align*}
\mathrm{Int}_r(U):=\lbrace p\in U\mid d(p,q)\geq r,\ \forall q\in X_z\setminus U\rbrace.
\end{align*} 
\end{enumerate}
\end{nota}

Vertical geodesics of $X$ can be understood as being normal to horospheres of $X$.

\begin{defn}[Projection on horospheres]\label{DefFlowHyper}\ 
\\Let $X$ be a Gromov hyperbolic, Busemann, proper, geodesically complete metric space. Then for all $A\subset X$ and all $z\leq h^-(A)$
\begin{equation}
\pi_z(A):=\left\lbrace x\in X_z |V_x\cap A\neq \emptyset\right\rbrace.
\end{equation}
\end{defn}

\begin{figure}
\begin{center}
\includegraphics[scale=0.9]{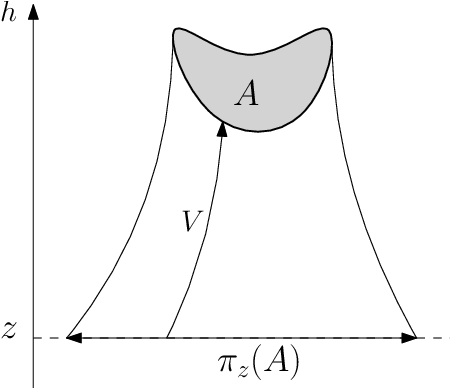} 
\end{center}
\caption{Projection of $A$ on $X_z$.}
\label{FigDefProj}
\end{figure}

The definition of this projection along the vertical flow is illustrated in Figure \ref{FigDefProj}.
The following Lemma shows that the projection of a disk on a horosphere is almost a disk, It will be used in further subsequent sections.

\begin{lemma}\label{LemmaDisqueInShadow}
Let $X$ be a Gromov hyperbolic, Busemann, proper, geodesically complete metric space. Let $z_0\in\mathbb{R}$ and $p\in X_{z_0}$. Then for $M\geq 288\delta$ we have that for all $z\leq z_0$ and for all $p_z\in\pi_z(\{p\})$
\begin{align*}
D_{2(z_0-z)-M}\big(p_z\big)\subset\pi_z\big(D_{M}(p)\big)\subset D_{2(z_0-z)+M}\big(p_z\big).
\end{align*}

\end{lemma}

\begin{figure}
\begin{center}
\includegraphics[scale=1.1]{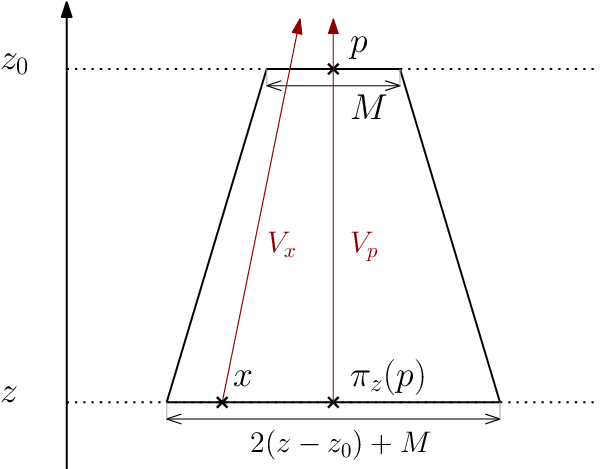} 
\end{center}
\caption{Proof of Lemma \ref{LemmaDisqueInShadow}}\label{FigProofLemmaDisqueInShadow}
\end{figure}

\begin{proof}
This Lemma is a corollary of Lemma \ref{LemmaBackward} and is illustrated in Figure \ref{FigProofLemmaDisqueInShadow}. Let $M=288\delta$ be the constant involved in Lemma \ref{LemmaBackward}.
\\Let us prove the first inclusion. Let $x\in D_{2(z_0-z)-M}\big(p_z\big)$, then $d(x,p_z)\leq 2(z_0-z)-M$. Let us denote $V_x$ a vertical geodesic containing $x$ and $V_p$ a vertical geodesic containing $p$ and $p_z$. We apply Lemma \ref{LemmaBackward} with $t_1=t_2=z$, $V_1=V_x$ and $V_2=V_p$, then $D=\frac{d(x,p_z)}{2}$. Moreover
\begin{align*}
z+D=z+\frac{d(x,p_z)}{2}\leq z +(z_0-z)-\frac{M}{2}\leq z_0.
\end{align*}
Therefore, by the Busemann convexity of $X$, the distance between vertical geodesic ray is convex and bounded, hence decreasing. Therefore
\begin{align*}
d(V_x(z_0),p)&=d(V_x(z_0),V_p(z_0))\leq d\big(V_x(z+D),V_p(z+D)\big),
\\&\leq M\quad,\text{ by Lemma \ref{LemmaBackward} used with }t=0,
\end{align*}
which means that $x\in \pi_z\big(D_M(p)\big)$.
\\Let us now prove the second inclusion, which is
\begin{equation}
\pi_z\big(D_{M}(p)\big)\subset D_{2(z_0-z)+M}\big(p_z\big).
\end{equation}
Let $x\in\pi_z\big(D_{M}(p)\big)$, then $d(V_x(z_0),V_p(z_0))\leq M$. Therefore by the triangle inequality
\begin{align*}
d(x,p_z)&=d\big(V_x(z),V_p(z)\big)\leq d\big(V_x(z),V_x(z_0)\big)+d\big(V_x(z_0),V_p(z_0)\big)+d\big(V_p(z_0),V_p(z)\big),
\\&\leq (z_0-z)+M+(z_0-z)=2(z_0-z)+M.
\end{align*}
Hence $x\in D_{2(z_0-z)+M}(p_z)$.
\end{proof}

Notations \ref{NotaBeginHyper} can be extended to horospherical products.

\begin{nota}\label{NotaBeginHoro}
Let $X$ and $Y$ be two proper, hyperbolic, geodesically complete, Busemann spaces. Then:
\begin{enumerate}
\item We denote the $r$-neighbourhood of $
U$, for all $U\subset X\bowtie Y$ and for all $r\geq 0$, by 
\begin{equation}
\mathcal{N}_r(U):=\lbrace p\in X\bowtie Y\mid d_{\bowtie}(p,U)\leq r\rbrace.
\end{equation}
\item The difference of height between two points $a,b\in X\bowtie Y$ is still denoted by $\Delta h\big(a,b\big):=\big|h(a)-h(b)\big|$.
\item We still denote, for all $z\in \mathbb{R}$ and $A\subset X\bowtie Y$, by $A_z:=A\cap h^{-1}(z)$ the "slice" of $A$ at the height $z$.
\item We still denote, for all $r\geq 0$ and $p\in X\bowtie Y$, by 
\begin{align*}
D_r(p):=\{x\in X\mid h(p)=h(x) \text{ and  } d_{\bowtie}(p,x)\leq r \}=B(p,r)\cap ( X\bowtie Y)_{h(p)},
\end{align*}
the ball of radius $r$ in the height level set containing $p$.
\end{enumerate}
\end{nota}

We recall other useful results of \cite{TF} that we will use later. First the fact that the height function is Lipschitz.

\begin{lemma}[Lemma 3.6 of \cite{TF}]\label{LemDistBigHaut}
Let $N$ be an admissible norm, and let $d_{\bowtie}$ the distance on $X\bowtie Y$ induced by $N$. Then the height function is $1$-Lipschitz with respect to the distance $d_{\bowtie}$, \textsl{i.e.},
\begin{equation}
\forall p,q \in X\bowtie Y ,\quad d_{\bowtie}(p,q)\geq \Delta h(p,q).
\end{equation}
\end{lemma}

Here is a description of the distance in Horospherical products.

\begin{thm}[Corollary 4.13 of \cite{TF}]\label{lengthGeod}
For all $p,q\in X\bowtie Y$
\begin{align*}
\left|d_{\bowtie}(p,q)-\left(d_Y\left(p^Y,q^Y\right)+d_X\left(p^X,q^X\right)-\Delta h(p,q)\right)\right|\preceq_{\bowtie}1.
\end{align*}
\end{thm}

Here is one central result of \cite{BH}, let us denote by $l(c)$ the length of a path $c$.

\begin{propo}[Proposition $1.6$, p400 of \cite{BH}]\label{LemmeBrid}
Let $X$ be a $\delta$-hyperbolic geodesic space. Let $c$ be a continuous path in X. If $[p,q]$ is a geodesic segment connecting the endpoints of $c$, then for every $x\in[p,q]$: $$d(x,\mathrm{im}(c))\leq \delta |\log_2 l(c)|+1.$$
\end{propo}

We also provide two more definitions that will be used in future sections. First a projection on level-sets of the height function.

\begin{defn}\label{DefProjHoro}
Let $z_0,z\in \mathbb{R}$ and let $U\subset (X\bowtie Y)_{z_0}$. Then we define the projection of U on $ (X\bowtie Y)_z$ by
\begin{align*}
\pi_{z}^{\bowtie}(U):=\Big\lbrace p\in (X\bowtie Y)_z\mid \exists V\text{ a vertical geodesic such that } p\in V\text{ and }V\cap U \neq \emptyset\Big\rbrace.
\end{align*}
\end{defn}

\begin{figure}[ht]
\begin{center}
\includegraphics[scale=0.6]{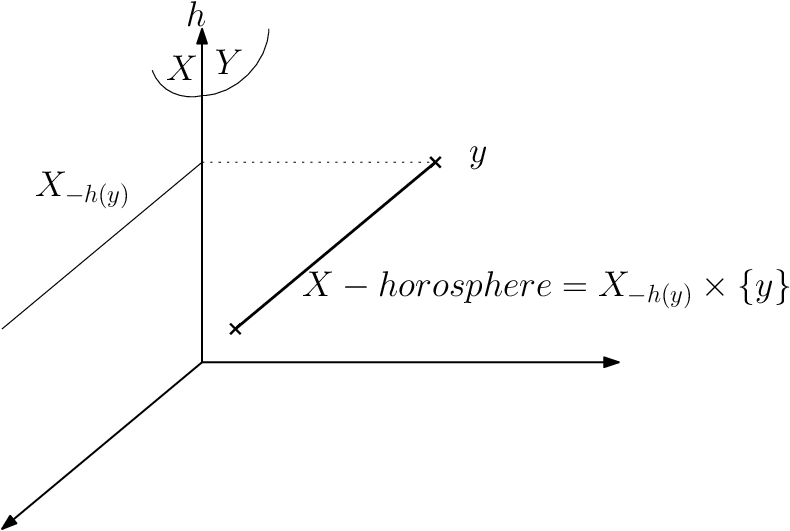} 
\end{center}
\caption{$X$-Horosphere in $X\bowtie Y$.}
\label{FigDefHorosph}
\end{figure}

Then we define $X$-horospheres and $Y$-horospheres as horospheres of hyperbolic spaces embedded in $X\bowtie Y$, illustrated in Figure \ref{FigDefHorosph}.

\begin{defn}\label{DefXYHorosphere}
The set $H\subset X\bowtie Y$ is called
\begin{enumerate}
\item an $X$-horosphere if there exists $y\in Y$ such that $H=X\bowtie\lbrace y\rbrace =X_{-h(y)}\times\{y\} $,
\item a $Y$-horosphere if there exists $x\in X$ such that $H=\lbrace x\rbrace \bowtie Y= \lbrace x\rbrace \times Y_{-h(x)}$.
\end{enumerate}  
\end{defn}

From now on, we will work in a horospherical product $X\bowtie Y$ of two proper, geodesically complete, $\delta$-hyperbolic and Busemann spaces.

\section{Metric aspects and metric tools in horospherical products}\label{ChapCoarsDifTetra}

Throughout this section we fix two constants $k\geq 1$ and $c\geq 0$. We recall the notions of quasi-isometry and quasi-geodesic.

\begin{defn}\label{DefQI}($(k,c)$-quasi-isometry)
\\Let $(E,d_E)$ and $(F,d_F)$ be two metric spaces. A map $\Phi:E\to F$ is called a $(k,c)$-quasi-isometry if and only if:
\begin{enumerate}
\item For all $x,x'\in E$, $k^{-1}d_E(x,x')-c\leq d_F\left(\Phi(x),\Phi(x')\right)\leq kd_E(x,x')+c$.
\item For all $y\in F$, there exists $x\in E$ such that $d(\Phi(x),y)\leq c$.
\end{enumerate}
A map satisfying Condition $1$ is called a quasi-isometric embedding of $E$.
\end{defn}

\begin{defn}($(k,c)$-quasigeodesic)
\\Let $E$ be a metric space. A $(k,c)$-quasigeodesic segment, respectively ray, line, of $E$ is a $(k,c)$-quasi-isometric embedding of a segment, respectively $[0;+\infty)$, $\mathbb{R}$, into $E$.
\end{defn}

In Lemma $2.1$ of \cite{GS}, Gouëzel and Shchur prove that any $(k,c)$-quasigeodesic segment is included in the $2c$-neighbourhood of a continuous $(k,4c)$-quasigeodesic segment sharing the same endpoints. Therefore, without loss of generality, we may consider that all quasi-geodesic segments are continuous.

This section gathers several geometric results on horospherical products, including the generalisation in our context of Lemmas 4.6, 3.1 and the coarse differentiation previously obtained by Eskin, Fisher and Whyte in \cite{EFW1}. Proposition \ref{LemmaQuatrePtSix}, Corollary \ref{LemmahIsQI} and Proposition \ref{LemmaTetrahedron} of this section will be especially useful in the following proofs. 

At first, a reader who is more interested in the rigidity result on horospherical product can take these propositions for granted and jump to the next sections.

When $A\asymp_{e} B$, and $e=( X\bowtie Y,d)$ is a horospherical product, we shall write $A\asymp_{\bowtie} B$ as a short-cut, and similarly $\preceq_{\bowtie}$, $\succeq_{\bowtie}$ and $M(\bowtie)$ for a constant depending only on the metric horospherical product $(X\bowtie Y,d_{\bowtie})$. 

\subsection{$\varepsilon$-monotonicity}
 
We introduce $\varepsilon$-monotone quasigeodesics, which happen to remain close to vertical geodesics. This fact plays a key role in our argument and will be proved later.

\begin{defn}($\varepsilon$-monotone quasigeodesic)
\\Let $\varepsilon\geq 0$ and let $\alpha:[0,R]\to X\bowtie Y$ be a quasigeodesic segment. Then $\alpha$ is called $\varepsilon$-monotone if and only if
\begin{equation}
\forall t_1,t_2\in[0,R],\ \ \Big(h\big(\alpha(t_1)\big)=h\big(\alpha(t_2)\big)\Big)\Rightarrow \Big(\big|t_1-t_2\big|\leq\varepsilon R \Big).
\end{equation}
\end{defn} 

Since $\alpha$ is assumed to be continuous, a $0$-monotone quasigeodesic has monotone height, $h\circ \alpha$ is either decreasing or increasing. We first show that in $X\bowtie Y$, the projections on $X$ and $Y$ of an $\varepsilon$-monotone quasigeodesic are also quasigeodesics.

\begin{thm}\label{LemmaProjAreQuasi}
Let $\varepsilon>0$, $R>\frac{1}{\varepsilon}$, and $\alpha=\left(\alpha^X,\alpha^Y\right):[0,R]\to X\bowtie Y$ be an $\varepsilon$-monotone $(k,c)$-quasigeodesic segment. Then there exists a constant $M(\bowtie,k,c)$ (depending only on $\bowtie$, $k$ and $c$) such that $\alpha^X$ and $\alpha^Y$ are $(4k,M\varepsilon R)$-quasigeodesics.
\end{thm} 

\begin{figure}
\begin{center}
\includegraphics[scale=1]{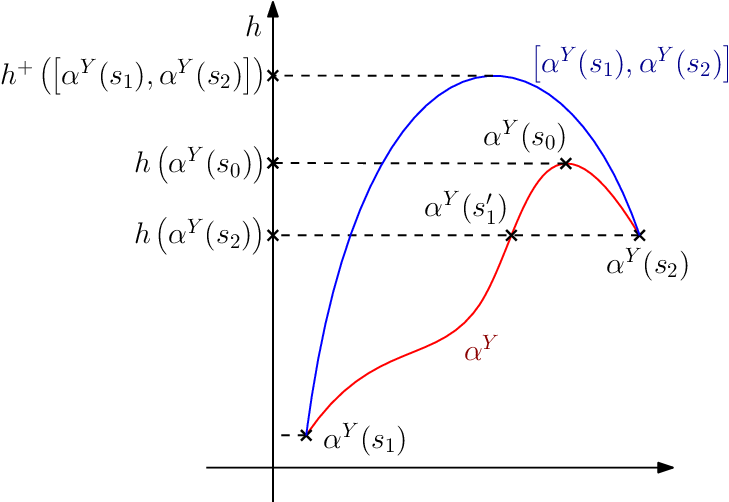} 
\end{center}
\caption{Proof of Theorem \ref{LemmaProjAreQuasi}.}
\label{FigProofProjAreQI}
\end{figure}

A portion of the proof of Theorem \ref{LemmaProjAreQuasi} is illustrated in figure \ref{FigProofProjAreQI}.

\begin{proof}
We know that $\forall p_1=\left(p_1^X,p_1^Y\right)$, $p_2=\left(p_2^X,p_2^Y\right)\in X\bowtie Y$ we have (this is the admissible assumption we made on the norm underneath the distance $d_{\bowtie}$)
\begin{equation}
d_{\bowtie}(p_1,p_2)\geq \frac{d_X\left(p_1^X,p_2^X\right)+d_Y\left(p_1^Y,p_2^Y\right)}{2}.
\end{equation} 
Therefore we have that $\alpha^X$ satisfies the upper-bound assumption of quasigeodesics
\begin{align*}
\forall s_1,s_2\in[0,R],\ d_X\left(\alpha^X(s_1),\alpha^X(s_2)\right)\leq 2 d_{\bowtie} \big(\alpha(s_1),\alpha(s_2)\big)\leq 2k |s_1-s_2|+2c.
\end{align*}
We want to find an appropriate $c'\geq c$ such that $\alpha^X$ satisfies the lower-bound condition of a $(4k,c')$-quasigeodesic. Let $c'\geq c$ and let us assume that $\alpha^X$ does not satisfy the lower-bound condition of a $(4k,c')$-quasigeodesic, we will show that this provides us with an upper-bound on $c'$. Indeed, consider $s_1,s_2\in[0,R]$ such that
\begin{align}
0\leq d_X\left(\alpha^X(s_1),\alpha^X(s_2)\right)&\leq \frac{1}{4k} |s_1-s_2|-c'.\label{InesAssumptionOns1s23}
\end{align} 
Therefore by the Lipschitz property of $h$
\begin{align}
\Delta h\left(\alpha^X(s_1),\alpha^X(s_2)\right)&\leq d_X\left(\alpha^X(s_1),\alpha^X(s_2)\right)\leq \frac{1}{4k} |s_1-s_2|-c',\nonumber
\\&\leq \frac{1}{4} d_{\bowtie}(\alpha(s_1),\alpha(s_2)\big)+\frac{c}{4}-c',\quad\text{ since }\alpha\text{ is a }(k,c)\text{- quasigeodesic}.\label{InesAssumptionOns1s2}
\end{align}
Theorem \ref{lengthGeod} gives us the existence of a constant $M(\bowtie)$ depending only on $X$, $Y$ and the underlying norm of $d_{\bowtie}$ such that
\begin{align}
&d_Y\left(\alpha^Y(s_1),\alpha^Y(s_2)\right)
\\\geq& d_{\bowtie}\big(\alpha(s_1),\alpha(s_2)\big)-d_X\left(\alpha^X(s_1),\alpha^X(s_2)\right)+\Delta h\big(\alpha(s_1),\alpha(s_2)\big)-M,\nonumber
\\\geq&  d_{\bowtie}\big(\alpha(s_1),\alpha(s_2)\big)- d_X\left(\alpha^X(s_1),\alpha^X(s_2)\right)-M\nonumber,
\\\geq&  d_{\bowtie}\big(\alpha(s_1),\alpha(s_2)\big)-\frac{1}{4k} |s_1-s_2|+c'-M,\quad\text{by assumption (\ref{InesAssumptionOns1s23}),}\nonumber
\\\geq&  d_{\bowtie}\big(\alpha(s_1),\alpha(s_2)\big)-\frac{1}{4}d_{\bowtie}\big(\alpha(s_1),\alpha(s_2)\big)-\frac{c}{4k}+c'-M,
\quad\text{since }\alpha\text{ is a }(k,c)\text{-quasigeodesic,}\nonumber
\\\geq&  \frac{1}{2}d_{\bowtie}\big(\alpha(s_1),\alpha(s_2)\big)-\frac{c}{4}+c'-M,\quad\text{since } k\geq 1.\label{IneqWithMBowtAndC}
\end{align} 
Without loss of generality, we may assume that $\max\Big(h\left(\alpha^Y(s_1)\right),h\left(\alpha^Y(s_2)\right)\Big)=h\left(\alpha^Y(s_2)\right)$. Applying Lemma \ref{LEM0} on the geodesic $\left[\alpha^Y(s_1),\alpha^Y(s_2)\right]$ of $Y$ gives us
\begin{align*}
h^+\left(\left[\alpha^Y(s_1),\alpha^Y(s_2)\right]\right)\geq h\left(\alpha^Y(s_2)\right)
+\frac{1}{2}\Big(d_Y\left(\alpha^Y(s_1),\alpha^Y(s_2)\right)-\Delta h\left(\alpha^Y(s_1),\alpha^Y(s_2)\right)\Big)-M(\bowtie).
\end{align*}
However $\alpha^Y$ is a continuous path between $\alpha^Y(s_1)$ and $\alpha^Y(s_2)$, then by Proposition \ref{LemmeBrid}, there exists $s_0\in[s_1,s_2]$ such that
\begin{align}
h\left(\alpha^Y(s_0)\right)\geq & h\left(\alpha^Y(s_2)\right)
+\frac{1}{2}\Big(d_Y\left(\alpha^Y(s_1),\alpha^Y(s_2)\right)-\Delta h\left(\alpha^Y(s_1),\alpha^Y(s_2)\right)\Big)\nonumber
\\&-\delta \log_2\big(d_Y\left(\alpha^Y(s_1),\alpha^Y(s_2)\right)\big)-M(\bowtie).\nonumber
\end{align}
Therefore by inequalities (\ref{InesAssumptionOns1s2}) and (\ref{IneqWithMBowtAndC})
\begin{align}
h\left(\alpha^Y(s_0)\right)\geq & h\left(\alpha^Y(s_2)\right)+\frac{1}{4}d_{\bowtie}\left(\alpha(s_1),\alpha(s_2)\right)-\frac{1}{8} d_{\bowtie}\left(\alpha(s_1),\alpha(s_2)\right)-\frac{c}{4}+c'-\frac{c}{8}+\frac{1}{2}c'\nonumber
\\&-\delta \log_2\big(d_Y\left(\alpha^Y(s_1),\alpha^Y(s_2)\right)\big)-\frac{M(\bowtie)}{2},\nonumber
\\\geq& h\left(\alpha^Y(s_2)\right)+\frac{1}{8}d_{\bowtie}\left(\alpha(s_1),\alpha(s_2)\right)-\delta \log_2\big(d_Y\left(\alpha^Y(s_1),\alpha^Y(s_2)\right)\big)+\frac{3}{2}c'-M(\bowtie,c).\nonumber
\end{align}
However $2d_{\bowtie}\geq d_X+d_Y\geq d_Y$, hence
\begin{align}
h\left(\alpha^Y(s_0)\right)\geq h\left(\alpha^Y(s_2)\right)+\frac{1}{8}d_{\bowtie}\left(\alpha(s_1),\alpha(s_2)\right)-\delta \log_2\big(d_{\bowtie}\left(\alpha(s_1),\alpha(s_2)\right)\big)+\frac{3}{2}c'-M(\bowtie,c).\label{Ineq1LemmeMonotVertCloseGeod}
\end{align}
Furthermore, there exists $r_0\in\mathbb{R}$ depending only on $\delta$ such that $\forall r\geq r_0$, $\frac{1}{8}r - \delta \log_2(r)> \frac{1}{10}r $ holds. Therefore, one of the two following statements holds:
\begin{enumerate}
\item[(a)] $d_{\bowtie}\left(\alpha(s_1),\alpha(s_2)\right)<r_0$.
\item[(b)] $\frac{1}{8}d_{\bowtie}\left(\alpha(s_1),\alpha(s_2)\right) - \delta \log_2\left(d_{\bowtie}\left(\alpha(s_1),\alpha(s_2)\right)\right)\geq \frac{1}{10}d_{\bowtie}\left(\alpha(s_1),\alpha(s_2)\right)$.
\end{enumerate}
We will deal with the first case $(a)$ at the end of the proof. Let us assume that $d_{\bowtie}\left(\alpha(s_1),\alpha(s_2)\right)\geq r_0$ hence $(b)$, then by inequality (\ref{Ineq1LemmeMonotVertCloseGeod})
\begin{align}
h\left(\alpha^Y(s_0)\right)\geq h\left(\alpha^Y(s_2)\right)+\frac{1}{10}d_{\bowtie}\left(\alpha(s_1),\alpha(s_2)\right)+\frac{3}{2}c'-M(\bowtie,c).\label{Ineq3LemmeMonotVertCloseGeod}
\end{align}
Then either $d_{\bowtie}\big(\alpha(s_1),\alpha(s_2)\big)\leq M(\bowtie,c)$ (up to multiplying by $10$ the constant $M$), or $h(\alpha^Y(s_0))\geq h(\alpha^Y(s_2))$. 
In the case $d_{\bowtie}\big(\alpha(s_1),\alpha(s_2)\big)\leq M(\bowtie,c)$, then $|s_1-s_2|\preceq_{k,c,\bowtie}1$ since $\alpha$ is a quasigeodesic, and therefore $c' \preceq_{k,c,\bowtie}1$ following assumption (\ref{InesAssumptionOns1s23}). In the other case we have $h(\alpha^Y(s_0))\geq h(\alpha^Y(s_2))$, therefore there exists $s_1'\in[s_1,s_0]$ such that $h\left(\alpha^Y(s_1')\right)=h\left(\alpha^Y(s_2)\right)$, since $\alpha$ is continuous. Hence
\begin{align}
d_{\bowtie}\big(\alpha(s_1'),\alpha(s_2)\big)\geq&\frac{1}{k}|s_1'-s_2|-c\geq\frac{1}{k}\big(|s_1'-s_0|+|s_0-s_2|\big)-M(c),\quad\text{since }\alpha\text{ is a quasigeodesic},\nonumber
\\\geq&\frac{1}{k^2}\Big(d_{\bowtie}\big(\alpha(s_1'),\alpha(s_0)\big)+d_{\bowtie}\big(\alpha(s_0),\alpha(s_2)\big)\Big)-M(k,c),\quad\text{since }\alpha\text{ is a quasigeodesic},\nonumber
\\\geq&\frac{1}{k^2}\Big(\Delta h\big(\alpha(s_1'),\alpha(s_0)\big)+\Delta h\big(\alpha(s_0),\alpha(s_2)\big)\Big)-M(k,c),\quad\text{by Lemma \ref{LemDistBigHaut}},\nonumber
\\\geq&\frac{2}{k^2}\Delta h\big(\alpha(s_0),\alpha(s_2)\big)-M(k,c),\quad\text{since }h\big(\alpha(s_1')\big)=h\big(\alpha(s_2)\big),\nonumber
\\\geq&\frac{1}{5k^2}d_{\bowtie}\left(\alpha(s_1),\alpha(s_2)\right)+\frac{3}{k^2}c'-M(k,c,\bowtie),\quad\text{by (\ref{Ineq3LemmeMonotVertCloseGeod})}.\label{Ineq2LemmeMonotVertCloseGeod}
\end{align} 
Moreover assumption (\ref{InesAssumptionOns1s23}) implies $|s_1-s_2|\geq 4kc'$. Then
\begin{align*}
d_{\bowtie}\big(\alpha(s_1),\alpha(s_2)\big)&\geq\frac{1}{k}|s_1-s_2|-c\geq 4c'-c.
\end{align*}
Combined with inequality (\ref{Ineq2LemmeMonotVertCloseGeod}) it gives us
\begin{align*}
d_{\bowtie}\big(\alpha(s_1'),\alpha(s_2)\big)\geq\frac{19}{5k^2}c'-M(k,c,\bowtie).
\end{align*}
Since $\alpha$ is $\varepsilon$-monotone and because $h\left(\alpha^Y(s_1')\right)=h\left(\alpha^Y(s_2)\right)$, we have
\begin{align*}
\varepsilon R\geq& |s_1'-s_2| \geq d_{\bowtie}\big(\alpha(s_1'),\alpha(s_2)\big)\geq\frac{19}{5k^2}c'-M(k,c,\bowtie).
\end{align*}
Hence
\begin{align*}
c'\leq M(k)\varepsilon R+M(k,c,\bowtie).
\end{align*}
We proved that if \(\alpha^X\) does not satisfy the lower bound inequality for being a \((4k, c')\)-quasigeodesic, then \(c' \leq M(k) \varepsilon R + M(k, c, \bowtie)\). 
Thus, when \(\varepsilon R \geq 1\), there exists a constant \(M(k, c, \bowtie)\) such that \(\alpha^X\) is a \((4k, M \varepsilon R)\)-quasigeodesic in both subcases of case \((b)\) under consideration. Similarly we show that $\alpha^Y$ is a $(4k,M\varepsilon R)$-quasigeodesic segment of $Y$. 
\\\\For case $(a)$, let us assume that each couple of times $(s_1,s_2)\in [0,R]^2$ that contradicts the lower-bound hypothesis of a $(4k,M\varepsilon R)$-quasigeodesic verifies that $d_{\bowtie}(\alpha(s_1),\alpha(s_2))<r_0$. Then $\alpha$ is a $(4k,r_0)$-quasigeodesic, with $r_0$ depending only on $\delta$. Therefore $\alpha$ is in both cases a  $(4k,M\varepsilon R)$-quasigeodesic, with $M$ depending only on $k,c$ and $X\bowtie Y$.
\end{proof}

In the sequel we denote by $d_{\mathrm{Hff}}$ the Hausdorff distance induced by  $d_\bowtie$. In the the proof of Lemma \ref{LemmaQuatrePtSix} we use a quantitative version of the quasigeodesic rigidity in a Gromov hyperbolic space, provided by the main theorem of \cite{GS}. 

\begin{thm}\label{TheoremQuantitativeEstimatQG}(\cite{GS})
\\Consider a $(k,C)$-quasigeodesic segment $\alpha$ in a $\delta$-hyperbolic space $X$, and $\gamma$ a geodesic segment between its endpoints. Then the Hausdorff distance $d_{\mathrm{Hff}}(\alpha,\gamma)$ between $\alpha$ and $\gamma$ satisfies 
\begin{align*}
d_{\mathrm{Hff}}(\alpha,\gamma)\leq 92k^2 (C+\delta).
\end{align*}
\end{thm}
This quantitative version allows us to have a linear control with respect to $C$ on the Hausdorff distance, which is mandatory in our cases since $C\asymp \varepsilon R$. Combining this rigidity with the fact that projections $\alpha^X$ and $\alpha^Y$ are also $\varepsilon$-monotone provides us with the existence of vertical geodesic segments close to $\alpha$.  

\begin{propo}\label{LemmaQuatrePtSix}
Let $\varepsilon>0$, $R>\frac{1}{\varepsilon}$, and $\alpha:[0,R]\to X\bowtie Y$ be an $\varepsilon$-monotone $(k,c)$-quasigeodesic segment. Then there exists a vertical geodesic segment $V:[0,R]\to X\bowtie Y$ such that
\begin{equation}
d_{\mathrm{Hff}}\big(\mathrm{im}(\alpha),\mathrm{im}(V)\big)\preceq_{k,c,\delta}\varepsilon R.
\end{equation}

This proposition corresponds to Lemma 4.6 in \cite{EFW1}.

\end{propo}

\begin{figure}
\begin{center}
\includegraphics[scale=1.4]{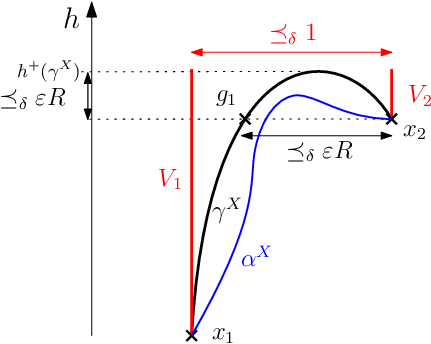} 
\end{center}
\caption{Proof of Proposition \ref{LemmaQuatrePtSix}}\label{FigProofLemmaQuatrePtSix}
\end{figure}

Figure \ref{FigProofLemmaQuatrePtSix} is an illustration of the proof.

\begin{proof}
By Theorem \ref{LemmaProjAreQuasi}, $\alpha^X$ is a $(4k,M\varepsilon R)$-quasi-geodesic in $X$ which is $\delta$-hyperbolic, hence by Theorem \ref{TheoremQuantitativeEstimatQG} there exists a geodesic $\gamma^X$ with the same endpoints as $\alpha^X$ such that
\begin{align*}
d_{\mathrm{Hff}}\Big(\mathrm{im}\left(\alpha^X\right),\mathrm{im}\left(\gamma^X\right)\Big)\preceq_{k,c,\delta}\varepsilon R.
\end{align*} 
Let us denote $x_1:=\alpha^X(0)$ and $x_2:=\alpha^X(R)$. The quasigeodesic $\alpha^X$ is also $\varepsilon$-monotone. Furthermore Proposition $2.2$ page 19 of \cite{Papa1} gives us that $\gamma^X$, which links $x_1$ to $x_2$, is included in the $24\delta$-neighbourhood of two vertical geodesic rays $V_1$ and $V_2$ such that $V_1(0)=x_1$ and $V_2(0)=x_2$. Let us denote $\tau:=h^+\left(\gamma^X\right)$, and let us recall that $\forall t_1,t_2\in\mathbb{R}^+$ and for $i\in\lbrace 1,2\rbrace$ we have $\Delta h \left(V_i(t_1),V_i(t_2)\right)=|t_1-t_2|$. Let us also denote by slight abuse $\gamma^X:=\mathrm{im}\left(\gamma^X\right)$, $\alpha^X:=\mathrm{im}\left(\alpha^X\right)$, $V_1:=\mathrm{im}\left(V_{1|[0,\tau-h(x_1)]}\right)$ and $V_2:=\mathrm{im}\left(V_{2|[0,\tau-h(x_2)]}\right)$. Since $ \tau=h^+\left(\gamma^X\right)=h^+\left(V_1\right)=h^+\left(V_2\right)$ we have
\begin{align*}
d_{\mathrm{Hff}}\Big(\gamma^X,V_1\cup V_2\Big)\preceq_{\delta} 1.
\end{align*}
Hence by the triangle inequality
\begin{align}
d_{\mathrm{Hff}}\Big(\alpha^X,V_1\cup V_2\Big)\preceq_{k,c,\delta} \varepsilon R.\label{IneqtwoGeodVertIntoOne}
\end{align}
Without loss of generality we can assume that $h(x_1)\leq h(x_2)$.
Furthermore $\gamma^X$ is continuous, therefore there exists a point of $\gamma^X$ close to both vertical geodesics (less than $24\delta$ apart). Furthermore $X$ is Busemann convex, hence the distance between the two vertical geodesics is decreasing. Therefore $d_X\Big(V_1\big(\tau-h(x_1)\big),V_2\big(\tau-h(x_2)\big)\Big)\preceq_{\delta}1$. We will use the $\varepsilon$-monotonicity of $\alpha^X$ to prove that $\tau\approx h(x_2)$. Let us denote by $x_1'$ a point of $\alpha^X$ such that $h(x_1')=h(x_2)$ and such that $d_X(x_1',\mathrm{V_1})\preceq_{k,c,\delta}\varepsilon R$. Since $\alpha^X$ is $\varepsilon$-monotone and a $(4k,M\varepsilon R)$-quasigeodesic we have that $d_X(x_1',x_2)\preceq_{k,c}\varepsilon R$, hence using the triangle inequality we have
\begin{align}
d_X\Big(V_1\big(h(x_2)-h(x_1)\big),x_2\Big)&\leq d_X\Big(V_1\big(h(x_2)-h(x_1)\big),x_1'\Big)+d_X\big(x_1',x_2\big),\nonumber
\\&\preceq_{k,c,\delta} \varepsilon R .\label{IneqMonotHoroCloseVertHoro}
\end{align} 
Let $g_1\in\mathrm{im}\left(\gamma^X\right)$ be the closest point to $x_1$ at height $h(x_2)$. Then we have:
\begin{enumerate}
\item$d_X\Big(g_1,V_1\big(h(x_2)-h(x_1)\big)\Big)\preceq_{\delta} 1$;
\item$d_X(g_1,x_2)\geq 2 \left(h^+\left(\gamma^X\right)-h(x_2)\right)$.
\end{enumerate}
We recall that $\tau=h^+\left(\gamma^X\right)$, then $d_X(g_1,x_2)\geq 2\tau-2h(x_2)\geq 0$, hence
\begin{align}
|\tau-h(x_2)|&\leq \frac{1}{2}d_X(g_1,x_2)\leq\frac{1}{2} d_X\Big(g_1,V_1\big(h(x_2)-h(x_1)\big)\Big)+\frac{1}{2}d_X\Big(V_1\big(h(x_2)-h(x_1)\big),x_2\Big),\nonumber
\\&\preceq_{k,c,\delta} \varepsilon R,\quad\text{by definition of }g_1\text{ and inequality }(\ref{IneqMonotHoroCloseVertHoro}).\nonumber
\end{align}
Hence $V_{2|[0,\tau-h(x_2)]}$ is a vertical geodesic segment of length $\preceq_{k,c,\delta} \varepsilon R$. Furthermore, $d_X(V_{1}(\tau-h(x_1)),V_{2}(\tau-h(x_2)))\leq_{\delta}$. Therefore by the triangle inequality, any point of $V_{2|[0,\tau-h(x_2)]}$ is (up to a multiplicative constant) $ \varepsilon R$-close  to $V_{1}(\tau-h(x_1))$. Therefore $d_{\mathrm{Hff}}(V_1\cup V_2,V_1)\preceq_{k,c,\delta} \varepsilon R$. Therefore, by the triangle inequality we can improve inequality (\ref{IneqtwoGeodVertIntoOne}) as follows
\begin{align*}
d_{\mathrm{Hff}}\Big(\alpha^X,V_1\Big)&\leq d_{\mathrm{Hff}}\Big(\alpha^X,V_1\cup V_2\Big)+d_{\mathrm{Hff}}\Big(V_1\cup V_2,V_1\Big),
\\&\preceq_{k,c,\delta} \varepsilon R,\quad\text{ by inequality (\ref{IneqtwoGeodVertIntoOne})}.
\end{align*}
We deduce similarly that $\alpha^Y$ is included in the $M\varepsilon R$-neighbourhood of a vertical geodesic segment $V_2'$. Therefore, $\alpha$ is included in the $M\varepsilon R$-neighbourhood of the vertical geodesic segment $(V_1,V_2')$.
\end{proof}

As a corollary, we show that the height function along an $\varepsilon$-monotone quasigeodesic is a quasi-isometry embedding of a segment into $\mathbb{R}$.

\begin{cor}\label{LemmahIsQI}
Let $\alpha:[0,R]\mapsto X\bowtie Y$ be an $\varepsilon$-monotone $(k,c)$-quasigeodesic segment. Then there exists a constant $M(k,c,\delta)$ such that the height function verifies $\forall t_1,t_2\in[0,R]$
\begin{align}
\frac{1}{k}|t_1-t_2|-M\varepsilon R\leq\Delta h\big(\alpha(t_1),\alpha(t_2)\big)\leq k|t_1-t_2|+M\varepsilon R.
\end{align}
\end{cor}

\begin{proof}
Let $t_1,t_2\in [0,R]$. The quasigeodesic upper-bound inequality is straightforward since $h$ is $1$-Lipschitz and $\alpha$ is a $(k,c)$-quasigeodesic.
\begin{align*}
\Delta h(\alpha(t_1),\alpha(t_2))\leq d_{\bowtie} (\alpha(t_1),\alpha(t_2))\leq k|t_1-t_2|+c.
\end{align*}
To achieve the lower-bound inequality we use Proposition \ref{LemmaQuatrePtSix}, hence there exists a vertical geodesic segment $V:[0,R]\to X\bowtie Y$ and a constant $M(k,c,\delta)$ such that
\begin{equation}
d_{\mathrm{Hff}}\big(\mathrm{im}(\alpha),\mathrm{im}(V)\big)\leq M\varepsilon R.
\end{equation}
For $i\in\{1,2\}$, let $s_i\in[0,R]$ be such that $d_{\bowtie}(\alpha(t_i),V(s_i))\leq M\varepsilon R$. Then by the triangle inequality
\begin{align*}
\Delta h(\alpha(t_1),\alpha(t_2))&\geq \Delta h(V(s_1),V(s_2))-2M\varepsilon R,
\\&=|s_1-s_2|-2M\varepsilon R,\quad\text{since V is vertical}.
\end{align*}
However we can achieve the lower-bound inequality on $|s_1-s_2|$
\begin{align*}
|s_1-s_2|= d_{\bowtie}((V(s_1),V(s_2))&\geq d_{\bowtie}(\alpha(t_1),\alpha(t_2))-2M\varepsilon R,\quad\text{by the triangle inequality,}
\\&\geq \frac{1}{k}|t_1-t_2|-c-2M\varepsilon R,\quad\text{since }\alpha\text{ is a quasigeodesic}.
\end{align*}
Which provides us with
\begin{align*}
\Delta h(\alpha(t_1),\alpha(t_2))&\geq |s_1-s_2|-2M\varepsilon R \geq  \frac{1}{k}|t_1-t_2|-5M\varepsilon R.
\end{align*}
\end{proof}

\subsection{Coarse differentiation of a quasigeodesic segment}\label{SecCoarseDiff}

The coarse differentiation of a quasigeodesic $\alpha$ consists in finding a scale $r>0$ such that a subdivision by pieces of length $r$ of $\alpha$ contains almost only $\varepsilon$-monotone components (which are therefore close to vertical geodesic segments).

Proposition \ref{LemmaExistGoodScale} provides us with the existence of such an appropriate scale .

\begin{lemma}\label{LemmaNonEpsilonMonoThenBigDeltaH}
Let $k\geq 1$, $c\geq 0$ and $\varepsilon> 0$. There exists $M(k,c,\bowtie,\varepsilon)$ such that for all $r\geq M$, $N\geq M$ and for all non $\varepsilon$-monotone, $(k,c)$-quasigeodesic segment $\alpha:[0,r]\to X\bowtie Y$ we have
\begin{equation}
\sum\limits_{j=0}^{N-1}\Delta h\left(\alpha\left(\dfrac{jr}{N}\right),\alpha\left(\dfrac{(j+1)r}{N}\right)\right)-\Delta h \big(\alpha(0),\alpha(r)\big)\succeq_{k,c\bowtie}\varepsilon r.
\end{equation}
\end{lemma}

This proposition corresponds to Lemma 4.7 in \cite{EFW1}.

\begin{figure}
\begin{center}
\includegraphics[scale=1.1]{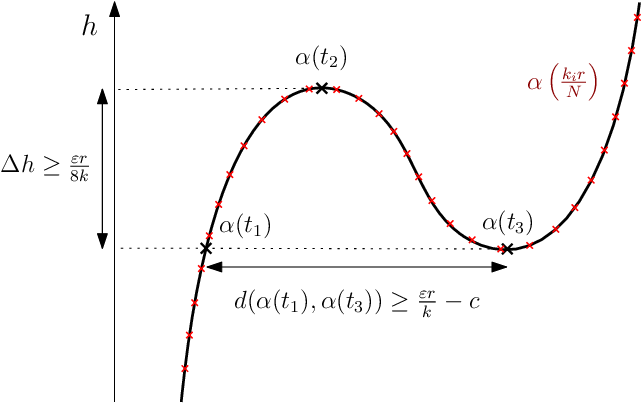} 
\end{center}
\caption{Subdivision of a quasi-geodesic.}
\end{figure}

\begin{proof}
Since $\alpha$ is non $\varepsilon$-monotone, there exist $t_1,t_3\in[0,r]$ such that
\begin{equation}
h(\alpha(t_1))=h(\alpha(t_3))\quad\text{and}\quad|t_1-t_3|>\varepsilon r.
\end{equation}
We can assume without loss of generality that $h(\alpha(0))\leq h(\alpha(t_1)) \leq h(\alpha(r))$ with $t_1<t_3$. Since $\alpha$ is a $(k,c)$-quasigeodesic we have $d_{\bowtie}\big(\alpha(t_1),\alpha(t_3)\big)\geq \dfrac{\varepsilon r}{k}-c$. By Corollary \ref{lengthGeod} of the first part of this manuscript, there exists $M(\bowtie)$ such that $d_{\bowtie}\leq d_X+d_Y +M$. Then at least one of the two following inequalities holds:
\begin{enumerate}
\item $d_{X}\left(\alpha^X(t_1),\alpha^X(t_3)\right)\succeq_{\bowtie} \frac{\varepsilon}{2k} r-M(\bowtie,c)$;
\item $d_{Y}\left(\alpha^Y(t_1),\alpha^Y(t_3)\right)\succeq_{\bowtie} \frac{\varepsilon}{2k} r-M(\bowtie,c)$.
\end{enumerate}
Let us assume that the first inequality is true. By Lemma \ref{LEM0} applied to the geodesic segment $\left[\alpha^X(t_1),\alpha^X(t_3)\right]$ we have
\begin{align*}
h^+\left(\left[\alpha^X(t_1),\alpha^X(t_3)\right]\right)\geq& \max\left(h\left(\alpha^X(t_1)\right),h\left(\alpha^X(t_3)\right)\right)
\\&+\frac{1}{2}\Big(d_X\left(\alpha^X(t_1),\alpha^X(t_3)\right)-\Delta h \left(\alpha^X(t_1),\alpha^X(t_3)\right)\Big)-96\delta,
\\=&h\left(\alpha^X(t_1)\right)+\frac{1}{2}d_X\left(\alpha^X(t_1),\alpha^X(t_3)\right)-96\delta.
\end{align*} 
Hence, there exists \(t_2 \in [t_1, t_3]\) such that the assumed inequality provides us with
\begin{align*}
\Delta h \big(\alpha(t_1),\alpha(t_2))&\succeq_{\bowtie}d_X\left(\alpha^X(t_2),\left[\alpha^X(t_1),\alpha^X(t_3)\right]\right)+\frac{\varepsilon r}{k}-M(\bowtie,c),
\\ &\succeq_{\bowtie} \frac{\varepsilon r}{k}-\delta\log_2\big(d_{\bowtie}(\alpha^X(t_1),\alpha^X(t_3))\big)-M(\bowtie,c),\quad\text{by Proposition }\ref{LemmeBrid},
\\&\succeq_{\bowtie}\frac{\varepsilon r}{k}-\delta\log_2(r)-M(\bowtie,c).
\end{align*}
Similarly, assuming the second inequality provides us with the same lower-bound on $\Delta h \big(\alpha(t_1),\alpha(t_2))$. Furthermore there exists $M(\varepsilon,\bowtie,c)$ such that for $r\geq M$ we have $\frac{1}{2}\varepsilon r\geq\delta\log_2(r)+M(\varepsilon,\bowtie,c)$, hence 
 \begin{align}
\Delta h \big(\alpha(t_1),\alpha(t_2))\succeq_{\bowtie} \frac{\varepsilon r}{2k}.\label{IneqDiffHautCoarsDiff}
\end{align}
\\Furthermore $\forall i\in\lbrace 1,2,3\rbrace$ there exists $n_i\in\lbrace 0,...,N-1\rbrace$ such that
\begin{equation}
\dfrac{n_i r}{N}\leq t_i \leq\dfrac{(n_i+1) r}{N}. \nonumber
\end{equation}
Computing the sum of the successive differences of heights provides us with
\begin{align*}
&\sum\limits_{j=0}^{N-1}\Delta h\left(\alpha\left(\dfrac{jr}{N}\right),\alpha\left(\dfrac{(j+1)r}{N}\right)\right)
\\&\geq \Delta h\left(\alpha\left(0\right),\alpha\left(\dfrac{n_1 r}{N}\right)\right)+\Delta h\left(\alpha\left(\dfrac{n_1r}{N}\right),\alpha\left(\dfrac{n_2r}{N}\right)\right)+\Delta h\left(\alpha\left(\dfrac{n_2r}{N}\right),\alpha\left(\dfrac{n_3r}{N}\right)\right)
\\&\quad+\Delta h\left(\alpha\left(\dfrac{n_3r}{N}\right),\alpha\left(r\right)\right),
\\&\geq  \Delta h\left(\alpha\left(0\right),\alpha\left(t_1\right)\right)+\Delta h\left(\alpha\left(t_1\right),\alpha\left(t_2\right)\right)+\Delta h\left(\alpha\left(t_2\right),\alpha\left(t_3\right)\right)+\Delta h\left(\alpha\left(t_3\right),\alpha\left(r\right)\right)
\\&\quad-6\left(\dfrac{kr}{N}+c\right),\quad\text{because h is Lipschitz, }\alpha\text{ is a quasigeodesic and by the triangle inequality},
\\&\geq\Delta h\big(\alpha(0),\alpha(r)\big)+2\Delta h\big(\alpha(t_1),\alpha(t_2)\big)-6\left(\dfrac{kr}{N}+c\right),\quad\text{since }h(\alpha(t_1))=h(\alpha(t_3)).
\end{align*}
Using inequality (\ref{IneqDiffHautCoarsDiff}) we have
\begin{align*}
\sum\limits_{j=0}^{N-1}\Delta h\left(\alpha\left(\dfrac{jr}{N}\right),\alpha\left(\dfrac{(j+1)r}{N}\right)\right)-\Delta h\big(\alpha(0),\alpha(r)\big)&\succeq_{\bowtie} \frac{\varepsilon r}{2k}-\frac{6kr}{N}-6c,
\\&\succeq_{k,c,\bowtie} \varepsilon r,\quad\text{since we assumed }N\geq M(k,c,\bowtie,\varepsilon).
\end{align*}
\end{proof}

The next lemma asserts that, at some scale, most segments of a quasigeodesic are $\varepsilon$-monotone.

\begin{propo}\label{LemmaExistGoodScale}
Let $k\geq 1$, $c\geq 0$, $\varepsilon>0$ and let $S$ be an integer. There exists $M(k,c,\bowtie,\varepsilon)$ such that for $r_0\geq M$ and $N\geq M$ the following occurs. Let us denote by $L=N^Sr_0$. Let $\alpha:[0,L]\to X\bowtie Y$ be a $(k,c)$-quasigeodesic segment. For all $s\in\lbrace 0,...,S\rbrace$ we cut $[0,L]$ into segments of length $N^s r_0$, and we denote by $A_s$ the set of these segment, that is
\begin{equation}
A_s := \left\lbrace \alpha\left(\left[kN^sr_0,(k+1)N^{s}r_0\right]\right)\vert k\in\lbrace 0,...,N^{S-s}-1\rbrace\right\rbrace, \nonumber
\end{equation} 
and let $\delta_s(\alpha)$ be the proportion of segments in $A_s$ which are not $\varepsilon$-monotone
\begin{equation}
\delta_s(\alpha):=\dfrac{\#\left\lbrace \beta\in A_s\vert\beta\text{ is not }\varepsilon\text{-monotone}\right\rbrace}{\# A_s}.
\end{equation}
Then
\begin{equation}
\sum\limits_{s=1}^{S}\delta_s(\alpha)\preceq_{k,c,\bowtie} \dfrac{1}{\varepsilon}.
\end{equation}
\end{propo}

\begin{proof}
The idea is to cut $[0,L]$ into $N$ segments of equal length, then to apply Lemma \ref{LemmaNonEpsilonMonoThenBigDeltaH} to the elements of this decomposition which are not $\varepsilon$-monotone. Afterwards we decompose every piece of this decomposition into $N$ segments of equal length to which we apply Lemma \ref{LemmaNonEpsilonMonoThenBigDeltaH} if they are not $\varepsilon$-monotone. The result follows by doing this sub-decomposition $S$ times in a row.
To begin with, we need to deal with $\alpha$ being $\varepsilon$-monotone or not. Hence $\delta_S(\alpha)=0$ or $1$ and in either case thanks to Lemma \ref{LemmaNonEpsilonMonoThenBigDeltaH} we have 
\begin{equation}
\sum\limits_{j=0}^{N-1}\Delta h\big(\alpha\left(jN^{S-1}r_0\right),\alpha\left((j+1)N^{S-1}r_0\right)\big)\succeq_{k,c,\bowtie}\Delta h \big(\alpha(0),\alpha(L)\big)+\delta_S(\alpha)\varepsilon L.
\end{equation}
Then for all $j\in\{0,...,N-1\}$ such that $\alpha\big([jN^{S-1}r_0,(j+1)N^{S-1}r_0]\big)$ is not $\varepsilon$-monotone
\begin{align*}
&\sum\limits_{k=0}^{N-1}\Delta h\big(\alpha\left(kN^{S-2}r_0+jN^{S-1}r_0\right),\alpha\left((k+1)N^{S-2}r_0+jN^{S-1}r_0\right)\big)\\\succeq_{k,c,\bowtie}&\Delta h \big(\alpha(jN^{S-1}r_0),\alpha((j+1)N^{S-1}r_0)\big)+\dfrac{\varepsilon L }{N},
\end{align*}
which happens $N\delta_{S-1}(\alpha)$ times. Therefore we have that
\begin{align*}
\sum\limits_{i=0}^{N^2-1}\Delta h\big(\alpha\left(iN^{S-2}r_0\right),\alpha\left((i+1)N^{S-2}r_0\right)\big)\succeq_{k,c,\bowtie}&\Delta h \big(\alpha(0),\alpha(r)\big)+\delta_S(\alpha)\varepsilon L+N\delta_{S-1}(\alpha)\dfrac{\varepsilon L}{N},
\\\succeq_{k,c,\bowtie}&\Delta h \big(\alpha(0),\alpha(r)\big)+\big(\delta_S(\alpha)+\delta_{S-1}(\alpha)\big)\varepsilon L.
\end{align*} 
By doing this another $S-2$ times we obtain
\begin{align*}
\sum\limits_{i=0}^{N^S-1}\Delta h\big(\alpha\left(ir_0\right),\alpha\left((i+1)r_0\right)\big)\succeq_{k,c,\bowtie}&\Delta h \big(\alpha(0),\alpha(r)\big)+\varepsilon L\sum\limits_{s=1}^S \delta_s(\alpha).
\end{align*}
Furthermore we have the following estimate using the Lipschitz property of $h$
\begin{align*}
\sum\limits_{i=0}^{N^S-1}\Delta h\big(\alpha\left(ir_0\right),\alpha\left((i+1)r_0\right)\big)&\leq \sum\limits_{i=0}^{N^S-1}d_{\bowtie}\big(\alpha\left(ir_0\right),\alpha\left((i+1)r_0\right)\big),
\\&\leq N^S(kr_0+c)\leq 2kL,\quad\text{with }r_0\geq\frac{c}{k}.
\end{align*}
Hence
\begin{equation}
\sum\limits_{s=1}^S \delta_s(\alpha)\preceq_{k,c,\bowtie} \dfrac{1}{\varepsilon L}2kL\preceq_{k,c,\bowtie}\dfrac{1}{\varepsilon}.
\end{equation}
\end{proof}

%

\subsection{Height respecting tetrahedric quadrilaterals}

In this subsection we show that a coarse tetrahedric quadrilateral whose sides are vertical geodesics, has two vertices on the same $X$-horosphere, and the other two on the same $Y$-horosphere (see \ref{DefXYHorosphere} for the definition of such horospheres). We call such a configuration a \textit{vertical quadrilateral}.

\begin{defn}(Orientation)
We define the orientation function on the paths of $X\bowtie Y$ as follows. For all $T>0$ and $\gamma:[0,T]\to X\bowtie Y$ we have
\begin{equation}
\mathrm{orientation}(\gamma) = \left\{
    \begin{array}{ll}
        \uparrow &\text{ if } h\big(\gamma(0)\big)< h\big(\gamma(T)\big),\quad \text{upward},\\
         \downarrow &\text{ if } h\big(\gamma(0)\big)> h\big(\gamma(T)\big),\quad \text{downward}.
    \end{array}
\right.
\end{equation}
\end{defn}

This lemma is strongly inspired by Lemma 3.1 of \cite{EFW1}, which establishes a similar result in the context of Diestel-Leader graphs and Sol geometries.

\begin{propo}\label{LemmaTetrahedron}(Vertical quadrilateral lemma)\\
Let $a_1$, $a_2$, $b_1$, $b_2\in X\bowtie Y$. Let $D>1$ and for $i,j\in\lbrace 1,2\rbrace$, let $V_{ij}:[0,l_{ij}]\rightarrow X\bowtie Y$ be vertical geodesic segments linking the $D$-neighbourhood of $a_i$ to the $D$-neighbourhood of $b_j$, and diverging quickly from each other. More specifically, we assume for all $i,j\in\lbrace 1,2\rbrace$:
\begin{enumerate}
\item[(a)] $d(V_{ij}(0),a_i)\leq D$;
\item[(b)] $d(V_{ij}(l_{ij}),b_j)\leq D$;
\item[(c)] $d(V_{i1}(t),\mathrm{im}(V_{i2}))\geq \dfrac{t}{10}-D,\ \forall t\in[0,l_{i1}]$;
\item[(d)] $d(V_{1j}(l_{1j}-t),\mathrm{im}(V_{2j}))\geq \dfrac{t}{10}-D,\ \forall t\in[0,l_{1j}]$.
\end{enumerate}
If for all $i,j\in\{1,2\}$, $l_{ij}> 2D$ and the vertical geodesic segments $V_{ij}$ share the same orientation, then there exists a constant $M(\bowtie)$ such that one of the two following statements holds:
\begin{enumerate}
\item The four vertical geodesics $V_{ij}$ are upward oriented and $a_2$ is in the $(MD)$-neighbourhood of the $X$-horosphere containing $a_1$, and $b_2$ is in the $(MD)$-neighbourhood of the $Y$-horosphere containing~$b_1$. Otherwise stated, we have $d_Y\left(a_1^Y,a_2^Y\right)\leq MD$ and $d_X\left(b_1^X,b_2^X\right)\leq MD$.
\item The four vertical geodesics $V_{ij}$ are downward oriented and $a_2$ is in the $(MD)$-neighbourhood of the $Y$-horosphere containing $a_1$, and $b_2$ is in the $(MD)$-neighbourhood of the $X$-horosphere containing~$b_1$. Otherwise stated, we have $d_X\left(a_1^X,a_2^X\right)\leq MD$ and $d_Y\left(b_1^Y,b_2^Y\right)\leq MD$.
\end{enumerate}
\end{propo}

\begin{figure}
\begin{center}
\includegraphics[scale=1]{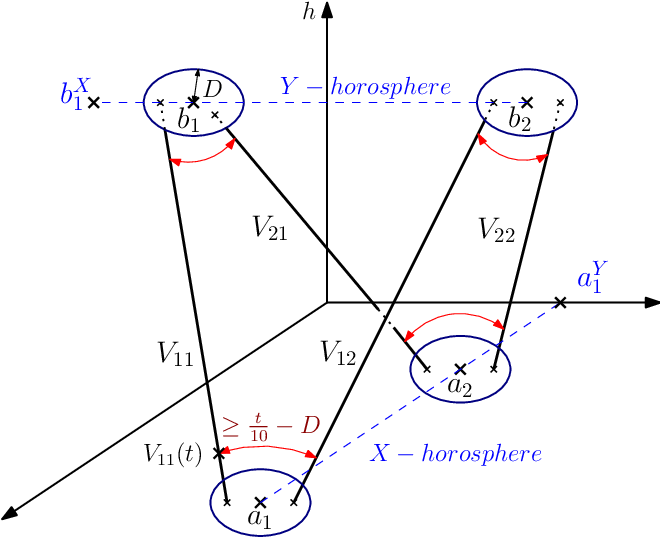} 
\end{center}
\caption{A coarse vertical quadrilateral of Proposition \ref{LemmaTetrahedron}.}\label{FigLemmaTetrahedron1}
\end{figure}

Proposition \ref{LemmaTetrahedron} is illustrated in Figure \ref{FigLemmaTetrahedron1}.

\begin{proof}\ 
\\For all $i,j\in\lbrace 1,2\rbrace$ let us denote by
\begin{equation}
a_i=\left(a_i^X,a_i^Y\right);\ b_j=\left(b_j^X,b_j^Y\right);\ V_{ij}=\left(V_{ij}^X,V_{ij}^Y\right).
\end{equation}
The hypothesis $(a)$ gives us
\begin{align}
d(V_{i1}(0),V_{i2}(0))&\leq d(V_{i1}(0),a_i)+d(a_i,V_{i2}(0))\leq 2D.\label{IneqStartPointsClose}
\end{align} 
By hypothesis $(b)$
\begin{equation}
\\d(V_{1j}(l_{1j}),V_{2j}(l_{2j}))\leq 2D.\nonumber
\end{equation}
Without loss of generality we can assume that for all $i,j\in \lbrace 1,2\rbrace$ $\mathrm{orientation}(V_{ij})=\uparrow$, which means that $h(a_i)\leq h(b_j)$. Then $\forall i,j\in\lbrace 1,2\rbrace$ and $t\in[0,l_{ij}]$ we have $h(V_{ij}(t))=t+h(V_{ij}(0))
$, hence
\begin{align*}
h(V_{ij}^X(t))&=t+h(V_{ij}(0));
\\h(V_{ij}^Y(t))&=-t-h(V_{ij}(0)).
\end{align*}
Since $X$ and $Y$ are Busemann convex spaces, $\forall i,j\in\lbrace 1,2\rbrace$ 
\begin{align*}
&t\mapsto d_Y\left(V_{i1}^Y(t),V_{i2}^Y(t)\right) \text{ is convex on }[0,\min(l_{i1},l_{i2})];
\\ &t\mapsto d_X\left(V_{1j}^X(l_{1j}-t),V_{2j}^X(l_{2j}-t)\right) \text{ is convex on }[0,\min(l_{1j},l_{2j})].
\end{align*}
Moreover, by the Busemann assumption, these maps remain convex up to a linear reparametrization. The chosen pair of vertical geodesics, whether in $X$ or $Y$, have endpoints separated by at most $2D$. Consequently, these vertical geodesics remain at a distance of at most $2D$ throughout the entire interval. Therefore
\begin{align}
&\ \forall t\in [0,\min(l_{i1},l_{i2})],\  d_Y\left(V_{i1}^Y(t),V_{i2}^Y(t)\right)\leq 2D;\label{IneqMiddle0Close}
\\&\ \forall t\in [0,\min(l_{1j},l_{2j})],\  d_X\left(V_{1j}^X(l_{1j}-t),V_{2j}^X(l_{2j}-t)\right)\leq 2D.\nonumber
\end{align}
We can assume without loss of generality that $l_{11}\leq l_{21}$ and that $l_{12}\leq l_{22}$. Then
\begin{align}
& d_X\left(V_{11}^X(0),V_{21}^X(l_{21}-l_{11})\right)\leq 2D;\label{IneqMiddle1Close}
\\& d_X\left(V_{12}^X(0),V_{22}^X(l_{22}-l_{12})\right)\leq 2D.\label{IneqMiddle2Close}
\end{align}
Let us denote $\Delta l_1=l_{21}-l_{11}$ and $\Delta l_2=l_{22}-l_{12}$, our goal is to show that these two real numbers are sufficiently close. We have $\forall i,j\in\lbrace 1,2\rbrace$
\begin{align*}
\Delta h (a_i,b_j)-2D&\leq l_{ij}\leq \Delta h (a_i,b_j)+2D.
\end{align*}
By subtracting these inequalities we get
\begin{align*}
-h(a_2)+h(a_1)-4D&\leq l_{21}-l_{11}\leq -h(a_2)+h(a_1)+4D;
\\-h(a_2)+h(a_1)-4D&\leq l_{22}-l_{12}\leq -h(a_2)+h(a_1)+4D.
\end{align*}
Then $|\Delta l_1-\Delta l_2|\leq 8D$. However
\begin{align*}
d_X\left(V_{21}^X(\Delta l_1),V_{22}^X(\Delta l_1)\right)\leq&d_X\left(V_{21}^X(\Delta l_1),V_{11}^X(0)\right)+d_X\left(V_{11}^X(0),V_{12}^X(0)\right)
\\&+d_X\left(V_{12}^X(0),V_{22}^X(\Delta l_2)\right)+d_X\left(V_{22}^X(\Delta l_2),V_{22}^X(\Delta l_1)\right).
\end{align*}
By the inequalities (\ref{IneqMiddle1Close}) and (\ref{IneqMiddle2Close}) we obtain
\begin{align}
d_X\left(V_{21}^X(\Delta l_1),V_{22}^X(\Delta l_1)\right)\leq& 2D+d_X\left(V_{11}^X(0),V_{12}^X(0)\right)+2D+|\Delta l_1-\Delta l_2|;\nonumber
\\\leq& 4D+2D+8D\leq 14D.\label{IneqMiddle3Close}
\end{align}
By using assumption $(c)$ and the characterisation of the distance on horospherical products we have
\begin{align*}
-D+\dfrac{\Delta l_1}{10}\leq& d_{\bowtie}\left(V_{21}(\Delta l_1),V_{22}(\Delta l_1)\right),
\\\leq& d_X\left(V_{21}^X(\Delta l_1),V_{22}^X(\Delta l_1)\right)+d_Y\left(V_{21}^Y(\Delta l_1),V_{22}^Y(\Delta l_1)\right)
\\&-\Delta h\left(V_{21}(\Delta l_1),V_{22}(\Delta l_1)\right)+M(\bowtie),\quad\text{ by Corollary \ref{lengthGeod},}
\\\leq& d_X\left(V_{21}^X(\Delta l_1),V_{22}^X(\Delta l_1)\right)+2D+M,\quad\text{by inequality (\ref{IneqMiddle0Close})}
\\\leq&16D+M,\quad\text{by inequality (\ref{IneqMiddle3Close})},
\end{align*}
which provides us with $\Delta l_1\leq 10(16D+M+D)=170D+10M$. We have
\begin{align*}
d_X\left(a_1^X,a_2^X\right)&\leq d_X\left(a_1^X,V_{11}^X(0)\right)+d_X\left(V_{11}^X(0),V_{21}^X(0)\right)+d_X\left(V_{21}^X(0),a_2^X\right),
\\&\leq d_X\left(V_{11}^X(0),V_{21}^X(\Delta l_1)\right)+d_X\left(V_{21}^X(\Delta l_1),V_{21}^X(0)\right)+2D,
\\&\leq 2D+170D+10M+2D\leq  174D+10M\quad,\text{ by inequality (\ref{IneqMiddle1Close})}.
\end{align*}
From this inequality we deduce that $|h(a_1)-h(a_2)|\leq 174D+10M\preceq_{\bowtie}D$. Similarly we deduce the following inequalities. 
\begin{align*}
&d_Y(b_1^Y,b_2^Y)\preceq_{\bowtie}D;
\\&|h(b_1)-h(b_2)|\preceq_{\bowtie}D.
\end{align*}
\end{proof}

Four points which satisfies the assumption of Proposition \ref{LemmaTetrahedron} are called a coarse vertical quadrilateral with nodes of scale $D$. 

\subsection{Orientation and tetrahedric quadrilaterals}\label{SecQuadLemma}

From now on we fix a $(k,c)$-quasi-isometry $\Phi:X\bowtie Y\to X\bowtie Y$. Let us consider a tetrahedric configuration consisting of two points on an \(X\)-horosphere, each connected by vertical geodesic segments to two points on a \(Y\)-horosphere, forming a total of four points and segments.
\\\\The following proposition \ref{LemmaSecondTetrahedronConfiguration} states that if two points on an $X$-horosphere are sufficiently far from each other, if two points on an $Y$-horosphere are sufficiently far from each other and if the vertical geodesic segments have $\varepsilon$-monotone images under a $(k,c)$-quasi-isometry $\Phi$, then all the images of the vertical geodesic segments by $\Phi$ share the same orientation.
\\We first show that their exists a constant $M(k,c,\bowtie)$ such that the concatenation of two consecutive $\varepsilon$-monotone quasigeodesic segments sharing the same orientation is an $M\varepsilon$-monotone quasigeodesic segment. This result will only be used in the proof of Proposition \ref{LemmaSecondTetrahedronConfiguration}.

\begin{lemma}\label{LemmaTwoUpwardIsUpward}
Let $k\geq 1$, $c\geq 0$, $D>0$, $\varepsilon>0$, $T\geq \frac{D+2c}{3\varepsilon}$ and let $\gamma:[0,T]\mapsto X\bowtie Y$ and $\gamma':[0,T]\mapsto X\bowtie Y$ be two $\varepsilon$-monotone, $(k,c)$-quasigeodesic segments such that:
\begin{enumerate}
\item $\mathrm{orientation}(\gamma)=\mathrm{orientation}(\gamma')$;
\item $d_{\bowtie}(\gamma(T),\gamma'(0))\leq D$.
\end{enumerate}
Let $\tilde{\gamma}:[0,2T]\to X\bowtie Y$ be the concatenation of $\gamma$ and $\gamma'$
\begin{equation}
\tilde{\gamma}(t) = \left\{
    \begin{array}{ll}
        \gamma(t) &\text{ if } t\in[0,T],\\
         \gamma'(t-T) &\text{ if } t\in]T,2T].
    \end{array}
\right.
\end{equation} 
Then there exists $M(k,c,D,\bowtie)$ such that $\tilde{\gamma}$ is an $M\varepsilon$-monotone, $(k,M\varepsilon T)$-quasigeodesic segment.
\end{lemma}

\begin{proof}
We can assume without loss of generality that $\gamma$ and $\gamma'$ are upward oriented, we first show that there exists $M(k,c,\bowtie)$ such that $\tilde{\gamma}$ is $M\varepsilon$-monotone. Let $t_1$, $t_2\in[0,2T]$ such that $h(\tilde{\gamma}(t_1))=h(\tilde{\gamma}(t_2))$. If both $t_1$ and $t_2$ are in $[0,T]$ or both are in $]T,2T]$, there is nothing to do since $\gamma$ and $\gamma'$ are $\varepsilon$-monotone. 
Then we can assume without loss of generality that $t_1\in[0,T]$ and $t_2\in]T,2T]$. Since $\gamma$ is upward oriented we have $h(\gamma(0))< h(\gamma(T))$, therefore, because $\gamma$ is $\varepsilon$-monotone and continuous, we have
\begin{equation}
h(\gamma(t_1))\leq h(\gamma(T))+k\varepsilon T +c\leq h(\gamma(T)) + 2k\varepsilon T,\label{UseIneqDeCont}
\end{equation}
otherwise, by continuity there exists $t_1'$ in $[0,t_1]$ such that $h(\gamma(t_1'))=h(\gamma(T))$ contradicting the $\varepsilon$-monotonicity. 
Two cases arise:
\begin{enumerate}
\item[(a)]$\Delta h\big(\gamma'(t_2-T),\gamma'(0)\big)\leq 2k\varepsilon T $;
\item[(b)]$\Delta h\big(\gamma'(t_2-T),\gamma'(0)\big)> 2k\varepsilon T $.
\end{enumerate} 
Let us consider the first case $(a)$. We know that $h(\gamma(t_1))=h(\tilde{\gamma}(t_1))=h(\tilde{\gamma}(t_2))=h(\gamma'(t_2-T))$ and that $\Delta h(\gamma(T),\gamma'(0))\leq D$, then by the triangle inequality we have
\begin{align*}
\Delta h\big(\gamma(t_1),\gamma(T)\big)= \Delta h\big(\gamma'(t_2-T),\gamma(T)\big)\leq \Delta h\big(\gamma'(t_2-T),\gamma'(0)\big)+\Delta h\big(\gamma'(0),\gamma(T)\big)\leq 2k\varepsilon T+D.
\end{align*}
According to Corollary \ref{LemmahIsQI}, $h$ is a $(k,M\varepsilon T)$-quasi-isometry along $\varepsilon$-monotone quasigeodesics. Hence
\begin{align*}
|t_1-T|&\leq k\Delta h\big(\gamma(t_1),\gamma(T)\big)+M\varepsilon T\leq (2k^2+M)\varepsilon T+kD \leq (4k^2+M)\varepsilon T,\quad\text{by assumption on }T;
\\|t_2-T|&\leq k\Delta h\big(\gamma'(t_2-T),\gamma'(0)\big)+M\varepsilon T\leq (2k^2+M)\varepsilon T .
\end{align*}
Therefore by the triangle inequality we obtain $|t_1-t_2|\leq (3k^2+M)\varepsilon (2T)$.
\\We consider now the second case $(b)$. By Corollary \ref{LemmahIsQI}, $h$ is a $(k,M\varepsilon T)$-quasi-isometry, therefore
\begin{align*}
\Delta h\big(\gamma'(t_2-T),\gamma'(0)\big)\geq \frac{1}{k}|t_2-T|-M\varepsilon T.
\end{align*}
Furthermore, $\gamma'$ is upward oriented, hence we have that $h(\gamma'(0))< h(\gamma'(t_2-T))$, otherwise, as for $\gamma$, by continuity one can construct $t_2'\in [t_2, T+T']$ contradicting the $\varepsilon$-monotonicity of $\gamma'$. 
Hence we have
\begin{align*}
h\big(\gamma'(t_2-T)\big)\geq h\big(\gamma'(0)\big) +\frac{1}{k}|t_2-T|-M\varepsilon T.
\end{align*}
In combination with inequality (\ref{UseIneqDeCont}) it provides us with
\begin{align*}
h(\gamma(t_1))&\leq  h(\gamma(T))+\varepsilon T \leq h(\gamma'(0))+D+\varepsilon T,
\\&\leq h(\gamma'(t_2-T))-\frac{1}{k}|t_2-T|+ D+(1+M)\varepsilon T.
\end{align*}
However $h(\gamma(t_1))=h(\gamma'(t_2-T))$ by definition of $t_1$ and $t_2$, therefore $0\leq-\frac{1}{k}|t_2-T|+ D+(1+M)\varepsilon T$, which gives
\begin{align}
|t_2-T|\leq (1+M)k\varepsilon T+kD\leq 3Mk\varepsilon T.\label{IneqLemmaOrient1}
\end{align} 
Hence
\begin{align*}
\Delta h(\gamma'(t_2-T),\gamma'(0))\leq d_{\bowtie}(\gamma'(t_2-T),\gamma'(0))\leq k|t_2-T|+c\leq(3Mk^2+1)\varepsilon T.
\end{align*}
Since $h(\gamma'(t_2-T))=h(\gamma(t_1))$, thanks to the triangle inequality we obtain
\begin{align}
\Delta h(\gamma(t_1),\gamma(T))&\leq \Delta h(\gamma(t_1),\gamma'(0))+\Delta h(\gamma'(0),\gamma(T)),\nonumber
\\&\leq(3Mk^2+1)\varepsilon T+D\leq(3Mk^2+2)\varepsilon T.\label{IneqLemmaOrient2}
\end{align}
Both inequalities (\ref{IneqLemmaOrient1}) and (\ref{IneqLemmaOrient2}) in combination with the fact that $h$ is a $(k,M\varepsilon T)$-quasigeodesic segment provide us with
\begin{align*}
|t_1-t_2|&= |t_1-T|+|T-t_2| \leq k (3Mk^2+2)\varepsilon T +M\varepsilon T+3Mk\varepsilon T ,
\\&\leq 9k^3 M\varepsilon T\leq \frac{9k^3 M}{2}\varepsilon(2T)\quad,\text{ since }k\geq 1,\ M\geq1.
\end{align*}
In the view of cases $(a)$ and $(b)$ we conclude that $\tilde{\gamma}$ is $ \frac{9k^3 M}{2}\varepsilon$-monotone. 
\\To prove that $\tilde{\gamma}$ is a $(k,3M\varepsilon T)$-quasigeodesic segment, we must check the upper-bound and lower bound required. Let $t_1,t_2\in[0,2T]$, as for the $\varepsilon$-monotonicity property, since $\gamma$ and $\gamma'$ are $(k,c)$-quasigeodesics, we can assume that $t_1\in[0,T]$ and $t_2\in]T,2T]$. By the triangle inequality, the upper-bound is straightforward.
\begin{align*}
d_{\bowtie}(\tilde{\gamma}(t_1),\tilde{\gamma}(t_2))&=d_{\bowtie}(\gamma(t_1),\gamma'(t_2-T)),
\\&\leq d_{\bowtie}(\gamma(t_1),\gamma(T))+d_{\bowtie}(\gamma(T),\gamma'(0))+d_{\bowtie}(\gamma'(0),\gamma'(t_2-T)),
\\&\leq k(T-t_1)+c+D+k(t_2-T)+c=k|t_2-t_1|+2c+D,
\\&\leq k|t_2-t_1|+3\varepsilon T\quad,\text{ by the assumed lower bound on }T.
\end{align*}
Last inequality holds because $\gamma$ and $\gamma'$ are $(k,c)$-quasigeodesics. To prove the lower-bound we will proceed similarly as for the $\varepsilon$-monotonicity. We have
\begin{align*}
d_{\bowtie}(\tilde{\gamma}(t_1),\tilde{\gamma}(t_2))&=d_{\bowtie}(\gamma(t_1),\gamma'(t_2-T)),
\\&\geq\Delta h\big(\gamma(t_1),\gamma'(t_2-T)\big),\quad\text{since }h\text{ is Lipschitz}.
\end{align*}
Similarly to inequality (\ref{UseIneqDeCont}) we have
\begin{align}
h(\gamma'(t_2-T))\geq h(\gamma'(0))-2k\varepsilon T.\label{UseIneqDeCont2}
\end{align}
Therefore
\begin{align*}
&\Delta h\big(\gamma(t_1),\gamma'(t_2-T)\big)\geq h(\gamma'(t_2-T))-h(\gamma(t_1))
\\=& \big(h(\gamma'(t_2-T))+\varepsilon T\big)-h(\gamma'(0))+h(\gamma'(0))-h(\gamma(T))+h(\gamma(T))-\big(h(\gamma(t_1))-\varepsilon T\big)-4k\varepsilon T,
\\=&\left|\big(h(\gamma'(t_2-T))+\varepsilon T\big)-h(\gamma'(0))\right|+\left|h(\gamma(T))-\big(h(\gamma(t_1))-\varepsilon T\big)\right|
\\&+h(\gamma'(0))-h(\gamma(T))-4k\varepsilon T\quad,\text{ by inequalities (\ref{UseIneqDeCont}) and (\ref{UseIneqDeCont2})},
\\\geq&\left|h(\gamma'(t_2-T))-h(\gamma'(0))\right|+\left|h(\gamma(T))-h(\gamma(t_1))\right|-D-8k\varepsilon T\quad,\text{ by the triangle inequality,}
\\\geq&\frac{1}{k}|t_2-T|-M\varepsilon T+\frac{1}{k}|T-t_1|-M\varepsilon T-D-8k\varepsilon T,\quad\text{because }h\text{ is a }(k,M\varepsilon T)\text{-quasigeodesic}.
\end{align*}
Hence
\begin{align*}
d_{\bowtie}(\tilde{\gamma}(t_1),\tilde{\gamma}(t_2))&\geq\Delta h\big(\gamma(t_1),\gamma'(t_2-T)\big),
\\&\geq \frac{1}{k}(t_2-t_1)-D-(2M+8k)\varepsilon T\geq \frac{1}{k}(t_2-t_1)-M'\varepsilon T,
\end{align*}
for $M'$ a constant depending on $k$, $c$, $D$ and $\bowtie$. This is the lower-bound we expected and proves that $\tilde{\gamma}$ is a $(k,M'\varepsilon T)$-quasigeodesic. 
\end{proof}

\begin{propo}\label{LemmaSecondTetrahedronConfiguration}
Let $h\in\mathbb{R}$ and let $k\geq 1$, $c\geq 0$ and $\varepsilon>0$. Let $\Phi:X\bowtie Y\to X'\bowtie Y' $ be a $(k,c)$-quasi-isometry. Let $D>1$ and $R>\frac{k2D+c}{\varepsilon}$. For $i,j\in\lbrace 1,2\rbrace$ let $a_i$, $b_j$ be four points of $X\bowtie Y$ verifying $d(a_1,a_2)>10kM\varepsilon R+2kc$ and $d(b_1,b_2)\geq  10kM\varepsilon R+2kc$, where $M$ is the constant involved in Lemma \ref{LemmaTwoUpwardIsUpward}, and let $V_{i,j}:[0,R]\to X\bowtie Y$ be four vertical geodesic segments linking the $D$-neighbourhood of $a_j$ to the $D$-neighbourhood of $b_i$, such that:
\begin{enumerate}
\item[$\bullet$] $h\big(V_{11}(0)\big)=h\big(V_{22}(0)\big)=h(a_1)=h(a_2)=h$;
\item[$\bullet$] $h\big(V_{11}(R)\big)=h\big(V_{22}(R)\big)=h(b_1)=h(b_2)=h+R$;
\item[$\bullet$] $h\big(V_{12}(0)\big)=h\big(V_{21}(0)\big)=h$;
\item[$\bullet$] $h\big(V_{12}(R)\big)=h\big(V_{21}(R)\big)=h+R$;
\item[$\bullet$] $\Phi\circ V_{i,j}$ is $\varepsilon$-monotone.
\end{enumerate}
Then the following statement holds:
\begin{align*}
&\mathrm{orientation}\big(\Phi\circ V_{11}\big)=\mathrm{orientation}\big(\Phi\circ V_{22}\big).
\end{align*} 
\end{propo} 

\begin{figure}
\begin{center}
\includegraphics[scale=0.7]{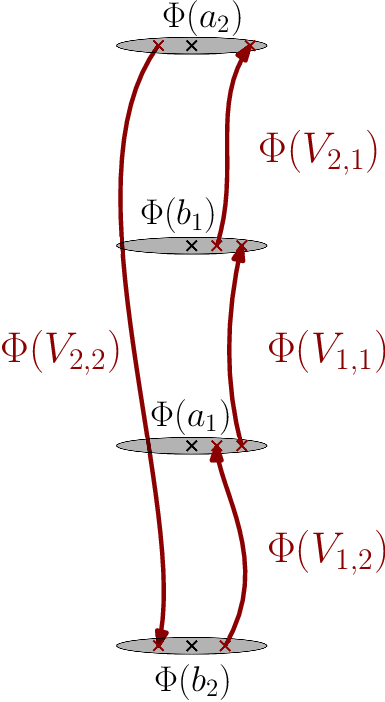} 
\end{center}
\caption{Case $(a)$ in proof of Proposition \ref{LemmaSecondTetrahedronConfiguration}.}
\label{FigSecondQuadConfig}
\end{figure}

\begin{proof}
Up to the additive constant $D$, one can consider $V_{1,1}\cup V_{2,1}\cup V_{2,2}\cup V_{1,2}$ as a coarse quadrilateral composed with $a_{i}$ and $b_{j}$ as its vertices, and with $V_{i,j}$ as its edges. To make the proof easier to follow, we shall use a vector of arrows to describe the orientations of the edges of the quadrilateral in play as follows:
\begin{align*}
\mathrm{orientation}\big(V_{1,1}, V_{2,1}, V_{2,2}, V_{1,2}\big)=(\uparrow,\downarrow,\uparrow,\downarrow).
\end{align*}
Similarly, we consider orientations of the image of $V_{1,1}\cup V_{2,1}\cup V_{2,2}\cup V_{1,2}$ by $\Phi$ as the successive orientations of the paths $\Phi\circ V_{i,j}$. We will proceed by contradiction to prove the lemma. Let us assume that $\mathrm{orientation}\big(\Phi\circ V_{1,1}\big)\neq\mathrm{orientation}\big(\Phi\circ V_{2,2}\big)$. We can assume without loss of generality that $\mathrm{orientation}\big(\Phi(V_{1,1})\big)=\uparrow$, therefore $\mathrm{orientation}\big(\Phi(V_{2,2})\big)=\downarrow$. Hence there are four possible orientations for $\Phi\big(V_{1,1}\cup V_{2,1}\cup V_{2,2}\cup V_{1,2}\big)$:
\begin{center}
\begin{table}[!ht]
\begin{tabular}{lllll}
&(a)$\quad(\uparrow,\uparrow,\downarrow,\uparrow)\quad$& (b)$\quad(\uparrow,\uparrow,\downarrow,\downarrow)\quad$& (c)$\quad(\uparrow,\downarrow,\downarrow,\uparrow)\quad$& (d)$\quad(\uparrow,\downarrow,\downarrow,\downarrow)\quad$
\end{tabular}
\end{table}
\end{center}
Let us consider the case $(a)$ (illustrated in Figure \ref{FigSecondQuadConfig}), we have  $\mathrm{orientation}\big(\Phi(V_{2,1})\big)=\uparrow$ and 
\newline$\mathrm{orientation}\big(\Phi(V_{1,2})\big)=\uparrow$. Hence we have
\begin{align*}
\mathrm{orientation}\big(\Phi(V_{1,2})\big)=\mathrm{orientation}\big(\Phi(V_{1,1})\big)=\mathrm{orientation}\big(\Phi(V_{2,1})\big).
\end{align*} 
Furthermore $\Phi$ is a $(k,c)$-quasi-isometry and both $V_{1,2}(R)$ and $V_{1,1}(0)$ are close to $a_1$, hence
\begin{align*}
d_{\bowtie'}\Big(\Phi(V_{1,2}(R)),\Phi(V_{1,1}(0))\Big)\leq k2D+c.
\end{align*}
Similarly we have
\begin{align*}
d_{\bowtie'}\Big(\Phi(V_{1,1}(R)),\Phi(V_{2,1}(0))\Big)\leq k2D+c.
\end{align*}
Then by Lemma \ref{LemmaTwoUpwardIsUpward}, there exists $M(k,c,\bowtie)$ such that the concatenation of $\Phi(V_{1,2})$, $\Phi(V_{1,1})$ and $\Phi(V_{2,1})$ is an $M\varepsilon$-monotone $(k,M\varepsilon T)$-quasigeodesic. Therefore by Proposition \ref{LemmaQuatrePtSix}, there exists a constant $M(k,c,\bowtie)$ and a vertical geodesic segment $\tilde{V}$ such that
\begin{align}
d_{\mathrm{Hff}}\big(\tilde{V},\Phi(V_{1,2})\cup\Phi(V_{1,1})\cup\Phi(V_{2,1})\big)\leq M \varepsilon R.
\end{align}
Furthermore, applying Proposition \ref{LemmaQuatrePtSix} on $\Phi(V_{2,2})$ provides us with the existence of a vertical geodesic segment $\tilde{V}'$ such that
\begin{align}
d_{\mathrm{Hff}}\big(\tilde{V}',\Phi(V_{2,2})\big)\leq M \varepsilon R.
\end{align}
Moreover $d_{\bowtie}\big(V_{2,2}(0),V_{2,1}(R)\big)\leq 2D$ (the two points are close to $a_2$) and $d_{\bowtie}\big(V_{2,2}(R),V_{1,2}(0)\big)\leq 2D$ (the two points are close to $b_2$), therefore $\tilde{V}$ and $\tilde{V}'$ are two vertical geodesics with endpoints $(k2D+c)+2M\varepsilon R$ close to $\Phi(a_2)$ and $\Phi(b_2)$. Thereby, these two vertical geodesic segments stay close to each other, we have
\begin{align*}
d_{\mathrm{Hff}}(\tilde{V},\tilde{V}')\leq (k2D+c)+2M\varepsilon R\leq 3M\varepsilon ,\quad\text{by assumption on }R.
\end{align*} 
Then, we show by the triangle inequality that $\Phi(a_1)$ is close to $\Phi(V_{2,2})$.
\begin{align}
d_{\bowtie'}\big(\Phi(a_1),\Phi(V_{2,2})\big)\leq d_{\bowtie'}\big(\Phi(a_1),\tilde{V}\big)+d_{\mathrm{Hff}}\big(\tilde{V},\tilde{V}'\big)+d_{\mathrm{Hff}}\big(\tilde{V}',\Phi(V_{2,2})\big)\leq 5M\varepsilon R.\label{IneqToContradict}
\end{align}
However, the assumption $d(a_1,a_2)>10kM\varepsilon R+2kc$ gives us that $a_1$ is sufficiently far from $V_{2,2}$
\begin{align*}
\forall t\in[0,R],\ &d_{\bowtie}\big(a_1,V_{2,2}(t)\big)\geq \Delta h\big(a_1,V_{2,2}(t)\big)=t,
\\\text{and},\ &d_{\bowtie}\big(a_1,V_{2,2}(t)\big)\geq d_{\bowtie}\big(a_1,a_2\big)-d_{\bowtie}\big(a_2,V_{2,2}(t)\big)> 10kM\varepsilon R+2kc-t.
\end{align*}
Therefore
\begin{align*}
\forall t\in[0,R],\ d_{\bowtie'}\big(\Phi(a_1),\Phi(V_{2,2}(t))\big)&\geq k^{-1}d_{\bowtie}\big(a_1,V_{2,2}(t)\big)-c,
\\&> \frac{t+10kM\varepsilon R+2kc-t}{2k}-c=5M\varepsilon R,
\end{align*}
which contradicts inequality (\ref{IneqToContradict}). Thereby, in case $(a)$, $\Phi\circ V_{1,1}$ and $\Phi\circ V_{2,2}$ share the same orientation.
\\The other three cases $(b)$, $(c)$ and $(d)$ are treated similarly. We first show that $\Phi\big(V_{1,1}\cup V_{2,1}\cup V_{2,2}\cup V_{1,2}\big)$ is in the $M\varepsilon R$-neighbourhood of two vertical geodesic segments which, depending on the case, have endpoints 
\begin{enumerate}
\item[(b)] close to $\Phi(a_1)$ and $\Phi(a_2)$;
\item[(c)] close to $\Phi(b_1)$ and $\Phi(b_2)$;
\item[(d)] close to $\Phi(a_1)$ and $\Phi(b_1)$.
\end{enumerate}
Which, depending on the case, contradicts the fact that:
\begin{enumerate}
\item[(b)] $d_{\bowtie}\big(b_1,V_{2,2}(t)\big)> 5M\varepsilon R$;
\item[(c)] $d_{\bowtie}\big(a_1,V_{2,2}(t)\big)> 5M\varepsilon R$;
\item[(d)] $d_{\bowtie}\big(b_2,V_{1,1}(t)\big)> 5M\varepsilon R$.
\end{enumerate}
\end{proof}

\section{Measure and Box-tiling}\label{ChapMeasure}

\subsection{Appropriate measure and horopointed admissible space}

In the setting of horospherical product, an important characteristics is that they are union of products of horospheres.

As such, if one wants to endow them with a measure, it makes sense that the measure should disintegrate along these horospherical product, and should be related somehow to the measures and the geometries of the initial spaces and its horospheres.

The properties we present are satisfied when our initial space are Riemannian manifolds for instance, or graphs of bounded geometry. We will also see in Section \ref{SecExample} that Heintze group are another set of spaces which satisfies them, making our requirements sound.

\begin{defn}\label{DefAdmissible2}(Admissible horopointed measured metric spaces.)
\\Let $(X,d)$ be a $\delta$-hyperbolic, Busemann, proper, geodesically complete, metric space, and let $a\in\partial X$ be a point on the Gromov boundary of $X$. A Borel measure $\mu^X$ on $X$ will be said $(X,a)$ horo-admissible if and only if $(E1)$, $(E2)$ and $(E3)$ are satisfied.
\begin{enumerate}
\item[(E1)] (There exists a direction $a\in\partial X$ such that) $\mu^X$ is disintegrable along the height function $h_a$, that is
\begin{align*}
\text{For all }z\in\mathbb{R},\text{ there exists a Borel measure }\mu^X_z\text{ on }X_z=h^{-1}(z)\text{ such that for any measurable set } A\subset X
\end{align*}
\begin{align*}
\mu^X(A)=\int\limits_{z\in\mathbb{R}} \mu^X_z\left(A_z\right)\mathrm{d}z. 
\end{align*}
\item[(E2)] Controllable geometry for the measures $\mu^X_z$ on horospheres, there exists $M_0\geq 288\delta$ such that
\begin{align*}
\forall x_1,x_2\in X,\text{ we have }\mu^X_{h(x_1)}\left(D_{M_0}(x_1)\right)\asymp_{X}\mu^X_{h(x_2)}\left(D_{M_0}(x_2)\right).
\end{align*}
\item[(E3)] There exists $m>0$ such that for all $z_0\in\mathbb{R}$, and for all measurable set $U\subset X_{z_0}$
\begin{align*}
\forall z \leq z_0,\ e^{m(z_0-z)}\mu^X_{z_0}(U)\asymp_{X}\mu^X_{z}(\pi_z(U)) .
\end{align*}
\end{enumerate} 
The space $(X,a,d,\mu^X)$ will be called a horo-pointed admissible metric measured space, or just admissible.
\end{defn}

The assumption $(E2)$, in combination with Lemma \ref{LemmaDisqueInShadow}, provides us with a uniform control on the measure of disks of any radius.

\begin{lemma}\label{LemmaMuDiskAreExp}
Let $r\geq M_0$. Then for all $x\in X$ we have
\begin{align*}
\mu_{h(x)}\left(D_r(x)\right)\asymp_{X} e^{m\frac{r}{2}}.
\end{align*}
\end{lemma}

\begin{figure}
\begin{center}
\includegraphics[scale=1.3]{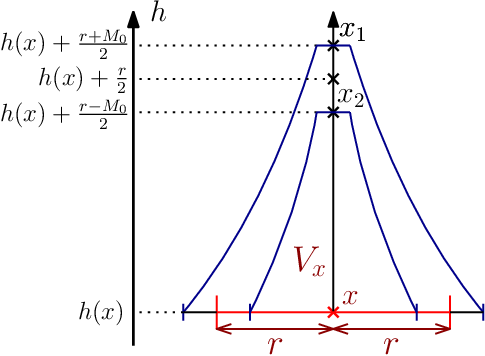} 
\end{center}
\caption{Proof of Lemma \ref{LemmaMuDiskAreExp}.}
\label{FigLemmaMuDiskIsExp}
\end{figure}

\begin{proof}
The proof is illustrated in Figure \ref{FigLemmaMuDiskIsExp}. Let $V_x$ be a vertical geodesic line containing $x$ and let $M_0\geq 288\delta$ be the constant involved in assumption $(E2)$. Let us denote $x_1$ the point of $V_x$ at the height $h(x)+\frac{r+M_0}{2}$ and let $x_2$ be the point of $V_x$ at the height $h(x)+\frac{r-M_0}{2}$. Applying Lemma \ref{LemmaDisqueInShadow} with $p=x_1$, $z_0=h(x)+\frac{r+M_0}{2}$ and $z=h(x)$ provides us with
\begin{align*}
D_r(x)=D_{2(z_0-z)-M_0}(x)\subset\pi_{h(x)}\left(D_{M_0}(x_1)\right).
\end{align*} 
Similarly, applying Lemma \ref{LemmaDisqueInShadow} with $p=x_2$, $z_0=h(x)+\frac{r-M_0}{2}$ and $z=h(x)$ provides us with $\pi_{h(x)}\left(D_{M_0}(x_2)\right)\subset D_r(x)$. Furthermore by assumption $(E3)$ then assumption $(E2)$ we have
\begin{align*}
\mu^X_{h(x)}\left(\pi_{h(x)}\left(D_{M_0}(x_1)\right)\right)\asymp_{X}e^{m(\frac{r+M_0}{2})}\mu^X_{h(x_1)}(D_{M_0}(x_1))\asymp_{X}e^{m\frac{r}{2}},
\end{align*}
since $M_0$ depends only on $X$. Similarly we have $\mu^X_{h(x)}\left(\pi_{h(x)}\left(D_{M_0}(x_2)\right)\right)\asymp_{X}e^{m\frac{r}{2}}$, therefore by the two previously obtained inclusions we have $\mu_{h(x)}\left(D_r(x)\right)\asymp_{X} e^{m\frac{r}{2}}$.
\end{proof} 

Heuristically, the next lemma asserts that the measure of the boundary of a disk is small in comparison to the measure of the disk.

\begin{lemma}\label{LemmaYetToBeWritten}
Let $M_0$ be the constant involved in assumption $(E_2)$ and let $M$ be the constant involved in Corollary \ref{CoroBackward}. Let $z_0\in\mathbb{R}$, $x_0\in X_{z_0}$ and $\mathcal{C}\subset X_{z_0}$ be a set containing $D_{M_0}(x_0)$ and contained in $D_{2M_0}(x_0)$. Then for all $z_1\leq z_0$, and for all $r\leq 2|z_1-z_0|-2M_0-M$ we have
\begin{align*}
\mu ^X_{z_1}\left(\mathrm{Int}_{r}\left(\pi^X_{z_1}(\mathcal{C})\right)\right)\asymp_{\bowtie}\mu ^X_{z_1}\left(\pi^X_{z_1}(\mathcal{C})\right).
\end{align*}
\end{lemma}

This Lemma might seem to contradict Lemma \ref{LemmaMuDiskAreExp}, however the $r$-interior of a disk of radius $R$ is very different from a disk of radius $R-r$ on horospheres, for $R$ sufficiently greater than $r$. 

\begin{proof}
Let us denote $J:= \mathrm{Int}_{r}\left(\pi^X_{z_1}(\mathcal{C})\right)$. By definition we have
\begin{align}
\pi^X_{z_1}(\mathcal{C})\setminus J:=\left\lbrace x\in\pi^X_{z_1}(\mathcal{C})\vert d_X\left(x,\pi^X_{z_1}(\mathcal{C})^c\right)< r  \right\rbrace.\label{charac1}
\end{align} 
At the height $z_1+\frac{r}{2}$, let $x_1\in \pi^X_{z_1+\frac{r}{2}}(\mathcal{C})\setminus \pi^X_{z_1+\frac{r}{2}}(J)$, then, at the height $z_1$, there exists $x_1'\in \pi^X_{z_1}(\mathcal{C})\setminus J$ such that $x_1\in V_{x_1'}$. Furthermore by the characterisation (\ref{charac1}), there exists $x_2'\in \pi^X_{z_1}(\mathcal{C})^c$ such that $d(x_1',x_2')\leq r$. Then by Lemma \ref{CoroBackward}, there exists $M(\delta)$ such that
\begin{align}
d_X\left(V_{x_2'}\left(z_1+\frac{r}{2}\right),V_{x_1'}\left(z_1+\frac{r}{2}\right)\right)=d_X\left(V_{x_2'}\left(z_1+\frac{r}{2}\right),x_1\right)\leq M\label{IneqJipep},
\end{align}  
with $V_{x_2'}\left(z_1+\frac{r}{2}\right)\in \pi^X_{z_1+\frac{r}{2}}(\mathcal{C})^c$. Therefore by the triangle inequality and Lemma \ref{LemmaDisqueInShadow}
\begin{align*}
d\left(x_1,\pi^X_{z_1+\frac{r}{2}}(x_0)\right)&\geq -d\left(x_1,V_{x_2'}\left(z_1+\frac{r}{2}\right)\right)+d\left(V_{x_2'}\left(z_1+\frac{r}{2}\right),\pi^X_{z_1+\frac{r}{2}}(x_0)\right),
\\&\geq 2|z_0-z_1|-r-M_0-M.
\end{align*}
Since last inequality holds for all $x_1\in  \pi^X_{z_1+\frac{r}{2}}(\mathcal{C})\setminus \pi^X_{z_1+\frac{r}{2}}(J)$, we have
\begin{align*}
D_{2|z_0-z_1|-r-M_0-M}(\pi^X_{z_1+\frac{r}{2}}(x_0))\subset \pi^X_{z_1+\frac{r}{2}}(J).
\end{align*}
Therefore by Lemma \ref{LemmaDisqueInShadow}
\begin{align*}
D_{2|z_0-z_1|-M_0-M}(\pi^X_{z_1}(x_0))\subset J.
\end{align*}
Moreover, $J\subset  \pi^X_{z_1}(\mathcal{C})\subset D_{2|z_0-z_1|+M_0}(\pi^X_{z_1}(x_0))$, hence by Lemma \ref{LemmaMuDiskAreExp}
\begin{align*}
\mu^X_{z_1}(J)\asymp_{X} e^{|z_0-z_1|m}\asymp_{X}\mu ^X_{z_1}\left(\pi^X_{z_1}(\mathcal{C})\right).
\end{align*}
\end{proof}

In order to achieve a rigidity result on horospherical products, we will need another measure $\lambda^X$ in the same measure class as $\mu^X$. 

\begin{defn}(measure $\lambda^X$ of $X$)\label{DefLambdaMeasure}
\\Let $X$ be an admissible horopointed space. The measure $\lambda^X$ on $X$ is defined from a set of weighted measure $\lambda_z^X$ on the level set $X_z$:
\begin{enumerate}
\item $\forall z\in\mathbb{R}$, $\lambda^X_z:=e^{mz}\mu^X_z$;
\item For all measurable set $A\subset X$, $\lambda^X(A):=\int\limits_{z\in\mathbb{R}} \lambda^X_z\left(A_z\right)\mathrm{d}z$,
\end{enumerate}
where $m$ is the constant involved in $(E3)$.
\end{defn}

For the Log model of the hyperbolic plane, this measure $\lambda^X$ turns out to be the Lebesgue measure on $\mathbb{R}^2$, and the measure $\mu^X$ is the Riemannian area. Up to a multiplicative constant, the measure $\lambda^X$ is constant along the projections. By assumption $(E3)$, the following property is immediate:

\begin{propr}\label{ProprLambdaSameOnLevels}
For all measurable set $ U\subset X$ we have 
\begin{equation}
\forall z_1,z_2 \leq h^-(U),\ \lambda^X_{z_1}(\pi_{z_1}(U))\asymp_{X}\lambda^X_{z_2}(\pi_{z_2}(U)).
\end{equation}
Otherwise stated we have the following relation between two push-forwards of the measure $\left(\pi_{ z_2}\right)_{\ast} \lambda^X_{z_2}\asymp_{X}\left(\pi_{ z_1}\right)_{\ast} \lambda^X_{z_1}$. They are push-forwards from a subset $U$ of $X$ to horospheres below $U$.
\end{propr}

Following the fact that height level sets of $X\bowtie Y$ are direct products of horospheres, we define disintegrable measures on the horospherical products from the disintegrable measures on $X$ and $Y$. We recall that $\forall z\in\mathbb{R}$
\begin{align*}
(X\bowtie Y)_z=X_z\times Y_{-z}
\end{align*}

\begin{defn}\label{DefProdMeasure}(Measure $\mu$ on $X\bowtie Y$)
\\Let $(X,\mu^X)$ and $(Y,\mu^Y)$ be two admissible spaces. Then for all measurable set $U\subset X\bowtie Y$, we define the measure $\mu$ on $X\bowtie Y$ by
\begin{align*}
\mu_{X\bowtie Y}(U):=\int\limits_{\mathbb{R}}\mu^X_z\otimes\mu^Y_{-z}\big(U_z\big) \mathrm{d}z.
\end{align*}
\end{defn} 

For all measurable set $U\subset X\bowtie Y$ we have
\begin{align*}
\mu_{X\bowtie Y}(U)=\int\limits_{\mathbb{R}}\left(~\int\limits_{y\in Y_{-z}}\mu^X_z\big(U_z^y\big) d\mu^Y_{-z}(y)\right)\mathrm{d}z,
\end{align*}
where $U_z^y:=\lbrace x\in X\mid (x,y)\in U_z \rbrace$. (This measure might be not well defined).
\begin{rem}
A couple $(X,Y)$ of horo-pointed admissible spaces is called admissible if the measure $\mu_{X\bowtie Y}$ of Definition \ref{DefProdMeasure} is well defined.
\end{rem}

From now on we fix four horo-pointed metric spaces $X$, $X'$, $Y$ and $Y'$, with $m>0$ (respectively $m'$, $n$, $n'$) the constant of assumption $(E3)$ for $X$ (respectively $X'$, $Y$, $Y'$). We will assume in Section \ref{SecThmOrient} and afterwards that $(X,Y)$ and $(X',Y')$ are two admissible couples with $m>n$ and $m'>n'$.

We define similarly a measure $\lambda_{X\bowtie Y}$ on $X\bowtie Y$.

\begin{defn}(Measure $\lambda$ on $X\bowtie Y$)
\\Let $(X,\mu^X)$ and $(Y,\mu^Y)$ be two admissible spaces. Then for all measurable subset $U\subset X\bowtie Y$
\begin{align*}
\lambda_{X\bowtie Y}(U):=\int\limits_{\mathbb{R}}\lambda^X_z\otimes\lambda^Y_{-z}\big(U_z\big) \mathrm{d}z=\int\limits_{\mathbb{R}}e^{(m-n)z}\mu^X_z\otimes\mu^Y_{-z}\big(U_z\big) \mathrm{d}z.
\end{align*}
\end{defn} 

For all measurable subset $U\subset X\bowtie Y$ we have
\begin{align*}
\lambda_{X\bowtie Y}(U)=\int\limits_{\mathbb{R}}\left(\int\limits_{y\in Y_{-z}}\lambda^X_z\big(U_z^y\big) d\lambda^Y_{-z}\right)\mathrm{d}z.
\end{align*}

From now on, we will simply denote by $\mu$ the measure $\mu_{X\bowtie Y}$ and by $\lambda$ the measure $\lambda_{X\bowtie Y}$.

\subsection{Box-tiling of X}\label{SecTill}

In this subsection we tile a proper, geodesically complete, Gromov hyperbolic and Busemann space $X$ with pieces called boxes. This is inspired from Lemma 3.4 of \cite{EFW1}, which constructs these tilings for trees and the hyperbolic space.

\begin{defn}(Box at scale $R$)
\\Let $X$ be admissible horo-pointed space. Let $M_0$ be the constant of $(E2)$, let $R>0$, let $x$ be a point of $X$ and let $\mathcal{C}(x)$ be a subset of $X_{h(x)}$ containing $D_{M_0}(x)$ and contained in $D_{2M_0}(x)$. Then, the box $\mathcal{B}(x,\mathcal{C}(x),R)$ is defined by
\begin{align*}
\mathcal{B}(x,\mathcal{C}(x),R):= \bigcup\limits_{z\in[h(x)-R,h(x)[}\pi_z\left(\mathcal{C}(x)\right).
\end{align*}
\end{defn}

We will often omit the parameter $\mathcal{C}(x)$ in the notation of a box. Later we depict an appropriate choice for these spaces $\mathcal{C}(x)$. The idea of the tiling is first to distinguish layers of thickness $R$, then to decompose each of these layers into disjoint boxes using a tiling of disjoint cells $\mathcal{C}(x)$ as the top of these boxes. In the Log model of the hyperbolic plane, when the cell $\mathcal{C}(x)$ is a segment of an horosphere, the associated box is a rectangle of $\mathbb{R}^2$. In \cite{EFW1}, Eskin, Fisher and Whyte tile the hyperbolic plane with translates of such a rectangle. However the space we consider might not be homogeneous, therefore we will tile Gromov hyperbolic spaces with boxes which are generically not the translate of one another. We recall that $\mathcal{N}_r$ refers to the $r$-neighbourhood of a subspace.

A subset of a metric space X is $k$-separated if and only if any two of its elements are at least at distance $k$. A maximal such set for the inclusion is called \textsl{maximal separating set}. We shall denote by $\mathcal{D}(X)$ such a set. The dependence of \(\mathcal{D}(X)\) on \(k\) should be indicated; however, for simplicity and by slight abuse of notation, we will omit this dependence.

One easily sees that a maximal separated set is then a $k$-covering. That is the union of the metric ball of radius k centred at the points of $\mathcal{D}(X)$ cover the whole space.


To construct a box tiling of $X$ we first fix a scale $R>0$. Let $M_0$ be the constant involved in assumption $(E2)$, then we chose a $2M_0$-maximal separating set $\mathcal{D}\left(X_{nR}\right)$ of the horospheres $X_{nR}$, with $n\in\mathbb{Z}$. Such maximal separating sets exist since $X$ is proper and so are $X_{nR}$. Let us call nuclei the points in these maximal separating sets. For every nucleus $x\in \mathcal{D}(X_{nR})$, we fix a cell $\mathcal{C}(x)$ such that $D_{M_0}(x)\subset \mathcal{C}(x)\subset D_{2M_0}(x)$. Therefore, given two different nuclei $x,x'\in\mathcal{D}\left(X_{nR}\right)$, we have $D_{M_0}(x)\cap D_{M_0}(x')=\emptyset$. We choose these cells such that they are $\mu_{nR}$ measurable and such that they tile their respective horospheres:
\begin{equation}
\forall n\in\mathbb{Z},\ \bigsqcup\limits_{x\in\mathcal{D}(X_{nR})} \mathcal{C}(x) =X_{nR}. \nonumber
\end{equation}
As an example, one can take Voronoi cells:
\begin{align*}
\mathcal{VC}(x):=\left\lbrace p\in X_{nR}|d(p,x)\leq d(p,x'), \text{ for all }x'\in \mathcal{D}\left(X_{nR}\right)\right\rbrace.
\end{align*} 
These cells might not be disjoint, but a point $p\in X_{nR}$ is contained in a finite number of Voronoi cells since $X$ is proper. Therefore, by choosing (for example thanks to an arbitrary order on $\mathcal{D}\left(X_{nR}\right)$) a unique cell containing $p$, and removing $p$ from the others, there exists a tiling of $X_{nR}$ by cells $\mathcal{C}(x)$. 
\\Now, for all $n\in\mathbb{Z}$ and for all $x\in \mathcal{D}\left(X_{nR}\right)$ we define the box $\mathcal{B}(x,R)$ at scale $R$ of nucleus $x$ by
\begin{equation}
\mathcal{B}(x,R):=\bigcup\limits_{z\in[(n-1)R;nR[}\pi_z\big(\mathcal{C}(x)\big).\nonumber
\end{equation}

\begin{figure}[ht]
\begin{center}
\includegraphics[scale=1]{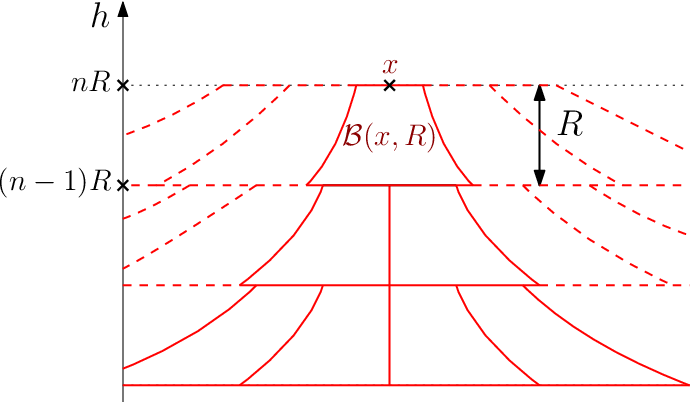} 
\end{center}
\caption{Box-tiling}
\end{figure}

In this definition, we chose $[(n-1)R;nR[$ for the boxes' heights. It is an arbitrary choice, one could prefer to use $](n-1)R;nR]$ as these heights intervals. Moreover, to construct the horospherical product of $X$ and $Y$, we will use intervals of the form $[\ldots;\ldots[$ for $X$ and $ ]\ldots;\ldots]$ for $Y$.
\\We recall that the cells $\mathcal{C}(x)$ tile the horospheres $X_{nR}$. Furthermore there exists a unique vertical geodesic ray leaving each point of $X$. Consequently we have a box tiling of $X$ at scale $R$:
\begin{equation}
X=\bigsqcup\limits_{n\in\mathbb{Z}}\bigsqcup\limits_{x\in \mathcal{D}\left(X_{nR}\right)}\mathcal{B}(x,R).
\end{equation}

The next lemma explains that any box contains and is contained in metric balls of similar scales.

\begin{lemma}\label{LemmaBoxAndBalls}
There exists a constant  $M(X)$ such that, for all $x\in X$ and $r>M$ there exist two boxes $\mathcal{B}\left(\frac{r}{2}\right)$ and $\mathcal{B}(3r)$ verifying  
\begin{align*}
\mathcal{B}\left(\frac{r}{2}\right)\subset B(x,r) \subset \mathcal{B}(3r).
\end{align*}
\end{lemma}

\begin{figure}[ht]
\begin{center}
\includegraphics[scale=1.2]{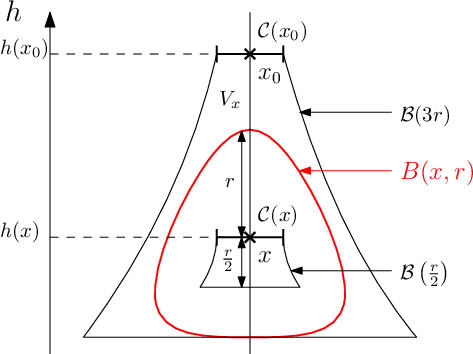} 
\end{center}
\caption{Proof of Lemma \ref{LemmaBoxAndBalls}}\label{FigProofLemmaBoxAndBalls}
\end{figure}

Proof is illustrated in Figure \ref{FigProofLemmaBoxAndBalls}. 

\begin{proof}
Let $\mathcal{C}(x)$ be a subset of $X_{h(x)}$ containing $D(x,M_0)$ and contained in $D(x,2M_0)$. Let us denote by $\mathcal{B}\left(\frac{r}{2}\right)$ the box at scale $\frac{r}{2}$ constructed from the cell $\mathcal{C}(x)$. For all $x'\in \mathcal{B}\left(\frac{r}{2}\right)$ let us denote by $x'':=V_{x'}(h(x))$ the point of $V_{x'}$ at the height $h(x)$, we have
\begin{align*}
d_X(x',x)\leq d_X\left(x',x''\right)+d_{\mathrm{X}}\left(x'',x\right)\leq \dfrac{r}{2} +2 M_0 \leq r,\quad\text{for }r\geq 4M_0,
\end{align*}
which gives us that $x'\in B(x,r)$. To prove the second inclusion, let us denote by $V_x$ the unique (since $X$ is Busemann convex) vertical geodesic ray leaving $x$. Let $x_0\in \mathrm{im}(V_x)$ such that $h(x_0)=h(x)+2r$ and $\mathcal{C}(x_0)$ be a subset of $X_{h(x_0)}$ containing $D(x_0,M)$ and contained in $D(x_0,2M)$. Then we claim that $B(x,r)$ is included in the box at scale $3r$ constructed from the cell $\mathcal{C}(x_0)$. Let $x'\in B(x,r)$, we recall that $d_r(x',x):=d_X(x',x)-\Delta h(x',x)$. By Lemma \ref{LEM0} we have that $d(V_x(h(x)+2r),V_{x'}(h(x)+2r))\leq 96\delta = M$ since $r\geq d_X(x',x)\geq \frac{1}{2} d_r(x',x)$ and since the distance between two vertical geodesics is decreasing in the upward direction.  Therefore $V_{x'}(h(x)+2r)\in \mathcal{C}(x_0)$. Furthermore $\Delta h (x_0,x')\leq \Delta h (x_0,x)+\Delta h (x,x')\leq 3r$, hence $x'\in \mathcal{B}(3r)$.
\end{proof}

\subsection{Tiling a big box by small boxes}

Let $R>0$ and $N\in\mathbb{N}$, next result shows that a box at scale $NR$ can be tiled with boxes at scale $R$. 

\begin{propo}\label{PropoTillBigBoxWithSmallBox}
Let $M_0$ be the constant of assumption $(E2)$. Let $R>0$ and $N\in\mathbb{N}$. Let $\mathcal{B}^X$ be a box at scale $NR$, and let us denote by $h^-:=h^-\left(\mathcal{B}^X\right)$ the lowest height of~$\mathcal{B}^X$. Then there exists a box tiling at scale $R$ of~$\mathcal{B}^X$. Otherwise stated for all $k\in\{1,\ldots,N\}$ there exists $\mathcal{D}_k\left(\mathcal{B}^X\right)\subset\mathcal{B}^X_{h^-+kR}$ such that:
\begin{enumerate}
\item For all $x\in \mathcal{D}_k\left(\mathcal{B}^X\right)$, there exists a cell $\mathcal{C}(x)$ such that $D_{M_0}(x)\subset \mathcal{C}(x)\subset D_{3M_0}(x)$.
\item We have $\bigsqcup\limits_{k=1}^N\bigsqcup\limits_{x\in\mathcal{D}_k\left(\mathcal{B}^X\right) }\mathcal{B}^X(x,\mathcal{C}(x),R)=\mathcal{B}^X$.
\end{enumerate} 
\end{propo}

\begin{proof}
To tile the box $\mathcal{B}^X$ we first tile by cells all of its level sets at height $h^-+kR$. Let $k\in\{1,\ldots,N\}$, and let $\mathcal{D}_k\left(\mathcal{B}^X\right)$ be an $2M_0$-maximal separating set of $\mathrm{Int}_{M_0}\left(\mathcal{B}^X_{h^-+kR}\right)$. Then:
\begin{enumerate}
\item For all $x,x'\in \mathcal{D}_k\left(\mathcal{B}^X\right)$ with $x\neq x'$ we have $D_{M_0}(x)\cap D_{M_0}(x')=\emptyset$.
\item $\mathrm{Int}_{M_0}\left(\mathcal{B}^X_{h^-+kR}\right) \subset \bigcup\limits_{x\in \mathcal{D}_k\left(\mathcal{B}^X\right)}D_{2M_0}(x) $.
\end{enumerate}
Furthermore $\mathcal{N}_{M_0}\left(\mathrm{Int}_{M_0}\left(\mathcal{B}^X_{h^-+kR}\right)\right)\subset\mathcal{B}^X_{h^-+kR}$, and for all $x\in \mathrm{Int}_{M_0}\left(\mathcal{B}^X_{h^-+kR}\right)$ we have $D_{M_0}(x)\subset \mathcal{B}^X_{h^-+kR}$. Therefore
\begin{align}
\bigsqcup\limits_{x\in \mathcal{D}_k\left(\mathcal{B}^X\right)}D_{M_0}(x)\subset \mathcal{B}^X_{h^-+kR} \subset \bigcup\limits_{x\in \mathcal{D}_k\left(\mathcal{B}^X\right)}D_{3M_0}(x).\label{SubsetTillgood}
\end{align}
For all $x\in \mathcal{D}_k\left(\mathcal{B}^X\right)$, we define
\begin{align*}
\mathcal{C}(x):=\left\lbrace p\in \mathcal{B}^X_{h^-+kR}\mid d(p,x)\leq d(p,x') \text{ for all }x'\in \mathcal{D}_k\left(\mathcal{B}^X\right) \right\rbrace.
\end{align*} 
As discussed at the beginning of Section \ref{SecTill}, these cells might intersect each other on their boundaries. However, a point contained in different cells can be removed in all of them except one, making them disjoint. The choice of cells on which we remove boundary points can be made thanks to an arbitrary order on the finite set $\mathcal{D}_k\left(\mathcal{B}^X\right)$.  
\\By the inclusions (\ref{SubsetTillgood}), for all $x\in \mathcal{D}_k\left(\mathcal{B}^X\right)$ we have $D_{M_0}(x)\subset \mathcal{C}(x)\subset D_{3M_0}(x)$ and
\begin{align*}
\bigsqcup\limits_{x\in\mathcal{D}_k\left(\mathcal{B}^X\right) }\mathcal{C}(x)=\mathcal{B}^X_{h^-+kR}.
\end{align*}
Furthermore, since vertical geodesic rays are uniquely determined by their starting point (because $X$ is Busemann), a tiling with cells provides us with a box tiling:
\begin{align*}
\bigsqcup\limits_{x\in\mathcal{D}_k\left(\mathcal{B}^X\right) }\mathcal{B}^X(x,\mathcal{C}(x),R)=\bigcup\limits_{z\in [h^-+(k-1)R;h^-+kR[}\mathcal{B}^X_z.
\end{align*}
Taking the union on $k\in\{1,\ldots,N\}$ provides us with the conclusion.
\end{proof}

\subsection{Box-tiling of $X\bowtie Y$}

\begin{figure}
\begin{center}
\includegraphics[scale=1]{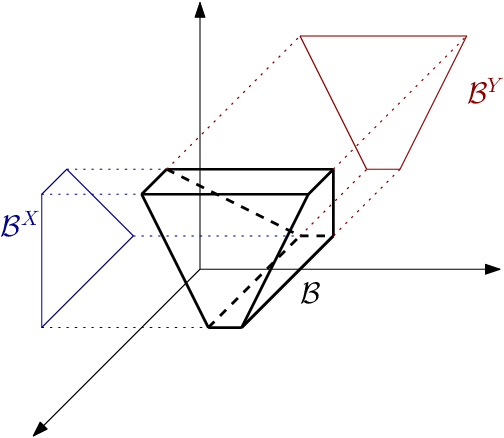}
\end{center}
\caption{Box in $X\bowtie Y$}
\label{FigBoxInHoro} 
\end{figure}

The boxes $\mathcal{B}$ of a horospherical product $X\bowtie Y$ are constructed as the horospherical products of boxes $\mathcal{B}^X\bowtie\mathcal{B}^Y$. Therefore they induce a tiling of $X\bowtie Y$. Such boxes are illustrated by Figure \ref{FigBoxInHoro}.

\begin{defn}(Box of $X\bowtie Y$ at scale $R$)
\\Let $X$ and $Y$ be two admissible spaces. A set $\mathcal{B}\subset X\bowtie Y$ is called box at scale $R$ of $X\bowtie Y$ if there exists $\mathcal{B}^X$ a box at scale $R$ of $X$ and $\mathcal{B}^Y$ a box at scale $R$ of $Y$ such that:
\begin{enumerate}
\item $h^{-}\left(\mathcal{B}^X\right)=-h^{+}\left(\mathcal{B}^Y\right)$;
\item $\mathcal{B}:=\mathcal{B}^X\bowtie \mathcal{B}^Y=\left\lbrace (x,y)\in\mathcal{B}^X\times \mathcal{B}^Y|h_X(x)=-h_Y(y)\right\rbrace$.
\end{enumerate}
\end{defn}

Let us point out that in the last definition, the box of $Y$ is in fact defined by

\begin{equation}
\mathcal{B}^Y(y,R):=\bigcup\limits_{z\in]-nR;(1-n)R]}\pi_z\big(\mathcal{C}(y)\big).
\end{equation}

This choice on the boundaries of the height intervals allows a precise match for the height inside the two boxes. Furthermore, one can see that given a box-tiling of $X$ and a box-tiling of $Y$, the natural subsequent tiling on $X\times Y$ provides the box tiling of $X\bowtie Y$ by restriction.

\begin{propo}(Box-tiling of $X\bowtie Y$ at scale $R$)
\\Let $X$ and $Y$ be two admissible spaces. Let $R$ be a positive number and let us consider the two following box tilings of $X$ and $Y$: 
\begin{align*}
X&=\bigsqcup\limits_{n\in\mathbb{Z}}\bigsqcup\limits_{x\in \mathcal{D}\left(X_{nR}\right)}\mathcal{B}^X(x,R);
\\Y&=\bigsqcup\limits_{n\in\mathbb{Z}}\bigsqcup\limits_{y\in \mathcal{D}\left(Y_{nR}\right)}\mathcal{B}^Y(y,R).
\end{align*}
Then the boxes of $X\bowtie Y$ constructed from boxes at opposite height in $X$ and $Y$ are a box tiling of $X\bowtie Y$. We have
\begin{align*}
X\bowtie Y=\bigsqcup\limits_{n\in\mathbb{Z}}\bigsqcup\limits_{(x,y)\in \mathcal{D}\left(X_{nR}\right)\times \mathcal{D}\left(Y_{(1-n)R}\right)}\mathcal{B}^X(x,R)\bowtie \mathcal{B}^Y(y,R).
\end{align*}
\end{propo}

\begin{proof}
Let us consider the box tilings of $X$ and $Y$:
\begin{align*}
X&=\bigsqcup\limits_{n\in\mathbb{Z}}\bigsqcup\limits_{x\in \mathcal{D}\left(X_{nR}\right)}\mathcal{B}^X(x,R);
\\Y&=\bigsqcup\limits_{n\in\mathbb{Z}}\bigsqcup\limits_{y\in \mathcal{D}\left(Y_{nR}\right)}\mathcal{B}^Y(y,R).
\end{align*}
We first show that the intersection of two distinct boxes is empty. Let $n_1,n_2\in\mathbb{R}$, $x_1\in \mathcal{D}\left(X_{n_1 R}\right)$, $x_2\in \mathcal{D}\left(X_{n_2 R}\right)$, $y_1\in \mathcal{D}\left(Y_{(1-n_1) R}\right)$ and $y_2\in \mathcal{D}\left(X_{(1-n_2) R}\right)$ such that $(x_1,y_1)\neq (x_2,y_2)$. Then we have either $x_1\neq x_2$ or $y_1\neq y_2$. Let us consider the case $x_1\neq x_2$, then $\mathcal{B}^X(x_1,R)\neq \mathcal{B}^X(x_2,R)$, and since they are two tiles of the box tiling of $X$, we have $\mathcal{B}^X(x_1,R)\cap \mathcal{B}^X(x_2,R)=\emptyset$. Therefore
\begin{align*}
\forall \left(p_1^X,p_1^Y\right) \in \mathcal{B}^X(x_1,R)\bowtie \mathcal{B}^Y(y_1,R),\ \forall \left(p_2^X,p_2^Y\right) \in \mathcal{B}^X(x_2,R)\bowtie \mathcal{B}^Y(y_2,R),\text{ we have } p_1^X\neq p_2^X.
\end{align*} 
Hence $(p_1^X,p_1^Y)\neq (p_2^X,p_2^Y)$, which gives us
\begin{align*}
\left(\mathcal{B}^X(x_1,R)\bowtie \mathcal{B}^Y(y_1,R) \right)\cap \left(\mathcal{B}^X(x_2,R)\bowtie \mathcal{B}^Y(y_2,R)\right)  =\emptyset.
\end{align*} 
The case when $y_1\neq y_2$ provide us with the same conclusion. Then we prove that the whole space $X\bowtie Y$ is covered by the horospherical product of boxes. Let $p=\left(p^X,p^Y\right)\in X\bowtie Y$. There exists $n\in\mathbb{Z}$ such that $(n-1)R\leq h(p) <nR$, hence there exist $x\in\mathcal{D}(X_{nR})$ and $y\in \mathcal{D}(Y_{(1-n)R})$ such that $p^X\in \mathcal{B}^X(x,R)$ and $p^Y\in\mathcal{B}^Y(y,R)$. Therefore $p\in\mathcal{B}^X(x,R)\bowtie \mathcal{B}^Y(y,R)$.
\end{proof}

\subsection{Measure of balls, boxes and neighbourhoods}

The results of this section focus on estimates on the measure $\mu$ of balls and boxes.

\begin{lemma}\label{LemmaBoxMeasure}
There exists $M(\bowtie)$ such that for all $r\geq M$ and all box $\mathcal{B}$ at scale $r$ of $X\bowtie Y$ we have
\begin{equation}
 \mu(\mathcal{B})\asymp_{\bowtie} e^{mr}.
\end{equation}
\end{lemma}
\begin{proof}
Without loss of generality we can assume that $h\left(\mathcal{B}\right)=[0;r[$. Let us denote by $\mathcal{C}^X$ the cell of $\mathcal{B}^X$ and $\mathcal{C}^Y$ the cell of $\mathcal{B}^Y$. Then
\begin{align*}
\mu(\mathcal{B})&=\int\limits_{0}^r \mu_z(\mathcal{B}_z)\mathrm{d}z=\int\limits_{0}^r \mu^X_z\left(\mathcal{B}^X_z\right)\mu^Y_{-z}\left(\mathcal{B}^Y_{-z}\right)\mathrm{d}z,\quad\text{since }\mathcal{B}_z=\mathcal{B}^X_z\times\mathcal{B}^Y_{-z},
\\&\asymp_{\bowtie} \int\limits_{0}^r e^{m(r-z)}\mu^X_r\left(\mathcal{C}^X\right)e^{nz}\mu^Y_0\left(\mathcal{C}^Y\right)\mathrm{d}z,\quad\text{by assumption  }(E3)\text{ and definition of boxes},
\\&\asymp_{\bowtie} e^{mr} \int\limits_{0}^r e^{(n-m)z}\mathrm{d}z,\quad\text{by Lemma \ref{LemmaMuDiskAreExp}},
\\&=\frac{e^{mr}-e^{nr}}{m-n}\preceq_{\bowtie} e^{mr}.
\end{align*}
However $m>n$, hence for $r\geq \frac{1}{m-n}$ we have $\frac{1}{2}e^{mr}\geq e^{nr}$. Therefore
\begin{align*}
\frac{e^{mr}-e^{nr}}{m-n}\geq\frac{e^{mr}}{2(m-n)}\succeq_{\bowtie}e^{mr}.
\end{align*}
\end{proof}

Combining Lemmas \ref{LemmaBoxAndBalls} and \ref{LemmaBoxMeasure} we get the next corollary.

\begin{cor}
There exists $M(\bowtie)$ such that for any $r\geq M$ and any $p\in X\bowtie Y$ we have
\begin{equation}
 e^{\frac{m}{2}r}\preceq_{\bowtie}\mu(B(p,r))\preceq_{\bowtie} e^{3mr}.
\end{equation}
\end{cor}

Therefore we have the following estimate between ball measures, corresponds to Lemma 5.2 in \cite{EFW1}.

\begin{cor}\label{CoroUniformBallMeasure}
There exists $M(\bowtie)$ such that for any $r_2>2r_1\geq M$ and for all $p_1,p_2\in X\bowtie Y$
\begin{align*}
\exp \big(\frac{1}{6}|r_2-r_1|m\big)\mu\big(B(p_1,r_1)\big)\leq\mu\big(B(p_2,r_2)\big)\leq \exp \big(6|r_2-r_1|m\big)\mu\big(B(p_1,r_1)\big).
\end{align*}
\end{cor}

\begin{cor}\label{CoroMeasureNeighbourhood}
There exists $M(\bowtie)$ such that for any $r_2>r_1\geq M$ and for all $A\subset X\bowtie Y$
\begin{align*}
\mu\left(\mathcal{N}_{r_2}(A)\right)\preceq_{\bowtie}e^{6|r_2-r_1|m}\mu\left(\mathcal{N}_{r_1}(A)\right).
\end{align*}
Furthermore, if there exists $z\in\mathbb{R}$ such that $A\subset X_z$ we have
\begin{align*}
\mu\left(\mathcal{N}_{M}(A)\right)\asymp_{\bowtie}\mu_z\left(\mathcal{N}_{M}(A)\cap X_z\right).
\end{align*}
In particular, for all $p\in (X\bowtie Y)_z$
\begin{align*}
\mu\left(B(p,M)\right)\asymp_{\bowtie}\mu_z\left(D_M(p)\right).
\end{align*}
\end{cor}
\begin{proof}
Since $X\bowtie Y$ is a proper metric space, by a covering lemma of \cite{HE}, there exists a set $Z\subset A$ such that:
\begin{enumerate}
\item The balls $ B(p,r_1)$ for $p\in Z$ are pairwise disjoint.
\item We have the following inclusions:
$$\bigsqcup\limits_{p\in Z} B(p,r_1)\subset \mathcal{N}_{r_1}(A)  \subset\bigcup\limits_{p\in Z} B(p,5r_1). $$
\end{enumerate} 
Therefore $\mathcal{N}_{r_2}(A)  \subset\bigcup\limits_{p\in Z} B(p,5r_1+(r_2-r_1))$.
\\\\Moreover, if $A\subset X_z$, for $r_1=M$ we have
\begin{align*}
\bigsqcup\limits_{p\in Z} D_M(p)\subset \mathcal{N}_{M}(A)\cap X_z  \subset\bigcup\limits_{p\in Z} D_{5M}(p),
\end{align*}
and for all $p\in Z$, $\mu_{z}(B(p,5M))\asymp_{\bowtie}1\asymp_{\bowtie}\mu_{z}(D_{5M}(p))$. Hence
\begin{align*}
\mu\left(\mathcal{N}_{M}(A)\right)&\asymp_{\bowtie}\sum\limits_{p\in Z} \mu(B(p,5M)),
\\&\asymp_{\bowtie}\sum\limits_{p\in Z} \mu_z(D_{5M}(p))\asymp_{\bowtie}\mu_z\left(\mathcal{N}_{M}(A)\cap X_z\right).
\end{align*}
\end{proof}

A $(k,c)$-quasi-isometry $\Phi:X\bowtie Y \to X'\bowtie Y'$ "quasi"-preserve the measure $\mu$.

\begin{lemma}\label{LemmaVolKeepByQI}
For all $(k,c)$-quasi-isometry $\Phi:X\bowtie Y \to X'\bowtie Y'$ and for all measurable subset $U\subset  X\bowtie Y$ we have
\begin{align*}
\mu\big(\mathcal{N}_{k(c+1)}(U)\big)\asymp_{k,c,\bowtie} \mu\Big(\mathcal{N}_{1}\big(\Phi(U)\big)\Big).
\end{align*}
\end{lemma}

\begin{proof}
Since $X\bowtie Y$ is a proper metric space, by a classical covering lemma of \cite{HE} there exists a set $Z\subset U$ such that:
\begin{enumerate}
\item The balls $ B(p,k(c+1))$ for $p\in Z$ are pairwise disjoint.
\item We have the following inclusions:
$$\bigsqcup\limits_{p\in Z} B(p,k(c+1))\subset \mathcal{N}_{k(c+1)}(U)  \subset\bigcup\limits_{p\in Z} B(p,5k(c+1)). $$
\end{enumerate} 
Since $\Phi$ is a $(k,c)$-quasi-isometry,  $\Phi(Z)$ verifies:
\begin{enumerate}
\item The balls $B(q,1)$ for $q\in \Phi(Z)$ are pairwise disjoint.
\item We have the following inclusions:
$$\bigsqcup\limits_{q\in \Phi(Z)} B(q,1)\subset \mathcal{N}_1\left(\Phi(U)\right)\subset \bigcup\limits_{q\in \Phi(Z)} B(q,5k^2(c+1)+c). $$
\end{enumerate} 
Furthermore, for all $p\in Z$ we have $$\mu(B(p,1))\asymp_{\bowtie}1\asymp_{\bowtie'}\mu(B(\Phi(p),1))\asymp_{k,c,\bowtie'}\mu(B(\Phi(p),5k^2(c+1)+c)),$$
hence $\mu\big(\mathcal{N}_{k(c+1)}(U)\big)\asymp_{k,c,\bowtie}\# Z \asymp_{k,c,\bowtie'} \mu\Big(\mathcal{N}_{1}\big(\Phi(U)\big)\Big)$.
\end{proof}

\subsection{Set of vertical geodesics}\label{SecVertSet}

Since $X$ is a Gromov hyperbolic, Busemann space, for any $x\in X$, there exists a unique vertical geodesic ray starting from $x$ in $X$, therefore, there is a one to one correspondence between portions of vertical geodesic rays in a box $\mathcal{B}^X$, and the points at the bottom of $\mathcal{B}^X$. A vertical geodesic segment of $\mathcal{B}^X$ is defined as the intersection of a vertical geodesic and $\mathcal{B}^X$. We recall that vertical geodesics are parametrised by arclength by their height.

Let $\mathcal{B}^X$ be a box at scale $R$ of $X$. Let us denote by $V\mathcal{B}^X$ the set of vertical geodesic segments of $\mathcal{B}$. A geodesic segment $v\in V\mathcal{B}^X$ intersects only in one point $x$ the bottom of $\mathcal{B}^X$, and $v$ is the only vertical geodesic segment of $ V\mathcal{B}^X$ intersecting $x$ by the Busemann assumption on $X$.

\begin{defn}(Measure $\eta$ on $V\mathcal{B}^X$)
\\Let $\mathcal{B}^X$ be a box at scale $R$ of $X$. The measure $\eta^X_{V\mathcal{B}^X}$ on $V\mathcal{B}^X$ is defined on all measurable subset $U\subset V\mathcal{B}^X$ by
\begin{equation}
\eta^X_{V\mathcal{B}^X}\big(U\big)=\lambda^X_{h^-(\mathcal{B}^X)} \Big(\left\lbrace \gamma\big(h^-(\mathcal{B}^X)\big)\mid \gamma \in U \right\rbrace\Big).
\end{equation}
\end{defn}

In particular, we say that $U$ is measurable if $\left\lbrace \gamma\big(h^-(\mathcal{B}^X)\big)\mid \gamma \in U \right\rbrace$ is measurable. Since the measure $\lambda$ is almost constant along projections, the measure on the set of vertical geodesic segment is related to the height of the boxes. Specifically we show that up to a multiplicative constant, the measure of a box is equal to the measure of its set of vertical geodesic segments multiplied by its height, as for rectangles in $\mathbb{R}^2$. In the sequel we might omit the index of the measure $\eta^X$.

\begin{propr}\label{PropMeasureGeodBoxHyper}
Let $\mathcal{B}^X$ be a box at scale $R$ of $X$ and let us denote $h^-:=h^-(\mathcal{B}^X)$ and $h^+:=h^+(\mathcal{B}^X)$. We have for all $z\in [h^-,h^+[$:
\begin{enumerate}
\item $\eta^X\big(V\mathcal{B}^X\big)\asymp_{X} \lambda_{z}^X\big(\mathcal{B}_z^X\big)\asymp_{X}e^{m h^+}$;
\item $\lambda^X\big(\mathcal{B}^X\big)\asymp_{X}R\lambda_{z}^X\big(\mathcal{B}_z^X\big)\asymp_{X} R \eta^X\big(V\mathcal{B}^X\big)\asymp_{X} R e^{m h^+}$.
\end{enumerate}

\end{propr}
\begin{proof}
Let $x\in X$ be such that $\mathcal{C}(x)$ is the cell of $\mathcal{B}^X$. We know that $D_{M_0}(x)\subset \mathcal{C}(x) \subset D_{2M_0}(x)$, hence by Lemma \ref{LemmaMuDiskAreExp} we have
\begin{equation}
\mu^X_{h(x)}(\mathcal{C}(x))\asymp_{X} 1.
\end{equation}
Then
\begin{align*}
\eta^X\big(V\mathcal{B}^X\big)&=\lambda^X_{h^-} \Big(\mathcal{B}^X\cap h^{-1}\big(h^-\big)\Big),\quad\text{by definition,}
\\&\asymp_{X}\lambda_{z}^X\big(\mathcal{B}_z^X\big)\asymp_{X}\lambda^X_{h^+} \big(\mathcal{C}(x)\big)\asymp_{X} e^{m h^+}\mu^X_{h^+} \big(\mathcal{C}(x)\big),\quad\text{by Property \ref{ProprLambdaSameOnLevels}}
\\&\asymp_{X} e^{m h^+},
\end{align*}
which proves the first point. The second point follows from the fact that the measures $\lambda_z$ are constant by projections on height level sets, up to the multiplicative constant $M(X)$. 
\begin{align*}
\lambda^X\big(\mathcal{B}^X\big)&=\int\limits_{h^-}^{h^+}\lambda^X_z\big(\mathcal{B}^X\cap h^{-1}(z)\big)\mathrm{d}z =\int\limits_{h^-}^{h^+}\lambda^X_z\big(\pi_z(\mathcal{C}(x))\big)\mathrm{d}z,
\\&\asymp_{X}\int\limits_{h^-}^{h^+}\lambda^X_{h^+}\big(\mathcal{C}(x)\big)\mathrm{d}z,\quad\text{ by Property \ref{ProprLambdaSameOnLevels}},
\\&\asymp_{X} R\lambda^X_{h^+}\big(\mathcal{C}(x)\big)\asymp_{X} R e^{m h^+}.
\end{align*}
\end{proof}

A vertical geodesic $V=(V^X,V^Y)\subset X\bowtie Y$ is a couple of vertical geodesics of $X$ and $Y$. Therefore, there is a bijection between the set of vertical geodesic segments $V\mathcal{B}$ of a box $\mathcal{B}:=\mathcal{B}^X\bowtie \mathcal{B}^Y$ and $V\mathcal{B}^X\times  V\mathcal{B}^Y$.

\begin{defn}\label{DefMeasureGeodVertBox}
Let $\mathcal{B}$ be a box at scale $R$ of $X\bowtie Y$. We define the measure $\eta_{V\mathcal{B}}$ on $V\mathcal{B}$ as
\begin{equation}
\eta_{V\mathcal{B}}:=\eta^X_{V\mathcal{B}^X}\otimes \eta^Y_{V\mathcal{B}^Y}.
\end{equation}
\end{defn} 

In the notation of measures on sets of vertical geodesic segments, we might omit the reference to the corresponding sets. The measures $\eta_{V\mathcal{B}}$, respectively $\eta_{V\mathcal{B}^X}^X$, $\eta_{V\mathcal{B}^Y}^Y$, will simply be denoted by $\eta$, respectively $\eta^X$, $\eta^Y$. 

\begin{propo}\label{PropMeasureGeodBoxHoro}
For each box $\mathcal{B}$ at scale $R$ of $X\bowtie Y$ we have for all $z_1,z_2\in[h^-,h^+[$:
\begin{enumerate}
\item $\eta\left(V\mathcal{B}\right)\asymp_{\bowtie }e^{m h^+}e^{-nh^-}\asymp_{\bowtie} \lambda_{z_1}^X\big(\mathcal{B}_{z_1}^X\big)\lambda_{-z_2}^Y\big(\mathcal{B}_{-z_2}^Y\big)$;
\item $\lambda\big(\mathcal{B}\big)\asymp_{\bowtie }R \eta\left(V\mathcal{B}\right)\asymp_{\bowtie} R\lambda_{z_1}^X\big(\mathcal{B}_{z_1}^X\big)\lambda_{-z_2}^Y\big(\mathcal{B}_{-z_2}^Y\big)$.
\end{enumerate}
\end{propo}
\begin{proof}
The first point follows from definition \ref{DefMeasureGeodVertBox} and Property \ref{PropMeasureGeodBoxHyper} applied on $\mathcal{B}^X$ and $\mathcal{B}^Y$. The proof of the second point is similar to the proof of Property \ref{PropMeasureGeodBoxHyper}
\begin{align*}
\lambda\big(\mathcal{B}\big)&=\int\limits_{h^-}^{h^+}\lambda^X_z\otimes\lambda^Y_{-z}\Big(\mathcal{B}^X_{z}\times\mathcal{B}^Y_{-z}\Big)\mathrm{d}z =\int\limits_{h^-}^{h^+}\lambda^X_z\Big(\mathcal{B}^X_{z}\Big)\lambda^Y_{-z}\Big(\mathcal{B}^Y_{-z}\Big)\mathrm{d}z,
\\&\asymp_{\bowtie}\int\limits_{h^-}^{h^+}\lambda^X_{h^-}\Big(\mathcal{B}^X_{h^-}\Big)\lambda^Y_{-h^+}\Big(\mathcal{B}^Y_{-h^+}\Big)\mathrm{d}z,\quad\text{ by Property \ref{ProprLambdaSameOnLevels}},
\\&\asymp_{\bowtie}\int\limits_{h^-}^{h^+}\eta^X \Big(V\mathcal{B}^X\Big)\eta^Y\Big(V\mathcal{B}^Y\Big)\mathrm{d}z,\quad\text{ by definition of }\eta,
\\&\asymp_{\bowtie}\eta\left(V\mathcal{B}\right)\int\limits_{h^-}^{h^+}1\mathrm{d}z = R\eta\left(V\mathcal{B}\right).
\end{align*}
Then applying twice Property \ref{PropMeasureGeodBoxHyper} provides us with the result.
\end{proof}

Let $\mathcal{B}$ be a box at scale R. Let $z\in \left[h^-\big(\mathcal{B}\big);h^+\big(\mathcal{B}\big)\right[$ and let $U\subset \mathcal{B}_z$. Then we denote $V_{\mathcal{B}}(U)$ the set of vertical geodesic segments of $V\mathcal{B}$ intersecting $U$, it is in bijection with
\begin{align*}
\lbrace(x,y)\in\mathcal{B}^X_{0}\times\mathcal{B}^Y_{-R} \big\mid \left(\pi^X_z(x),\pi^Y_{-z}(y)\right)\in U\big\rbrace.
\end{align*}

We need the following property stating that the measure of a given subfamily of vertical geodesics can be computed on any level of our box.

\begin{propr}\label{PropMeasureGeodVertSameMeasureSpace}
Let $\mathcal{B}$ be a box at scale $R$ of $X\bowtie Y$. Then for all $z\in \left[h^-\big(\mathcal{B}\big);h^+\big(\mathcal{B}\big)\right[$ and for all measurable subset $U_z\subset \mathcal{B}_z$
\begin{align*}
\eta\big(V_{\mathcal{B}}(U_z)\big)\asymp_{\bowtie} \lambda_z\big(U_z\big).
\end{align*}
\end{propr}

\begin{proof}
Without loss of generality we can assume that $ \left[h^-\big(\mathcal{B}\big);h^+\big(\mathcal{B}\big)\right[=[0:R[$. By definition we have
\begin{align*}
\eta\big(V_{\mathcal{B}}(U_z)\big)&:=\int\limits_{x_0\in\mathcal{B}^X_{0}}\int\limits_{y_0\in\mathcal{B}^Y_{-R}}\mathbbm{1}_{\big\lbrace(x,y)\in\mathcal{B}^X_0\times\mathcal{B}^Y_{-R} \big\mid \left(\pi^X_z(x),\pi^X_{-z}(y)\right)\in U_z\big\rbrace}\big(x_0,y_0\big) d\lambda^Y_{-R} d\lambda^X_{0},
\\&=\int\limits_{x_0\in\mathcal{B}^X_0}\int\limits_{y_0\in\mathcal{B}^Y_{-R}}\mathbbm{1}_{U_z}\Big(\pi_z^X(x_0),\pi_{-z}^Y(y_0)\Big) d\lambda^Y_{-R} d\lambda^X_{0},
\\&=\int\limits_{x_0\in\mathcal{B}^X_{0}}\left(\int\limits_{y\in\mathcal{B}^Y_{-z}}\mathbbm{1}_{U_z}\Big(\pi_z^X(x_0),y\Big) d\left(\left(\pi^Y_{-z}\right)_{*}\lambda^Y_{-R}\right) \right)d\lambda^X_{0},\quad\text{with a pushforward of }\lambda^Y_{-R}\text{ by }\pi^Y_{-z},
\\&=\int\limits_{y\in\mathcal{B}^Y_{-z}}\left(\int\limits_{x_0\in\mathcal{B}^X_{0}}\mathbbm{1}_{U_z}\Big(\pi_z^X(x_0),y\Big) d\lambda^X_{0} \right)d\left(\left(\pi^Y_{-z}\right)_{*}\lambda^Y_{-R}\right),\quad\text{by Fubini's Theorem,}
\\&=\int\limits_{y\in\mathcal{B}^Y_{-z}}\left(\int\limits_{x\in\mathcal{B}^X_{z}}\mathbbm{1}_{U_z}\big(x,y\big) d\left(\left(\pi^Y_{z}\right)_{*}\lambda^X_{0}\right) \right)d\left(\left(\pi^Y_{-z}\right)_{*}\lambda^Y_{-R}\right),\quad\text{with a pushforward of }\lambda^X_{0}\text{ by }\pi^X_{z},
\\&\asymp_{\bowtie}\int\limits_{y\in\mathcal{B}^Y_{-z}}
\int\limits_{x\in\mathcal{B}^X_{z}}\mathbbm{1}_{U_z}\big(x,y\big) d\lambda^X_{z}
d\lambda^Y_{-z},\quad\text{by using Property \ref{ProprLambdaSameOnLevels} twice,}
\\&\asymp_{\bowtie}\lambda_z\big( U_z \big).
\end{align*}
\end{proof}

\subsection{Projections of set of almost full measure}

Let us denote by $p^X:X\bowtie Y\to X~;~(x,y)\mapsto x$ and by $p^Y:X\bowtie Y\to X~;~(x,y)\mapsto y$ the projections on the two coordinates of $X\bowtie Y$. We also denote by slight abuse the projection on a set of vertical geodesic segments $p^X:V\mathcal{B}\to V\mathcal{B}^X~;~(v^X,v^Y)\mapsto v^X$ and $p^X:V\mathcal{B}\to V\mathcal{B}^Y~;~(v^X,v^Y)\mapsto v^Y$. Given a subset $U\subset \mathcal{B}$, we might simply denote by $U^X$, respectively $U^Y$, its projection on $X$, respectively on $Y$, and similarly for subsets of $V\mathcal{B}$.
\\\\In this section, we show that if a subset of a box has almost full measure, then most of the fibers with respect to these projections also have almost full measure.
 \\\\Let $0<\alpha\leq 1$, let $V_1\subset V\mathcal{B}$ be a measurable subset (it will be chosen later as a subset of small measure, containing "bad" vertical geodesics). Let us denote for all $v^X\in V\mathcal{B}^X$
\begin{align*}
G^Y\left(v^X\right):=&\left\lbrace v^Y\in V\mathcal{B}^Y\mid \left(v^X,v^Y\right)\in V_0\right\rbrace = p^{Y}\left((p^X)^{-1}\left(v^X\right)\cap (V\mathcal{B}\setminus V_1)\right).
\\G^X:=&\left\lbrace v^X\in V\mathcal{B}^X\mid \eta^Y \left(G^Y\left(v^X\right)\right)\geq \big(1-\sqrt{\alpha}\big)\eta^Y\left(V_1^Y\right)\right\rbrace.
\end{align*}
The set $G^X$ is the set of vertical geodesics in $V\mathcal{B}^X$ whose fibers have almost full intersection with $V\mathcal{B}\setminus V_1$.

The following lemma asserts that almost all fibers have almost full intersection with $V\mathcal{B}\setminus V_1$.

\begin{lemma}\label{LemmaMarkov1}
Let $0<\alpha \leq 1$ and let $V_1\subset V\mathcal{B}$ be a measurable subset such that $\eta(V_1)\leq \alpha \eta(V\mathcal{B})$, then
\begin{align*}
\eta^X\left(G^X\right)\geq \big(1-\sqrt{\alpha}\big)\eta^X\left(V\mathcal{B}^X\right).
\end{align*}
\end{lemma}

\begin{proof}
By construction we have 
\begin{align*}
\bigcup _{v^X\in V\mathcal{B}^X} G^Y\left(v^X\right) =& \left(V\mathcal{B}\setminus V_1\right)^Y.
\end{align*}
To prove the Lemma we proceed by contradiction. Let us assume that $\eta^X\left(G^X\right)< \big(1-\sqrt{\alpha}\big)\eta^X\left(V\mathcal{B}^X\right)$, then $\eta^X\left(V\mathcal{B}^X\setminus G^X\right)> \sqrt{\alpha}\eta^X\left(V\mathcal{B}^X\right)$. Therefore
\begin{align*}
\eta\left(V_1\right)=&\int\limits_{V\mathcal{B}}\mathbbm{1}_{V_1}(v)\mathrm{d}\eta(v),
\\=&\int\limits_{V\mathcal{B}^X}\int\limits_{V\mathcal{B}^Y}\mathbbm{1}_{V_1}\left(v^X,v^Y\right)\mathrm{d}\eta^Y(v^Y)\mathrm{d}\eta^X(v^X),\quad\text{ by definition of }\eta,
\\=&\int\limits_{V\mathcal{B}^X}\int\limits_{V\mathcal{B}^Y}\mathbbm{1}_{V\mathcal{B}^Y\setminus G^Y\left(v^X\right)}\left(v^Y\right)\mathrm{d}\eta^Y(v^Y)\mathrm{d}\eta^X(v^X),\quad\text{ by definition of }G^Y\left(v^X\right),
\\=&\int\limits_{V\mathcal{B}^X}\eta^Y\left(V\mathcal{B}^Y\setminus G^Y\left(v^X\right)\right)\mathrm{d}\eta^X(v^X)\geq \int\limits_{V\mathcal{B}^X\setminus G^X}\eta^Y\left(V\mathcal{B}^Y\setminus G^Y\left(v^X\right)\right)\mathrm{d}\eta^X(v^X).
\end{align*}
Furthermore, when $v^X\in V\mathcal{B}^X\setminus G^X$ we have that $\eta^Y\left(G^Y\left(v^X\right)\right)<(1-\sqrt{\alpha})\eta^Y\left(V\mathcal{B}^Y\right)$, hence $\eta^Y\left(V\mathcal{B}^Y\setminus G^Y\left(v^X\right)\right)\geq \sqrt{\alpha}\eta^Y\left(V\mathcal{B}^Y\right)$. Therefore
\begin{align*}
\eta\left(V_1\right)\geq& \int\limits_{V\mathcal{B}^X\setminus G^X}\sqrt{\alpha} \eta^Y\left(V\mathcal{B}^Y\right)\mathrm{d}\eta^X(v^X),
\\\geq&\sqrt{\alpha} \eta^Y\left(V\mathcal{B}^Y\right)\eta^X\left(V\mathcal{B}^X\setminus G^X\right),
\\\geq&\sqrt{\alpha}\sqrt{\alpha} \eta^Y\left(V\mathcal{B}^Y\right)\eta^X\left(V\mathcal{B}^X\right),\quad\text{by the contradiction assumption},
\\&>\alpha \eta(V\mathcal{B}),\quad\text{since }V\mathcal{B}\text{ is a product},
\end{align*}
which contradicts $\eta\left(V_1\right)\leq\alpha\eta\left(V\mathcal{B}\right)$.
\end{proof}

In the previous Lemma we only used the fact that the set of vertical geodesic segments $V\mathcal{B}$ was the product of its projections endowed with a product measure $\eta$. We will use it once again on the product of two measured spaces endowed with a product measure in the proof of Proposition \ref{PropoOrientPreserving}.
\\\\We recall that for any $U\subset X\bowtie Y$ we denote $V\mathcal{B}(U):=\lbrace v\in V\mathcal{B}\mid \text{im}(v)\cap U\neq\emptyset\rbrace$. Similarly for all $V_1\subset V\mathcal{B}$ we denote $V_1(U):=\lbrace v\in V_1\mid \text{im}(v)\cap U\neq\emptyset\rbrace$.
\\\\The next Lemma is a local version of Lemma \ref{LemmaMarkov1}. Let $V_1\subset V\mathcal{B}$. Let $M>0$ be a constant, let $a\in \mathcal{B}$ and let us denote $VD:=V\mathcal{B}\left(D_{M}(a)\right)$ and $V_1 D:=V_1\left(D_{M}(a)\right)$. For all $v=(v^X,v^Y)\in V_{\mathcal{B}}$, let us denote by
\begin{align*}
E^Y\left(v^X\right):=&\left\lbrace v^Y\in VD^Y\mid \left(v^X,v^Y\right)\in V_1 D\right\rbrace = \left(p^{Y}\right)^{-1}\left(p^X\left(v^X\right)\cap  V_1D\right);
\\F^X:=&\left\lbrace v^X\in VD^X \mid \eta^Y\big(E^Y(v^X)\big)\geq \sqrt{\alpha}\eta^Y\left(VD^Y\right)\right\rbrace.
\end{align*} 

\begin{lemma}\label{LemmaMarkov2}
Let $0<\alpha\leq 1$. If $\eta(V_1D)\leq \alpha\eta(VD)$ then
\begin{equation}
\eta^X\left(F^X\right)\leq\sqrt{\alpha}\eta^X\left(VD^X\right).
\end{equation}
\end{lemma}

\begin{proof}
Let us proceed by contradiction. We assume that
\begin{equation}
\eta^X\left(F_i^X\right)>\sqrt{\alpha}\eta^X\left(VD_i^X\right).
\end{equation}
Then we have
\begin{align*}
\eta\left(V_1 D\right)=&\int\limits_{v^X\in V D^X}\int\limits_{v^Y\in V D^Y}\mathbbm{1}_{V_1 D}(v^X,v^Y)\mathrm{d}\eta^Y\mathrm{d}\eta^X,
\\=&\int\limits_{v^X\in V D^X}\int\limits_{v^Y\in V D^Y}\mathbbm{1}_{E^Y(v^X)}(v^Y)\mathrm{d}\eta^Y\mathrm{d}\eta^X,
\\=&\int\limits_{v^X\in V D^X}\eta^Y\left(E^Y(v^X)\right)\mathrm{d}\eta^X,\quad\text{by the definition of }E^Y(v^X),
\\\geq&\int\limits_{v^X\in F^X}\eta^Y\left(E^Y(v^X)\right)\mathrm{d}\eta^X,\quad\text{ since }F^X\subset VD^X,
\\>&\sqrt{\alpha}\eta^X\left(VD^X\right)\sqrt{\alpha}\eta^Y(VD^Y)>\alpha\eta(VD),
\end{align*}
which contradicts assumption on $VD$. Hence $\eta^X(F^X)\leq\sqrt{\alpha}\eta^X(VD^X)$.
\end{proof}

The following lemma asserts that for almost all points of the box, almost all vertical geodesics passing through the disc $D_{M_0}(x)$ do not belong to $V_1$. 
 
\begin{lemma}\label{LemmaCloseGeodAreGood}
There exists a constant $0<\alpha (\bowtie )\leq 1$ such that for all $0<\alpha\leq \alpha (\bowtie )$ the following statement holds. Let $M_0$ be the constant involved in assumption $(E2)$ and let $\mathcal{B}$ be a box at scale $R$. If $V_1\subset V\mathcal{B}$ satisfies $\eta\left(V_1\right)\leq \alpha  \eta\left(V\mathcal{B}\right)$, then
\begin{equation}
\lambda\left(\left\lbrace x\in\mathcal{B}\mid \dfrac{\eta\left(V_1\big(D_{M_0}(x)\big)\right)}{\eta\left(V\mathcal{B}\big(D_{M_0}(x)\big)\right)}>\alpha^{\frac{1}{4}}\right\rbrace\right)\leq\alpha^{\frac{1}{4}}\lambda\left(\mathcal{B}\right).
\end{equation}
\end{lemma}
\begin{proof}
Without loss of generality we may assume that $h(\mathcal{B})=[0;R[$. Let us denote
\begin{equation}
U=\left\lbrace x\in\mathcal{B}\mid \dfrac{\eta\left(V_1\big(D_{M_0}(x)\big)\right)}{\eta\left(V\mathcal{B}\big(D_{M_0}(x)\big)\right)}>\alpha^{\frac{1}{4}}\right\rbrace.
\end{equation}
We proceed by contradiction, let us assume that $\lambda(U)>\alpha^{\frac{1}{4}}\lambda (\mathcal{B})$. In this case there exists $z\in [0;R[$ such that $\lambda_z(U_z)> \alpha^{\frac{1}{4}} \lambda_z (\mathcal{B}_z)$. Let $U_z'\subset U_z$ be a $2M_0$ maximal separating set of $U_z$. We have that $\bigsqcup\limits_{x\in U_z'}D_{M_0}(x)$ is a disjoint union and that $U_z \subset\bigcup\limits_{x\in U_z'}D_{2M_0}(x)$. Then we have
\begin{align}
\lambda_z\left(\bigsqcup\limits_{x\in U_z'}D_{M_0}(x)\right)&= \sum\limits_{x\in U_z'}\lambda_z\left(D_{M_0}(x)\right) =\sum\limits_{x\in U_z'}\lambda_z\left(D_{2M_0}(x)\right)\dfrac{\lambda_z\left(D_{M_0}(x)\right)}{\lambda_z\left(D_{2M_0}(x)\right)},\nonumber
\\&\asymp_{\bowtie} \sum\limits_{x\in U_z'}\lambda_z\left(D_{2M_0}(x)\right),\quad\text{by Lemma \ref{LemmaMuDiskAreExp}}\nonumber,
\\&\geq \lambda_z\left(\bigcup\limits_{x\in U_z'}D_{2M_0}(x)\right)\geq \lambda_z\left(U_z\right),\nonumber
\\&\succeq_{\bowtie} \alpha^{\frac{1}{4}}\lambda_z\left(\mathcal{B_z}\right),\quad\text{by assumption on }U_z.\label{IneqDisqueCoverLevel}
\end{align}
However $\forall x\in U_z'$ we have $\eta\left(V_1\big(D_{M_0}(x)\big)\right)>\alpha^{\frac{1}{4}}\eta\left(V\mathcal{B}\big(D_{M_0}(x)\big)\right)$, therefore
\begin{align*}
\eta\left(V_1\left(\bigcup\limits_{x\in U_z'}D_{M_0}(x)\right)\right)&>\alpha^{\frac{1}{4}}\eta\left(V\mathcal{B}\left(\bigcup\limits_{x\in U_z'}D_{M_0}(x)\right)\right),
\\&\asymp_{\bowtie}\alpha^{\frac{1}{4}}\lambda_z\left(\bigcup\limits_{x\in U_z'}D_{M_0}(x)\right),\quad\text{by Lemma \ref{PropMeasureGeodVertSameMeasureSpace}},
\\&\geq\alpha^{\frac{1}{4}}\alpha^{\frac{1}{4}}\lambda_z\left(\mathcal{B}\right)=\sqrt{\alpha}\lambda_z\left(\mathcal{B}\right),\quad\text{by inequality (\ref{IneqDisqueCoverLevel})},
\\&\asymp_{\bowtie}\sqrt{\alpha}\eta\left(V\mathcal{B}\right),\quad\text{by Lemma \ref{PropMeasureGeodVertSameMeasureSpace}}.
\end{align*}
Since $\eta\left(V_1\right)\geq\eta\left(V_1\left(\bigcup\limits_{x\in U_z'}D_{M_0}(x)\right)\right)$ and since $\sqrt{\alpha }> M(\bowtie)\alpha $ for $\alpha< \frac{1}{M^2}$, it contradicts the assumptions of the lemma.
\end{proof}

Let us point out that in this Lemma, we first showed that on a fixed level-set, most of its point were surrounded by almost only of vertical geodesic not in $V_1$. This remark will be relevant in the proof of Proposition \ref{PropoOrientPreserving}.
\\\\The next three lemmas are estimates on the quantity of $Y$-horospheres verifying specific properties. They are used in section \ref{SecPartII}. Let $\mathcal{B}$ be a box, $x\in\mathcal{B}^X$ let $U\subset \mathcal{B}$ and let us denote by
\begin{align*}
H_x:=\{ x\}\bowtie \mathcal{B}^Y=\{(x,y)\mid y\in\mathcal{B}^Y,\ h(y)=-h(x)\}=(p^X)^{-1}(x),
\end{align*}
a $Y$-horosphere of $\mathcal{B}$. Let us denote by
\begin{align*}
E^Y(x):=&\left\lbrace y \in \mathcal{B}^{Y} \mid (x,y) \in U^c\right\rbrace = p^Y(p^{X -1}(x) \cap U^c)=(H_x\cap U^c)^Y;
\\E^X:=&\left\lbrace x\in \mathcal{B}^X\mid \lambda^Y_{-h(x)}\left(E^Y(x)\right)>\sqrt{\alpha}\lambda^Y\left(H^Y_x\right)\text{ and }h(x)\geq h^-(\mathcal{B}^X)+\frac{R}{2}\right\rbrace.
\end{align*}

The set $E^X$ is in bijection with the "bad" $Y$-horospheres $H$ above the middle of $\mathcal{B}$, the ones which have more than $\sqrt{\alpha}$ fraction of their measure $\lambda^Y$ in $U^c$.

The following lemma asserts that almost all $Y$-horospheres in the upper half of the box are good $Y$-horospheres. 

\begin{lemma}
If $\lambda(U)\geq \left(1-\alpha\right)\lambda (\mathcal{B})$ with $0<\alpha< 1$, then we have
\begin{align*}
\lambda^X(E^X)<\sqrt{\alpha}\lambda^X\big(\mathcal{B}^X\big).
\end{align*}
\end{lemma}

\begin{proof}
Without loss of generality we can assume that $h(\mathcal{B})=[0;R[$. We proceed by contradiction, let us assume that $\lambda^X\left(E^X\right)\geq \sqrt{\alpha}\lambda^X\left(\mathcal{B}^X\right)$. Then we compute the measure of $U^c$:
\begin{align*}
\lambda\left(U^c\right)&=\int\limits_{0}^R\lambda_z^X\otimes\lambda_{-z}^Y\left(U_z^c\right)\mathrm{d}z=\int\limits_{0}^R\int\limits_{\mathcal{B}^X_z}\lambda_{-z}^Y\big(\lbrace y\in Y_{-z}\mid (x,y)\in U^c_z\rbrace\big)\mathrm{d}\lambda_z^X(x)\mathrm{d}z,\quad\text{by definition,}
\\&=\int\limits_{0}^R\int\limits_{\mathcal{B}^X_z}\lambda^Y_{-z}\left((H_x\cap U^c)^Y\right)\mathrm{d}\lambda_z^X(x)\mathrm{d}z,
\\&\geq \int\limits_{0}^R\int\limits_{E^X_z}\lambda^Y_{-z}\left((H_x\cap U^c)^Y\right)\mathrm{d}\lambda_z^X(x)\mathrm{d}z,\quad\text{since }E_z^X\subset \mathcal{B}_z^X,
\\&>\sqrt{\alpha}\int\limits_{0}^R\left[\int\limits_{E^X_z}\lambda^Y_{-z}\left(H^Y_x\right)\mathrm{d}\lambda_z^X(x)\right]\mathrm{d}z ,\quad\text{by the definition of }E^X,
\\&=\sqrt{\alpha}\int\limits_{0}^R\left[\lambda^Y_{-z}\left(\mathcal{B}^Y_{-z}\right)\lambda^X_{z}\left(E^X_z\right)\right]\mathrm{d}z,\quad\text{by the definition of }H_x,
\\&\geq \sqrt{\alpha}\sqrt{\alpha}\int\limits_{0}^R\lambda^Y_{-z}\left(\mathcal{B}^Y_{-z}\right)\lambda^X_{z}\left(\mathcal{B}^X_{z}\right)\mathrm{d}z\geq \alpha \lambda (\mathcal{B}),\quad\text{by assumption on }E^X,
\end{align*}
which contradicts the assumption on $U$.
\end{proof}

For all $U\subset\mathcal{B}$ we denote $Sh(U)$ and call shadow of $U$ the set of points of $\mathcal{B}$ below $U$ such that
\begin{align*}
Sh(U):=\{p\in \mathcal{B}\mid  \exists V\in V\mathcal{B}\text{ containing }p\text{ and intersecting }U\text{ on a point }p'\text{ such that }h(p')\geq h(p)\}.
\end{align*}

For $S$ a subset of X, we shall call \textsl{large $Y$-horosphere} the subset $H_S$ defined by
\begin{align*}
H_S:=S\bowtie Y= (p^X)^{-1}(S).
\end{align*}

\begin{figure}
\begin{center}
\includegraphics[scale=1]{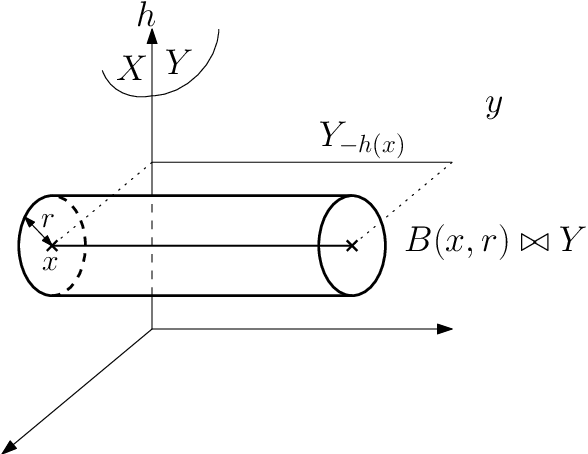} 
\end{center}
\caption{Large $X$-Horosphere in $X\bowtie Y$.}
\end{figure}

Let $M_0$ be the constant involved in assumption $(E2)$. Let us denote by $F^X\subset \mathcal{B}^X$ the subset

\begin{align*}
F^X:=\left\lbrace x\in \mathcal{B}^X\mid \lambda\left(Sh(H_{D_{M_0}(x)})\cap U^c\right)>\alpha^{\frac{1}{4}}\lambda\left(Sh(H_{D_{M_0}(x)})\right)\text{ and }h(x)\geq h^-(\mathcal{B}^X)+\frac{R}{2}\right\rbrace.
\end{align*}

The set $F^X$ is in bijection with the "bad" $Y$-horospheres $H$ that are above the middle of the box $\mathcal{B}$. By "bad" we mean the ones which have more than $\alpha^{\frac{1}{4}}$ fraction of the measure $\lambda$ of their shadow in $U^c$.

In the following lemma, we show that the shadow of almost all the $Y$-horospheres in the upper half of the box have almost full measure.

\begin{lemma}\label{LemmaGoodShadowForHorocycle}
There exists a constant $0<\alpha (\bowtie )\leq 1$ such that for all $0<\alpha\leq \alpha (\bowtie )$ the following statement holds. If $\lambda(U)\geq \left(1-\alpha\right)\lambda (\mathcal{B})$, then we have
\begin{align*}
\lambda^X(F^X)<\alpha^{\frac{1}{4}}\lambda^X\big(\mathcal{B}^X\big).
\end{align*}
\end{lemma}
\begin{proof}
Without loss of generality we can assume that $h(\mathcal{B})=[0;R[$. We proceed by contradiction, let us assume that $\lambda^X(F^X)\geq\alpha^{\frac{1}{4}}\lambda^X\big(\mathcal{B}^X\big)$. Therefore, there exists $z_0\in\left[\frac{R}{2},R\right[$ such that
\begin{align*}
\lambda^X_{z_0}(F_{z_0}^X)\geq \alpha^{\frac{1}{4}} \lambda^X_{z_0}\big(\mathcal{B}_{z_0}^X\big).
\end{align*}
Let $Z$ be a $2M_0$-maximal separating subset of $F_{z_0}^X$. Then we have
\begin{align*}
\lambda(U^c)\geq& \lambda \left(Sh\left(\bigsqcup\limits_{x\in Z}H_{D_{M_0}(x)}\right)\cap U^c\right)=\sum\limits_{x\in Z}\lambda \big(Sh\left(H_{D_{M_0}(x)}\right)\cap U^c\big),\quad\text{since this is a disjoint union},
\\\geq&\alpha^{\frac{1}{4}}\sum\limits_{x\in Z}\lambda \big(Sh\left(H_{D_{M_0}(x)}\right)\big)\asymp_{\bowtie} \alpha^{\frac{1}{4}}\sum\limits_{x\in Z} z_0\lambda_{z_0}\big(H_{D_{M_0}(x)}\big),\quad\text{ by definition of }F_{z_0}^X\text{ and Proposition \ref{ProprLambdaSameOnLevels}}.
\end{align*}
However $\lambda_{z_0}\big(H_{D_{M_0}(x)}\big)=\lambda_{z_0}^X(D_{M_0}(x))\lambda_{-z_0}^Y(\mathcal{B}^Y_{-z_0})$ since $H_{D_{M_0}(x)}=D_{M_0}(x)\times\mathcal{B}^Y_{-z_0}$, hence
\begin{align*}
\lambda(U^c)\succeq_{\bowtie}&\alpha^{\frac{1}{4}} z_0\sum\limits_{x\in Z}\lambda_{z_0}^X(D_{M_0}(x))\lambda_{-z_0}^Y(\mathcal{B}^Y_{-z_0}),
\\\asymp_{\bowtie}&\alpha^{\frac{1}{4}} z_0\lambda_{-z_0}^Y(\mathcal{B}^Y_{-z_0})\sum\limits_{x\in Z} \lambda^X_{z_0}\big(D_{2M_0}(x)\big),\quad\text{by Lemma \ref{LemmaMuDiskAreExp}},
\\\geq& \alpha^{\frac{1}{4}} z_0\lambda_{-z_0}^Y(\mathcal{B}^Y_{-z_0})\lambda^X_{z_0}\left(\bigcup\limits_{x\in Z}D_{2M_0}(x)\right)\geq \alpha^{\frac{1}{4}} z_0\lambda_{-z_0}^Y(\mathcal{B}^Y_{-z_0})\lambda^X_{z_0}\big(F_{z_0}^X\big) ,\quad\text{by definition of }Z,
\\\geq& \alpha^{\frac{1}{4}}\alpha^{\frac{1}{4}} z_0\lambda_{-z_0}^Y(\mathcal{B}^Y_{-z_0}) \lambda^X_{z_0}\big(\mathcal{B}_{z_0}^X\big) ,\quad\text{by assumption on }F_{z_0}^X,
\\\geq&\sqrt{\alpha}\frac{R}{2} \lambda_{-z_0}^Y(\mathcal{B}^Y_{-z_0})\lambda^X_{z_0}\big(\mathcal{B}_{z_0}^X\big)\asymp_{\bowtie} \frac{1}{2}\sqrt{\alpha}\lambda\big(\mathcal{B}\big),\quad\text{since }z_0\geq\frac{R}{2}\text{ and by Proposition \ref{PropMeasureGeodBoxHoro}},
\end{align*} 
which contradicts the assumptions on $U$ for $\alpha <\frac{1}{M(\bowtie)^2}$.
\end{proof}

The following lemma asserts that the projection on a level-set of almost all the $Y$-horospheres have almost full measure. 

\begin{lemma}\label{LemmaGoodPlanForHorocycle}
If $\lambda(U)\geq \left(1-\alpha\right)\lambda (\mathcal{B})$, then there exists a constant $M(\bowtie)$ such that for any large $Y$-horosphere $H_{D_{M_0}(x)}$ with $x\in \mathcal{B_x}\setminus F_X$ as in Lemma \ref{LemmaGoodShadowForHorocycle}, and for $1\geq M\rho\geq M^2\alpha^{\frac{1}{4}}>0$, there exists $P$ a level set of the height function in $\mathcal{B}$, such that
\begin{align*}
\lambda_{h(P)}(P\cap \mathrm{Sh}(H_{D_{M_0}(x)})\cap U^c)\preceq_{\bowtie}\alpha^{\frac{1}{4}}\lambda_{h(P)}(P\cap \mathrm{Sh}(H_{D_{M_0}(x)})).
\end{align*}
Furthermore, $P$ can be chosen such that $\rho R< d_{\bowtie}(P,H)<2\rho R$.
\end{lemma}
\begin{proof}
We proceed by contradiction, let us assume that such a plane $P$ does not exist. Computing the measure $\lambda$ of $\mathrm{Sh}(H_{D_{M_0}(x)})\cap U^c\cap \mathcal{B}_{[h(H)-2\rho R;h(H)-\rho R]}$ contradicts the fact that 
\begin{align*}
\lambda(\mathrm{Sh}(H_{D_{M_0}(x)})\cap U^c)\leq \alpha^{\frac{1}{4}}\lambda(\mathrm{Sh}(H_{D_{M_0}(x)})).
\end{align*}
Indeed, we show the contradiction using Lemma \ref{ProprLambdaSameOnLevels} and because we integrate on a sufficiently large portion of $[0,R]$ ($\rho \geq M \alpha ^{\frac{1}{4}}$).
\end{proof}

In the following lemma we show that almost all level-sets admit a point with large $X$-horospheres and $Y$-horospheres.  
 
\begin{lemma}\label{LemmaJustMarkov}
There exists a constant $0<\alpha (\bowtie )\leq 1$ such that for all $0<\alpha\leq \alpha (\bowtie )$ the following statement holds. Let $U\subset \mathcal{B}$ be such that $\lambda (U)\geq (1-\alpha) \lambda(\mathcal{B})$. Then there exists $U'\subset U$ such that:
\begin{enumerate}
\item $\lambda \big(U'\big)\geq \left(1-\alpha^{\frac{1}{4}}\right)\lambda\big(\mathcal{B}\big)$;
\item For all $z\in h(U')$ there exists $(x_{0,z},y_{0,z})\in  U'_z$ such that for all $(x_1,y_1)\in U'_z$, we have $(x_1,y_{0,z})\in U'_z$ and $(x_{0,z},y_1)\in U'_z$.
\end{enumerate}
\end{lemma}

\begin{proof}
We may assume without loss of generality that $h(\mathcal{B})=[0,R[$. Let us denote by
\begin{align*}
H_U:=\left\lbrace z\in [0,R[\ \vert\ \lambda_z\left(U_z\right)\geq \left(1-\alpha^{\frac{1}{4}}\right)\lambda_z\left(\mathcal{B}_{z}\right)  \right\rbrace.
\end{align*}
Then we claim that $\mathrm{Leb}(H_U)\geq \left(1-\alpha^{\frac{1}{4}}\right)R$. To prove this claim we proceed by contradiction. Let us assume that $\mathrm{Leb}(H_U)<\left(1-\alpha^{\frac{1}{4}}\right)R$, then $\mathrm{Leb}([0,R[\setminus H_U)\geq\alpha^{\frac{1}{4}}R$.
Furthermore, for all $z\in [0,R[\setminus H_U$ we have $\lambda_z\left(U_z\right)< \left(1-\alpha^{\frac{1}{4}}\right)\lambda_z\left(\mathcal{B}_{z}\right) $, hence
\begin{align}
\lambda_z\left(\mathcal{B}_z\setminus U_z\right)\geq \alpha^{\frac{1}{4}}\lambda_z\left(\mathcal{B}_{z}\right).\label{osk}
\end{align}
Therefore, by computing the measure of $\mathcal{B}\setminus U$ we have
\begin{align*}
\lambda \left( \mathcal{B}\setminus U\right)=&\int\limits_{z\in [0,R[}\lambda_z\left(\mathcal{B}_z\setminus U_z\right)\mathrm{d}z\geq \int\limits_{z\in ([0,R[\setminus H_U)}\lambda_z\left(\mathcal{B}_z\setminus U_z\right)\mathrm{d}z,
\\\geq& \int\limits_{z\in ([0,R[\setminus H_U)}\alpha^{\frac{1}{4}}\lambda_z\left(\mathcal{B}_{z}\right)\mathrm{d}z ,\quad\text{by inequality (\ref{osk})},
\\\succeq_{X}& \alpha^{\frac{1}{2}} \lambda (\mathcal{B}),\quad\text{by the contradiction assumption and Property \ref{ProprLambdaSameOnLevels}},
\end{align*}
which contradicts the assumption on $U$ for $\alpha$ small enough. Hence $\mathrm{Leb}(H_U)\geq \left(1-\alpha^{\frac{1}{4}}\right)R$.
\\Let us denote for $z\in[0;R[$
\begin{align*}
U^y&:=\left\lbrace x\in \mathcal{B}^X_{z} \ \vert\ (x,y)\in U\right\rbrace;
\\H&:=\left\lbrace z\in [0,R[\ \vert\ \exists y\in \mathcal{B}^Y_{-z}\ ,\  \lambda_z^X\left(U^y\right)\geq \left(1-\alpha^{\frac{1}{4}}\right)\lambda^X_z\left(\mathcal{B}^X_{z}\right)  \right\rbrace.
\end{align*}
In particular, for all $y\in \mathcal{B}^Y_{-z}$ we have  $U^y \subset U^X_{z}$, and by the definition of $\lambda$
\begin{align*}
\lambda(U)=\int\limits_{z\in [0,R[}\int\limits_{y\in \mathcal{B}^Y_{-z}}\lambda^X_z(U^y).
\end{align*}
We claim that $\mathrm{Leb}(H)\geq \left(1-\alpha^{\frac{1}{4}}\right)R$. To prove this claim, we also proceed by contradiction. Let us assume that $\mathrm{Leb}(H)<\left(1-\alpha^{\frac{1}{4}}\right)R$, then $\mathrm{Leb}([0,R[\setminus H)\geq \alpha^{\frac{1}{4}}R$. Furthermore for all $z\in [0,R[\setminus H$ we have that
\begin{align*}
\forall y\in \mathcal{B}^Y_{-z}\ ,\  \lambda_z^X\left(U^y\right)< \left(1-\alpha^{\frac{1}{4}}\right)\lambda^X_z\left(\mathcal{B}^X_{z}\right).
\end{align*}
Therefore, by the definition of $U_y$ we have that $\forall y\in \mathcal{B}^Y_{-z}$
\begin{align}
\lambda_z^X\left(\left\lbrace x\in \mathcal{B}^X_{z} \vert(x,y)\notin U\right\rbrace\right) \geq \alpha^{\frac{1}{4}} \lambda^X_z\left(\mathcal{B}^X_{z}\right).\label{PleaseLastLemma}
\end{align} 
Hence, by computing the measure of $\mathcal{B}\setminus U$ we have
\begin{align*}
\lambda \left( \mathcal{B}\setminus U\right)=&\int\limits_{z\in [0,R[}\int\limits_{y\in \mathcal{B}^Y_{-z}}\lambda^X_z\left(\left\lbrace x\in U^X_{z} \vert(x,y)\notin U\right\rbrace\right)\mathrm{d}\lambda^Y_{-z}\mathrm{d}z,
\\\geq& \int\limits_{z\in \big([0,R[\setminus H\big)}\int\limits_{y\in \mathcal{B}^Y_{-z}}\lambda^X_z\left(\left\lbrace x\in U^X_{z} \vert(x,y)\notin U\right\rbrace\right)\mathrm{d}\lambda^Y_{-z}\mathrm{d}z,
\\\geq& \int\limits_{z\in \big([0,R[\setminus H\big)}\int\limits_{y\in \mathcal{B}^Y_{-z}}\alpha^{\frac{1}{4}} \lambda^X_z\left(\mathcal{B}^X_{z}\right)\mathrm{d}\lambda^Y_{-z}\mathrm{d}z,\quad\text{by inequality (\ref{PleaseLastLemma})},
\\=& \alpha^{\frac{1}{4}}\int\limits_{z\in \big([0,R[\setminus H\big)}\lambda^Y_{-z}\left(\mathcal{B}^Y_{-z}\right) \lambda^X_z\left(\mathcal{B}^X_{z}\right)\mathrm{d}z,
\\\succeq_{X}& \alpha^{\frac{1}{4}}\alpha^{\frac{1}{4}}\lambda\left(\mathcal{B}\right)= \alpha^{\frac{1}{2}}\lambda\left(\mathcal{B}\right),\quad\text{by the contradiction assumption and Property \ref{ProprLambdaSameOnLevels}},
\end{align*}
which contradicts the assumption $\lambda \left( \mathcal{B}\setminus U\right)<\alpha\lambda\left(\mathcal{B}\right)$, for $\alpha< \frac{1}{M(\bowtie)^2}$. Let us denote for all $x\in \mathcal{B}^X_{z}$
\begin{align*}
U^x&:=\left\lbrace y\in \mathcal{B}^Y_{-z} \ \vert\ (x,y)\in U\right\rbrace;
\\H'&:=\left\lbrace z\in [0,R[\ \vert\ \exists x\in \mathcal{B}^X_{z}\ ,\  \lambda_{-z}^Y\left(U^x\right)\geq \left(1-\alpha^{\frac{1}{4}}\right)\lambda^Y_{-z}\left(\mathcal{B}^Y_{-z}\right)  \right\rbrace.
\end{align*}
We show similarly that $\mathrm{Leb}(H')\geq\left(1-\alpha^{\frac{1}{4}}\right)R$, therefore $\mathrm{Leb}(H\cap H'\cap H_U)\geq\left(1-3\alpha^{\frac{1}{4}}\right)R$.
\\\\For all $z\in H\cap H'$ there exists $(x_{0,z},y_{0,z})\in \mathcal{B}_z$ such that
\begin{align}
 \lambda_z^X\left(U^{y_0}\right)&\geq \left(1-\alpha^{\frac{1}{4}}\right)\lambda^X_z\left(\mathcal{B}^X_{z}\right) \label{osk2};
\\\lambda_{-z}^Y\left(U^{x_0}\right)&\geq \left(1-\alpha^{\frac{1}{4}}\right)\lambda^Y_{-z}\left(\mathcal{B}^Y_{-z}\right)\label{osk3}.
\end{align}
Let us define for all $z\in H_U\cap H\cap H'$, $U'_z:=\left(U^{x_{0,z}}\times U^{y_{0,z}}\right)$. Then we have:
\begin{enumerate}
\item $U'\subset U$.
\item $\lambda_z (U'_z)=\lambda_z \left(\left(U^{x_{0,z}}\times U^{y_{0,z}}\right)\cap U_z\right)\geq \left(1-3\alpha^{\frac{1}{4}}\right)\lambda_z(\mathcal{B})$ by inequalities (\ref{osk2}), (\ref{osk3}) and by the definition of $H_U$.
\item For all $(x_1,y_1)\in U'_z$ we have $(x_1,y_{0,z})\in U'_z$ and $(x_{0,z},y_1)\in U'_z$
\end{enumerate}
Let $(x_1,y_1)\in U'_z$, then  $(x_1,y_{0,z})\in U'$ hence $(x_{0,z},y_{0,z})\in U'$. Furthermore we have that $\mathrm{Leb}(H_U\cap H\cap H')\geq (1-3\alpha^{\frac{1}{4}})R$, hence $\mathrm{Leb}([0,R[\setminus (H_U\cap H\cap H'))\leq 3\alpha^{\frac{1}{4}}R$. Therefore
\begin{align*}
\lambda \left(\mathcal{B}\setminus U'\right)=&\int\limits_{z\in [0,R[}\lambda_z\left((\mathcal{B}\setminus U')_z\right)\mathrm{d}z ,
\\=& \int\limits_{z\in ([0,R[\setminus (H_U\cap H\cap H'))}\lambda_z\left(\mathcal{B}_z\setminus\left(U^{x_{0,z}}\times U^{y_{0,z}}\right)\right)\mathrm{d}z,
\\\leq& \int\limits_{z\in ([0,R[\setminus (H_U\cap H\cap H'))}\left(3\alpha^{\frac{1}{4}}\right)\lambda_z(\mathcal{B}_z)\mathrm{d}z ,\quad\text{by construction of }U'_z,
\\\preceq_{X}& 9\alpha^{\frac{1}{2}}\lambda (\mathcal{B}),\quad\text{ by the measure of }[0,R[\setminus (H_U\cap H\cap H')\text{ and by Property \ref{ProprLambdaSameOnLevels}}.
\end{align*}
Hence $\lambda (U')\geq \left(1-\alpha^{\frac{1}{4}}\right)\lambda(\mathcal{B})$, since $\alpha^{\frac{1}{4}}>9M(X)\alpha^{\frac{1}{2}}$ ($\alpha$ small enough in comparison to a constant depending only on $X$).
\end{proof}

These points $(x_{0,z},y_{0,z})$ will play a key role in the definition of the product map close to a given quasi-isometry in Theorems \ref{ThmPartI}.

\subsection{Divergence}

Two distinct vertical geodesics in a $\delta$-hyperbolic and Busemann space diverge quickly from each other. However this statement, based on Corollary 6.0.3, depends on the pair of geodesics. The next lemma aims at making this more precise for $X$ an admissible horo-pointed space.
More specifically we are going to look at a point $x$ and at all the vertical geodesics passing by a point of the disc centred at $x$ of radius $M_0$ (the $(E2)$ constant) along the horosphere at height $h(x)$, that is $VD_{M_0}(x)$. Let $V_0$ be a geodesic containing $x$, we want to quantify the vertical geodesics in $VD_{M_0}(x)$  which start diverging from the vertical geodesic $V_0$ between the heights  $h(x)-l$ and $h(x)+l$. We shall denote this set by $\mathrm{Div}(V_0)$:
\begin{align*}
\mathrm{Div}(V_0):=\{V\in VD_{M_0}(x)\mid|h_{\mathrm{Div}}(V_0,V)-h(x)|\leq l \}.
\end{align*}

\begin{lemma}\label{LemGeodDivergeOften}
With the above notations we have
\begin{align*}
\eta^X \big(VD_{M_0}(x)\backslash\mathrm{Div}(V_0)\big)\preceq_{X} e^{-ml}\eta^X\Big(VD_{M_0}(x)\Big).
\end{align*}
\end{lemma}

\begin{proof}

By slight abuse of notation, we may intersect a set of vertical geodesic segments \(E \subset V\mathcal{B}\) with a subset \(F \subset \mathcal{B}\). By this, we mean the intersection of \(F\) with the union of the images of the vertical geodesics of \(E\), where each image is a vertical geodesic in \(\mathcal{B}\). For example:
\begin{align*}
VD_{M_0}(x)\cap \mathcal{B}_{h(x)}=D_{M_0}(x).
\end{align*}
Any vertical geodesic segment $V\in VD_{M_0}(x)$ did not start to diverge from the vertical geodesic $V_0$ at the height $h(x)$, we have $h_{\mathrm{Div}}(V,V_0)\leq h(x)$. Therefore, all the vertical geodesic segments which did not start to diverge at the height $h(x)-l$, denoted by $VD_{M_0}(x)\backslash\mathrm{Div}(V_0)$, are still $M_0$-close to $\pi_{h(x)-l}(x)$:
\begin{align}
\Big(VD_{M_0}(x)\backslash\mathrm{Div}(V_0)\Big)\cap \mathcal{B}_{h(x)-l}\subset D_{M_0}(\pi_{h(x)-l}(x)).\label{IncluHere}
\end{align}
We use Lemma \ref{LemmaDisqueInShadow} with $z_0=h(x)$ and $z=h(x)-l$, which gives
\begin{align}
D_{2l-M_0}\big(\pi_{h(x)-l}\left(x\right)\big)\subset \pi_{h(x)-l}\left(D_{M_0}(x)\right)=VD_{M_0}(x)\cap \mathcal{B}_{h(x)-l}. \label{IncluHere2}
\end{align}
Therefore
\begin{align*}
\frac{\eta^X\big(VD_{M_0}(x)\backslash\mathrm{Div}(V_0)\big)}{\eta^X\big(VD_{M_0}(x)\big)}\asymp_X& \frac{\lambda^X_{h(x)-l}\big(VD_{M_0}(x)\backslash\mathrm{Div}(V_0)\cap \mathcal{B}_{h(x)-l}\big)}{\lambda^X_{h(x)-l}\big(VD_{M_0}(x)\cap \mathcal{B}_{h(x)-l}\big)},\quad\text{by Property \ref{PropMeasureGeodVertSameMeasureSpace}},
\\ \leq&\frac{\lambda^X_{h(x)-l}\left( D_{M_0}(\pi_{h(x)-l}(x))\right)}{\lambda^X_{h(x)-l}\big(VD_{M_0}(x)\cap \mathcal{B}_{h(x)-l}\big)},\quad\text{by inequality \ref{IncluHere}},
\\ \leq&\frac{\lambda^X_{h(x)-l}\left( D_{M_0}(\pi_{h(x)-l}(x))\right)}{\lambda^X_{h(x)-l}\left(D_{2l-M_0}\big(\pi_{h(x)-l}\left(x\right)\big)\right)},\quad\text{by inequality \ref{IncluHere2}}.
\end{align*}
Moreover by the definition of $\lambda^X$ and Lemma \ref{LemmaMuDiskAreExp}
\begin{align}
\frac{\lambda^X_{h(x)-l}\left( D_{M_0}(\pi_{h(x)-l}(x))\right)}{\lambda^X_{h(x)-l}\left(D_{2l-M_0}\big(\pi_{h(x)-l}\left(x\right)\big)\right)}=\frac{\mu^X_{h(x)-l}\left( D_{M_0}(\pi_{h(x)-l}(x))\right)}{\mu^X_{h(x)-l}\left(D_{2l-M_0}\big(\pi_{h(x)-l}\left(x\right)\big)\right)}\preceq_{X} e^{-ml}.
\end{align}
Therefore
\begin{align*}
\frac{\eta^X\big(VD_{M_0}(x)\backslash\mathrm{Div}(V_0)\big)}{\eta^X\big(VD_{M_0}(x)\big)}\preceq_{X} e^{-ml}.
\end{align*}
\end{proof}

Heuristically, the previous lemma asserts that most of the vertical geodesics segments passing close to a point $x$, start diverging from each other close to the height $h(x)$.

We now provide an estimate on the exponential contraction of the measure $\mu$ along the vertical direction. 

\begin{lemma}\label{LemmaProjCroitExp}
There exists $M(\bowtie)$ such that the following holds. Let $h_0\in \mathbb{R}$, let $U\subset(X\bowtie Y)_{h_0}$ be a measurable subset. Let $\Delta>M$ and let $A\subset (X\bowtie Y)_{h_0-\Delta}$ be a measurable subset. Suppose also that all vertical rays $V$ intersecting $U$ intersect $A$. Then
\begin{align*}
\mu_{h_0-\Delta}\big(A\big)\succeq_{\bowtie} e^{(m-n)\Delta}\mu_{h_0}\big(U\big).
\end{align*}
\end{lemma}

\begin{proof}
Since $\pi^{\bowtie}_{h_0-\Delta}(U)\subset A$ we have
\begin{align*}
\mu_{h_0-\Delta}\left(\pi^{\bowtie}_{h_0-\Delta}(U)\right)\leq \mu_{h_0-\Delta}\left(A\right),
\end{align*} 
where $\pi^{\bowtie}$ is defined in Notations \ref{DefProjHoro}. We recall that for all $x\in X$, $U_x^Y:=\lbrace y\in Y\mid (x,y)\in U\rbrace$. By definition
\begin{align}
\mu_{h_0}\big(U\big)&=\mu_{h_0}^X\otimes\mu_{-h_0}^Y\left(U\right)=\int\limits_{X_{h_0}}\mu_{-h_0}^Y\big(U^Y_x\big)\mathrm{d}\mu_{h_0}^X(x).\label{IneqUseInLemmaPommeIci}
\end{align}
For all $x\in U^X$ let us denote $U_x:=\lbrace (x,y)\in U \mid y \in U^Y\rbrace$, then
\begin{align*}
(U_x)^Y=U^Y_x:=\lbrace y\in Y\mid (x,y)\in U\rbrace.
\end{align*} 
Furthermore $U^Y_x\subset \pi_{-h_0}^Y \left[\pi^{Y}_{\Delta-h_0}\left(U^Y_x\right)\right]$, hence
\begin{align*}
\mu^Y_{-h_0}\left(U^Y_x\right)&\leq \mu^Y_{-h_0}\left(\pi_{-h_0}^Y \left[\pi^{Y}_{\Delta-h_0}\left(U^Y_x\right)\right]\right)\asymp_{\bowtie} e^{n\Delta} \mu_{\Delta-h_0}^Y\left[\pi^{Y}_{\Delta-h_0}\left(U^Y_x\right)\right],\quad\text{ by assumption (E3),}
\end{align*}
which gives us,
\begin{align}
\mu_{h_0}\big(U\big)&\preceq_{\bowtie} e^{n\Delta}\int\limits_{U^X}\mu_{\Delta-h_0}^Y\left[\pi^{Y}_{\Delta-h_0}\left(U^Y_x\right)\right]\mathrm{d}\mu_{h_0}^X(x),\quad\text{by definition of }\mu_{h_0}.\label{IneqUseInLemmaPommeIci3}
\end{align}
However we have
\begin{align}
\pi^Y_{\Delta-h_0}(U_x^Y)&=\left(\pi^{\bowtie}_{h_0-\Delta}(U_x)\right)^Y=\left(\pi^{\bowtie}_{h_0-\Delta}(U)\right)^Y_{\pi^X_{h_0-\Delta}(x)},\label{IneqUseInLemmaPommeIci2}
\\&=\left\lbrace y\in \left(\pi^{\bowtie}_{h_0-\Delta}(U)\right)^Y\mid(\pi^{X}_{h_0-\Delta}(x),y)\in \pi^{\bowtie}_{h_0-\Delta}(U)\right\rbrace.\nonumber
\end{align}
Hence
\begin{align*}
\mu_{h_0}\big(U\big)&\preceq_{\bowtie}  e^{n\Delta}\int\limits_{U^X}\mu_{\Delta-h_0}^Y\left[\left(\pi^{\bowtie}_{h_0-\Delta}(U)\right)^Y_{\pi^X_{h_0-\Delta}(x)}\right]\mathrm{d}\mu_{h_0}^X(x),\quad\text{by (\ref{IneqUseInLemmaPommeIci3}) and (\ref{IneqUseInLemmaPommeIci2})},
\\&= e^{n\Delta} \int\limits_{\pi^{X}_{h_0-\Delta}\left(U^X\right)}\mu_{\Delta-h_0}^Y\left[\left(\pi^{\bowtie}_{h_0-\Delta}(U)\right)^Y_{x'}\right]\mathrm{d}\pi^{X}_{h_0-\Delta}\ast \mu_{h_0}^X(x'),
\\&\asymp_{\bowtie}e^{n\Delta}e^{-m\Delta} \int\limits_{\pi^{X}_{h_0-\Delta}\left(U^X\right)}\mu_{\Delta-h_0}^Y\left[\left(\pi^{\bowtie}_{h_0-\Delta}(U)\right)^Y_{x'}\right]\mathrm{d}\mu_{h_0-\Delta}^X(x'),\quad\text{ by assumption (E3),}
\\&= e^{(n-m)\Delta} \mu_{h_0-\Delta}\left(\pi^{\bowtie}_{h_0-\Delta}(U)\right).
\end{align*}
Furthermore, as said at the beginning we have $\mu_{h_0-\Delta}\left(\pi^{\bowtie}_{h_0-\Delta}(U)\right)\leq \mu_{h_0-\Delta}\left(A\right)$, therefore
\begin{align*}
\mu_{h_0-\Delta}\big(A\big)\succeq_{\bowtie} e^{(m-n)\Delta}\mu_{h_0}\big(U\big).
\end{align*} 
\end{proof}

In the next Lemma we transfer a control on the measure $\mu$ to a control on the measure $\eta$.

\begin{lemma}\label{LemmaNoVertIfSmall}
Let $M_0$ be the constant involved in assumption $(E2)$, $\mathcal{B}$ be a box and $z\in h(\mathcal{B})$. Let $A\subset (\mathcal{B})_{z}$ and let $E\subset \mathcal{B}$ such that $h^+(E)\leq h(A)$. Then, if there exists $Q\geq 1$ such that $\mu \left(\mathcal{N}_{M_0}(E)\right)\leq Q^{-1}\mu \left(\mathcal{N}_{M_0}(A)\right)$, we have that
\begin{align*}
\eta\left(V\mathcal{N}_{M_0}(E)\right)\preceq_{\bowtie} Q^{-1}\eta\left(V\mathcal{N}_{M_0}(A)\right).
\end{align*} 
\end{lemma}
\begin{proof}
Let $Z\subset E$ be a $2M_0$-maximal separating set, we have:
\begin{enumerate}
\item The balls $ B(p,M_0)$ for $p\in Z$ are pairwise disjoint.
\item We have the following inclusions:
$$\bigsqcup\limits_{p\in Z} B(p,M_0)\subset \mathcal{N}_{M_0}(E)  \subset\bigcup\limits_{p\in Z} B(p,3M_0). $$
\end{enumerate} 
The radius $3M_0$ is required since we cover all $\mathcal{N}_{M_0}(E) $ and not only $E$. Furthermore, all balls and disks of radius $M_0$ have comparable measure $\mu$ by assumption $(E2)$ and Corollary (\ref{CoroMeasureNeighbourhood}), therefore
\begin{align}
\mu(\mathcal{N}_{M_0}(E))\asymp_{\bowtie} \#Z\asymp_{\bowtie}\sum\limits_{p\in Z}\mu(B(p,M_0))\asymp_{\bowtie}\sum\limits_{p\in Z}\mu_{h(p)}(D_{M_0}(p)).\label{AiSeLaFine2}
\end{align} 
Moreover, for all $v\in VE$, there exists $p\in Z$ such that $v\cap D_{3M_0}(p)\neq \emptyset$. Consequently we have $V\mathcal{N}_{M_0}(E)\subset \bigcup\limits_{p\in Z}VD_{3M_0}(p)$, hence
\begin{align*}
\eta\left(V\mathcal{N}_{M_0}(E)\right)&\leq \sum\limits_{p\in Z}\eta\left(VD_{3M_0}(p)\right)\asymp_{X}\sum\limits_{p\in Z}\lambda_{h(p)}\left(D_{3M_0}(p)\right),\quad\text{by Property \ref{PropMeasureGeodVertSameMeasureSpace}},
\\&\leq \sum\limits_{p\in Z}\lambda^X_{h(p)}\left(D_{6M_0}(p^X)\right)\lambda^Y_{-h(p)}\left(D_{6M_0}(p^Y)\right).
\end{align*}
Furthermore, disks of radius $r$ are included in rectangles of width $2r$, hence
\begin{align*}
\eta\left(V\mathcal{N}_{M_0}(E)\right)&\preceq_{\bowtie}\sum\limits_{p\in Z}e^{h(p)(m-n)}\mu_{h(p)}\left(D_{6M_0}(p)\right),\quad\text{by the definition of }\lambda_{h(p)},
\\&\leq e^{h(a)(m-n)}\sum\limits_{p\in Z}\mu_{h(p)}\left(D_{6M_0}(p)\right),\quad\text{because }h^+(E)\leq h(A),
\\&\preceq_{\bowtie} e^{h(a)(m-n)}\mu\left(\mathcal{N}_{M_0}(E)\right),\quad\text{by inequalities (\ref{AiSeLaFine2})}.
\end{align*}
Using similar arguments we obtain
\begin{align*}
\eta\left(V\mathcal{N}_{M_0}(A)\right)\asymp_{\bowtie} \lambda_{h(a)}\left(V\mathcal{N}_{M_0}(A)\right)\asymp_{\bowtie} e^{h(a)(m-n)}\mu\left(V\mathcal{N}_{M_0}(A)\right).
\end{align*}
Combined with the assumption $\mu \left(\mathcal{N}_{M_0}(E)\right)\leq Q^{-1}\mu \left(\mathcal{N}_{M_0}(A)\right)$ we have
\begin{align*}
\eta\left(V\mathcal{N}_{M_0}(A)\right)\succeq_{\bowtie} e^{h(a)(m-n)} Q\mu \left(\mathcal{N}_{M_0}(E)\right)\succeq_{\bowtie}  Q\eta\left(V\mathcal{N}_{M_0}(E)\right).
\end{align*}
\end{proof}

Heuristically, if a set $E$ is sufficiently small and below a set $A$, then the set of vertical geodesic segments intersecting $E$ will also be small. 

\section{Proof of the geometric rigidity}\label{ChapBigProof}

The aim of this chapter is to present a proof of our key result. Let $(X,Y)$ and $(X',Y')$ be two horopointed admissible couples, of parameter respectively $(m,n)$ and $(m',n')$. Let us assume that $m>n$ and $m'> n'$. 

\begin{thm}\label{ThmProd}
Let $\Phi:X\bowtie Y\to X'\bowtie Y'$ be a $(k,c)$ quasi-isometry, then there exist two quasi-isometries $\Phi^X:X\to X'$ and $\Phi^Y:Y\to Y'$ such that
\begin{align*}
d_{\bowtie}\Big(\Phi,\left(\Phi^X,\Phi^Y\right)\Big)\preceq_{k,c,\bowtie}1.
\end{align*}
\end{thm}

Although this statement is similar to the statement in the case of Sol and Diestel-Leader, our broader setting of admissible spaces requires additional key arguments, such as lemma \ref{LemmaYetToBeWritten}, and therefore relies heavily on the previous sections.
\\\\To make the exposition of the various statements in this chapter smoother, we made the following abuse of notation. In a statement, when a parameter, say $\theta$, needs to be sufficiently small, we will write it by "For $\theta\preceq_{\bowtie} 1$ we have ..." instead of "There exists a constant $M(\bowtie)$ such that if $\theta \leq \frac{1}{M}$, then ...". 
\\From now until the end of this chapter we consider $\Phi:X\bowtie Y\to X'\bowtie Y'$ a $(k,c)$-quasi-isometry with fixed constants $k\geq 1$ and $c\geq 0$.

\subsection{Vertical geodesics with $\varepsilon$-monotone image}

In order to construct a product map, the key idea is to use the quadrilateral lemmas of Section \ref{SecQuadLemma} on the image by the quasi-isometry $\Phi$ of a quadrilateral in $X\bowtie Y$. To do so we need to locate which vertical geodesic segments are sent close to vertical geodesic segments. Thanks to Proposition \ref{LemmaQuatrePtSix} it is sufficient to look for vertical geodesic segments with an $\varepsilon$-monotone image under $\Phi$, where $0\leq\varepsilon<1$ is a parameter to be determined later (depending on $\bowtie$, $k$ and $c$). We call \textit{good} these vertical geodesic segments.

\begin{nota}
We recall that we denote $V\mathcal{B}$ the set of vertical geodesic segments of the box $\mathcal{B}$. Let us denote by $V^g\mathcal{B}$ the set of good vertical geodesic segments and $V^b\mathcal{B}$ the set of bad vertical geodesic segments, that is
\begin{align*}
&V^g\mathcal{B}:=\big\lbrace \gamma \in V\mathcal{B}\mid \Phi\circ\gamma \text{ is }\varepsilon\text{-monotone}\big\rbrace;
\\&V^b\mathcal{B}:=\big\lbrace \gamma \in V\mathcal{B}\mid \Phi\circ\gamma \text{ is not }\varepsilon\text{-monotone}\big\rbrace=V\mathcal{B}\setminus V^g\mathcal{B}.
\end{align*}
\end{nota}

In the following Proposition, we prove the existence of an appropriate scale on which almost all boxes possess almost only good vertical geodesics. We shall denote by $\eta:=\eta_{V\mathcal{B}}$, $\eta^X:=\eta^X_{V\mathcal{B}^X}$ and $\eta^Y:=\eta^Y_{V\mathcal{B}^Y}$.

\begin{propo}\label{Lemma411}
For $0<\theta\preceq_{\bowtie}1$, there exist two positive constants $M(k,c,\bowtie,\varepsilon)$ and $M'(k,c,\bowtie)$ such that for all $r_0\geq M$, $N\geq \frac{M'}{\varepsilon}$ and $S\geq\frac{M'}{\varepsilon\theta^3}$ and boxes $\mathcal{B}$ at scale $L:=N^Sr_0$, there exist $k_0\in\lbrace 1,...,S\rbrace$, a box tiling $\bigsqcup\limits_{i\in I}\mathcal{B}_i= \mathcal{B}$ at scale $R=N^{k_0}r_0$ and $I_g\subset I$ such that:
\begin{enumerate}
\item $\lambda\left(\bigcup\limits_{i\in I_g}\mathcal{B}_i\right)\geq (1-\theta) \lambda\left(\mathcal{B}\right)$ (Boxes indexed by $I_g$ cover almost all $\mathcal{B}$);
\item $\forall i\in I_g$, $\dfrac{\eta_i\big(V^b\mathcal{B}_i\big)}{\eta_i\big(V\mathcal{B}_i\big)}\leq \theta$ (almost all vertical geodesic segments in $\mathcal{B}_i$ have $\varepsilon$-monotone image);
\end{enumerate}
where $\eta_i:=\eta_{V\mathcal{B}_i}$.
\end{propo}

\begin{proof}
We recall from Lemma \ref{LemmaExistGoodScale} the definition of $\delta_s(\alpha)$ for a quasi-geodesic segment $\alpha$.
\begin{equation}
A_s := \left\lbrace \alpha\left(\left[kN^sr_0,(k+1)N^{s}r_0\right]\right)\vert k\in\lbrace 0,...,N^{S-s}-1\rbrace\right\rbrace. \nonumber
\end{equation} 
Then $\delta_s(\alpha)$ is the proportion of segments in $A_s$ which are not $\varepsilon$-monotone:
\begin{equation}
\delta_s(\alpha):=\dfrac{\#\left\lbrace \beta\in A_s\vert\beta\text{ is not }\varepsilon\text{-monotone}\right\rbrace}{\# A_s}.
\end{equation}
Using Proposition \ref{LemmaExistGoodScale} on every vertical geodesic segment in $\mathcal{B}$ we have that $\forall \alpha \in V\mathcal{B}$
\begin{equation}
\sum\limits_{s=1}^{S}\delta_s(\alpha)\preceq_{\bowtie,k,c} \dfrac{1}{\varepsilon}.\label{RefDeCos}
\end{equation}
We now integrate the inequality (\ref{RefDeCos}) with respect to $\eta$ over $V\mathcal{B}$ to get
\begin{align*}
\dfrac{1}{\varepsilon}&\succeq_{\bowtie,k,c}\dfrac{1}{\eta\big(V\mathcal{B}\big)}\int\limits_{\alpha\in V\mathcal{B}}\left(\sum\limits_{s=1}^{S}\delta_s(\alpha)\right)\mathrm{d}\eta =\sum\limits_{s=1}^{S}\left( \dfrac{1}{\eta\big(V\mathcal{B}\big)}\int\limits_{\alpha\in V\mathcal{B}}\delta_s(\alpha)\mathrm{d}\eta\right).
\end{align*}
Consequently there exists $k_0\in\lbrace1,...,S\rbrace$ such that
\begin{equation}
\dfrac{1}{\eta\big(V\mathcal{B}\big)}\int\limits_{\alpha\in V\mathcal{B}}\delta_{k_0}(\alpha)d\eta\preceq_{\bowtie,k,c} \dfrac{1}{S\varepsilon}\preceq_{\bowtie} \theta^3,\quad\text{ by assumption on }S.\label{IneqTheta1ETNonTheta2Ici}
\end{equation}
From now on we denote $R:=N^{k_0}r_0$. There are $\frac{L}{R}$ layers of boxes at scale $R$ in  $\mathcal{B}$. We average $\delta_{k_0}(\alpha)$ along all $\alpha\in V\mathcal{B}$:
\begin{align}
\dfrac{1}{\eta\big(V\mathcal{B}\big)}\int\limits_{\alpha\in V\mathcal{B}}\delta_{k_0}(\alpha)\mathrm{d}\eta\nonumber=&\dfrac{1}{\eta\big(V\mathcal{B}\big)}\int\limits_{\alpha\in V\mathcal{B}}\dfrac{R}{L}\sum\limits_{k=0}^{\frac{L}{R}-1}\delta_{k_0}\big(\alpha([kR;(k+1)R])\big)\mathrm{d}\eta,\nonumber
\\=&\dfrac{1}{\eta\big(V\mathcal{B}\big)}\dfrac{R}{L}\sum\limits_{k=0}^{\frac{L}{R}-1}\int\limits_{\alpha\in V\mathcal{B}}\delta_{k_0}\big(\alpha([kR;(k+1)R])\big)\mathrm{d}\eta.\label{IneqTheta1ETNonTheta7Ici}
\end{align}
Let us denote by $\mathcal{B}_{[k]}:=\mathcal{B}\cap h^{-1}\big([kR; (k+1)R[\big)$ the $k$-th layer of $\mathcal{B}$. Since vertical geodesic segments of $X\bowtie Y$ are couples of vertical geodesic segments, $V\mathcal{B}_{[k]}$ is in bijection with $V\mathcal{B}_{[k]}^X\times V\mathcal{B}_{[k]}^Y$ which is itself in bijection with $\mathcal{B}_{kR}^X\times \mathcal{B}_{-(k+1)R}^Y $ as explained in Section \ref{SecVertSet}. Let us denote by $f$ this bijection.
\begin{align*}
f: \mathcal{B}_{[k]}&\to \mathcal{B}_{kR}^X\times \mathcal{B}_{-(k+1)R}^Y
\\\alpha &\mapsto \left(\alpha^X(kR),\alpha^Y(-(k+1)R)\right).
\end{align*}  
For all $\alpha\in V\mathcal{B}$ and for all $k\in \left\lbrace 0,...,\frac{L}{R}-1\right\rbrace$ we have $\delta_{k_0}\big(\alpha([kR;(k+1)R])\big)=0$ or $1$, hence
\begin{align*}
\delta_{k_0}\big(\alpha([kR;(k+1)R])\big)&= \mathbbm{1}_{V^b\mathcal{B}_{[k]}}\Big(\alpha([kR;(k+1)R])\Big),
\\&=\mathbbm{1}_{f\left(V^b\mathcal{B}_{[k]}\right)}\Big(\alpha_X(kR),\alpha_Y(-(k+1)R)\Big).
\end{align*}
Therefore
\begin{align}
&\int\limits_{\alpha\in V\mathcal{B}}\delta_{k_0}\big(\alpha([kR;(k+1)R]\big)\mathrm{d}\eta\nonumber
\\=&\int\limits_{\left(\alpha^X,\alpha^Y\right)\in V\mathcal{B}^X\times  V\mathcal{B}^Y}\mathbbm{1}_{f\left(V^b\mathcal{B}_{[k]}\right)}\Big(\alpha_X(kR),\alpha_Y(-(k+1)R)\Big)\mathrm{d}\eta^X\mathrm{d}\eta^Y,\nonumber
\\=&\int\limits_{\left(x,y\right)\in \mathcal{B}_0^X\times  \mathcal{B}_{-L}^Y}\mathbbm{1}_{f\left(V^b\mathcal{B}_{[k]}\right)}\Big(\pi^X_{kR}(x),\pi^Y_{-(k+1)R}(y)\Big)\mathrm{d}\lambda^X_{0}\mathrm{d}\lambda^Y_{-L},\quad\text{by definition }\eta^X\text{ and }\eta^Y,\nonumber
\\\asymp_{\bowtie}&\int\limits_{\left(x',y'\right)\in \mathcal{B}_{kR}^X\times  \mathcal{B}_{-(k+1)R}^Y}\mathbbm{1}_{f\left(V^b\mathcal{B}_{[k]}\right)}\Big(x',y'\Big)\mathrm{d}\lambda^X_{kR}\mathrm{d}\lambda^Y_{-(k+1)R},\quad\text{by Property \ref{ProprLambdaSameOnLevels}}.\label{UseIneqMeasurLayer1}
\end{align}
Let $\bigsqcup_{i\in I} \mathcal{B}_i$ be the box tiling at scale $R$ as in Proposition \ref{PropoTillBigBoxWithSmallBox}, and for all $k\in\{0,\ldots,\frac{L}{R}-1\}$ let us denote by $I_k\subset I$ the indices of the boxes $\mathcal{B}_i$ which tile $\mathcal{B}_{[k]}$. Then we have $V\mathcal{B}_{[k]}=\bigsqcup_{i\in I_k}V\mathcal{B}_i$ and $V^b\mathcal{B}_{[k]}=\bigsqcup_{i\in I_k}V^b\mathcal{B}_i$. Therefore for all $\left(x,y\right)\in \mathcal{B}_{kR}^X\times  \mathcal{B}_{-(k+1)R}^Y$
\begin{align*}
\mathbbm{1}_{f\left(V^b\mathcal{B}_{[k]}\right)}\big(x,y\big)=\mathbbm{1}_{f\left(\bigsqcup_{i\in I_k}V^b\mathcal{B}_i\right)}\big(x,y\big)=\sum\limits_{i\in I_k}\mathbbm{1}_{f\left(V^b\mathcal{B}_i\right)}\big(x,y\big).
\end{align*}
Hence from inequality (\ref{UseIneqMeasurLayer1}) we have
\begin{align*}
\int\limits_{\alpha\in V\mathcal{B}}\delta_{k_0}\big(\alpha([kR;(k+1)R]\big)\mathrm{d}\eta\asymp_{\bowtie}&\int\limits_{\left(x,y\right)\in \mathcal{B}_{kR}^X\times  \mathcal{B}_{-(k+1)R}^Y}\sum\limits_{i\in I_k}\mathbbm{1}_{f\left(V^b\mathcal{B}_i\right)}\big(x,y\big)\mathrm{d}\lambda^X_{kR}\mathrm{d}\lambda^Y_{-(k+1)R},
\\=&\sum\limits_{i\in I_k}\int\limits_{\left(x,y\right)\in \mathcal{B}_{kR}^X\times  \mathcal{B}_{-(k+1)R}^Y}\mathbbm{1}_{f\left(V^b\mathcal{B}_i\right)}\big(x,y\big)\mathrm{d}\lambda^X_{kR}\mathrm{d}\lambda^Y_{-(k+1)R},
\\=&\sum\limits_{i\in I_k}\int\limits_{\alpha\in V\mathcal{B}_i}\mathbbm{1}_{V^b\mathcal{B}_i}\big(\alpha\big)\mathrm{d}\eta_i =\sum\limits_{i\in I_k}\eta_i\left(V^b\mathcal{B}_i\right).
\end{align*}
In combination with equality (\ref{IneqTheta1ETNonTheta7Ici}) we have
\begin{align*}
\dfrac{1}{\eta\big(V\mathcal{B}\big)}\int\limits_{\alpha\in V\mathcal{B}}\delta_{k_0}(\alpha)\mathrm{d}\eta&\succeq_{\bowtie}  \dfrac{1}{\eta\big(V\mathcal{B}\big)}\dfrac{R}{L}\sum\limits_{k=0}^{\frac{L}{R}-1}\sum\limits_{i\in I_k}\eta_i\big(V^b\mathcal{B}_{i}\big),
\\&\succeq_{\bowtie} \sum\limits_{i\in I}\dfrac{R\eta_i\big(V\mathcal{B}_{i}\big)}{L\eta\big(V\mathcal{B}\big)} \dfrac{\eta_i\big(V^b\mathcal{B}_{i}\big)}{\eta_i\big(V\mathcal{B}_{i}\big)},
\\&\succeq_{\bowtie}\sum\limits_{i\in I}\dfrac{\lambda\big(\mathcal{B}_i\big)}{\lambda\big(\mathcal{B}\big)} \dfrac{\eta_i\big(V^b\mathcal{B}_{i}\big)}{\eta_i\big(V\mathcal{B}_{i}\big)},\quad\text{by Proposition \ref{PropMeasureGeodBoxHoro}}.
\end{align*}
Let us denote by $I_b$ the set of indices $i$ of boxes $\mathcal{B}_i$ such that $\dfrac{\eta_i\big(V^b\mathcal{B}_{i}\big)}{\eta_i\big(V\mathcal{B}_i \big)}\geq \theta$, and $I_g:=I\setminus I_b$. By definition, $I_g$ satisfies the second part of our proposition, we are left with proving that it also satisfies the first part. To do so we assume by contradiction that $\lambda\left(\bigcup\limits_{i\in I_b}\mathcal{B}_i\right)\geq \theta \lambda\left(\mathcal{B}\right)$, then
\begin{align*}
\dfrac{1}{\eta\big(V\mathcal{B}\big)}\int\limits_{\alpha\in V\mathcal{B}}\delta_{k_0}(\alpha)d\eta&\succeq_{\bowtie}\sum\limits_{i\in I_b}\dfrac{\lambda\big(\mathcal{B}_i\big)}{\lambda\big(\mathcal{B}\big)} \dfrac{\eta_i\big(V^b\mathcal{B}_{i}\big)}{\eta_i\big(V\mathcal{B}_{i}\big)},\quad\text{since }I_b\subset I,
\\&\succeq_{\bowtie} \theta\frac{\sum\limits_{i\in I_b}\lambda\big(\mathcal{B}_i\big)}{\lambda\big(\mathcal{B}\big)},\quad\text{by the definition of }I_b,
\\&\succeq_{\bowtie}\theta^2,\quad\text{by the contradiction assumption},
\end{align*}
which contradicts inequality (\ref{IneqTheta1ETNonTheta2Ici}) for $\theta\preceq_{\bowtie}1$. Therefore $\lambda\left(\bigcup\limits_{i\in I_b}\mathcal{B}_i\right)< \theta \lambda\left(\mathcal{B}\right)$, hence $\lambda\left(\bigcup\limits_{i\in I_g}\mathcal{B}_i\right)\geq (1-\theta) \lambda\left(\mathcal{B}\right)$.
\end{proof}
 
Let $\mathcal{B}$ be a box at scale $R$. Let us denote the upward and downward oriented vertical geodesic segments by
\begin{align*}
V^{\uparrow}\mathcal{B}:=&\left\lbrace V\in V^g\mathcal{B}\mid h\big(\Phi\circ V(0)\big)\leq h\big(\Phi\circ V(R)\big)\right\rbrace;
\\V^{\downarrow}\mathcal{B}:=&\left\lbrace V\in V^g\mathcal{B}\mid h\big(\Phi\circ V(0)\big)\geq h\big(\Phi\circ V(R)\big)\right\rbrace.
\end{align*}

We are now going to show that in a given box $\mathcal{B}_i$ with $i\in I_g$, almost all vertical geodesic segments share the same orientation. 
 
\begin{lemma}\label{LemmaDominantOrientation}
For $0<\varepsilon^2\preceq_{k,c,\bowtie}\theta\preceq_{k,c,\bowtie}1$, and for $R\succeq_{k,c,\bowtie} \frac{1}{\varepsilon}$ we have that if $\mathcal{B}$ is a box at scale $R$ such that $\eta\left(V^b\mathcal{B}\right)\leq\theta\eta\left(V\mathcal{B}\right)$, then one of the two following statements holds:
\begin{enumerate}
\item$\eta\left(V^{\uparrow}\mathcal{B}\cap V^g\mathcal{B}\right)\geq(1-3\sqrt{\theta})\eta\left(V\mathcal{B}\right)$;
\item$\eta\left(V^{\downarrow}\mathcal{B}\cap V^g\mathcal{B}\right)\geq(1-3\sqrt{\theta})\eta\left(V\mathcal{B}\right)$.
\end{enumerate}
\end{lemma}

In the proof, we first characterise a set of vertical geodesic segment whose images share the same orientation, then we show that this set has almost full measure.

\begin{proof}
Without loss of generality we can assume that $h(\mathcal{B})=[0,R[$. Let us denote by
\begin{align*}
G^Y\left(v^X\right):=&\left\lbrace v^Y\in V\mathcal{B}^Y\mid \left(v^X,v^Y\right)\in V^g\mathcal{B}\right\rbrace;
\\G^X:=&\left\lbrace v^X\in V\mathcal{B}^X\mid \eta^Y \left(G^Y\left(v^X\right)\right)\geq \big(1-\sqrt{\theta}\big)\eta^Y\left(V\mathcal{B}^Y\right)\right\rbrace.
\end{align*}
By construction we have 
\begin{align*}
\bigcup _{v^X\in V\mathcal{B}^X} G^Y\left(v^X\right) =& \left(V^g\mathcal{B}\right)^Y.
\end{align*}
Applying Lemma \ref{LemmaMarkov1} with $V_1:=V^g(\mathcal{B})$ and $\alpha:=\theta$ we get
\begin{align}
\eta^X\left(G^X\right)\geq \big(1-\sqrt{\theta}\big)\eta^X\left(V\mathcal{B}^X\right).\label{IneqPP1PP2PP12}
\end{align}
Let $v_1^X:[0,R]\to X$ and $v_2^X:[0,R]\to X$ be two vertical geodesic segments of $G^X$, then
\begin{align*}
\eta^Y \left(G^Y\left(v_1^X\right)\right)\geq \big(1-\sqrt{\theta}\big)\eta^Y\left(V\mathcal{B}^Y\right);
\\\eta^Y \left(G^Y\left(v_2^X\right)\right)\geq \big(1-\sqrt{\theta}\big)\eta^Y\left(V\mathcal{B}^Y\right).
\end{align*}
Hence 
\begin{align}
\eta^Y\left(G^Y\big(v_1^X\big)\cap G^Y\big(v_2^X\big)\right)\geq \big(1-2\sqrt{\theta}\big)\eta^Y\left(V\mathcal{B}^Y\right).\label{bobob}
\end{align}
Let $v_1^Y,v_2^Y\in G^Y\left(v_1^X\right)\cap G^Y\left(v_2^X\right)$ and let us denote by $V_{i,j}:=(v_i^X,v_j^Y)$ with $i,j=1,2$. By definition of $v_1^Y$ and $v_2^Y$, the quasigeodesic segments $\Phi\left(V_{i,j}\right)$ are $\varepsilon$-monotone. 

Two cases occur. As a first case let us assume that
\begin{align*}
d_X\left(v_1^X(0),v_2^X(0)\right)&>\sqrt{\theta}R;
\\d_Y\left(v_1^Y(0),v_2^Y(0)\right)&>\sqrt{\theta} R.
\end{align*} 
Let $M$ be the constant involved in Proposition \ref{LemmaSecondTetrahedronConfiguration}. For $R\geq 4kc$ and $\varepsilon\leq \frac{\sqrt{\theta}}{20kM}$ we have that $\sqrt{\theta}R\geq 10kM\varepsilon R +2kc$, hence we can apply Proposition \ref{LemmaSecondTetrahedronConfiguration} on $V_{1,1}$ and $V_{2,2}$, which gives us that they share the same orientation.
\\\\The second case, that is when either $d_X\left(v_1^X(0),v_2^X(0)\right)\leq \sqrt{\theta} R$ or $d_Y\left(v_1^Y(0),v_2^Y(0)\right)\leq \sqrt{\theta} R$, is treated thanks to an auxiliary geodesic segment. Hence without loss of generality we focus on the case $d_X\left(v_1^X(0),v_2^X(0)\right)\leq \sqrt{\theta} R$ and consider a geodesic segment $v_3^X\in G^X$ verifying $d_X\left(v_1^X(0),v_3^X(0)\right)> \sqrt{\theta} R$ and $d_X\left(v_2^X(0),v_3^X(0)\right)> \sqrt{\theta} R$. 
To prove its existence, we consider the measure of
\begin{equation}
G^X\setminus V_{\mathcal{B}^X}\Big(D_{\sqrt{\theta} R}\Big(v_1^X(0)\Big)\cup D_{\sqrt{\theta} R}\Big(v_2^X(0)\Big)\Big).
\end{equation} 
Let $M_0$ be the constant of assumption $(E2)$. By Lemma \ref{LemmaMuDiskAreExp} we have for all $r_1\geq r_2>M_0$ and for all $x\in X_0$ that $\mu_0(D_{r_1}(x))\asymp_{\bowtie}e^{m\frac{|r_1-r_2|}{2}}\mu_0(D_{r_2}(x))$, therefore
\begin{align}
\lambda_{0}\left(D_{\sqrt{\theta} R}\Big(v_1^X(0)\Big)\right)\preceq_{\bowtie}& e^{m\frac{\sqrt{\theta }R-R}{2}}\lambda_{0}\left(D_{R}\Big(v_1^X(0)\Big)\right)\leq e^{-m\frac{R}{4}}\lambda_{0}\left(D_{R}\Big(v_1^X(0)\Big)\right),\quad\text{since }\theta\leq\frac{1}{4}.\label{bobibobibob}
\end{align}
Furthermore, by Lemma \ref{LemmaDisqueInShadow} the bottom of $\mathcal{B}$ contains a disk of radius $2R-M_0$, hence by Lemma \ref{LemmaMuDiskAreExp} we have $\eta^X\left(V\mathcal{B}^X\right)\asymp_{X}\lambda_0\left(D_{2R}(v_1^X(0))\right)$. Combined with inequality (\ref{bobibobibob}) we have
\begin{align*}
\lambda_{0}\left(D_{\sqrt{\theta} R}\Big(v_1^X(0)\Big)\right)\preceq_{\bowtie}e^{-m\frac{R}{4}}\eta^X\left(V\mathcal{B}^X\right).
\end{align*}
The same formula holds for $v_2^X$ instead of $v_1^X$. By inequality (\ref{IneqPP1PP2PP12}) we have that
\begin{align*}
\eta^X\left(G^X\right)\geq(1-\sqrt{\theta})\eta^X\left(V\mathcal{B}^X\right)\geq \frac{1}{2}\eta^X\left(V\mathcal{B}^X\right),
\end{align*}
hence there exists $M(\bowtie)$ such that
\begin{align*}
\eta^X\left(G^X\setminus V_{\mathcal{B}^X}\left(D_{\sqrt{\theta} R}\Big(v_1^X(0)\Big)\cup D_{\sqrt{\theta} R}\Big(v_2^X(0)\Big)\right)\right)\geq&\left(\frac{1}{2}-2Me^{-m\frac{R}{4}}\right)\eta^X\left(V\mathcal{B}^X\right),
\\>&0,\quad\text{for }R\geq \frac{4}{m}\ln(4M+1).
\end{align*} 
Therefore there exists $v_3^X\in G^X$ such that
\begin{align*}
d_X\left(v_1^X(0),v_3^X(0)\right)>& \sqrt{\theta} R;
\\d_X\left(v_2^X(0),v_3^X(0)\right)>& \sqrt{\theta} R.
\end{align*} 
Applying twice Lemma \ref{LemmaSecondTetrahedronConfiguration}, first on $V_{1,1}$ and $V_{3,3}$, then on $V_{2,2}$ and $V_{3,3}$, we get that the  $\Phi\big(V_{1,1}\big)$ has the same orientation as $\Phi\big(V_{3,3}\big)$ which has the same orientation as $\Phi\big(V_{2,2}\big)$. Therefore $\Phi\big(V_{1,1}\big)$ and $\Phi\big(V_{2,2}\big)$ share the same orientation.  
\\\\Let us fix $v_0^X\in G^X$ and $v_0^Y\in G^Y\left(v_0^X\right)$. Then the image of every vertical geodesic segment $V\in\bigcup\limits_{v^X\in G^X}\left\lbrace v^X\right\rbrace\times \left(G^Y\left(v_0^X\right)\cap G^Y\left(v^X\right)\right)$ shares the same orientation as the image of $\left(v_0^X,v_0^Y\right)$. Furthermore
\begin{align*}
\eta \left(\bigcup\limits_{v^X\in G^X}\left\lbrace v^X\right\rbrace\times \left(G^Y\left(v_1^X\right)\cap G^Y\left(v^X\right)\right)\right)=&\int\limits_{v^X\in G^X}\eta^Y\left(G^Y\left(v_1^X\right)\cap G^Y\left(v^X\right)\right)d\eta^X,
\\\geq&\int\limits_{v^X\in G^X}\big(1-2\sqrt{\theta}\big)\eta^Y\left(V\mathcal{B}^Y\right)d\eta^X,\quad\text{by inequality (\ref{bobob})},
\\=&\big(1-2\sqrt{\theta}\big)\eta^Y\left(V\mathcal{B}^Y\right)\eta^X\left(G^X\right)
\\\geq&\big(1-2\sqrt{\theta}\big)\eta^Y\left(V\mathcal{B}^Y\right)\big(1-\sqrt{\theta}\big)\eta^X\left(V\mathcal{B}^X\right),\quad\text{by inequality (\ref{IneqPP1PP2PP12})},
\\\geq&\big(1-3\sqrt{\theta}\big)\eta\left(V\mathcal{B}\right),
\end{align*}
which proves the lemma.
\end{proof}

\subsection{Factorisation of a quasi-isometry in small boxes}\label{SecPartI}

Proposition \ref{Lemma411} gives us two scales $R$ and $L$ such that all boxes at scale $L$ can be tiled with boxes at scale $R$. Moreover, almost all of them, that is the $\mathcal{B}_i$ for $i\in I_g$, contained almost only vertical geodesic segments with $\varepsilon$-monotone image under $\Phi$. 
\\\\A map $f:X\bowtie Y\to X'\bowtie Y'$ is called a \textbf{product map} if there exist two maps $f^X$ and $f^Y$ such that one of the two following holds:
\begin{enumerate}
\item We have $f^X:X\to X'$, $f^Y:Y\to Y'$ and $\forall p=(p^X,p^Y)\in X\bowtie Y$, $f(p)=\left(f^X\left(p^X\right),f^Y\left(p^Y\right)\right)$.
\item We have $f^X:X\to Y'$, $f^Y:Y\to X'$ and $\forall p=(p^X,p^Y)\in X\bowtie Y$, $f(p)=\left(f^Y\left(p^Y\right),f^X\left(p^X\right)\right)$.
\end{enumerate}  
In particular, when we denote by $\left(f^X,f^Y\right)$ a product map on a horospherical product, it implies that when $h(x)+h(y)=0$, we have $h\left(f^X(x)\right)+h\left(f^Y(y)\right)=0$. Therefore a product map is height respecting.

\begin{thm}\label{ThmPartI}
For $0<\theta\leq\varepsilon\preceq_{\bowtie} 1$, $r_0\succeq_{\bowtie}\frac{\sqrt[\varepsilon]{2}}{\varepsilon}$, $N\succeq_{\bowtie} 1$ and for $S\succeq_{\bowtie} \dfrac{1}{\varepsilon \theta^2}$, we have that for any $i\in I_g$, there exists a product map $\hat{\Phi}_i$, and $U_i'\subset \mathcal{B}_i$ such that:
\begin{enumerate}
\item $\lambda(U_i')\geq \left(1-\theta^{\frac{1}{8}}\right)\lambda\left(\mathcal{B}_i\right)$;
\item For all $(x,y)\in U_i'$, $d_{\bowtie'}\left(\Phi(x,y),\hat{\Phi}_i(x,y)\right)\preceq_{k,c,\bowtie} \varepsilon R$.
\end{enumerate}
In particular we have $\Delta h\left(\Phi(x,y),\hat{\Phi}_i(x,y)\right)\preceq_{k,c,\bowtie}\varepsilon R$.
\end{thm}

This proposition corresponds to Proposition 4.14 of \cite{EFW1}.

Since almost all the points in a good box are surrounded by almost only good vertical geodesic segment (Lemma \ref{LemmaCloseGeodAreGood}), we show that given two points sharing the same $X$ coordinates, we can almost always construct a quadrilateral verifying the hypotheses of Proposition \ref{LemmaTetrahedron}.

\begin{lemma}\label{LemmaExistsPointForTetra}
Let $M_0$ be the constant of assumption $(E2)$. For $0<\theta\preceq_{\bowtie}1$ and for $R\succeq_{\bowtie}\frac{1}{\theta}$, let $\mathcal{B}$ be a box at scale $R$ of $X\bowtie Y$. Let us assume the existence of a subset $U$ of $\mathcal{B}$ such that:
\begin{enumerate}
\item[$(a)$] $\lambda(U)\geq  (1-\theta)\lambda(\mathcal{B})$;
\item[$(b)$] For all $x\in U$, $\eta\left(V^b_{\mathcal{B}}\big(D_{M_0}(x)\big)\right)\leq \sqrt{\theta}\eta\left(V_{\mathcal{B}}\big(D_{M_0}(x)\big)\right)$.
\end{enumerate}
Then we have:
\begin{enumerate}
\item For all $a_1,a_2\in U$ such that $a_1^X=a_2^X$, there exist $b_1,b_2\in\mathcal{B}$ and four vertical geodesic segments $\gamma_{i,j}$ linking $a_i$ to $b_j$ such that $a_1$, $a_2$, $b_1$ and $b_2$ form a coarse vertical quadrilateral with nodes of scale $D=\theta R$, meaning that the configuration verifies the assumptions of Proposition \ref{LemmaTetrahedron}.
\item For $i,j\in\{1,2\}$, $\Phi(\gamma_{i,j})$ has $\varepsilon$-monotone image under $\Phi$. 
\end{enumerate} 
\end{lemma}

By Lemma \ref{LemmaCloseGeodAreGood}, the boxes $\mathcal{B}_i$, with $i\in I_g$, verify the assumptions of this Lemma. Moreover, we recall that a vertical quadrilateral satisfy the assumptions of Proposition \ref{LemmaTetrahedron}. 

\begin{proof}[Proof of Lemma \ref{LemmaExistsPointForTetra}]
Let $M_0$ be the constant of assumption $(E2)$. Let $a_1,a_2\in U$. For $i\in\{1,2\}$ let us denote $VD_i:=V_{\mathcal{B}}\left(D_{M_0}(a_i)\right)$ and $V^b D_i:=V_{\mathcal{B}}^b\left(D_{M_0}(a_i)\right)$. 
For all $v=(v^X,v^Y)\in V_{\mathcal{B}}$ and all $i\in\{1,2\}$ let us denote by:
\begin{enumerate}
\item $E_i^Y(v^X):=\{v^Y\in VD_i^Y\mid (v^X,v^Y)\in V^b D_i\}$;
\item $F_i^X:=\left\lbrace v^X\in VD_i^X \mid \eta^Y\big(E_i^Y(v^X)\big)\geq \theta^{\frac{1}{4}}\eta^Y(VD_i^Y)\right\rbrace$ .
\end{enumerate}
Thanks to Lemma \ref{LemmaMarkov2}, applied with $V_1:= V^b\mathcal{B}$, $\alpha:=\sqrt{\theta}$ and $a=a_i$, we have that
\begin{equation}
\eta^X\left(F_i^X\right)<\theta^{\frac{1}{4}}\eta^X\left(VD_i^X\right).
\end{equation}
Let us take $a_1$ and $a_2$ in $U$ such that $a_1^X=a_2^X$, then $VD^X_1=VD_2^X$:
\begin{enumerate}
\item $\eta^X(VD_i^X\backslash (F_1^X\cup F_2^X))\geq (1-2\theta^{\frac{1}{4}})\eta^X(VD_i^X)$;
\item For all $v^X\in VD_i^X\backslash (F_1^X\cup F_2^X)$ and $i\in\{1,2\}$ we have $\eta^Y\big(E_i^Y(v^X)\big)< \theta^{\frac{1}{4}}\eta^Y(VD_i^Y)$.
\end{enumerate}
The sets $VD_i^X\backslash (F_1^X\cup F_2^X)$ enclose the vertical geodesic segments in $\mathcal{B}^X$ passing close to $a_1^X=a_2^X$ such that almost all the induced vertical geodesic segments around $a_1$ and $a_2$ in $\mathcal{B}$ are good (ie. have $\varepsilon$-monotone images under the quasi-isometry $\Phi$). 
\\Since we have a sufficient proportion of good vertical geodesic segments, we will be able to find several of them that intersect the same neighbourhood in two different points sufficiently far from each other. If $h(a_1^X)< \theta R$, the construction of the quadrilateral of Proposition \ref{LemmaTetrahedron} with $D=\theta R$ is straightforward since the four points $a_1$, $a_2$, $b_1$ and $b_2$ would be $\theta R$ close, hence without loss of generality we may assume that $h(a_1^X)\geq \theta R$. Moreover, as we did before we can also suppose that $h(\mathcal{B})=[0,R[$.
\\We apply Lemma \ref{LemmaDisqueInShadow} with $z_0=h(a_1)$ and $z=h(a_1)-\theta R$ to get the following inclusions:
\begin{align}
D^X_{2\theta R-M_0}\big(\pi_{h(a_1)-\theta  R}\left(a_1^X\right)\big)\subset \pi_{h(a_1)-\theta  R}\left(D_{M_0}\left(a_1^X\right)\right) \subset D^X_{2\theta R+M_0}\big(\pi_{h(a_1)-\theta  R}\left(a_1^X\right)\big).\label{Chuiout}
\end{align}
We now suppose by contradiction that any couple of good vertical geodesic segments does not diverge quickly. This means that they stay $M_0$-close until they attain a height lower than $h\left(a_1^X\right)-\theta R$. Therefore
\begin{align*}
\pi_{h(a_1)-\theta R}\left(VD_i^X\backslash (F_1^X\cup F_2^X)\right)&\subset D^X_{M_0}\left(\pi_{h(a_1)-\theta R}\left(a_1^X\right)\right).
\end{align*} 
Thanks to the inclusions (\ref{Chuiout}) we have $VD^X_{2\theta R-M_0}(\pi_{h(a_1)-\theta  R}(a_1^X))\subset VD_1^X $, hence, combined with Property \ref{PropMeasureGeodVertSameMeasureSpace} we obtain
\begin{align*}
\frac{\eta^X\left(VD_1^X\backslash (F_1^X\cup F_2^X)\right)}{\eta^X(VD_1^X)}\preceq_{\bowtie}&\frac{\lambda^X_{h(a_1)-\theta R}\left(D_{M_0}\left(\pi_{h(a_1)-\theta R}\left(a_1^X\right)\right)\right)}{\lambda^X_{h(a_1)-\theta R}\left(D_{2\theta R}\left(\pi_{h(a_1)-\theta R}\left(a_1^X\right)\right)\right)},
\\\preceq_{\bowtie}&e^{\frac{m(M_0-2\theta R)}{2}},\quad\text{by Lemma \ref{LemmaMuDiskAreExp}},
\end{align*}
which, for $R$ large enough in comparison to $\frac{1}{\theta}$, contradicts the fact that $\eta^X(VD_1^X\backslash (F_1^X\cup F_2^X))\geq (1-2\theta^{\frac{1}{4}})\eta^X(VD_1^X)$, the first conclusion of the previously used Lemma \ref{LemmaMarkov2}. Hence there exists a couple of vertical geodesic segments $V^X_1$ and $V^X_2$ of $VD_i^X\backslash (F_1^X\cup F_2^X)$ diverging quickly from each other. Furthermore we have $\eta^Y\big(E_i^Y(v^X)\big)< \theta^{\frac{1}{4}}\eta^Y(VD_i^Y)$, hence there exists segments $V_1^Y$ and $V_2^Y$ such that $(V_1^X,V_1^Y)\in V^g_{\mathcal{B}}(D_M(a_1))$ and $(V_2^X,V_2^Y)\in V^g_{\mathcal{B}}(D_M(a_2))$. 
\\Let us define $b_i^X=V^X_i\left(h(a_1)-\frac{1}{2}d\left(a_1^X,a_2^X\right)\right)$, so that $b_1^X$ and $b_2^X$ are at the height where $V_1^X$ and $V_2^X$ diverge. Similarly, let us define $b_1^Y=b_2^Y=V^Y_1\left(-h(a_1)+\frac{1}{2}d\left(a_1^X,a_2^X\right)\right)$ such that $V_1^Y$ and $V_2^Y$ diverge, and $\gamma_{ij}=\left(V_{i}^X,V_{j}^Y\right)$ to ensure that the vertical geodesic segments of the quadrilateral $\gamma_{11}\cup\gamma_{12}\cup\gamma_{22}\cup\gamma_{21}$ have close endpoints. Furthermore by construction, they diverge from each other and have $\varepsilon$-monotone image under $\Phi$. 
\end{proof}

In the next proofs, we will be using Proposition \ref{LemmaQuatrePtSix} on each of the four images $\Phi(\gamma_{ij})$, which will provide us with a new quadrilateral $(\varepsilon+\theta)R$ close to $\Phi\left(\gamma_{11}\cup\gamma_{12}\cup\gamma_{22}\cup\gamma_{21}\right)$ on which the assumptions of Lemma \ref{LemmaTetrahedron} are verified.
\\\\Finally we deduce that on a good box, the quasi-isometry $\Phi$ is close to a product map.

\begin{proof}[Proof of Theorem \ref{ThmPartI}]
Let $i\in I_g$ and $\mathcal{B}_i$ a good box (defined in Lemma \ref{Lemma411}). Then following Lemma \ref{Lemma411}, we have $\eta_i(V^b\mathcal{B}_i)\leq \theta\eta_i(V\mathcal{B}_i)$. Therefore by Lemma \ref{LemmaDominantOrientation}, one of the two following statements hold:
\begin{enumerate}
\item$\eta\left(V^{\uparrow}\mathcal{B}\cap V^g\mathcal{B}\right)\geq(1-3\sqrt{\theta})\eta\left(V\mathcal{B}\right)$;
\item$\eta\left(V^{\downarrow}\mathcal{B}\cap V^g\mathcal{B}\right)\geq(1-3\sqrt{\theta})\eta\left(V\mathcal{B}\right)$.
\end{enumerate}
Let us first assume that the dominant orientation is upward. Let us choose $V_1=V\mathcal{B}\setminus \left(V^{\uparrow}\mathcal{B}\cap V^g\mathcal{B}\right)$, the vertical geodesics which have neither dominant orientation nor $\varepsilon$-monotone image by $\Phi$. By Lemma $\ref{LemmaCloseGeodAreGood}$, used with $\alpha:=\theta^2$, we have that there exists $U_i\subset\mathcal{B}_i$ such that:
\begin{enumerate}
\item $\lambda(U_i)\geq (1-\sqrt{\theta})\lambda(\mathcal{B}_i)$;
\item For $p\in U_i$ we have $\eta\left(V_1\big(D_{M_0}(x)\big)\right)<\eta\left(V\mathcal{B}\big(D_{M_0}(x)\big)\right)\sqrt{\theta}$.
\end{enumerate}
Let us apply Lemma \ref{LemmaJustMarkov}, with $U:=U_i$ and $\alpha:=\sqrt{\theta}$, then there exists $U'\subset U_i$ of almost full measure such that $\forall z\in h(U')$, $\exists (x_{0,z},y_{0,z})\in U'_z$ such that $\forall(x_1,y_1)\in U'_z$, we have $(x_1,y_{0,z})\in U'$ and $(x_{0,z},y_1)\in U'$.
Let $a,a_0\in U'$ such that $a^X=a_0^X$. By Lemma \ref{LemmaExistsPointForTetra} applied on $a_0$ and $a$, there exist $b_1,b_2\in\mathcal{B}_i$ and four vertical geodesics $V_{ij}$ in $V^{\uparrow}\mathcal{B}\cap V^g\mathcal{B}$ such that $b_1$ and $b_2$ form a coarse vertical quadrilateral $T$ with $a_0$ and $a$, where $V_{ij}$ are the edges of $T$. Proposition \ref{LemmaQuatrePtSix} gives a constant $M(k,c,\bowtie)$ and four vertical geodesic segments $M \varepsilon R$-close to the four sides of $\Phi(T)$. Furthermore we assumed that the dominant orientation is upward, hence the images of the four sides are all upward oriented. Hence thanks to Proposition \ref{LemmaTetrahedron} we get
\begin{align*}
d_{X'}\left(\Phi(a_0)^{X'},\Phi(a)^{X'}\right)\preceq_{k,c,\bowtie}\varepsilon R.
\end{align*}
Then for all $a\in U'$ such that $a^X=a_0^X$
\begin{align}
d_{X'}\left(\Phi(a_0)^{X'},\Phi(a)^{X'}\right)\preceq_{k,c,\bowtie}\varepsilon R.\label{IneqFinThm1}
\end{align}
We show similarly that for all $a\in U'$ such that $a^Y=a_0^Y$ we have
\begin{align}
d_{Y'}\left(\Phi(a_0)^{Y'},\Phi(a)^{Y'}\right)\preceq_{k,c,\bowtie}\varepsilon R.\label{IneqFinThm2}
\end{align}
Let us define the product map  $\hat{\Phi}_i:=\left(\hat{\Phi}_i^X,\hat{\Phi}_i^Y\right):X\bowtie Y\to X'\bowtie Y'$. For all $z\in h(U')$, let $\left(x_{0,z},y_{0,z}\right)\in U'_z$ be the points involved in Lemma \ref{LemmaJustMarkov}, and for all $z\in [0,R[\setminus h(U')$, let us fix an arbitrary point $(x_{0,z},y_{0,z})\in (\mathcal{B}_i)_z$. We can therefore define for all $x\in X$
\begin{align*}
\hat{\Phi}_i^X(x):=V^{X'}_{\Phi(x,y_{0,z})}\big(h\circ\Phi(x_{0,z},y_{0,z})\big).
\end{align*}
Then for all $(x,y)\in U'$ the triangle inequality gives
\begin{align}
&d_{X'}\left(\hat{\Phi}_i^X(x),\Phi(x,y)^{X'}\right)= d_{X'}\left(V^{X'}_{\Phi(x,y_{0,z})}\big(h\circ\Phi(x_{0,z},y_{0,z})\big),\Phi(x,y)^{X'}\right),\nonumber
\\\leq&  d_{X'}\left(V^{X'}_{\Phi(x,y_{0,z})}\big(h\circ\Phi(x_{0,z},y_{0,z})\big),\Phi(x,y_{0,z})^{X'}\right)+d_{X'}\left(\Phi(x,y_{0,z})^{X'},\Phi(x,y)^{X'}\right).\label{Jpep45}
\end{align}
Furthermore, as the distance between two points of the same vertical geodesics is equal to their difference of height, we can write the following equality
\begin{align*}
d_{X'}\left(V^{X'}_{\Phi(x,y_{0,z})}\big(h\circ\Phi(x_{0,z},y_{0,z})\big),\Phi(x,y_{0,z})^{X'}\right)&=\Delta h\left(\Phi(x,y_{0,z})^{X'},\Phi(x_{0,z},y_{0,z})^{X'}\right),\nonumber
\\&=\Delta h\left(\Phi(x,y_{0,z})^{Y'},\Phi(x_{0,z},y_{0,z})^{Y'}\right).
\end{align*}
We combine it with inequality (\ref{Jpep45}), and then use the Lipschitz Property of $h$ to get
\begin{align*}
d_{X'}\left(\hat{\Phi}_i^X(x),\Phi(x,y)^{X'}\right)&\leq  \Delta h\left(\Phi(x,y_{0,z})^{Y'},\Phi(x_{0,z},y_{0,z})^{Y'}\right)+d_{X'}\left(\Phi(x,y_{0,z})^{X'},\Phi(x,y)^{X'}\right),
\\&\leq d_{Y'}\left(\Phi(x,y_{0,z})^{Y'},\Phi(x_{0,z},y_{0,z})^{Y'}\right)+d_{X'}\left(\Phi(x,y_{0,z})^{X'},\Phi(x,y)^{X'}\right),
\\&\preceq_{k,c,\bowtie} 2\varepsilon R,\quad\text{ by inequalities (\ref{IneqFinThm1}) and (\ref{IneqFinThm2})}.
\end{align*}
Similarly, we define $\hat{\Phi}_i^Y(y)$ by
\begin{align*}
\hat{\Phi}_i^Y(y):=V^{Y'}_{\Phi(x_{0,z},y)}\big(h\circ\Phi(x_{0,z},y_{0,z})\big).
\end{align*}
and we show that $d_{Y'}\left(\hat{\Phi}_i^Y(y),\Phi(x,y)^{Y'}\right)\preceq_{k,c,\bowtie} \varepsilon R$. Furthermore for all $(x,y)\in U_i$ we have $h\left(\hat{\Phi}_i^X(x)\right)=-h\left(\hat{\Phi}_i^Y(y)\right)$, hence $\hat{\Phi}_i:=(\hat{\Phi}_i^X,\hat{\Phi}_i^Y): X\bowtie Y\to X'\bowtie Y'$ is a well defined product map. Then we chose $U_i':=U'$ to conclude the proof.
\\The downward orientation case is dealt in the same way by switching the definitions of $\hat{\Phi}^X_i$ and $\hat{\Phi}^Y_i$.
\end{proof}

\subsection{Shadows and orientation}\label{SecThmOrient}

We use the fact that $m> n$ to prove that $\Phi$ is orientation preserving, hence the upward orientation is dominant, on each good box at scale $R$.

\begin{propo}\label{PropoOrientPreserving}
Assume that $m>n$ and that $m'>n'$. For $R\succeq_{\bowtie}\frac{1}{\theta}$ the product map $\hat{\Phi}_i$ of Theorem \ref{ThmPartI} is orientation preserving for each $i\in I_g$.
\end{propo}

We recall that given a box $\mathcal{B}$, the shadow of a subset $U\subset\mathcal{B}$, we denote by $Sh(U)$, the set of points of $\mathcal{B}$ below $U$ in the following sens:
\begin{align*}
Sh(U):=\{p\in \mathcal{B}\mid  \exists V\in V\mathcal{B}\text{ containing }p\text{ and intersecting }U\text{ on a point }p'\text{ such that }h(p')\geq h(p)\}.
\end{align*}
And we remind the reader that given a subset $S\subset X$, the large $Y$-horosphere given by $S$ and denote by $H_S\subset X\bowtie Y$, is the set
\begin{align*}
H_S:=S\bowtie Y.
\end{align*} 
Let us denote $\mathcal{B}=\mathcal{B}_i$ for $i\in I_g$. Thanks to Theorem \ref{ThmPartI}, there exist $U=U_i$ with $\lambda(U)\geq\left(1-\theta^{\frac{1}{4}}\right)\lambda(\mathcal{B})$ such that $\Phi$ is close to a product map on $U$. We consider two parameters $\rho_1$ and $\rho_2$ with $1\succeq_{\bowtie} \rho_2\succeq_{\bowtie} \rho_1\succeq_{\bowtie}\theta^{\frac{1}{16}}$. The relations between them will be specified later. Hence Lemma \ref{LemmaGoodShadowForHorocycle} applies with $\alpha=\theta^{\frac{1}{4}}$, and it gives us a $Y$-horosphere $H_{x_0}$ such that
\begin{align*}
\lambda\left(Sh(H_{D_{M_0}(x_0)})\cap U^c\right)<\theta^{\frac{1}{16}}\lambda\left(Sh(H_{D_{M_0}(x0)})\right).
\end{align*}
Then we apply twice Lemma \ref{LemmaGoodPlanForHorocycle} with $\alpha=\theta^{\frac{1}{4}}$, and $\rho=\rho_i$ for $i\in\{1,2\}$ to get two level sets of $h$ in $\mathcal{B}$, $P_1$ and $P_2$, such that for $i\in\{1,2\}$
\begin{align*}
\lambda_{h(P_i)}(P_i\cap \mathrm{Sh}(H_{D_{M_0}(x_0)})\cap U^c)\preceq_{\bowtie}\theta^{\frac{1}{16}}\lambda_{h(P_i)}(P_i\cap \mathrm{Sh}(H_{D_{M_0}(x_0)})),
\end{align*}
and such that $\rho_i R< \Delta h(P_i,H_{x_0})<2\rho_i R$.

The next lemma will gives us the existence of two subsets below a $Y$-horosphere $H$, which are sufficiently big (for the measure $\mu$ in comparison to the horosphere) and sufficiently apart from each other so that any path linking them must get close to $H$. 

This lemma is strongly inspired from Lemma 5.9 of \cite{EFW1}.

\begin{lemma}\label{LemmaPathBetweenS1S2}
There exist a constant $M_1(k,c,\bowtie)$ a constant depending on $k,c$ and on the metric measured spaces $X\bowtie Y$ with the following property. In the settings above, for $R\succeq_{\bowtie}\frac{1}{\rho_2}$, there exist $S_1$ and $S_2$, two subsets of $P_2\cap \mathcal{B}$ such that for $j\in\{1,2\}$ we have:
\begin{enumerate}
\item $\forall s_1\in S_1$, $s_2\in S_2$, $d_X(s_1^X,s_2^X)\geq \rho_2 R$;
\item $\lambda_{h(P_2)}(S_j\cap U^c)\preceq_{\bowtie}\theta^{\frac{1}{32}}\lambda_{h(P_2)}(S_j)$;
\item $\mu_{h(P_2)}\left(S_j\right)\succeq_{\bowtie}\exp\left(\frac{m-n}{2}\rho_2 R\right) \mu_{h(H)}\big(\mathcal{N}_{M_0}(H)\big)$;
\item Any path $\gamma$ joining $S_1$ and $S_2$ of length $l(\gamma)\leq M_1\rho_2 R$ intersects $\mathcal{N}_{6 \rho_1 R}(H)$.
\end{enumerate}
\end{lemma}

\begin{figure}
\begin{center}
\includegraphics[scale=0.75]{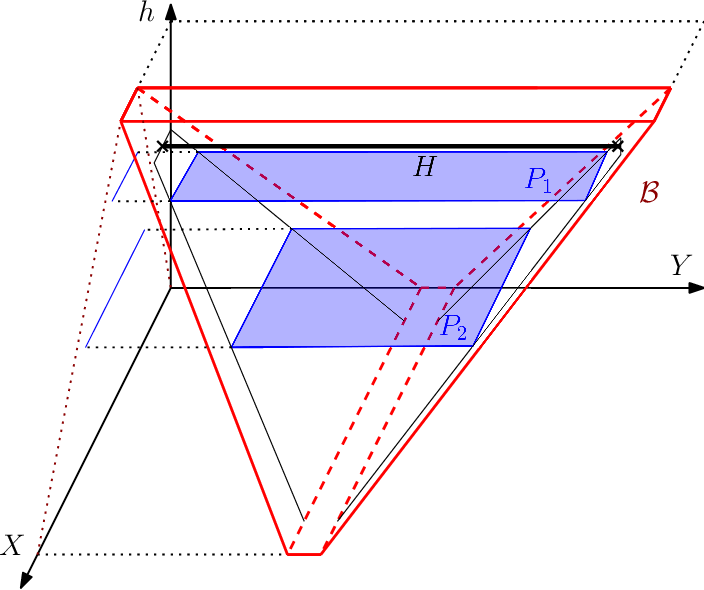} 
\end{center}
\caption{Configuration of Lemma \ref{LemmaPathBetweenS1S2}}
\label{FigProofExistSHoroProd}
\end{figure}

\begin{figure}
\begin{center}
\includegraphics[scale=1.4]{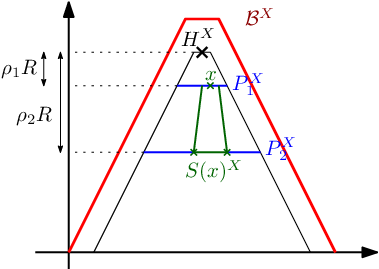} 
\end{center}
\caption{Construction of $S(x)^X$ in Lemma \ref{LemmaPathBetweenS1S2}}
\label{FigProofExistSDansX}
\end{figure}

\begin{proof}
For $j\in\{1,2\}$, let us denote by $Q_j:= P_j\cap \mathrm{Sh}(H_{D_{M_0}(x_0)})$. We tile $Q_1^X$ with the top of boxes as in a box tiling. More precisely, let $M_0$ be the constant involved in assumption $(E2)$, and let $Z\subset Q_1^X$ be an $2M_0$-maximal separating set of $Q_1^X$. Then there exists a set of disjoint cells $\{\mathcal{C}(x)\mid x\in Z\}$ such that:
\begin{enumerate}
\item $\forall x \in Z$, $D(x,M_0)\subset\mathcal{C}(x)\subset D(x,2M_0)$;
\item $Q_1^X =\bigcup_{x\in Z}\mathcal{C}(x)$.
\end{enumerate} 
Thanks to this tessellation, we tile $Q_1$ with the large horosphere $H_{\mathcal{C}(x)}:=\mathcal{C}(x)\times \mathcal{B}_{-h(P_1)}^Y=\mathcal{C}(x)\times Q_1^Y$. Furthermore for any two points $x_1,x_2\in Z$
\begin{align*}
\lambda_{h(P_1)}\Big(H_{\mathcal{C}(x_1)}\Big)&=\lambda^X_{h(P_1)}\big(\mathcal{C}(x_1)\big)\lambda^Y_{-h(P_1)}\left(\mathcal{B}^Y_{-h(P_1)}\right),
\\&\asymp_{\bowtie}\lambda^X_{h(P_1)}\big(\mathcal{C}(x_2)\big)\lambda^Y_{-h(P_1)}\left(\mathcal{B}^Y_{-h(P_1)}\right) ,\quad\text{by Lemma \ref{LemmaMuDiskAreExp}},
\\&=\lambda_{h(P_1)}\Big(H_{\mathcal{C}(x_2)}\Big).
\end{align*}
Therefore $\lambda_{h(P_1)}\left(Q_1\right)\asymp_{\bowtie} \lambda_{-h(P_1)}^Y (Q_1^Y)\# Z$. We tile $Q_2$ by projections of the tessellation of $Q_1$, these projections look like stripes on $Q_2$ 
\begin{align}
Q_2=\bigsqcup\limits_{x\in Z}\pi^X_{h(P_2)}\big(\mathcal{C}(x)\big)\times  \mathcal{B}_{-h(P_2)}^Y.\label{Tess1}
\end{align}
Let us denote these stripes by $S(x):=\pi^X_{h(P_2)}(\mathcal{C}(x))\times  \mathcal{B}_{-h(P_2)}^Y$ for all $ x\in Z$. For all $x_1,x_2\in Z$, $d_X(x_1,x_2)\geq M_0$, hence by Lemma \ref{LemmaBackward} $\forall \left(s_1^X,s_1^Y\right)\in \mathrm{Int}_{M_0} (S(x_1))$ and $\forall \left(s_2^X,s_2^Y\right)\in \mathrm{Int}_{M_0}(S(x_2))$ we have
\begin{align}
d_X\left(s_1^X,s_2^X\right)&\geq 2\Delta h(P_1,P_2)-M_0=2\rho_2 R-2\rho_1R-M_0-M,
\\&\geq 2(\rho_2-2\rho_1)R,\quad\text{for }R\geq \frac{2(M_0+M)}{\rho_1}.\label{IneqDistS1S2}
\end{align}
Furthermore we have by construction that
\begin{align*}
\lambda^X_{h(P_2)}\left(\pi^X_{h(P_2)}(\mathcal{C}(x_1))\right)\asymp_{\bowtie} \lambda^X_{h(P_2)}\left(\pi^X_{h(P_2)}(\mathcal{C}(x_2))\right).
\end{align*}
Hence, combined with Lemma \ref{LemmaYetToBeWritten}, we get
\begin{align*}
\lambda_{h(P_2)}(\mathrm{Int}_{M_0}(S(x_1)))\asymp_{\bowtie}\lambda_{h(P_2)}(S(x_1))\asymp_{\bowtie} \lambda_{h(P_2)}(S(x_2))\asymp_{\bowtie}\lambda_{h(P_2)}(\mathrm{Int}_{M_0}(S(x_2)))
\end{align*}
Therefore, by the tessellation (\ref{Tess1}), $\lambda_{h(P_2)}\left(Q_2\right)\asymp_{\bowtie} \lambda_{-h(P_2)}^Y (Q_2^Y)\# Z$. By Lemma \ref{LemmaGoodPlanForHorocycle}, used with $\alpha:=\theta^{\frac{1}{4}}$, we get
\begin{equation}
\lambda_{h(P_2)}\left(Q_2\cap U^c\right)\preceq_{\bowtie}\theta^{\frac{1}{16}}\lambda_{h(P_2)}\left(Q_2\right).\nonumber
\end{equation}
Moreover, for all $x_1,x_2 \in Z$ we have $\lambda_{h(P_2)}(S(x_1))\asymp_{\bowtie} \lambda_{h(P_2)}(S(x_2))$ and the set of stripes $S(x)$ for $x\in Z$ tile the set $Q_2$. Therefore there exists $Z'\subset Z$ such that $\#Z'\geq\left(1-\theta^{\frac{1}{32}}\right)\#Z$ and such that for all $x\in Z'$ we have $\lambda_{h(P_2)}(S(x)\cap U^c)\leq \theta^{\frac{1}{32}}\lambda_{h(P_2)} (S(x))$. \\We are now able to define $S_1$ and $S_2$. Let $x_1,x_2\in Z$ be distinct and, for $j\in\{1,2\}$, let us denote by $S_j$ the following subset of $S(x_j)$
\begin{align}
S_j:= \pi^X_{h(P_2)}(\mathcal{C}(x_j))\times\mathrm{Int}_{M\rho_2 R}\left(  \mathcal{B}_{-h(P_2)}^Y\right).\label{DefSi}
\end{align} 
By Lemma \ref{LemmaYetToBeWritten}, applied with $r=M\rho_2 R$, $z_0=-h^-(\mathcal{B})$ and $z_1= -h(P_2)$, we have $\mu^Y_{h(P_2)}\left(\mathcal{B}_{-h(P_2)}^Y\right)\asymp_{\bowtie}\mu^Y_{h(P_2)}\left(\mathrm{Int}_{M\rho_2 R}\left(  \mathcal{B}_{-h(P_2)}^Y\right)\right)$, therefore
\begin{align}
\mu_{h(P_2)}(S_j)\asymp_{\bowtie}\mu_{h(P_2)}(S(x_j)). \label{redemption}
\end{align} 
The first point of the Lemma holds by inequality (\ref{IneqDistS1S2}), and the second point holds because we choose $x_1$ and $x_2$ in $Z'$.  
\\\\Let us now prove the third point. Let $y_0\in Y$ be the nucleus of the cell of $\mathcal{B}^Y$, we have $\mathcal{B}^Y_{-z}:=\pi^Y_{-z}(\mathcal{C}(y_0))$. LEt us denote by $h^-:=h^-(\mathcal{B})$. By Lemma \ref{LemmaDisqueInShadow} applied with $p=y_0$,  $z_0=-h^-$ and $z=-(h(H)-\rho_2R)=-h(P_2)$ we have
\begin{align*}
D_{2|h^- -h(P_2)|-M_0}\left(\pi^Y_{-h(P_2)}(y_0)\right)\subset \mathcal{B}^Y_{-h(P_2)}\subset D_{2|h^- -h(P_2)|+M_0}\left(\pi^Y_{-h(P_2)}(y_0)\right).
\end{align*}
It follows that, for $x\in Z$,
\begin{align*}
&\pi^X_{h(P_2)}(\mathcal{C}(x))\times D^Y_{2\left(|h^- -h(H)|+\rho_2 R\right)-M_0}(\pi^Y_{-h(P_2)}(y_0))\subset S(x),
\\\subset~& \pi^X_{h(P_2)}(\mathcal{C}(x))\times D^Y_{2\left(|h^- -h(H)|+\rho_2 R\right)+M_0}(\pi^Y_{-h(P_2)}(y_0)).
\end{align*}
By Lemma \ref{LemmaDisqueInShadow}, $\pi_{h(P_2)}^X(\mathcal{C}(x))$ resembles a disk of radius $2|h(P_1)-h(P_2)|\pm M_0=2(\rho_2-\rho_1)R\pm M_0$. Lemma \ref{LemmaMuDiskAreExp} gives $\mu^X_{h(P_2)}\left(\pi^X_{h(P_2)}(\mathcal{C}(x))\right)\asymp e^{m(\rho_2-\rho_1)R}$. Again by Lemma \ref{LemmaMuDiskAreExp} applied on 
\begin{align*}
D^Y_{2\left(|h^--h(H)|+\rho_2 R\right)+M_0}\left(\pi^Y_{-h(P_2)}(y_0)\right),
\end{align*}
we have
\begin{align*}
\mu_{h(P_2)}(S(x))\asymp_{\bowtie}e^{m(\rho_2-\rho_1)R}e^{n\left(|h^--h(H)|+\rho_2R\right)}.
\end{align*}
Similarly $Q_2$ resembles a product $D_{2\rho_2R\pm M_0}\times B^Y_{-h(P_2)}$, hence
\begin{align*}
\mu_{h(P_2)}\left(Q_2\right)\asymp_{\bowtie}e^{m\rho_2R}e^{n\left(|h^--h(H)|+\rho_2R\right)}.
\end{align*} 
Therefore we obtain an estimate of $\# Z$
\begin{align}
\frac{\mu_{h(P_2)}\left(Q_2\right)}{\mu_{h(P_2)}(S(x))}\asymp_{\bowtie} e^{m\rho_1 R}.\label{ControlMaSWP}
\end{align}
Applying Lemma \ref{LemmaProjCroitExp} with $A=Q_2$, $U=\mathcal{N}_{M_0}(H)$ and $\Delta=\rho_2 R$ gives
\begin{align*}
\mu_{h(P_2)}\left(Q_2\right)\succeq_{\bowtie}\exp\left((m-n)\rho_2R)\right)\mu_{h(H)}\left(\mathcal{N}_{M_0}(H)\right).
\end{align*} 
In combination with inequalities (\ref{redemption}) and (\ref{ControlMaSWP}) we have for $j\in \{1,2\}$
\begin{align*}
\mu_{h(P_2)}\left(S_j\right)\succeq_{\bowtie}&\exp((m-n)\rho_2 R-m\rho_1R) \mu_{h(H)}\big(\mathcal{N}_{M_0}(H)\big),
\\\succeq_{\bowtie}&\exp\left(\frac{m-n}{2}\rho_2 R\right) \mu_{h(H)}\big(\mathcal{N}_{M_0}(H)\big),
\end{align*}
where the last inequality holds since $(m-n)\rho_2 -m\rho_1\geq\frac{m-n}{2}\rho_2 $ when $\rho_1\leq \frac{m-n}{m}\rho_2$. Therefore the third conclusion of this Lemma holds.
\\\\Let us prove the fourth conclusion. Let $\gamma$ be a path joining $s_1\in S_1$ and $s_2\in S_2$ such that $l(\gamma)\leq M\rho_2 R$. By inequality (\ref{IneqDistS1S2}), $d_X\left(s_1^X,s_2^X\right)\geq 2\rho_2R-4\rho_1R$. By Lemma \ref{LEM0} there exists a constant $M'(\delta)$ such that the geodesic segment $\left[s_1^X,s_2^X\right]$ contains a point $s_3^X$ within $4\rho_1 R-M'(\delta)\leq 5\rho_1 R$ of $H^X=\{x_0\}$, for $R\geq \frac{M'(\delta)}{\rho_1}$. Therefore by Proposition \ref{LemmeBrid}
\begin{align*}
l(\gamma^X)\geq 2^{\delta d_X\left(\gamma^X,s_3^X\right)}.
\end{align*}
However, every $\delta$-hyperbolic space with $\delta\leq 1$ is also $1$-hyperbolic. Therefore we can assume without loss of generality that $\delta \geq 1$. Then we have
\begin{align*}
 l(\gamma^X)\geq  2^{d_X\left(\gamma^X,s_3^X\right)}\geq 2^{ d_X\left(\gamma^X,H^X\right)-5\rho_1 R}. 
\end{align*}
Hence $\log_2(M\rho_2 R)\geq d\left(\gamma^X,H^X\right)-5\rho_1 R$. Furthermore, there exists $M'(k,c,\bowtie)$ such that for $R\geq \frac{M'}{\rho_2}$ we have $\log_2(M\rho_2 R)\leq \rho_1 R$. In this case
\begin{align*}
&d\left(\gamma^X,H^X\right)\leq 6 \rho_1 R .
\end{align*}
Therefore there exists $t\in\mathbb{R}$ such that $\Delta h(\gamma(t),H)\leq 6 \rho_1 R$. Let us now look at $\gamma ^Y$. Two cases arise, we have either $\gamma ^Y(t)\in \mathrm{Sh}\left(\mathcal{B}_{-h(P_2)}^Y\right)$ or $\gamma ^Y(t)\notin \mathrm{Sh}\left(\mathcal{B}_{-h(P_2)}^Y\right)$. 
\\In the first case, there exists $y\in H^Y$ such that $\gamma ^Y(t)\in V_y$. Furthermore $\Delta h(\gamma(t),H)\leq 6 \rho_1 R$, hence $d_Y\left(\gamma^Y(t),H^Y\right)=\Delta h\left(\gamma^Y(t),H^Y\right)\leq 6 \rho_1 R$ and consequently $d_Y\left(\gamma^Y,H^Y\right)\leq 6 \rho_1 R$. Which proves $d\left(\gamma,H\right)\leq 6 \rho_1 R$. 
\\In the second case, when $\gamma ^Y(t)\notin \mathcal{B}_{-h(P_2)}^Y$, by our claim (\ref{DefSi}) we have that the vertical geodesic ray $V_{\gamma ^Y(t)}$ starting at $\gamma ^Y(t)$ intersect $Y_{-h\left(P_2\right)}$ in a point $y$ such that $d_Y(y,S_1^Y\cup S_2^Y)>M \rho_2 R$. Therefore
\begin{align*}
M\rho_2 R\geq l(\gamma)&\geq \frac{1}{2}l\left(\gamma^Y\right)\geq \frac{1}{2}\left( d(s_1,\gamma(t))+d(\gamma(t),s_2) \right),
\\&>\frac{2M\rho_2 R}{2}>M\rho_2 R,
\end{align*}
which is absurd, hence the second case when $\gamma ^Y(t)\notin \mathcal{B}_{-h(P_2)}^Y$ does not occur. Therefore we always have that $\gamma$ intersect the $6 \rho_1 R$-neighbourhood of $H$.
\end{proof}

\begin{proof}[Proof of Proposition \ref{PropoOrientPreserving}]
Let us be in the settings above. Let us assume by contradiction that $\hat{\Phi}$ is orientation reversing, which means that there exist $\hat{\Phi}^X:X\to Y'$ and $\hat{\Phi}^X:Y\to X'$ such that for all $(x,y)\in \mathcal{B}$ we have $\hat{\Phi}(x,y)=(\hat{\Phi}^Y(y),\hat{\Phi}^X(x))$.
\\For all $p\in X'\bowtie Y'$ such that $d_{\bowtie'}\left(p,\hat{\Phi}(H\cap U)\right)\leq \rho_1 R$ there exists $q\in H\cap U$ such that $d_{\bowtie'}\left(p,\hat{\Phi}(q)\right)\leq \rho_1 R$. Therefore by the triangle inequality
\begin{align*}
d_{\bowtie'}\left(p,\Phi(q)\right)& \leq d_{\bowtie'}\left(p,\hat{\Phi}(q)\right)+d_{\bowtie'}\left(\hat{\Phi}(q),\Phi(q)\right)\preceq_{k,c,\bowtie} \rho_1 R+\varepsilon R,\quad\text{by Theoreom \ref{ThmPartI} since }q\in U ,
\\&\preceq_{k,c,\bowtie} \rho_1 R,\quad\text{since }\varepsilon\leq \rho_1.
\end{align*}
Hence there exists $M(k,c,\bowtie)$ such that $\mathcal{N}_{\rho_1 R}(\hat{\Phi}(H\cap U))\subset \mathcal{N}_{M\rho_1 R}(\Phi(H\cap U))$. We show similarly that for $j\in\{1,2\}$
\begin{align}
\mathcal{N}_{\rho_1 R}(\Phi(S_j\cap U))\subset \mathcal{N}_{M\rho_1 R}(\hat{\Phi}(S_j\cap U)).\label{IneqNeighCloseNeigh}
\end{align} 
Let $M'(\bowtie)$ be the constant involved in Corollary \ref{CoroMeasureNeighbourhood}. Then
\begin{align*}
\mu\left(\mathcal{N}_{8k\rho_1R}(\Phi(H))\right)&\preceq_{k,c,\bowtie}e^{48k\rho_1 R m'}\mu\left(\mathcal{N}_{kc+c}(\Phi(H))\right),\quad\text{by Corollary \ref{CoroUniformBallMeasure}},
\\&\preceq_{k,c,\bowtie}e^{48k\rho_1 R m'} \mu\left(\mathcal{N}_{1}(H)\right),\quad\text{by Lemma \ref{LemmaVolKeepByQI}},
\\&\leq e^{48k\rho_1 R m'} \mu\left(\mathcal{N}_{M'}(H)\right),
\\&\asymp_{\bowtie}e^{48k\rho_1 R m'} \mu_{h(H)}\left(\mathcal{N}_{M'}(H)\right),\quad\text{by the second part of Corollary \ref{CoroMeasureNeighbourhood}},
\\&\preceq_{\bowtie}e^{48k\rho_1 R m'} \mu_{h(H)}\left(\mathcal{N}_{M_0}(H)\right),\quad\text{by the first part of Corollary \ref{CoroMeasureNeighbourhood}}.
\end{align*}
Combined with point 3 of Lemma \ref{LemmaPathBetweenS1S2} we have
\begin{align*}
\mu\left(\mathcal{N}_{8k\rho_1R}(\Phi(H))\right)&\preceq_{\bowtie}e^{-(m-n)\frac{\rho_2}{2} R }e^{48k\rho_1 R m'} \mu_{h(P_2)}\left(S_j\right),
\\&\preceq_{\bowtie}e^{-(m-n)\frac{\rho_2}{2} R }e^{48k\rho_1 R m'} \mu_{h(P_2)}\left(S_j\cap U\right),\quad\text{ thanks to 2. of Lemma  \ref{LemmaPathBetweenS1S2}},
\\&\preceq_{\bowtie}e^{-(m-n)\frac{\rho_2}{4} R} \mu_{h(P_2)}\left(\mathcal{N}_{1}(S_j\cap U)\right),\quad\text{since }\rho_1\leq\frac{m-n}{96k m'}\rho_2,
\\&\asymp_{\bowtie}e^{-(m-n)\frac{\rho_2}{4} R} \mu\left(\mathcal{N}_{M'}(S_j\cap U)\right),\quad\text{by Corollary \ref{CoroMeasureNeighbourhood}}.
\\&\leq e^{-(m-n)\frac{\rho_2 }{4}R} \mu\left(\mathcal{N}_{M'+kc+c}(S_j\cap U)\right).
\end{align*}
Hence using Lemma \ref{LemmaVolKeepByQI} on $\mathcal{N}_{M'}(S_j\cap U)$
\begin{align}
\mu\left(\mathcal{N}_{8k\rho_1R}(\Phi(H))\right)&\preceq_{k,c,\bowtie} e^{-(m-n)\frac{\rho_2 }{4}R} \mu\left(\mathcal{N}_{M'+1}(\Phi(S_j\cap U))\right),\nonumber
\\&\leq e^{-(m-n)\frac{\rho_2 }{4}R} \mu\left(\mathcal{N}_{\rho_1 R}(\Phi(S_j\cap U))\right),\quad\text{for }R\geq\frac{M'}{\rho_1},\nonumber
\\&\leq e^{-(m-n)\frac{\rho_2 }{4}R} \mu\left(\mathcal{N}_{M\rho_1 R}(\hat{\Phi}(S_j\cap U))\right),\quad\text{by inequality (\ref{IneqNeighCloseNeigh})},\nonumber
\\&\preceq_{k,c,\bowtie} e^{-(m-n)\frac{\rho_2 }{4}R} e^{6M\rho_1 R m'} \mu\left(\mathcal{N}_{M'}(\hat{\Phi}(S_j\cap U))\right),\quad\text{by Corollary \ref{CoroMeasureNeighbourhood}},\nonumber
\\&\preceq_{k,c,\bowtie} e^{-(m-n)\frac{\rho_2 }{8}R} \mu\left(\mathcal{N}_{M'}(\hat{\Phi}(S_j\cap U))\right),\quad\text{since }\rho_1\leq\frac{m-n}{48M m'}\rho_2,\nonumber
\\&\asymp_{k,c,\bowtie} e^{-(m-n)\frac{\rho_2 }{8}R} \mu_{\hat{z_0}}\left(\mathcal{N}_{M'}(\hat{\Phi}(S_j\cap U))\cap X'_{\hat{z_0}}\right),\quad\text{by the second part of Corollary \ref{CoroMeasureNeighbourhood}}.\nonumber
\end{align}
where $\hat{z_0}:=\hat{\Phi}(P_2)$. Since $\hat{\Phi}$ is orientation reversing, we can now apply Lemma \ref{LemmaNoVertIfSmall} with $A_j=\hat{\Phi}(S_j\cap U)$, $E=\mathcal{N}_{8k\rho_1R}(\Phi(H))$ and $Q=e^{(m-n)\frac{\rho_2 }{8}R}$ we have that
\begin{align*}
\eta\left(V\mathcal{N}_{M_0}\left(\hat{\Phi}(S_j\cap U)\right)\right)\succeq_{k,c,\bowtie} e^{(m-n)\frac{\rho_2 }{8}R}\eta\left(V\mathcal{N}_{M_0}\left(E\right)\right).
\end{align*}
Then, as pointed out below Lemma \ref{LemmaCloseGeodAreGood}, we can apply it on a $A_j$ with $V_1= VE$. Hence let us take $U_{A_j}\subset A_j$ maximal for the inclusion such that:
\begin{enumerate}
\item[$\cdot$]$\lambda_{\hat{z_0}}(U_{A_j})\geq (1- e^{(m-n)\frac{\rho_2 }{8}R})\lambda_{\hat{z_0}}(A_j)$.
\item[$\cdot$]For all $p\in U_{A_j}$, most of the vertical geodesic in $D_{M_0}(p)$ do not intersect $E$.
\end{enumerate} 
By Property \ref{PropMeasureGeodVertSameMeasureSpace} we have
\begin{align*}
\lambda_{\hat{z_0}}\left(\mathcal{N}_{M_0}\left(\hat{\Phi}(S_j\cap U)\right)\right)\succeq_{k,c,\bowtie} e^{(m-n)\frac{\rho_2 }{8}R}\lambda_{\hat{z_0}}\left(\pi^{\bowtie}_{\hat{z_0}}\left(\mathcal{N}_{M_0}\left(E\right)\right)\right).
\end{align*}
Hence by the definition of $\lambda_{\hat{z_0}}$
\begin{align}
\mu_{\hat{z_0}}\left(\mathcal{N}_{M_0}\left(\hat{\Phi}(S_j\cap U)\right)\right)\succeq_{k,c,\bowtie}e^{(m-n)\frac{\rho_2 }{8}R} \mu_{\hat{z_0}}\left(\pi^{\bowtie}_{\hat{z_0}}\left(\mathcal{N}_{M_0}\left(E\right)\right)\right).\label{IneqDesMuPrSj}
\end{align}
Let us denote $E':=\mathcal{N}_{M_0}(\hat{\Phi}(S_j\cap U)\setminus U_{A_j})$. By the construction of $U_{A_j}$, $E'\cap X'_{z_0}$ is of almost full measure in $\pi^{\bowtie}_{\hat{z_0}}\left(\mathcal{N}_{M_0}\left(E\right)\right)$. Furthermore, by Theorem $\ref{ThmPartI}$ $\hat{\Phi}$, is $M\varepsilon R$-close to $\Phi$ on $U$, hence we have (similarly as in inequality (\ref{IneqNeighCloseNeigh})) that
\begin{align*}
\mathcal{N}_{\rho_1 R}\left(\hat{\Phi}^{-1}\left(E'\right)\right)\subset \mathcal{N}_{M\rho_1 R}\left(\Phi^{-1}\left(E'\right)\right).
\end{align*}
Therefore
\begin{align*}
\mu\left(\mathcal{N}_{\rho_1 R}\left(\hat{\Phi}^{-1}\left(E'\right)\right)\right)&\preceq_{k,c,\bowtie} \mu\left(\mathcal{N}_{M\rho_1 R}\left(\Phi^{-1}\left(E'\right)\right)\right),
\\&\preceq_{k,c,\bowtie } e^{6M\rho_1 R m} \mu\left(\mathcal{N}_{kc+c}\left(\Phi^{-1}\left(E'\right)\right)\right),\quad\text{by the first part of Corollary \ref{CoroMeasureNeighbourhood}},
\\&\asymp_{k,c,\bowtie } e^{6M\rho_1 R m} \mu\left(\mathcal{N}_{1}\left(E'\right)\right),\quad\text{by Lemma \ref{LemmaVolKeepByQI}},
\\&\asymp_{k,c,\bowtie } e^{6M\rho_1 R m} \mu_{\hat{z_0}}\left(\mathcal{N}_{M_0}\left(E'\right)\right),\quad\text{by the second part of Corollary \ref{CoroMeasureNeighbourhood}},
\\&\asymp_{k,c,\bowtie } e^{6M\rho_1 R m} \mu_{\hat{z_0}}\left(\pi^{\bowtie}_{\hat{z_0}}\left(\mathcal{N}_{M_0}\left(E\right)\right)\right),
\\&\preceq_{k,c,\bowtie }e^{-(m-n)\frac{\rho_2 }{8}R} e^{6M\rho_1 R m}\mu_{\hat{z_0}}\left(\mathcal{N}_{M_0}\left(\hat{\Phi}(S_j\cap U)\right)\right),\quad\text{by the definition of }U_{A_j},
\\&\preceq_{k,c,\bowtie } e^{-(m-n)\frac{\rho_2 }{16}R}\mu_{h(P_2)}\left(\mathcal{N}_{M_0}\left(S_j\cap U\right)\right),\quad\text{since }\rho_1\leq\frac{\rho_2}{M},
\\&\leq e^{-(m-n)\frac{\rho_2 }{16}R}\mu_{h(P_2)}\left(\mathcal{N}_{M_0}\left(S_j\right)\right).
\end{align*}
Following the second conclusion of Lemma \ref{LemmaPathBetweenS1S2}, there exists a constant $M(\bowtie)$ such that $\lambda_{h(P_2)}(S_j \cap U^c)\leq M\theta^{\frac{1}{32}}\lambda_{h(P_2)}(S_j)$. 
\\We apply twice Lemma \ref{LemmaMarkov1} for $j=1,2$ with $(V_1,\eta)= \left(\mathcal{N}_{M_0}\left(S_j^X\right)\times \mathcal{N}_{M_0}\left(S_j^Y\right),\mu_{h(P_2)} \right)$, $V_0 =  U^{c}\cap\mathcal{N}_{\rho_1 R}\left(\hat{\Phi}^{-1}\left(E'\right)\right)$ and $\alpha:=e^{(m-n)\frac{\rho_2 }{16}R}\mu_{h(P_2)}+M\theta^{\frac{1}{32}}$. Let us denote $G^Y\left(p^X\right):=\left\lbrace p^Y\in V_1^Y\mid \left(p^X,p^Y\right)\in V_0\right\rbrace$, we have that
\begin{align*}
\mu^X_{h(P_2)}\left(\left\lbrace p^X\in V_1^X\mid \mu_{-h(P_2)}^Y \left(G^Y\left(p^X\right)\right)\right\rbrace\right)\geq \left(1-e^{-(m-n)\frac{\rho_2 }{32}R}\right)\mu_{h(P_2)}^Y\left(V_1^Y\right).
\end{align*}
Since $e^{-(m-n)\frac{\rho_2 }{32}R}+M\theta^{\frac{1}{32}}<\frac{1}{2}$, there exists $s_1\in (S_1\cap U)\setminus\hat{\Phi}^{-1}\left(E'\right) $ and $s_2\in  (S_2\cap U)\setminus\hat{\Phi}^{-1}\left(E'\right)$ such that $s_1^Y=s_2^Y$. 
\\\\Let us denote by $\hat{s}_j:=\hat{\Phi}(s_j)$ for $j\in\{1,2\}$. By construction we have $\hat{s}_j\in A_j$, then $VD_{M_0}(\hat{s_j})$ contains almost only vertical geodesic segments which do not intersect $E$. Since $\hat{s}_1^{X'}=\hat{s}_2^{X'}$, and by Lemma \ref{LemmaMarkov2}, we can find two vertical geodesics $v_1\in VD_{M_0}(\hat{s_1}) $ and $v_2\in VD_{M_0}(\hat{s_2})$ which do not intersect $E=\mathcal{N}_{8k\rho_1R}(\Phi(H))$, and such that $v_1^X=v_2^X$. Since $v_1^Y$ and $v_2^Y$ meet (up to an additive constant) at the height $-\hat{z}_0+\frac{1}{2}d_{Y'}(\hat{s}_1^{Y'},\hat{s}_2^{Y'})$, there exist $M(\delta)$ such that the concatenation of $v_1$ and $v_2$ is $(1,M(\delta))$-quasigeodesic linking $\hat{s}_1$ to $\hat{s}_2$. 
\\Let us denote by $\gamma:=\Phi^{-1}(v_1\cup v_2)$, then $\gamma$ is a $(k,c+M)$-quasigeodesic. By Lemma 2.1 of \cite{GS}, there exists a $2k$-Lipschitz, $(k,4(M+c))$-quasi-geodesic $\gamma'$ in the $2(M+c)$-neighbourhood of $\gamma$, linking $\Phi^{-1}(\hat{s}_1)$ to $\Phi^{-1}(\hat{s}_2)$. Let us denote $s_1'=\Phi^{-1}(\hat{s}_1)$ and $s_2'=\Phi^{-1}(\hat{s}_2)$. Because $\gamma'$ is $2k$-Lipschitz, and since $\Phi^{-1}$ is a $(k,c)$-quasi-isometry we have
\begin{align}
l(\gamma')&\leq 2k d_{\bowtie'}(\hat{s}_1,\hat{s}_2) \leq k^2 d_{\bowtie}(s_1',s_2')+c.\label{cafcaf1}
\end{align}
Furthermore, $\gamma'$ does not intersect the $\frac{1}{k}(7k\rho_1R-2c)-c$-neighbourhood of $H$ since $\Phi^{-1}$ is a quasi-isometry. Moreover $s_j'$ and $s_j$ are $\varepsilon$R close to each other, that is
\begin{align}
d_{\bowtie}(s_j',s_j)&=d_{\bowtie}(\Phi^{-1}(\hat{\Phi}(s_j)),s_j),\nonumber
\\&\leq k d_{\bowtie'}(\hat{\Phi}(s_j),\Phi(s_j))\preceq_{k,c,\bowtie}\varepsilon R,\quad\text{since }s_j\in U.\label{ClF}
\end{align}
Consequently by the triangle inequality we get
\begin{align}
d_{\bowtie}(s_1',s_2') &\leq d_{\bowtie}(s_1',s_1)+d_{\bowtie}(s_1,s_2)+d_{\bowtie}(s_2,s_2'),\nonumber
\\&\preceq_{k,c,\bowtie}\varepsilon R + d_{\bowtie}(s_1,s_2),\quad\text{ since }\hat{\Phi}^{-1}(s_j)\in U.\label{cafcaf2}
\end{align}
Furthermore $s_1^Y=s_2^Y$, therefore by Corollary \ref{lengthGeod}, with $M=15C_0$ we obtain
\begin{align*}
d_{\bowtie}(s_1,s_2)\leq d_X\left(s_1^X,s_2^X\right)+M\leq 2\rho_2 R+M,\quad\text{by the first point of Lemma \ref{LemmaPathBetweenS1S2}}.
\end{align*}
Combined with inequalities (\ref{cafcaf1}) and (\ref{cafcaf2}) we get
\begin{align*}
l(\gamma')\preceq_{k,c,\bowtie} 2k^2 (2\rho_2 R+M+2\varepsilon R)+c\preceq_{k,c,\bowtie}\rho_2 R,\quad\text{for }R\geq \frac{M+c}{\rho_2}.
\end{align*}
For $j\in\{1,2\}$, let $\gamma_j:=[s_j,s_j']$, by inequality (\ref{ClF}) we have $l(\gamma_j)\preceq_{k,c,\bowtie}\varepsilon R$. Hence the path $\gamma''$, constructed as the concatenation of $\gamma_1$, $\gamma'$ and $\gamma_2$, is a path linking $s_1\in S_1$ to $s_2\in S_2$, of length $l(\gamma)\preceq_{k,c,\bowtie}\rho_2 R$ since $\varepsilon\leq \rho_2$. Furthermore, by construction, $\gamma''$ does not intersect the $7\rho_1R-3c-2M\varepsilon R>6\rho_1 R$-neighbourhood of $H$. This contradicts the fourth point of Lemma \ref{LemmaPathBetweenS1S2}, therefore $\Phi$ is orientation preserving.
\end{proof}


\subsection{Factorisation of a quasi-isometry in big boxes}\label{SecPartII}

In Section \ref{SecPartI} we proved that for all $i\in I_g$, $\Phi_{|\mathcal{B}_i}$ is close to a quasi-isometry product $\hat{\Phi}_i=\left(\hat{\Phi}_i^X,\hat{\Phi}_i^Y\right)$ on a set of almost full measure $U_i\subset\mathcal{B}_i$. In this section we prove that $\Phi$ is close to $\hat{\Phi}$ on all boxes at scale $L$ on a set of almost full measure. This is a step-forward since this is true on all boxes at scale $L$ and not only a significant number of them.
 
\begin{thm}\label{ThmPartII}
For $0<\theta \preceq_{k,c,\bowtie}1$ there exists $L_0(k,c,\bowtie,\theta)>0$ such that for all $L\geq L_0$ and for all box $\mathcal{B}$ at scale $L$, there exists $M(k,c,\bowtie)$, $U\subset\mathcal{B}$ and a $(k,M\sqrt{\theta} L)$-quasi-isometry product map $\hat{\Phi}=\left(\hat{\Phi}^X,\hat{\Phi}^Y\right)$, with $\hat{\Phi}^X:X\to X'$ and $\hat{\Phi}^Y:Y\to Y'$, such that:
\begin{enumerate}
\item $\lambda(U)\geq \left(1-\theta\right)\lambda(\mathcal{B})$;
\item $d_{\bowtie'}\left(\Phi_{\mid U},\hat{\Phi}_{\mid U}\right)\preceq_{\bowtie}\theta L$.
\end{enumerate}

\end{thm}

Let $\mathcal{B}$ be a box at scale $L$, let $i\in I_g$ and for all $i\in I_g$ let $U_i\subset\mathcal{B}_i$ be as in Theorem \ref{ThmPartI}, where $U_i$ is the subset of $\mathcal{B}_i$ on which $\Phi$ is close to a product map $\hat{\Phi}_i$. Let us denote by $W\subset \mathcal{B}$ the "good" set of $\mathcal{B}$ 
\begin{align*}
W:=\bigsqcup\limits_{i\in I_g} U_i.
\end{align*}
where "good" means the set on which $\Phi$ is close to a product map on boxes at scale $R$.
We introduce the function	 $P$ which quantifies the portion of a geodesic segment which is not in~$W$.

\begin{defn}
Let $\gamma:[0,L]\to X\bowtie Y$ be a vertical geodesic segment of $mathcal{B}$. We denote the measure of points in $\gamma\cap W^c$ by
\begin{equation}
P(\gamma):=\mathrm{Leb}\left(\gamma^{-1}(W^c)\right).
\end{equation}
\end{defn}

The value of $P(\gamma)$ is related to $\gamma$ being $\varepsilon$-monotone. The following lemma is mostly inspired from Lemma 5.10 of \cite{EFW1}.

\begin{lemma}\label{LemmaPimplisMonot}
For $0\leq \varepsilon\preceq_{k,c,\bowtie}  \sqrt{\theta} \preceq_{k,c,\bowtie} 1$, there exists $M(\bowtie,k,c)$ such that for all vertical geodesic segments $\gamma:[0,L]\to X\bowtie Y$ we have
\begin{align*}
P(\gamma)\leq \sqrt{\theta} L\Rightarrow \Phi\circ \gamma\text{ is }M\sqrt{\theta}\text{-monotone}. 
\end{align*}
\end{lemma}

\begin{proof}
Let $t_1,t_2\in[0,L]$ such that $h(\Phi(\gamma(t_1)))=h(\Phi(\gamma(t_2)))$ and such that $t_2\geq t_1$. Let us decompose $[t_1,t_2]$ into segments of length $\sqrt{\theta} R$. Without loss of generality we can assume that $t_2-t_1\geq \sqrt{\theta} L$. Let us denote $N:=\left\lfloor\frac{t_2-t_1}{ \sqrt{\theta} R}\right\rfloor $, $I_{i}:=[t_1+i \sqrt{\theta} R,t_1+(i+1) \sqrt{\theta} R[$ for any $i\in\left\{0,...,N-1\right\}$ and $I_N:=[t_1+(N-1) \sqrt{\theta} R,t_2]$. We have
\begin{align*}
[t_1,t_2]:=\bigsqcup\limits_{i=0}^N I_i.
\end{align*}
Then for all $i\in\{0,...,N\}$ let us choose $s_i\in I_i$ such that $\gamma(s_i)\in W$ if possible, and any $s_i \in I_i$ otherwise. Let us denote by $J$ the set of odd indexes in $\lbrace 0,...,N\}$, we split $J$ into the following sets:
\begin{align*}
J_0:=&\left\{ j\in J\mid \gamma(s_j)\text{ and }\gamma(s_{j+2})\text{ are both in the same box and in }W\right\};\\
J_1:=&\left\{ j\in J\mid \gamma(s_j)\text{ and }\gamma(s_{j+1})\text{ are in different boxes}\right\};\\
J_1':=&\left\{ j\in J\mid \gamma(s_{j+1})\text{ and }\gamma(s_{j+2})\text{ are in different boxes}\right\};\\
J_2:=&\left\{ j\in J\mid I_j \subset W^c\right\};\\
J_2':=&\left\{ j\in J\mid I_{j+2} \subset W^c\right\}.
\end{align*}
We claim that
\begin{align*}
J=J_0\sqcup \big(J_1\cup J_1'\cup J_2\cup J_2'\big).
\end{align*}
To prove it, one can see that two cases arise when an odd index $j$ is not in $J_0$. The first case is when $\gamma(s_j)$ and $\gamma(s_j+2)$ are not in the same box, which leads to the fact that either $j\in J_1$ or $j\in J_1'$. The second case happens when $\gamma(s_j)$ or $\gamma(s_j+2)$ are not in $W$, which leads to either $I_j\subset W^c$ or $I_{j+2}\subset W^c$. Therefore, we proved that an odd index is either in $J_0$ or in $J_1\cup J_1'\cup J_2\cup J_2'$.

We have that $P(\gamma)\leq \sqrt{\theta} L$, hence $\# J_2\leq \frac{\sqrt{\theta}L}{ \sqrt{\theta} R}=\frac{L}{R}$ and similarly $\# J_2'\leq \frac{L}{R}$.  Furthermore there are less than $\frac{L}{R}$ boxes intersecting $\gamma$, therefore $\# J_1\leq \frac{t_2-t_1}{R}\leq \frac{L}{R}$ and $\# J_1'\leq \frac{L}{ R}$, hence
\begin{align*}
&\#\big(J_1\cup J_1'\cup J_2\cup J_2'\big)\leq 4\frac{L}{ R};\\
&\# J_0=\# J-\#\big(J_1\cup J_1'\cup J_2\cup J_2'\big)\geq \frac{t_2-t_1}{2\sqrt{\theta}R}-4\frac{L}{R}.
\end{align*}
We see that the "good" indexes are in majority compared to the "bad" indexes. We now use that fact to prove that $|t_2-t_1|$ is smaller than $\sqrt{\theta}L$. Let us denote $q(t):=h\circ\Phi\circ\gamma(t)$ for all $t\in[0,L]$. We assume that $N$ is odd, the case where $N$ is even is treated identically. By assumption $q(t_1)=q(t_2)$ therefore
\begin{align}
0=&q(t_2)-q(t_1)=q(t_2)-q(s_N)+\sum\limits_{i\in J}\left(q(s_{i+2})-q(s_{i})\right)+q(s_1)-q(t_1),\nonumber
\\=&q(t_2)-q(s_N)+\sum\limits_{i\in J_0}\left(q(s_{i+2})-q(s_i)\right)+\sum\limits_{i\in J\setminus J_0}\left(q(s_{i+2})-q(s_i)\right)+q(s_1)-q(t_1).\label{IneqRatsd}
\end{align}
However we proved that $\#J_0$ is much bigger than $\#(J\setminus J_0)$, and for any $i\in J_0$, $q(s_{i+2})-q(s_i)$ is a positive number by the upward orientation of the quasi-isometry on $W$. Therefore we will show that $|t_1-t_2|$ must be small for this equality to hold. First, we have to consider that $\forall i\in \{0,...,N\}$
\begin{align*}
&l(I_{i+1})\leq|s_i-s_{i+2}|\leq l(I_{i})+l(I_{i+1})+l(I_{i+2}),
\\\Rightarrow & \sqrt{\theta}R\leq|s_i-s_{i+2}|\leq 3 \sqrt{\theta}R,
\\\Rightarrow & |q(s_i)-q(s_{i+2})|\preceq_{k,c,\bowtie}  \sqrt{\theta}R.
\end{align*}
Hence for all $i\in J\setminus J_0$ we have $q(s_{i+2})-q(s_i)\succeq_{k,c,\bowtie} -\sqrt{\theta}R $. Furthermore for all $i\in J_0$, $s_i$ and $s_{i+2}$ are in the same box and in $W$, therefore by Corollary \ref{LemmahIsQI}, there exists $M(k,c,\bowtie)$ such that
\begin{align*}
q(s_{i+2})-q(s_i)\geq \frac{1}{k}|s_i-s_{i+2}|-M\varepsilon R\succeq_{k,c,\bowtie}\sqrt{\theta}R;\quad\text{since }\sqrt{\theta}\geq 2M\varepsilon.
\end{align*}
Combined with equality (\ref{IneqRatsd})
\begin{align*}
0\succeq_{k,c,\bowtie} \sqrt{\theta}R \# J_0  -\sqrt{\theta}R \#\big(J_1\cup J_1'\cup J_2\cup J_2'\big) \succeq |t_2-t_1|-\sqrt{\theta}L.
\end{align*}
Hence $|t_2-t_1|\preceq_{k,c,\bowtie} \sqrt{\theta}L$, which proves that there exists $M(k,c,\bowtie)$ such that $\gamma$ is $M\sqrt{\theta}$-monotone.
\end{proof}


Let $M$ be the constant involved in Lemma \ref{LemmaPimplisMonot}, let $\theta'=\theta^{\frac{1}{16}}$ and let $\varepsilon':= 2M\sqrt{\theta'}$. We now show that almost all vertical geodesic segments of boxes at scale $L$ have $\varepsilon'$-monotone images under $\Phi$.

Let us denote by $V^g\mathcal{B}\subset V\mathcal{B}$ the set of vertical geodesic segments of $V\mathcal{B}$ whose image by $\Phi$ are $\varepsilon'$-monotone.

\begin{lemma}\label{LemmaGoodGeodScaleL}
For $L\succeq_{k,c,\bowtie} \frac{1}{\theta}$ and for any box $\mathcal{B}$ at scale $L$ we have that
\begin{align}
\eta \left(V^g\mathcal{B}\right)\geq (1-\theta^{\frac{1}{32}})\eta \left(V\mathcal{B}\right).
\end{align} 
\end{lemma}
\begin{proof}
Lemma \ref{LemmaPimplisMonot} tells us that $P(\gamma)\geq \sqrt{\theta} L$ for all $\gamma \in V^b\mathcal{B}$. Computing the measure $\lambda$ of $W^c$ we have 
\begin{align}
\lambda(W^c)&=\int\limits_0^L\lambda_z\left(W^c_z\right)\mathrm{d}z\asymp_{\bowtie}\int\limits_0^L\eta\left(V_{\mathcal{B}}\left(W_z^c\right)\right)\mathrm{d}z,\quad\text{by Proposition \ref{PropMeasureGeodVertSameMeasureSpace},}\nonumber
\\&\asymp_{\bowtie}\int\limits_0^L\int\limits_{V\mathcal{B}}\mathbbm{1}_{V_{\mathcal{B}}\left(W_z^c\right)}(\gamma)\mathrm{d}\eta(\gamma)\mathrm{d}z\asymp_{\bowtie}\int\limits_{V\mathcal{B}}\int\limits_0^L\mathbbm{1}_{V_{\mathcal{B}}\left(W_z^c\right)}(\gamma)\mathrm{d}z\mathrm{d}\eta(\gamma),\quad\text{by Fubini Theorem.}\label{IneqInEtaVg}
\end{align}
However we have
\begin{equation}
\mathbbm{1}_{V_{\mathcal{B}}\left(W_z^c\right)}(\gamma)=\left\lbrace \begin{array}{ll}
       0 &\text{ if } z\in \gamma^{-1}(W), \\
       1 &\text{ if } z\in \gamma^{-1}(W^c).
   \end{array}
   \right.
\end{equation}
Therefore $\mathbbm{1}_{V_{\mathcal{B}}\left(W_z^c\right)}(\gamma)=\mathbbm{1}_{\gamma^{-1}(W^c)}(z)$. With inequality (\ref{IneqInEtaVg}) it gives us
\begin{align}
\lambda(W^c)&\asymp_{\bowtie}\int\limits_{V\mathcal{B}}\int\limits_0^L\mathbbm{1}_{\gamma^{-1}(W^c)}(z)\mathrm{d}z\mathrm{d}\eta(\gamma)\geq\int\limits_{V^b\mathcal{B}} \int\limits_0^L\mathbbm{1}_{\gamma^{-1}(W^c)}(z)\mathrm{d}z\mathrm{d}\eta(\gamma),\quad\text{since }V^b\mathcal{B}\subset V\mathcal{B},\nonumber
\\&\geq \int\limits_{V^b\mathcal{B}}\mathrm{Leb}\left(\gamma^{-1}(W^c)\right)\mathrm{d}\eta(\gamma)= \int\limits_{V^b\mathcal{B}}P(\gamma)\mathrm{d}\eta(\gamma).\label{IneqJsaiplusLesNoms}
\end{align}
Let us assume by contradiction that $\eta \left(V^g\mathcal{B}\right)< (1-\sqrt{\theta'})\eta \left(V\mathcal{B}\right)$, hence we have  $\eta \left(V^b\mathcal{B}\right)>\sqrt{\theta'}\eta \left(V\mathcal{B}\right)$. Therefore by inequality (\ref{IneqJsaiplusLesNoms})
\begin{align*}
\lambda(W^c)&\succeq_{\bowtie} \eta\left(V^b\mathcal{B}\right)\sqrt{\theta'}L \geq\sqrt{\theta'} \eta\left(V\mathcal{B}\right)\sqrt{\theta'}L ,
\\&\asymp_{\bowtie} \theta' \lambda (\mathcal{B}),
\end{align*}
which contradicts the first conclusion of Theorem \ref{ThmPartII} for $\theta\preceq_{k,c,\bowtie}1$.
\end{proof}

As in Section \ref{SecPartI}, we deduce that, in boxes which have almost only vertical geodesic segment with $2M\sqrt{\theta'}$-monotone image, $\Phi$ is close to a product map. Let us denote $\varepsilon':=2M\theta^{\frac{1}{16}}$ and $\theta':=2M\theta^{\frac{1}{16}}$, then for $0<\theta'\preceq_{k,c,\bowtie}1$ we have that $\theta'\leq \varepsilon'\leq \sqrt{\theta'}$.

\begin{proof}[Proof of Theorem \ref{ThmPartII}]
The proof is similar to Theorem \ref{ThmPartI}. The Lemma \ref{LemmaGoodGeodScaleL} plays the role of the second conclusion of Lemma \ref{Lemma411}, with $\varepsilon'$ instead of $\varepsilon$. In a box at scale $L$, almost all vertical geodesic segment have $\varepsilon'$-monotone image by $\Phi$. 

Then, because $\varepsilon'\preceq_{k,c,\bowtie}\sqrt{\theta'}$, Lemma \ref{LemmaDominantOrientation} provides us with a dominant orientation. In combination with Lemma \ref{LemmaCloseGeodAreGood}, we get Lemma \ref{LemmaExistsPointForTetra}, which provides us with the vertical quadrilateral. 

Afterwards, we make use of them, as in the proof of Theorem \ref{ThmPartI}, to construct the quasi-isometry product $\hat{\Phi}$. In a box at scale $R$, the upper-bound $\varepsilon R$ on the distance between $\Phi$ and $\hat{\Phi}$ is achieved since $\theta'\leq \varepsilon$, and in our box at scale $L$, it is achieved since $\theta'\leq \varepsilon'$.

Finally, the exponents on $\theta$ of Theorem \ref{ThmPartII} can be removed since we can fix $\theta$, then do the proof with a parameter $\tilde{\theta}=\theta^8$, then replace $\tilde{\theta}$ by $\theta^8$.
\end{proof}

This is a step forward since now, Theorem \ref{ThmPartII} holds for all boxes at scale $L$, and not only a significant proportion of boxes at scale $R$.


\subsection{A quasi-isometry quasi-respects the height}

Let $p,q\in X\bowtie Y$ be such that $h(p)=h(q)$. In this section we are aiming the following theorem, which estimates the difference of height between the images of $p$ and $q$ under $\Phi$.

\begin{thm}\label{ThmControlonHeight}
For $0<\theta\preceq_{k,c,\bowtie}1$, there exists $M(k,c,\bowtie,\theta)$ (here $M$ depends also on $\theta$) such that for all $p$ and $q$ in $X\bowtie Y$ with $h(p)=h(q)$ we have
\begin{align}
\Delta h\big(\Phi(p),\Phi(q)\big)\leq \theta d_{\bowtie}(p,q)+M.
\end{align}
\end{thm}

By the previous section, we know that in a box of a sufficiently large scale, the quasi-isometry $\Phi$ is (on a set of almost full measure) close to a product map. We first show that this product map is coarsely an homothety along the height function. 
\\Let $L_0$ be the constant of Theorem \ref{ThmPartII}, let $L\geq L_0$ and let $\mathcal{B}$ be a box at scale $L$. Let us denote by $h^+:=\sup\{h(p)|p\in\mathcal{B}\}$ and by $h^-:=\inf\{h(p)|p\in\mathcal{B}\}$. Let $\hat{\Phi}:=( \hat{\Phi}^X,\hat{\Phi}^Y):X\bowtie Y \to X'\bowtie Y'$ be the corresponding product map of Theorem \ref{ThmPartII}.

\begin{lemma}\label{LemmaQiProdHomotHeight}
Let $a\in\mathcal{B}_{h^+}$ and $b\in\mathcal{B}_{h^-}$ be two points of $\mathcal{B}$, one on its top part and one on its bottom part. Then we have both:
\begin{align*}
\Big|\Delta h(\hat{\Phi}(a),\hat{\Phi}(b))-\frac{m}{m'}L\Big|&\preceq_{k,c,\bowtie}\theta^{\frac{1}{2}} L;
\\\Big|\Delta h(\hat{\Phi}(a),\hat{\Phi}(b))-\frac{n}{n'}L\Big|&\preceq_{k,c,\bowtie}\theta^{\frac{1}{2}} L.
\end{align*}
\end{lemma}  

\begin{proof}
Let $U\subset \mathcal{B}$ be the set involved in Theorem \ref{ThmPartII}, we recall that $\lambda(U)\geq(1-\theta)\lambda(\mathcal{B}) $ and that for all $p\in U$, $d_{\bowtie'}(\Phi(p),\hat{\Phi}(p))\preceq_{\bowtie,k,c}\theta L$. Since the measure $\lambda$ identically weights the level sets of $\mathcal{B}$, by a Markov inequality there exists $z^+\in [h^+-\theta^{\frac{1}{2}} L, h^+]$ and $z^-\in [h^-, h^-+\theta^{\frac{1}{2}} L]$ such that
\begin{align*}
\lambda_{z^+}(U_{z^+})&\geq (1-\theta^{\frac{1}{2}})\lambda_{z^+}(\mathcal{B}_{z^+});
\\\lambda_{z^-}(U_{z^-})&\geq (1-\theta^{\frac{1}{2}})\lambda_{z^-}(\mathcal{B}_{z^-}).
\end{align*}
By the definition of $\lambda_{z^+}$ we have that $\frac{1}{2}\mu_{z^+}(\mathcal{B}_{z^+})\leq\mu_{z^+}(U_{z^+})\leq \mu_{z^+}(\mathcal{B}_{z^+})$. Furthermore $\mu_{z^+}(\mathcal{B}_{z^+})\asymp_{\bowtie}e^{nL}e^{(m-n)|h^+-z^+|}$ since we went down by a height $|h^+-z^+|$ in the box. Therefore:
\begin{align*}
e^{nL}e^{(m-n)\theta^{\frac{1}{2}}L}\preceq_{k,c,\bowtie}\mu_{z^+}\big(\mathcal{N}_{kc+c}(U_{z^+})\big).
\end{align*}
Furthermore, $B_{z^+}$ resembles a rectangle of width $2|h^+-z^+|$ in $X$ and $2(L-|h^+-z^+|)$ in $Y$, hence we have:
\begin{align*}
\mu^Y_{z^+}\big(\mathcal{N}_{kc+c}(U^Y_{z^+})\big)\succeq_{k,c,\bowtie}e^{nL}e^{(m-n)\theta^{\frac{1}{2}}L}\frac{1}{\mu^X_{z^+}\big(\mathcal{B}^X_{z^+})}\succeq_{k,c,\bowtie}e^{nL}e^{2(m-n)\theta^{\frac{1}{2}}L}.
\end{align*}
By Lemma \ref{LemmaVolKeepByQI}, and since $\hat{\Phi}$ is close to $\Phi$ on $U$, we deduce:
\begin{align}
e^{nL}e^{2(m-n)\theta^{\frac{1}{2}}L}\preceq_{k,c,\bowtie}\mu^Y_{z^+}\big(\mathcal{N}_{1}(\hat{\Phi}(U^Y_{z^+}))\big).\label{IneqLastRepa1}
\end{align}
Let $\Delta\geq 0$ be $|h(\hat{\Phi}(U_{z^+}))-h(\hat{\Phi}(U_{z^-}))|$. For all $p\in U_{z^+}$ there exists a vertical geodesic $V_p$ of $\theta$-monotone image under $\Phi$ passing close to $p$. Furthermore, $d_Y\left(V_p^Y(z^-),U^Y_{z^{-}}\right)\leq 2\theta^\frac{1}{2}L$ since $\mathcal{B}^Y_{z^-}$ has a relatively small diameter. Therefore, all vertical geodesics starting at $\mathcal{N}_1\left(\hat{\Phi}^Y(U^Y_{z^+})\right)$ intersect $\mathcal{N}_{M\theta^{\frac{1}{2}} L}\left(\hat{\Phi}^Y(U^Y_{z^-})\right)$. Hence we have that
\begin{align*}
\mathcal{N}_1\left(\hat{\Phi}^Y(U^Y_{z^+})\right)\subset\pi_{h(\hat{\Phi}(U_{z^+}))}\Big(\mathcal{N}_{M\theta^{\frac{1}{2}} L}\left(\hat{\Phi}^Y(U^Y_{z^-})\right)\Big).
\end{align*}
Therefore:
\begin{align*}
\mu^{Y'}\big(\mathcal{N}_1(\hat{\Phi}^Y(U^Y_{z^+}))\big)&\leq \mu^{Y'}\Big(\pi_{h(\hat{\Phi}(U_{z^+}))}\Big(\mathcal{N}_{M\theta^{\frac{1}{2}} L}\left(\hat{\Phi}^Y(U^Y_{z^-})\right)\Big)\Big),
\\&\preceq_{k,c,\bowtie} e^{n'\Delta} \mu^{Y'}\Big(\mathcal{N}_{M\theta^{\frac{1}{2}} L}(\hat{\Phi}^Y(U^Y_{z^-}))\Big),
\\&\preceq_{k,c,\bowtie} e^{n'\Delta}e^{-n'M\theta^{\frac{1}{2}} L} \mu^{Y'}\Big(\mathcal{N}_{1}(\hat{\Phi}^Y(U^Y_{z^-}))\Big),\quad\text{by Corollary \ref{CoroUniformBallMeasure}},
\\&\preceq_{k,c,\bowtie} e^{n'\Delta}e^{n'M\theta^{\frac{1}{2}} L} e^{\theta^{\frac{1}{2}} L},\quad\text{because  }\mathcal{B}_{z^-}^Y\text{has small }\mu\text{ measure},
\\&=e^{n'\Delta}e^{n'(M+1)\theta^{\frac{1}{2}} L}.
\end{align*}
Combined with inequality (\ref{IneqLastRepa1}) we obtain
\begin{align*}
e^{nL}e^{2(m-n)\theta^{\frac{1}{2}}L}\preceq_{k,c,\bowtie}e^{n'\Delta}e^{n'(M+1)\theta^{\frac{1}{2}} L},
\end{align*}
which provides us with $e^{nL}\preceq_{k,c,\bowtie}e^{n'\Delta}e^{M'\theta^{\frac{1}{2}} L}$, where $M'$ is a constant depending on $k$, $c$, $\bowtie$ and $\bowtie'$. Then there exists $M''(k,c,\bowtie,\bowtie')$ such that by taking the logarithm we get
\begin{align*}
nL\leq n'\Delta+M''\theta^{\frac{1}{2}} L.
\end{align*}
Similarly, we do the same proof on $\Phi^{-1}$, on the box of height $\Delta$ containing $\hat{\Phi}(U_{z^+}\cup U_{z^-})$ which provides us with 
\begin{align*}
n'\Delta\leq nL+M''\theta^{\frac{1}{2}} L.
\end{align*}
Therefore $|\Delta -\frac{n}{n'}L|\preceq_{k,c,\bowtie}\theta^{\frac{1}{2}} L$. To obtain the same results with the constants $m$ and $m'$, we focus on the sets $U^X_{z^+}$ and $U^X_{z^-}$ instead of $U^Y_{z^+}$ and $U^Y_{z^-}$.
\end{proof}

As a corollary we obtain a first quasi-isometry invariant for horospherical products.

\begin{propo}\label{PropoQIinvariant}
If $X\bowtie Y$ and $X'\bowtie Y'$ are quasi-isometric, then $\frac{m}{n}=\frac{m'}{n'}$.
\end{propo}
\begin{proof}
By Lemma \ref{LemmaQiProdHomotHeight}, and by the triangle inequality we have that $|\frac{m}{m'}-\frac{n}{n'}|\preceq_{k,c,\bowtie} \frac{1}{L}$ for all $L\geq L_0$. Therefore, $\frac{m}{m'}=\frac{n}{n'}$, hence $\frac{m}{n}=\frac{m'}{n'}$.
\end{proof}

\begin{lemma}
Let $0<\theta\preceq_{k,c,\bowtie} 1$. Let $p:=(p^X,p^Y),q:=(q^X;q^Y)\in X\bowtie Y$ such that $d_{\bowtie}(p,q)\geq L_0^2$ and such that $p^Y=q^Y$ (hence $h(p)=h(q)$). Then we have:
\begin{align*}
\Delta h (\Phi(p),\Phi(q))\preceq_{k,c,\bowtie} \theta^{\frac{1}{2}}d_{\bowtie}(p,q).
\end{align*} 
\end{lemma}

\begin{proof}
Let $\mathcal{B}$ be a box of scale $L=d_{\bowtie}(p,q)$, such that $p$ and $q$ are contained in its bottom part. Let $V_p^X\in V\mathcal{B}^X$ be the vertical geodesic segment of $X$ of length $L$ starting at $p$. We apply Proposition \ref{Lemma411} on $V_p^X\bowtie \mathcal{B}^Y$ (as a box of an embedded copy of $\mathbb{R}\bowtie Y$ inside $X\bowtie Y$) with $r_0=L_0$ and $L\geq L_0^2$, we obtain that there exists $R\geq L_0$, a box tilling $\mathcal{B}\bigcup\limits_{i\in I} \mathcal{B}_i$ of boxes at scale $R$ and $I_g\subset I$ such that:
\begin{enumerate}
\item $\lambda\left(\bigcup\limits_{i\in I_g}\mathcal{B}_i\right)\geq (1-\theta) \lambda\left(\mathcal{B}\right)$ (Boxes indexed by $I_g$ cover almost all $\mathcal{B}$);
\item $\forall i\in I_g$, $\dfrac{\eta_i\big(V^b\mathcal{B}_i\big)}{\eta_i\big(V\mathcal{B}_i\big)}\leq \theta$ (almost all vertical geodesic segments in $\mathcal{B}_i$ have $\varepsilon$-monotone image),
\end{enumerate}
where $\eta_i:=\eta_{V\mathcal{B}_i}$. In this setting, we have that $V\mathcal{B}_i:=\big\lbrace(V^X_p,V^Y)|V^Y\in V\mathcal{B}_i^Y \big\rbrace$, hence most vertical geodesics in $\mathcal{B}_i^Y$ are a good vertical geodesic of $\mathcal{B}_i$ when coupled with a portion of $V^X_p$. 

Let us denote by $J:=\{0,...,\frac{L}{R}-1\}$ and for all $j\in J$ let us denote by $p^X_j:=V^X_p(jR)$, then $V^X_p(jR):=\bigcup\limits_{j\in J}[p^X_j;p^X_{j+1}]$.

Since the measure of the good boxes cover almost all $\mathcal{B}$, and because the measure $\lambda$ equally weights the level sets, by a Markov inequality argument there exists $J_g\subset J$ such that for all $j\in J_g$, $\mathcal{B}_{[jR;(j+1)R]}$ is almost entirely covered by boxes of $I_g$. Therefore, again by Markov inequality argument, there exists $W_p\subset V^Y\mathcal{B}$ such that:
\begin{enumerate}
\item $\eta(W_p)\geq (1-\theta^{\frac{1}{2}})\eta(V^Y\mathcal{B})$;
\item $\forall V^Y\in W_p$ and $\forall j\in J_g$ we have  $V^Y([-(j+1)R; -jR])\in\bigcup\limits_{i\in I_g}V^g\mathcal{B}^Y_i$.
\end{enumerate}
Let $V^Y\in W_p$, for all $j\in J$ let us denote by $p_i:=\left(V^X_p(jR),V^Y(-jR)\right)$. By Lemma \ref{LemmaQiProdHomotHeight}, for all $j\in J$ we have
\begin{align}
|\Delta h \big(\hat{\Phi}(p_j),\hat{\Phi}(p_{j+1})\big)-\frac{m}{m'}R|\preceq_{k,c,\bowtie}\theta^{\frac{1}{2}} R.\label{IneqClaFin1}
\end{align}
For all $j\in J_g$, let us denote by $\mathcal{B}_j$ the box at scale $R$ containing $[p_j;p_{j+1}]$. By the choice of $R$, we have that most vertical geodesic segments of $\mathcal{B}_j^Y$ have $\theta$-monotone image when coupled with $\left[p_j^X;p_{j+1}^X\right]$. 

Furthermore $\mathcal{B}_j$ contains almost only good vertical geodesic segments, therefore, there exists $v\in V^g\mathcal{B}_j$ such that $\left(\left[p_j^X;p_{j+1}^X\right],v^Y\right)\in V^g\mathcal{B}_j$ and such that $\left(v^X,\left[p_j^Y;p_{j+1}^Y\right]\right)\in V^g\mathcal{B}_j$. 
Therefore there exists a good coarse vertical quadrilateral containing $p_j$ and $p_{j+1}$, hence $d\left(\Phi(p_j),\hat{\Phi}(p_j)\right)\preceq_{k,c,\bowtie}\theta R$. Similarly we have $d\left(\Phi(p_{j+1}),\hat{\Phi}(p_{j+1})\right)\preceq_{k,c,\bowtie}\theta R$. Hence combined with inequality (\ref{IneqClaFin1}) we get
\begin{align*}
|\Delta h \big(\Phi(p_j),\Phi(p_{j+1})\big)-\frac{m}{m'}R|\preceq_{k,c,\bowtie}\theta^{\frac{1}{2}} R.
\end{align*} 
Therefore by the triangle inequality, there exists $M(k,c,\bowtie)$ and $M'(k,c,\bowtie)$ such that:
\begin{align*}
\Delta h(\Phi(p_0),\Phi(p_{\frac{L}{R}-1})&\leq \sum\limits_{j=0}^{\frac{L}{R}-1}\Delta h(\Phi(p_j),\Phi(p_{j+1})),
\\&\leq \sum\limits_{j\in J_g}\Delta h(\Phi(p_j),\Phi(p_{j+1}))+ \sum\limits_{j\in J\setminus J_g}\Delta h(\Phi(p_j),\Phi(p_{j+1})),
\\&\leq \#J_g\left(\frac{m}{m'}R+M\theta^{\frac{1}{2}} R\right) + \#(J\setminus J_g)(kR+c),
\\&\leq \frac{L}{R}\left(\frac{m}{m'}R+M\theta^{\frac{1}{2}} R\right) + \theta\frac{L}{R}(kR+c),
\\&\leq \frac{m}{m'}L+M'\theta^{\frac{1}{2}} L.
\end{align*}
Similarly we have $\Delta h(\Phi(p_0),\Phi(p_{\frac{L}{R}-1})\geq \frac{m}{m'}L-M'\theta^{\frac{1}{2}} L$. By doing the same reasoning on $q$ we have that for all $V^Y\in W_q$, $|\Delta h(\Phi(q_0),\Phi(q_{\frac{L}{R}-1})-\frac{m}{m'}L|\preceq_{k,c,\bowtie}\theta^{\frac{1}{2}} L$, where $q_j:=(V^X_q(jR),V^Y(-jR))$. Furthermore $W_p\cap W_q$ is non empty for $\theta^{\frac{1}{2}}\preceq_{k,c,\bowtie} 1$, then let $V^Y\in W_p\cap W_q$. Without loss of generality we can assume that $\Phi(p)\geq\Phi(q)$, we have:
\begin{align*}
\Delta h (\Phi(p),\Phi(q))&= h(\Phi(p))-h(\Phi(p_0))+h(\Phi(p_0))-h\left(\Phi\left(p_{\frac{L}{R}-1}\right)\right)+h\left(\Phi\left(p_{\frac{L}{R}-1}\right)\right)
\\&-h\left(\Phi\left(q_{\frac{L}{R}-1}\right)\right)+h\left(\Phi\left(q_{\frac{L}{R}-1}\right)\right)-h(\Phi(q_0))+h(\Phi(q_0))-h(\Phi(q)),
\\&\preceq_{k,c,\bowtie} d_{\bowtie}(p,p_0)-\frac{m}{m'}L+M\theta^{\frac{1}{2}} L+ d_{\bowtie}\left(p_{\frac{L}{R}-1},q_{\frac{L}{R}-1}\right)+\frac{m}{m'}L+M\theta^{\frac{1}{2}} L+ d_{\bowtie}(q,q_0),
\\&\preceq_{k,c,\bowtie} d_{\bowtie}(p,p_0)+ d_{\bowtie}\left(p_{\frac{L}{R}-1},q_{\frac{L}{R}-1}\right)+ d_{\bowtie}(q,q_0)+2\theta^{\frac{1}{2}} L.
\end{align*}
However $d_{\bowtie}(p,p_0)\leq M_0$ since they share the same $X$ coordinate and because the top part of $\mathcal{B}^Y$ as a diameter of at most $M_0$, similarly $d_{\bowtie}(q,q_0)\leq M_0$. By construction $p_{\frac{L}{R}-1}^Y=q_{\frac{L}{R}-1}^Y$, furthermore the top part of $\mathcal{B}^X$ has a diameter of at most $M_0$, hence $d_{\bowtie}\left(p_{\frac{L}{R}-1},q_{\frac{L}{R}-1}\right)\leq M_0$. Finally we obtain:
\begin{align*}
\Delta h (\Phi(p),\Phi(q))\preceq_{k,c,\bowtie} \theta^{\frac{1}{2}} L=\theta^{\frac{1}{2}}d_{\bowtie}(p,q).
\end{align*} 
\end{proof}

\begin{cor}\label{CoroAllVertRayAreGood}
Any vertical geodesic ray $V$ of $X\bowtie Y$ satisfies, for all $t_1,t_2\in\mathbb{R}$
\begin{align*}
h\left(\Phi\circ V(t_1)\right)=h\left(\Phi\circ V(t_2)\right)\quad\Rightarrow\quad |t_1-t_2|\preceq_{k,c,\bowtie} 1.
\end{align*}
\end{cor}
\begin{proof}
Suppose $V$ is a vertical geodesic segment parametrised by arclength. Suppose $0<t_1<t_2$ are such that $h\big(\Phi(V(t_1))\big)=h\big(\Phi(V(t_2))\big)$. We apply Theorem \ref{ThmControlonHeight} on $\Phi^{-1}$ with $p=\Phi(V(t_1))$, $q=\Phi(V(t_2))$, where $\theta$ is here fixed and depends only on $k,c$ and the metric measured space $(X\bowtie Y,d_{\bowtie})$. Hence there exists $M(k,c,\bowtie)>0$ such that 
\begin{align}
\Delta h\big(V(t_1),V(t_2)\big)\preceq_{k,c,\bowtie} \theta^{\frac{1}{2}}|t_1-t_2|+M.
\end{align}
However $\Delta h\big(V(t_1),V(t_2)\big)=|t_1-t_2|$, hence
\begin{align*}
\left(1-\theta^{\frac{1}{2}}\right)\big|t_1-t_2\big|\preceq_{k,c,\bowtie}1.
\end{align*}
Therefore $\big|t_1-t_2\big|\preceq_{k,c,\bowtie}1$ since $\theta^{\frac{1}{2}}\leq \frac{1}{2}$.
\end{proof}

This is stronger than being $\varepsilon$-monotone since it holds for all $t_1$, $t_2\in\mathbb{R}$.

\subsection{Factorisation of a quasi-isometry on the whole space}

Finally, we provide the proof of the Theorem \ref{ThmProd}, which states that $\Phi$ is close to a product map $\hat{\Phi}$ on the whole space $X\bowtie Y$.

\begin{proof}[Proof of Theorem \ref{ThmProd}]
We first pick an arbitrary vertical geodesic $V_0^X$ of $X$ and an arbitrary vertical geodesic $V_0^Y$ of $Y$. Then we work with the two embedded copies $X_0:= X\bowtie V_0^Y$ and $Y_0:= V_0^X\bowtie Y$ of $X$ and $Y$ in $X\bowtie Y$. Let $p\in X\bowtie Y$, there exist a unique $a\in X_0$ and a unique $b\in Y_0$ such that $p^X=a^X$ and $p^Y =b^Y$. We can construct a coarse vertical quadrilateral $Q$ containing $p$ and $a$ as in Lemma \ref{LemmaExistsPointForTetra}. Thanks to Corollary \ref{CoroAllVertRayAreGood}, we know that $\Phi(Q)$ is in the $M(k,c,\bowtie)$-neighbourhood of a coarse vertical quadrilateral $Q'$ on which we use Proposition \ref{LemmaTetrahedron}. This gives us
\begin{align}
d_{X'}\left(\Phi(p)^{X'},\Phi(a)^{X'}\right)&\preceq_{k,c,\bowtie}1;\label{IneqFin1}
\\\Delta h\left(\Phi(p)^{X'},\Phi(a)^{X'}\right)&\preceq_{k,c,\bowtie}1.\label{IneqFin2}
\end{align}
Similarly we have $d_{Y'}\left(\Phi(p)^{Y'},\Phi(b)^{Y'}\right)\preceq_{k,c,\bowtie}1$. Let us denote
\begin{align*}
\hat{\Phi}^X:X&\rightarrow X'
\\ x&\mapsto \Phi\left(x,V_0^Y(-h(x))\right)^{X'}.
\end{align*}
By rewriting inequality (\ref{IneqFin1}) we have
\begin{align*}
d_{X'}\left(\Phi(p)^{X'},\hat{\Phi}^{X}\left(p^X\right)\right)=&d_{X'}\left(\Phi(p)^{X'},\hat{\Phi}^X\left(a^X\right)\right)=d_{X'}\left(\Phi(p)^{X'},\Phi\left(a^X,V_0^Y(-h(a^X))\right)^{X'}\right),
\\=&d_{X'}\left(\Phi(p)^{X'},\Phi\left(a\right)^{X'}\right)\preceq_{k,c,\bowtie}1.
\end{align*}
Similarly by denoting $\hat{\Phi}^Y:=\Phi\left(V_0^X(-h(y)),y\right)^{Y'}$ for all $y\in Y$, we have
\begin{equation}
d_Y\left(\Phi(p)^Y,\hat{\Phi}^Y\left(p^Y\right)\right)\preceq_{k,c,\bowtie}1.
\end{equation} 
The last problem is that given a point $p$, the heights of $\hat{\Phi}^X\left(p^X\right)$ and $\hat{\Phi}^Y\left(p^Y\right)$ may differ. As in the proof of Theorem \ref{ThmPartI}, inequality (\ref{IneqFin2}) guaranties that they are sufficiently close, which allows us to chose $\hat{\Phi}^X$ and $\hat{\Phi}^Y$ such that $\hat{\Phi}:=(\hat{\Phi}^X,\hat{\Phi}^Y)$ is a well defined product map on $X\bowtie Y$. Then we have
\begin{align*}
d_{\bowtie'}\left(\Phi(p),\hat{\Phi}(p)\right)&\preceq_{k,c,\bowtie}1;
\\\Delta h\left(\Phi(p),\hat{\Phi}(p)\right)&\preceq_{k,c,\bowtie}1.
\end{align*}
We now prove that $\hat{\Phi}^X$ and $\hat{\Phi}^X$ are quasi-isometries. Let $x,x'\in X$, then
\begin{align*}
d_{X'}\left(\hat{\Phi}^X(x),\hat{\Phi}^X(x')\right)&\preceq_{k,c,\bowtie} d_{X'}\left(\Phi\left(x,V_0^Y(-h(x))\right)^{X'},\Phi\left(x',V_0^Y(-h(x'))\right)^{X'}\right) ,
\\&\leq d_{\bowtie'}\left(\Phi\left(x,V_0^Y(-h(x))\right),\Phi\left(x',V_0^Y(-h(x'))\right)\right) ,
\\&\leq kd_{\bowtie}\left(\left(x,V_0^Y(-h(x))\right),\left(x',V_0^Y(-h(x'))\right)\right) +c,
\\&\leq kd_X(x,x')+kd_Y(V_0^Y(-h(x)),V_0^Y(-h(x'))) +c+M(k,c,\bowtie),\quad\text{by Corollary \ref{lengthGeod}},
\\&\leq kd_X(x,x')+\Delta h(x,x') +c+M\leq (k+1)d_X(x,x')+c+M.
\end{align*}
Similarly, and because $d_{\bowtie'}\geq \frac{d_{X'}+d_{Y'}}{2}$,
\begin{align*}
&d_{X'}\left(\hat{\Phi}^X(x),\hat{\Phi}^X(x')\right)
\\&= d_{X'}\left(\Phi\left(x,V_0^Y(-h(x))\right)^{X'},\Phi\left(x',V_0^Y(-h(x'))\right)^{X'}\right) ,
\\&\geq 2d_{\bowtie'}\left(\Phi\left(x,V_0^Y(-h(x))\right),\Phi\left(x',V_0^Y(-h(x'))\right)\right) -d_{Y'}\left(\Phi\left(x,V_0^Y(-h(x))\right)^{Y'},\Phi\left(x',V_0^Y(-h(x'))\right)^{Y'}\right) ,	
\\&\geq \frac{1}{k} d_{X}\left(x,x'\right)-c-d_{Y'}\left(\hat{\Phi}^Y\left(V_0^Y(-h(x))\right),\hat{\Phi}^Y\left(V_0^Y(-h(x))\right)\right) -2M,\quad\text{by the triangle inequality},
\\&\geq \frac{1}{k} d_{X}\left(x,x'\right)-c-2M.
\end{align*}
The proof that $\hat{\Phi}^Y$ is a quasi-isometry is similar.
\end{proof}

\section{Some solvable Lie groups as horospherical products}\label{SecExample}

In this chapter, we provide a characterisation of the quasi-isometry group of the horospherical product of two Heintze groups. See Theorem \ref{ThmBilip} for the precise description. 

\subsection{Admissibility of Heintze groups}

In this section we show that a Heintze group satisfies the conditions required to apply our main rigidity result \ref{ThmProd}. 
\begin{defn}(Heintze group)
\\A Heintze group is a solvable Lie group $S=N\rtimes_A \mathbb{R}$ where $N$ is a connected, simply connected, nilpotent Lie group, and $A$ is a derivation of $\mathrm{Lie}(N)$ whose eigenvalues all have positive real parts. 
\end{defn}

Heintze obtained in his work \cite{Heintze} that any negatively curved homogeneous manifold is isometric to a Heintze group.

\begin{rem}
A Heintze group admits a left-invariant metric with strictly negative sectional curvature, see \cite{Heintze} for further details. From now on we fix $g$ a left-invariant metric on $N\rtimes_A \mathbb{R}$ with maximal sectional curvature $-1$. Since $N\rtimes_A \mathbb{R}$ is simply connected, it is a $CAT(-1)$-space. 
\end{rem}

From now on we fix the metric $g$ such that $S=N\rtimes_A\mathbb{R}$ is a $\mathrm{CAT}(-1)$ space. Therefore S is a $\delta$-hyperbolic, Busemann, proper, geodesically complete metric space. Moreover, we show that S satisfies all three assumptions of Definition \ref{DefAdmissible2}. 
The assumption $(E1)$ holds thanks to the decomposition $S=N\rtimes_A\mathbb{R}$. We have for all $(n,z)\in N\rtimes_A\mathbb{R}$, $g_{(n,z)}=\exp(-zA)(g_{N})_n\exp(-zA)^{t}\oplus\mathrm{d}z^2$, where $g_N$ is the restriction of $g$ to the Lie algebra of $N$. 
Let us denote by $g_z:=\exp(-zA)g_{N}\exp(-zA)^{t}$ a left invariant metric on $N$, then let us denote by $\mu:=\mu_g$ the measure on $S$ induced by $g$ and by $\mu_z:=\mu_{g_z}$ the measure on $N$ induced by $g_z$. Then for all measurable subset $U\subset S$ we have
\begin{align*}
\mu(U)&:=\int\limits_{S}\mathbbm{1}_U(n,z)\mathrm{d}\mu_g(n,z)=\int\limits_{\mathbb{R}}\int\limits_{N}\mathbbm{1}_U(n,z)\mathrm{d}\mu_{g_{z}}(n)\mathrm{d}z,
\\&=\int\limits_{\mathbb{R}}\mu_z(U_z)\mathrm{d}z,
\end{align*}
where $U_z:=\{n\in N\vert (n,z)\in U\}$.  
Assumption $(E2)$ holds with constant $M_0=1$ since $g_{n,z}$ is left-invariant, and assumption $(E3)$ arises from the fact that $\mathrm{det}(g_z)=\exp(-2z\cdot\mathrm{tr}(A))\mathrm{det}(g)$. Therefore, any Heintze group is an admissible horo-pointed space.  Let us denote $S_1:=N_1\rtimes_{A_1}\mathbb{R}$ and $S_2:=N_2\rtimes_{A_2}\mathbb{R}$, then
\begin{align*}
S_1\bowtie S_2=(N_1\times N_2)\rtimes_A\mathbb{R},
\end{align*}
with $A$ the matrix $\mathrm{diag}(A_1,-A_2)$. Similarly let us denote by $S_1':=N_1'\rtimes_{A_1'}\mathbb{R}$ and $S_2':=N_2'\rtimes_{A_2'}\mathbb{R}$ two Heintze groups, with $N_1'$, $N_2'$ being two simply connected nilpotent Lie groups and $A_1'$, $A_2'$ being two derivations.

\subsection{Precision on the components of the product map}

We first refine Theorem \ref{ThmProd} for Heintze groups.

\begin{rem}
For any vertical geodesics $V$ of $(N_1\times N_2)\rtimes_A\mathbb{R}$ there exist $n_1\in N_1$, $n_2\in N_2$ and an arclength parametrisation of $V$ such that $V(t)=(n_1,n_2,t)$. 
\end{rem}

Let $\Phi:(N_1\times N_2)\rtimes_A\mathbb{R}\to (N_1'\times N_2')\rtimes_{A'}\mathbb{R}$ be a $(k,c)$-quasi-isometry. Let us assume that $\mathrm{tr}(A_1)>\mathrm{tr}(A_2)$ and that $\mathrm{tr}(A_1')>\mathrm{tr}(A_2')$. By Theorem \ref{ThmProd} there exist $\hat{\Phi}_1:S_1\to S_1'$ and $\hat{\Phi}_2:S_2\to S_2'$ such that
\begin{align*}
d_{\bowtie}\big(\Phi,(\hat{\Phi}_1,\hat{\Phi}_2))\preceq_{k,c,\bowtie}1.
\end{align*}

\begin{lemma}\label{LemmaVertsendonVertLie}
Let $i\in\{1,2\}$, then for any vertical geodesic $V\in S_i$, there exists a vertical geodesic $V'\in S_j'$ such that
\begin{align*}
d_{\mathrm{Hff}}\left(\hat{\Phi}_i(V),V'\right)\preceq_{k,c,\bowtie}1.
\end{align*}
\end{lemma}

This lemma also holds for any horospherical product where our main result, the geometric rigidity, applies.

\begin{proof}

Since $S_i=N_i\rtimes_{A_i}\mathbb{R}$ is a Gromov hyperbolic space, there exists $M(k,c,\bowtie)$ such that the image of a vertical geodesic by $\hat{\Phi}_i$ is in a $M$-neighbourhood of a geodesic $\gamma$ of $S_i'$. By Corollary \ref{CoroAllVertRayAreGood} $\gamma$ is a vertical geodesic, hence for $V':=\gamma$ we have $d_{\mathrm{Hff}}\left(\hat{\Phi}_i(V),V'\right)\preceq_{k,c,\bowtie}1$.
\end{proof}

Let $n\in N_i$ and let us denote by $V_n$ the vertical geodesic $V_n~:~\mathbb{R}\to S_i~;~t\mapsto(n,t)$. By Lemma \ref{LemmaVertsendonVertLie} there exists a vertical geodesic $V_n'$ such that
\begin{align}
d_{\mathrm{Hff}}(\hat{\Phi}_i(V_n),V_n')\preceq_{k,c,\bowtie} 1\label{IneqEsayIci0}
\end{align}
Furthermore $V_n'$ is unique since it is an infinite geodesic of the Heintze group $S_i$. We define a map $\Psi_i:N_i\to N_i'$ as the following
\begin{equation}
\text{For all }n\in N_i\ ,\ \Psi_i(n)= P\left(V_n'(0)\right),
\end{equation}  
where $P:N_i'\rtimes_{A_i} \mathbb{R}\to N_i'$ is the natural projection on $N_i$. 
\\\\The goal of this subsection is to prove the following theorem.

\begin{thm}\label{ThmProduct}
There exists $t_0\in\mathbb{R}$ such that for the aforementioned $\Psi_i$ we have
\begin{align*}
d_{\bowtie}\left(\Phi,\left(\Psi_1,\Psi_2,\frac{\mathrm{tr}(A_1)}{\mathrm{tr}(A_1')} id_{\mathbb{R}}+t_0\right)\right)\preceq_{k,c,\bowtie}1.
\end{align*}
\end{thm} 
We can replace $\frac{\mathrm{tr}(A_1)}{\mathrm{tr}(A_1')}$ by $\frac{\mathrm{tr}(A_2)}{\mathrm{tr}(A_2')}$ thanks to Proposition \ref{PropoQIinvariant}. We first show $\hat{\Phi}_i$ and $\Psi_i$ are related.

\begin{lemma}\label{LemmaQiProd}
Let $i\in\{1,2\}$. There exists $f_i:\mathbb{R}\to\mathbb{R}$ such that for all $(n,t)\in S_i$
\begin{align*}
d_{S_i}\Big(\hat{\Phi}_i(n,t),\big(\Psi_i(n),f_i(t)\big)\Big)\preceq_{k,c,\bowtie}1.
\end{align*}  
\end{lemma}
\begin{proof}
Let $f_i:\mathbb{R}\to\mathbb{R};t \mapsto h(\hat{\Phi}_i(e_{N_i},t))$. Then by Theorem \ref{ThmProd} we have that $h\big(\hat{\Phi}_i(n,t)\big)=f_i(t)$ for all $n\in N_i$. Therefore by the definition of $\Psi_i$ we have $\left(\Psi_i(n),f_i(t)\right)=V_n'(f_i(t))$. Hence
\begin{align}
d_{S_i'}\left(\hat{\Phi}_i(n,t),(\Psi_i(n),f_i(t))\right)&= d_{S_i'}\left(\hat{\Phi}_i(n,t),V_n'(f_i(t))\right).\label{IneqEsayIci1}
\end{align}
However by inequality (\ref{IneqEsayIci0}), there exists $s_t\in\mathbb{R}$ such that 
\begin{equation}
d_{S_i'}\left(\hat{\Phi}_i(n,t),V_n'(s_t)\right)\preceq_{k,c,\delta}.1\label{IneqEsayIci3}
\end{equation}
Furthermore we know that
\begin{align}
1\succeq_{k,c,\bowtie }d_{S_i'}\left(\hat{\Phi}_i(n,t),V_n'(s_t)\right)\geq \Delta h\left(\hat{\Phi}_i(n,t),V_n'(s_t)\right) =|f_i(t)-s_t|.\label{IneqEsayIci2}
\end{align}
Therefore
\begin{align*}
d_{S_i'}\left(\hat{\Phi}_i(n,t),V_n'(f_i(t))\right)&\leq d_{S_i'}\left(\hat{\Phi}_i(n,t),V_n'(s_t)\right)+d_{S_i'}\left(V_n'(s_t),V_n'(f_i(t))\right)\quad,\text{by the triangle inequality,}
\\&= d_{S_i'}\left(\hat{\Phi}_i(n,t),V_n'(s_t)\right)+|f_i(t)-s_t|\preceq_{k,c,\bowtie} 1\quad,\text{by inequalities (\ref{IneqEsayIci3}) and (\ref{IneqEsayIci2})}.
\end{align*}
Combined with equality (\ref{IneqEsayIci1}) it provides us with $d_{S_i'}\left(\hat{\Phi}_i(n,t),(\Psi_i(n),f_i(t))\right)\preceq_{k,c,\bowtie}1$.
\end{proof}

\begin{cor}(Quasi-isometries quasi-preserve the horosphere volume)\label{LemmaQIpreserveHoroVol}
\\Let $t\in\mathbb{R}$, $r>0$ and $n\in N_i$. Then the map $\tilde{\Phi}_i:=(\Psi_i,f_i)$ quasi-preserves the volume of any disk $D:=D_r(n,t)$
\begin{align*}
\mu^{S_i}_t \left(D\right) \asymp_{k,c,\bowtie}\mu^{S_i'}_{f_i(t)} \left(\mathcal{N}_1\left(\tilde{\Phi}_i(D)\right)\right).
\end{align*}
\end{cor}

\begin{proof}
By Lemma \ref{LemmaQiProd}, there exists $M(k,c,\bowtie)$ such that $\tilde{\Phi}_i$ is $M$-close to $\hat{\Phi}_i$. Therefore, there exists $k',c'$ depending only on $k,c$ and $S_1\bowtie S_2$ such that $\tilde{\Phi}_i$ is a $(k',c')$-quasi-isometry.
\\We first pick a $2k'(c' + 1)$-maximal separating set $Z$ of D. Then $\tilde{\Phi}_i(Z)$ verifies:
\begin{enumerate}
\item The disks $D_1(p)$ with $p\in \tilde{\Phi}_i(Z)$ are pairwise disjoints.
\item $\bigcup\limits_{p\in \tilde{\Phi}_i(Z)} D_{1}(p)\subset\mathcal{N}	_1\left(\tilde{\Phi}_i(D)\right)\subset \bigcup\limits_{p\in \tilde{\Phi}_i(Z)} D_{2k'\cdot k'(c'+1)+c'+1}(p)$.
\end{enumerate}
Furthermore by Lemma \ref{LemmaMuDiskAreExp}, we have $\forall (n,t)\in Z$
\begin{align*}
\mu^{S_i}_t\big(D_{k'(c'+1)}(n,t)\big)\asymp_{k,c,\bowtie}1\asymp_{k,c,\bowtie}\mu^{S_i}_t\big(D_{2k'(c'+1)}(n,t)\big).
\end{align*}
Hence $\mu^{S_i}_t \left(D\right) \asymp_{k,c,\bowtie}\# Z$. Furthermore, by Lemma \ref{LemmaMuDiskAreExp} we also have $\forall (n,t)\in Z$
\begin{align*}
\mu^{S_i'}_{f_i(t)}\big(D_{1}\big(\Phi_i(n,t)\big)\big)\asymp_{k,c,\bowtie}1\asymp_{k,c,\bowtie}\mu^{S_i'}_{f_i(t)}\big(D_{2k'\cdot k'(c'+1)+c'+1}\big(\Phi_i(n,t)\big)\big).
\end{align*}
Therefore
\begin{align*}
\mu^{S_i}_t \left(D\right) \asymp_{k,c,\bowtie}\# Z\asymp_{k,c,\bowtie}\mu^{S_i'}_{f_i(t)} \left(\mathcal{N}_1\left(\tilde{\tilde{\Phi}}_i(D)\right)\right).
\end{align*}
\end{proof}

\begin{lemma}\label{LemmaTimeNotDistorded}(Quasi-isometries quasi-translate the height)
\\Let $f_i:\mathbb{R}\to \mathbb{R}$ be the function involved in Lemma \ref{LemmaQiProd}. Then for all $t\in\mathbb{R}$
\begin{align*}
\left|\frac{\mathrm{tr}(A_1)}{\mathrm{tr}(A_1')}t-(f_i(t)-f_i(0))\right|\preceq_{k,c,\bowtie}1.
\end{align*}
\end{lemma}

\begin{proof}
We recall that for all $t\in\mathbb{R}$, $f_i(t):=h\left(\hat{\Phi}_i(e_{N_i},t)\right)$. Let $n\in N_i$, $r>0$, $t\in\mathbb{R}$, and let us denote $U\subset N_i$ such that $D_r(n,0)=(U,0)$. Then we have
\begin{align}
\mu^{S_i}_0(U,0)=e^{2\mathrm{tr}(A_i)t}\mu^{S_i}_t(U,t).\label{EqComparU}
\end{align}
However $\tilde{\Phi}_i(U,0)=\big(\Psi_i(U),f_i(0)\big)$ and $\tilde{\Phi}_i(U,t)=\big(\Psi_i(U),f_i(t)\big)$, therefore
\begin{align}
\mu^{S_i'}_{f_i(0)}\left(\mathcal{N}_1\big(\tilde{\Phi}_i(U,0)\big)\right)&=\mu^{S_i'}_{f_i(0)}\left(\mathcal{N}_1\big(\Psi_i(U),f_i(0)\big)\right),\nonumber
\\&=e^{2\mathrm{tr}(A_i')(f_i(t)-f_i(0))}\mu^{S_i'}_{f_i(t)}\left(\mathcal{N}_1\big(\Psi_i(U),f_i(t)\big)\right),\nonumber
\\&=e^{2\mathrm{tr}(A_i')(f_i(t)-f_i(0))}\mu^{S_i'}_{f_i(t)}\left(\mathcal{N}_1\big(\tilde{\Phi}_i(U,t)\big)\right).\label{EqComparPhiU}
\end{align}
Furthermore by Lemma \ref{LemmaQIpreserveHoroVol} we have
\begin{align*}
\mu^{S_i}_0(U,0)&\asymp_{k,c,\bowtie}\mu^{S_i'}_{f_i(0)}\left(\mathcal{N}_1\big(\tilde{\Phi}_i(U,0)\big)\right);
\\\mu^{S_i}_t(U,t)&\asymp_{k,c,\bowtie}\mu^{S_i'}_{f_i(t)}\left(\mathcal{N}_1\big(\tilde{\Phi}_i(U,t)\big)\right).
\end{align*}
In combination with equalities (\ref{EqComparU}) and (\ref{EqComparPhiU}), it provides us with
\begin{align*}
\mu^{S_i}_0(U,0)&=e^{2\mathrm{tr}(A_i)t}\mu_t(U,t)\asymp_{k,c,\bowtie}e^{2\mathrm{tr}(A_i)t}\mu^{S_i'}_{f_i(t)}\left(\mathcal{N}_1\big(\tilde{\Phi}_i(U,t)\big)\right),
\\&=e^{2\mathrm{tr}(A_i)t} e^{2\mathrm{tr}(A_i')(f_i(0)-f_i(t))}\mu^{S_i'}_{f_i(0)}\left(\mathcal{N}_1\big(\tilde{\Phi}_i(U,0)\big)\right),
\\&\asymp_{k,c,\bowtie}e^{2\mathrm{tr}(A_i)t} e^{2\mathrm{tr}(A_i)(f_i(0)-f_i(t))}\mu^{S_i}_0(U,0).
\end{align*}
Hence we have $e^{2\mathrm{tr}(A_i)t}\asymp_{k,c,\bowtie}e^{2\mathrm{tr}(A_i')(f_i(t)-f_i(0))}$, which, composed with the logarithm, gives us
\begin{equation}
\left|\frac{\mathrm{tr}(A_1)}{\mathrm{tr}(A_1')}t-(f_i(t)-f_i(0))\right|\preceq_{k,c,\bowtie}1.
\end{equation}
\end{proof}

\begin{cor}\label{CorFin}
There exists $t_0\in\mathbb{R}$ such that for $i\in\{1,2\}$ and for all $(n,t)\in N_i\times \mathbb{R}$
\begin{align*}
d_{S_i}\left(\hat{\Phi}_i(n,t),\left(\Psi_i(n),\frac{\mathrm{tr}(A_1)}{\mathrm{tr}(A_1')}t+t_0\right)\right)\preceq_{k,c,\bowtie}1.
\end{align*}  

\end{cor}
\begin{proof}
The proof is a direct application of Lemmas \ref{LemmaQiProd} and \ref{LemmaTimeNotDistorded} by taking $t_0:=f_i(0)$.
\end{proof}

In this corollary $t_0$ depends on $\Phi$. 

\begin{proof}[Proof of Theorem \ref{ThmProduct}]
Using Corollary \ref{CorFin} on $N_1$ and $N_2$ provides us with Theorem \ref{ThmProduct}.
\end{proof}

\subsection{Hamenstädt distance and product maps of bilipschitz maps.}

As presented in section 5.3 of \cite{CKLNO}, the parabolic visual boundary of $N_i\rtimes \mathbb{R}$ may be identified with the Lie group $N_i$ endowed with the following $A_i$-homogeneous Hamenstädt distance.

\begin{defn}\label{DefHamen}(Hamenstädt distance)
For any $n,m\in N_i$, we define their Hamenstädt distance as
\begin{align*}
d_{N_i,A_i,H}(n,m):=\exp\left(-\frac{1}{2}\lim\limits_{s\rightarrow+\infty}\Big(2s-d_{N_i\rtimes_{A_i}\mathbb{R}}\big((n,-s),(m,-s)\big)\Big)\right).
\end{align*}
\end{defn}
We might omit $A_i$ and $N_i$ in the notation. We denote $\mathrm{Bilip}\big(N\big)$ the group of bilipschitz maps of $N$ for the Hamenstädt distance.
\begin{align*}
\mathrm{Bilip}(N_i):=\left\lbrace \Psi :(N_i,d_H)\to (N_i,d_H) ~|~ \exists k\geq 1, \Psi\text{ is a }(k,0)\text{-quasi-isometry}\right\rbrace.
\end{align*}
This is indeed a distance when the left invariant metric $g$ is normalized so that $\mathbb{R}\ltimes_{A_i} N_i$ is a $\mathrm{CAT}(-1)$ space.

Two quasi-isometries $\Phi$ and $\Phi'$ are said to be equivalent when they are at finite distance from each other.
\begin{align*}
\Phi\sim \Phi'\quad\Leftrightarrow\quad \sup\limits_{x}d_{\bowtie}\big(\Phi(x),\Phi'(x)\big)<+\infty.
\end{align*}
In this section we prove the following characterisation of the quasi-isometry group of $S_1\bowtie S_2=(N_1\times N_2)\rtimes_A\mathbb{R}$.

\begin{thm}\label{ThmQISolvable}
Let $N_1\rtimes_{A_1}\mathbb{R}$ and $N_2\rtimes_{A_2}\mathbb{R}$ be two Heintze group such that $\mathrm{tr}(A_1)\neq \mathrm{tr}(A_1)$, let $\Phi\in \mathrm{QI}\big((N_1\times N_2)\rtimes_A\mathbb{R}\big)$ and let $\Psi_1$, $\Psi_2$ be as in Theorem \ref{ThmProduct}. The we have the following isomorphism.
\begin{align*}
f:\mathrm{QI}\big((N_1\times N_2)\rtimes_A\mathbb{R}\big)\big/\mathord{\sim}&\to\mathrm{Bilip}\big(N_1\big)\times\mathrm{Bilip}\big(N_2\big)
\\\Phi&\mapsto\big(\Psi_1,\Psi_2\big).
\end{align*} 
\end{thm}

This distance is related to the height divergence of vertical geodesics in the following way.

\begin{lemma}\label{LemmaBackExtended}(Extended Backward Lemma)
Let $n,m\in N_i$, let $V:t\mapsto(n,t)$ and let $W:t\mapsto(m,t)$, then
\begin{align*}
d_H(n,m)\asymp_{k,c,\bowtie}\exp\left(h_{\mathrm{Div}}(V,W)\right).
\end{align*}
\end{lemma}

See Corollary \ref{CoroBackward} for the definition of $h_{\mathrm{Div}}(V,W)$.

\begin{proof}
By the Corollary \ref{CoroBackward} there exists a height $h_{\mathrm{Div}}(V,W)\in\mathbb{R}$ such that $V$ and $W$ diverge from each other at the height $h_{\mathrm{Div}}(V,W)$. Hence there exists $M(k,c,\bowtie)$ such that for all $s_1\leq s_2\leq h_{\mathrm{Div}}(V,W)$
\begin{align*}
d\big(V(s_2),W(s_2)\big)-M&\leq d_{S_i}\big(V(s_1),W(s_1)\big)+2|s_2-s_1|\leq d_{S_i}\big(V(s_2),W(s_2)\big)+M.
\end{align*}
Therefore
\begin{align}
\exp\left(d_{S_i}\big(V(s_1),W(s_1)\big)+2|s_2-s_1|\right)&\asymp_{k,c,\bowtie} \exp\left(d_{S_i}\big(V(s_2),W(s_2)\big)\right).\label{IneqUse1LemmaFinExemple}
\end{align}
Let us denote $h_0:=h_{\mathrm{Div}}(V,W)$. Then we can compute the Hamenstädt distance $d_H(n,m)$
\begin{align*}
d_H(n,m)&=\exp\left(-\frac{1}{2}\lim\limits_{s\rightarrow+\infty}\Big(2s-d_{S_i}\big(V(-s), W(-s)\big)\Big)\right),
\\&\asymp_{k,c,\bowtie}\exp\left(-\frac{1}{2}\lim\limits_{s\rightarrow+\infty}\Big(2s-d_{S_i}\big( V(h_0), W(h_0)\big)-(2h_0+2s)\Big)\right),\quad\text{by inequality (\ref{IneqUse1LemmaFinExemple})},
\\&\asymp_{k,c,\bowtie}\exp\left(-\frac{1}{2}\lim\limits_{s\rightarrow+\infty}\Big(-d_{S_i}\big(V(h_0), W(h_0)\big)-2h_0\Big)\right),
\\&=\exp\left(\frac{d_{S_i}\big(V(h_0), W(h_0)\big)}{2}+h_0\right)=\exp\left(\frac{d_{S_i}\big(V(h_0), W(h_0)\big)}{2}\right)\exp\left(h_0\right),
\\&\asymp_{k,c,\bowtie}\exp\left(h_0\right),\quad\text{by definition of }h_{\mathrm{Div}}(V,W).
\end{align*}
\end{proof}
%

We show that the aforementioned maps $\Psi_i$ are bilipschitz.

\begin{thm}\label{ThmBilip}
Let $\Psi_i$ be the map of Theorem \ref{ThmProduct}. Then $\Psi_i$ is a bilipschitz homeomorphism either from $(N_i,d_H)$ to $\left(N_i',(d_H)^{\frac{\mathrm{tr}(A_1)}{\mathrm{tr}(A_1')}}\right)$ or from $\left(N_i,(d_H)^{\frac{\mathrm{tr}(A_1')}{\mathrm{tr}(A_1)}}\right)$ to $(N_i',d_H)$.
\end{thm}
\begin{proof}
Let $n,m\in N_i$ and let $V:t\mapsto (n,t)$ and $W:t\mapsto (m,t)$ be two vertical geodesics of $N_i\rtimes_{A_i} \mathbb{R}$. Let us denote by $\lambda_0:=\frac{\mathrm{tr}(A_1)}{\mathrm{tr}(A_1')}$. By Lemma \ref{LemmaBackExtended} we have 
\begin{align*}
d_H(n,m)\asymp_{k,c,\bowtie}\exp\left(h_{\mathrm{Div}}(V,W)\right).
\end{align*}  
Since $\title{\Phi}_i:=(\Psi_i,\lambda_0\mathrm{id}_{\mathbb{R}}+t_0)$ is a $(k',c')$-quasi-isometry, we have:
\begin{enumerate}
\item$d_{S_i}\big((\Psi_i(n),\lambda_0 h_{\mathrm{Div}}(V,W)+t_0),(\Psi_i(m),\lambda_0 h_{\mathrm{Div}}(V,W)+t_0)\big)\asymp_{k,c,\bowtie} 1$;
\item$\forall s\geq h_{\mathrm{Div}}(V,W)$, $d_{S_i}\big((\Psi_i(n),\lambda_0 s+t_0),(\Psi_i(m),\lambda_0 s+t_0)\big)\preceq_{k,c,\bowtie} 1$.
\end{enumerate}  
Furthermore, for all $n\in N_i$, $\tilde{\Phi}_i(V_n)=V_{\Psi_i(n)}$ hence $\tilde{\Phi}_i(V_n)$ is a vertical geodesics of $S_i'$. Then there exists $M(k,c,\bowtie)$ such that
\begin{align*}
\left(\lambda_0 h_{\mathrm{Div}}(V,W)+t_0\right)-M\leq h_{\mathrm{Div}}\left(\tilde{\Phi}_i(V),\tilde{\Phi}_i(W)\right)\leq \left(\lambda_0 h_{\mathrm{Div}}(V,W)+t_0\right)+M.
\end{align*}
Consequently Lemma \ref{LemmaBackExtended} provides us with
\begin{align*}
d_H\left(\Psi_i(n),\Psi_i(m)\right)&\asymp_{k,c,\bowtie}\exp\left(h_{\mathrm{Div}}(V_{\Psi_i(n)},W_{\Psi_i(m)})\right)=\exp\left(h_{\mathrm{Div}}\left(\tilde{\Phi}_i(V),\tilde{\Phi}_i(W)\right)\right),
\\&\asymp_{k,c,\bowtie}\exp(t_0)\exp\left(\lambda_0 h_{\mathrm{Div}}\left(V,W\right)\right),
\\&\asymp_{k,c,\bowtie}\exp(t_0)\left(d_H(n,m)\right)^{\lambda_0},\quad\text{ by Lemma \ref{LemmaBackExtended}}.
\end{align*} 
Here $t_0$ depends only on $\Phi$. Furthermore, if $\lambda_0\leq 1$, $(d_H)^{\lambda_0}$ is still a distance by concavity. Hence, depending on the value of $\lambda_0$, either $\Psi_i:\big(N_i,d_H\big)\to\big(N_i',(d_H)^{\lambda_0}\big)$ or $\Psi_i:\big(N_i,(d_H)^{\lambda_0}\big)\to\big(N_i',d_H\big)$ is a bilipschitz map.
\end{proof}

We now focuses on self quasi-isometries of $(N_1\times N_2)\rtimes_A\mathbb{R}$.

\begin{proof}[Proof of Theorem \ref{ThmQISolvable}:]
Let $\Psi_1$, $\Psi_2$ be as in Theorem \ref{ThmProduct}, and let $f$ be the map
\begin{align*}
f:\mathrm{QI}\big((N_1\times N_2)\rtimes_A\mathbb{R}\big)\big/\mathord{\sim}&\to\mathrm{Bilip}\big(N_1\big)\times\mathrm{Bilip}\big(N_2\big)
\\\Phi&\mapsto\big(\Psi_1,\Psi_2\big).
\end{align*}
We first show that this application is well defined. Let $\Phi,\Phi'\in\mathrm{QI}\big((N_1\times N_2)\rtimes_A\mathbb{R}\big)$ be such that $\Phi\sim\Phi'$, which means that $d_{\bowtie}(\Phi,\Phi')\preceq_{k,c,\bowtie} 1$. 

By Theorems \ref{ThmProduct} and \ref{ThmBilip}, there exist $\Psi_i,\Psi'_i\in \mathrm{Bilip}(N_i)$ such that:
\begin{enumerate}
\item$d\big(\Phi,(\Psi_1,\Psi_2,\mathrm{id}_{\mathbb{R}})\big)\preceq_{k,c,\bowtie}1$;
\item$f(\Phi)=(\Psi_1,\Psi_2)$;
\item$d\big(\Phi',(\Psi'_1,\Psi'_2,\mathrm{id}_{\mathbb{R}})\big)\preceq_{k,c,\bowtie}1$;
\item$f(\Phi')=(\Psi'_1,\Psi'_2)$.
\end{enumerate} 
By the definition of $\Psi_i$ and $\Psi_i'$ , for all $n\in N$ we have
\begin{align*}
\Psi_i(n)=& P\left(V_n'(0)\right);
\\\Psi_i'(n)=& P\left(V_n''(0)\right),
\end{align*} 
where $V_n'$ is the unique vertical geodesic close to $\hat{\Phi}_i(V_n)$ and $V_n''$ the unique vertical geodesic close to $\hat{\Phi}'_i(V_n)$. However $\Phi\sim\Phi'$, then $\hat{\Phi}_i(V_n)$ and $\hat{\Phi}_i'(V_n)$ are $M$-close to each other for some $M(k,c,\bowtie)$, therefore $d_{\mathrm{Hff}}(V_n',V_n'')\preceq_{k,c,\bowtie} 1$. However these vertical geodesics are unique, then $V_n'=V_n''$. Consequently, $\Psi_i(n)=\Psi'_i(n)$, hence $\Psi_i=\Psi_i'$, therefore $f$ is well defined.
\\\\Let us now prove that $f$ is injective. Let $\Phi$ and $\Phi'$ be two quasi-isometries of $(N_1\times N_2)\rtimes_A\mathbb{R}$ such that $f(\Phi)=f(\Phi')$. Then by Theorem \ref{ThmProduct} and by the triangle inequality
\begin{align*}
d_{\bowtie}\left(\Phi,\Phi'\right)\leq d_{\bowtie}\left(\Phi,(\Psi_1, \Psi_2,\mathrm{id}_{\mathbb{R}})\right)+d_{\bowtie}\left((\Psi_1, \Psi_2,\mathrm{id}_{\mathbb{R}}),\Phi'\right)\preceq_{k,c,\bowtie,\Phi,\Phi'}1.
\end{align*}
Hence $\Phi\sim \Phi'$, which proves that $f$ is injective.
\\\\Let $\Psi_i\in\mathrm{Bilip}\big(N_i,d_H\big)$, our goal is to show that $(\Psi_i,\mathrm{id}_{\mathbb{R}})$ is a quasi-isometry of $(N_i\rtimes_{A} \mathbb{R},d_{S_i})$. Let $(n,t_n),(m,t_m)\in S_i$.  By Lemma \ref{LemmaBackExtended} applied on $n$ and $m$, there exists a constant $M(k,c,\bowtie)$ such that
\begin{align}
\ln\left(d_H\big(n,m\big)\right)-M\leq h_{\mathrm{Div}}\left(V_{n},V_{m}\right)\leq\ln\left(d_H\big(n,m\big)\right)+M.\label{Ineq559}
\end{align}
Similarly, by Lemma \ref{LemmaBackExtended} applied on $\Psi_i(n)$ and $\Psi_i(m)$
\begin{align}
\ln\left(d_H\big(\Psi_i(n),\Psi_i(m)\big)\right)-M\leq h_{\mathrm{Div}}\left(V_{\Psi_i(n)},V_{\Psi_i(m)}\right)\leq\ln\left(d_H\big(\Psi_i(n),\Psi_i(m)\big)\right)+M.\label{Ineq459}
\end{align}
We know that $\Psi_i\in\mathrm{Bilip}\big(N_i,d_H\big)$ hence $d_H(n,m)\asymp d_H\big(\Psi_i(n),\Psi_i(m)\big)$. Therefore by inequalities (\ref{Ineq559}) and (\ref{Ineq459}) we have
\begin{align}
\left|h_{\mathrm{Div}}\left(V_{n},V_{m}\right)-h_{\mathrm{Div}}\left(V_{\Psi_i(n)},V_{\Psi_i(m)}\right)\right|\preceq 1.\label{Ineq659}
\end{align}
Moreover by Lemma \ref{LemmaBackward} we can characterise the distance between two points thanks to the height of divergence of their associated vertical geodesics. Let us denote $h_0=h_{\mathrm{Div}}\left(V_{n},V_{m}\right)$. By inequality (\ref{Ineq659}) and by Lemma \ref{LemmaBackward}, if $h_0\geq \max(t_n,t_m)$ we have both:
\begin{align*}
&\left|d_{S_i}\big((n,t_n),(m,t_m)\big)-\Big(|t_m-h_0|+|t_n-h_0|\Big)\right|\preceq_{\delta} 1;
\\&\left|d_{S_i}\big((\Psi_i(n),t_n),(\Psi_i(m),t_m)\big)-\Big(|t_m-h_0|+|t_n-h_0|\Big)\right|\preceq_{\delta} 1.
\end{align*}
Consequently by the triangle inequality there exists $M(\delta)$ such that
\begin{align*}
d_{S_i}\big((n,t_n),(m,t_m)\big)-M\leq d_{S_i}\big((\Psi_i(n),t_n),(\Psi_i(m),t_m)\big)\leq d_{S_i}\big((n,t_n),(m,t_m)\big)+M.
\end{align*} 
Similarly, if $h_0\leq \max(t_n,t_m)$ we have both:
\begin{align*}
&\left|d_{S_i}\big((n,t_n),(m,t_m)\big)-\Big(|t_m-t_n|\Big)\right|\preceq_{\delta} 1;
\\&\left|d_{S_i}\big((\Psi_i(n),t_n),(\Psi_i(m),t_m)\big)-\Big(|t_m-t_n|\Big)\right|\preceq_{\delta} 1.
\end{align*}
Hence again
\begin{align*}
d_{S_i}\big((n,t_n),(m,t_m)\big)-M\leq d_{S_i}\big((\Psi_i(n),t_n),(\Psi_i(m),t_m)\big)\leq d_{S_i}\big((n,t_n),(m,t_m)\big)+M.
\end{align*} 
Therefore  $(\Psi_i,\mathrm{id}_{\mathbb{R}})$ is a $(1,M)$-quasi-isometry of $N_i\rtimes \mathbb{R}$, hence $(\Psi_1,\Psi_2,\mathrm{id}_{\mathbb{R}})$ is also a $(1,M)$-quasi-isometry, which provides us with $f(\Psi_1,\Psi_2,\mathrm{id}_{\mathbb{R}})=(\Psi_1,\Psi_2)$. Hence $f$ is surjective, and finally bijective.
\\\\Let us now prove that $f$ is a morphism. Let $\Phi$,$\Phi'\in\mathrm{QI}\big((N_1\times N_2)\rtimes_A\mathbb{R}\big)$.
Furthermore, $d_{\bowtie}\left(\Phi',(\Psi_1',\Psi_2',\mathrm{id}_{\mathbb{R}})\right)\preceq 1$, hence $d_{\bowtie}\left(\Phi\circ\Phi',\Phi\circ(\Psi_1',\Psi_2',\mathrm{id}_{\mathbb{R}})\right)\preceq 1$ since $\Phi$ is a quasi-isometry. Moreover, $d_{\bowtie}\left(\Phi,(\Psi_1,\Psi_2,\mathrm{id}_{\mathbb{R}})\right)\preceq 1$, therefore by the triangle inequality
\begin{align*}
d_{\bowtie}\left(\Phi\circ\Phi',(\Psi_1,\Psi_2,\mathrm{id}_{\mathbb{R}})\circ(\Psi_1',\Psi_2',\mathrm{id}_{\mathbb{R}})\right)\preceq 1.
\end{align*}
However
\begin{align*}
(\Psi_1,\Psi_2,\mathrm{id}_{\mathbb{R}})\circ (\Psi_1',\Psi_2',\mathrm{id}_{\mathbb{R}})=(\Psi_1\circ \Psi_1',\Psi_2\circ\Psi_2',\mathrm{id}_{\mathbb{R}}),
\end{align*}
which provides us with
\begin{align*}
d_{\bowtie}\left(\Phi\circ\Phi',(\Psi_1\circ \Psi_1',\Psi_2\circ\Psi_2',\mathrm{id}_{\mathbb{R}})\right)\preceq 1.
\end{align*}
Consequently $f(\Phi\circ\Phi')=(\Psi_1\circ \Psi_1',\Psi_2\circ\Psi_2')$.
\end{proof}

In this proof we showed that $\Phi\sim(\Psi_1,\Psi_2,\mathrm{id}_{\mathbb{R}})$, therefore any quasi-isometry is in the equivalence class of an $(1,M)$-quasi-isometry. 

\subsection{Quasi-isometric classification and necessary conditions to being quasi-isometric}

Thanks to Proposition \ref{PropoQIinvariant} and Theorem \ref{ThmBilip} we are able to provide necessary conditions and quasi-isometric classifications for families of solvable Lie groups of the form $\mathbb{R}\ltimes_{\mathrm{Diag}(A_1,A_2)} (N_1\times N_2)$. 

Let us recall two consequences implied by being quasi-isometric in the Lie group setting. For $i\in\{1,2\}$, let $N_i$, $N_i'$ be two simply connected, nilpotent Lie groups and let $A_i$, $A_i'$ be two matrices whom eigenvalues have positive real parts, acting by derivation on the corresponding Lie algebra. Let us assume that $\mathrm{tr}(A_1)>\mathrm{tr}(A_2)$ and $\mathrm{tr}(A_1')>\mathrm{tr}(A_2')$. If $\mathbb{R}\ltimes_{\mathrm{Diag}(A_1,A_2)} (N_1\times N_2)$ and $\mathbb{R}\ltimes_{\mathrm{Diag}(A_1',A_2')} (N_1'\times N_2')$ are quasi-isometric then:
\begin{enumerate}
\item $\frac{\mathrm{tr}(A_1)}{\mathrm{tr}(A_2)}=\frac{\mathrm{tr}(A_1')}{\mathrm{tr}(A_2')}$ (Proposition \ref{PropoQIinvariant});
\item For $i\in\{1,2\}$, $N_i$ and $N_i'$ are bilipschitz. (Theorem \ref{ThmBilip}).
\end{enumerate}
Let us denote by
\begin{align*}
S_{N_1,N_2}:=\mathbb{R}\ltimes_{\mathrm{Diag}(A_1,-A_2)} (N_1\times N_2).
\end{align*}

Combining Lemma $4.1$ of paper \cite{PS} and Theorem \ref{ThmBilip} we obtain the following statement.

\begin{propo}
Let us assume that $\mathrm{tr}(A_1)>\mathrm{tr}(A_2)$ and $\mathrm{tr}(A_1')>\mathrm{tr}(A_2')$. If $S_{N_1,N_2}$ and $S_{N_1',N_2'}$ are quasi-isometric, then we have that for $i\in\{1,2\}$, $A_i$ and $\frac{\mathrm{tr}(A_1)}{\mathrm{tr}(A_1')}A_i'$ share the same characteristic polynomial.
\end{propo}

A \textit{Carnot} group $N$ is a simply connected, nilpotent Lie group with a Lie algebra $\mathrm{Lie}(N)$ which admits a grading: there exists a family of subspaces $V_i$ with $i\in\{1,...,r\}$ for some $r\geq 1$ such that $V_{i+1} =[V_1,V_i]$ for $i<r$ and such that  
\begin{align*}
\mathrm{Lie}(N)=\bigoplus\limits_{i=1}^r V_i.
\end{align*}
A Carnot group is equipped with a 1-parameter family of automorphisms called \textit{dilations} on $N$ and defined for $t\in\mathbb{R}$ by $\delta_t:=\exp(tD)$, with $D$ a Lie derivation on $\mathrm{Lie}(N)$ verifying that $Dv=iv$ for $v\in V_i$ and $i\in\{1,...,r\}$. Such a derivation is called a \textit{Carnot derivation}. A Lie group $S(N_1,N_2)$ is \textit{Carnot-Sol type} if $N_1$ and $N_2$ are Carnot groups and if their respective derivations $A_1$ and $A_2$ are Carnot derivations. Combining Theorem \ref{ThmBilip} and Theorem 2 of \cite{Pansu}, we get the following necessary condition 

\begin{propo}\label{QIClass2}
Let $S(N_1,N_2)$ and $S(N_1',N_2')$ be two Carnot-Sol type Lie groups and assume that $\mathrm{tr}(A_1)>\mathrm{tr}(A_2)$ and that $\mathrm{tr}(A_1')>\mathrm{tr}(A_2')$. Then: 
\begin{align*}
S_{N_1,N_2}\text{ and }S_{N_1',N_2'}\text{ are quasi-isometric}\quad\Rightarrow\quad\text{For }i\in\{1,2\},\  N_i\text{ and }N_i'\text{ are isomorphic.}
\end{align*}
\end{propo}

Furthermore, for a given Carnot derivation $A$ on a Carnot group $N$, there exists a positive real $\alpha > 0$ such that $\mathbb{R}\ltimes_{A} N = S_{\alpha}$ where $S_{\alpha}:=\mathbb{R}\ltimes_{\alpha} N$ is  the group defined by the action of $\mathbb{R}$ via the dilation $(\delta_{\alpha t})_{t\in\mathbb{R}}$ on $N$. Let $N_1$ and $N_2$ be two Carnot groups and for any two positive reals $\alpha,\beta>0$, let $G_{\alpha,\beta}:=\mathbb{R}\ltimes_{\alpha,-\beta}(N_1\times N_2)$ be the group defined by the action of $\mathbb{R}$ on $N_1\times N_2$,
\begin{align*}
\mathbb{R}\to\mathrm{Aut}(N\times N),\quad t\mapsto\left(\delta_{\alpha t},\delta_{-\beta t}\right).
\end{align*}
Note that $G_{\alpha,\beta}=S_{\alpha}\bowtie S_{\beta}$. Thanks to the quasi-isometry invariant of Proposition \ref{PropoQIinvariant}, we obtain the quasi-isometry classification for Carnot-Sol type Lie groups. 

\begin{propo}\label{QIClass3}
Let $(\alpha,\beta)$ and $(\sigma,\tau)$ be two pairs of positive reals with $\alpha>\beta$ and $\sigma>\tau$, then
\begin{align*}
G_{\alpha,\beta} \text{ quasi-isometric to } G_{\sigma,\tau}\quad \Leftrightarrow \quad\frac{\alpha}{\beta}=\frac{\sigma}{\tau}\quad \Leftrightarrow \quad G_{\alpha,\beta} \text{ isomorphic to } G_{\sigma,\tau}.
\end{align*}
\end{propo}

\begin{proof}
If $\frac{\alpha}{\beta}=\frac{\sigma}{\tau}$ , then $G_{\alpha,\beta}$ and $G_{\sigma,\tau}$ are isomorphic and thus in particular quasi-isometric (or even bilipschitz) with respect to any left-invariant Riemannian metrics on the groups. Indeed, the map 
\begin{align*}
G_{\alpha,\beta}\to G_{\sigma,\tau},\quad (x,y,t)\mapsto(x,y,\lambda t).
\end{align*}
is an isomorphism. For $(x_i,y_i,t_i)\in G_{\alpha,\beta}$ for $i\in\{1,2\}$, we have in $G_{\sigma,\tau}$
\begin{align*}
(x_1,y_1,\lambda t_1)\cdot(x_2,y_2,\lambda t_2)&=(x_1\cdot \delta_{\sigma\lambda t_1}x_2,y_1\cdot\delta_{-\tau\lambda t_1}y_2,\lambda (t_1+t_2)),
\\&=(x_1\cdot \delta_{\alpha t_1}x_2,y_1\cdot\delta_{-\beta t_1}y_2,\lambda (t_1+t_2)),
\end{align*}
which is the image of $(x_1,y_1, t_1)\cdot(x_2,y_2, t_2)$. Proposition \ref{PropoQIinvariant} conclude the proof since the ratios of traces of the respective derivations are $\frac{\alpha}{\beta}$ and $\frac{\sigma}{\tau}$.  
\end{proof}

Otherwise stated, two non-unimodular Carnot-Sol type solvable Lie groups are quasi-isometric if and only if they are isomorphic.

\end{document}